\documentclass[a4paper,12pt]{report}
\usepackage{amsmath}
\usepackage{amsthm}
\usepackage{amsfonts}
\usepackage{amssymb}
\usepackage{graphicx}
\usepackage[all]{xy}
\usepackage{a4wide}
\usepackage{enumerate}
\usepackage{mathrsfs}
\newtheorem{theorem}{Theorem}[section]
\newtheorem{proposition}[theorem]{Proposition}
\newtheorem{corollary}[theorem]{Corollary}
\newtheorem{lemma}[theorem]{Lemma}
\theoremstyle{definition}

\newtheorem{definition}[theorem]{Definition}
\newtheorem{definitions}[theorem]{Definitions}
\newtheorem{example}[theorem]{Example}
\newtheorem{examples}[theorem]{Examples}

\newtheorem{remark}[theorem]{Remark}
\newtheorem{remarks}[theorem]{Remarks}

\newcommand{\bcomod}[2]{{}^{#1}\mathcal{M}^{#2}}
\newcommand{\cat}[1]{\mathbf{#1}}
\newcommand{\cohom}[3]{\mathrm{h}_{#1}(#2,#3)}
\newcommand{\coring}[1]{\mathfrak{#1}}
\newcommand{\cotensor}[1]{\square_{#1}}
\renewcommand{\hom}[3]{\mathrm{Hom}_{#1}(#2,#3)}
\newcommand{\lcomod}[1]{{}^{#1}\mathcal{M}}
\newcommand{\ldual}[1]{{ }^*#1}
\newcommand{\rdual}[1]{#1^*}
\newcommand{\lDual}{{ }^*(-)}
\newcommand{\rDual}{(-)^*}
\newcommand{\rcomod}[1]{\mathcal{M}^{#1}}
\newcommand{\rmod}[1]{\mathcal{M}_{#1}}
\newcommand{\lmod}[1]{{}_{#1}\mathcal{M}}
\renewcommand{\bmod}[2]{{}_{#1}\mathcal{M}_{#2}}
\newcommand{\tensor}[1]{\otimes_{#1}}
\newcommand{\nat}[2]{\mathrm{Nat}(#1,#2)}
\newcommand{\h}[2]{\mathrm{H}(#1,#2)}
\newcommand{\e}[2]{\mathrm{e}_{#1}(#2)}
\newcommand{\pic}[1]{\mathrm{Pic}_k(#1)}
\newcommand{\rpic}[1]{\mathrm{Pic}_k^r{(#1)}}
\newcommand{\lpic}[1]{\mathrm{Pic}_k^l{(#1)}}
\newcommand{\aut}[1]{\mathrm{Aut}_k(#1)}
\newcommand{\inn}[1]{\mathrm{Inn}_k(#1)}
\newcommand{\rinn}[1]{\mathrm{Inn}_k^r{(#1)}}
\newcommand{\out}[1]{\mathrm{Out}_k(#1)}
\newcommand{\rout}[1]{\mathrm{Out}_k^r{(#1)}}
\newcommand{\oper}[2]{\operatorname{#1}(#2)}

\begin{document}

\thispagestyle{empty}

\begin{center}
{\Large University of Granada \\ Faculty of
Sciences\\
Department of Algebra}
\end{center}

\vspace{4cm}
\begin{center}
\textbf{Ph.D. Thesis}

\vspace{3cm} {\LARGE \textbf{Adjoint and Frobenius Pairs of
Functors, Equivalences and the Picard Group for Corings}}

\vspace{3cm}
{\Large Mohssin Zarouali Darkaoui}\\

 \vspace{4cm} Granada, 2007

\end{center}

\newpage
\thispagestyle{empty}

\begin{center}
{\Large Adjoint and Frobenius Pairs of Functors, Equivalences and
the Picard Group for Corings}

\vspace{4cm}
\textbf{Ph.D. Thesis realized by Mohssin Zarouali Darkaoui
under the supervision of Prof. Dr. Jos\'e G\'omez Torrecillas}
(University of Granada)

\vspace{5cm}
Thesis defended on  May 7, 2007, at the Faculty of
Sciences of the University of Granada, whose Jury members are

\vspace{2cm}
\textbf{President:} Prof. Dr. Jos\'e Luis G\'omez Pardo
(University of Santiago, Spain)\\
Prof. Dr. Tomasz Brzezi\'nski (University of Wales Swansea, U.K)\\
Prof. Dr.  Constantin N\u{a}st\u{a}sescu (University of Bucharest,
Romania)\\
Prof. Dr. Blas Torrecillas Jover (University of Almeria, Spain)\\
\textbf{Secretary:} Prof. Dr.  Jos\'e Luis Bueso Montero
(University of Granada, Spain)
\end{center}

\vspace{2cm}
Degree obtained with the highest distinction (``Sobresaliente Cum
Laude'') (by unanimity).








\tableofcontents

\chapter*{Introduction}

The aim of this work is:
\begin{enumerate}[1.]
\item to study adjoint and Frobenius pairs of functors between comodule
categories over corings;

\item to extend the comatrix coring defined in
\cite{ElKaoutit/Gomez:2003,Brzezinski/Gomez:2003} to a quasi-finite
comodule over a coring and to study equivalences between of comodule
categories over corings. We think that our generalized comatrix
coring gives more concrete characterizations of these equivalences;

\item to extend the Picard group from algebras \cite{Bass:1967,Bass:1968}
and coalgebras over fields \cite{Torrecillas/Zhang:1996} to corings
and extend some of their properties;

\item and to apply the obtained results to induction functors over
different base rings and to corings coming from entwined structures
and graded ring theory, in order to give more concrete and new
results. For example, we give a characterization of adjoint and
Frobenius pairs of functors between categories of graded modules
over $G$-sets.
\end{enumerate}

We conclude that the notion of a coring unifies several important
algebraic structures of current interest and is a powerful tool
which allows to obtain more concrete and simple results than the
already existing ones.

Corings were introduced by M. Sweedler in \cite{Sweedler:1975} as a
generalization of coalgebras over commutative rings to the case of
non-commutative rings, to give a formulation of a predual of the
Jacobson-Bourbaki's theorem for intermediate extensions of division
ring extensions. Thus, a coring over an associative ring with unit
$A$ is a comonoid in the monoidal category of all $A$-bimodules.
Further study of corings are given by several authors, for example,
see \cite{Guzman:1985,Guzman:1989}. It is well known that algebras
and coalgebras over fields (see
\cite{Sweedler:1969,Abe:1977,Dascalescu/Nastasescu/Raianu: 2001}),
or more generally over rings (see \cite{Brzezinski/Wisbauer:2003}),
are special cases of corings.

Recently, corings were intensively studied. The main motivation of
these studies is an observation by Takeuchi, namely that an
entwining structure (resp. an entwined module) can be viewed as a
suitable coring (resp. as a comodule over a suitable coring) (see
Theorem \ref{Takeuchi}). T. Brzezi\'{n}ski has given in
\cite{Brzezinski:2002} some new examples and general properties of
corings. Among them, a study of Frobenius corings is developed,
extending previous results on entwining structures
\cite{Brzezinski:2000} and relative Hopf modules
\cite{Caenepeel/Militaru/Zhu:1997}.

A pair of functors $(F,G)$ is said to be a \emph{Frobenius pair}
\cite{Castano/Gomez/Nastasescu:1999}, if $G$ is at the same time a
left and right adjoint to $F$. That is a standard name which we use
instead of Morita's original ``\emph{strongly adjoint pairs}''
\cite{Morita:1965}. The functors $F$ and $G$ are known as
\emph{Frobenius functors} \cite{Caenepeel/Militaru/Zhu:1997}. The
study of Frobenius functors was motivated by a paper of K. Morita,
where he proved \cite[Theorem 5.1]{Morita:1965} that given a ring
extension $i:A\to B$, the induction functor
$-\tensor{A}B:\rmod{A}\to \rmod{B}$ is a Frobenius functor if and
only if the morphism $i$ is Frobenius in the sense of
\cite{Kasch:1961} (see also \cite{Nakayama/Tsuzuku:1960}): $_AB$ is
finitely generated projective and $_BB_A\simeq\hom{A}{_AB}{{}_AA}$
as $(B,A)$-bimodules. There is another interesting result relating
Forbenius extensions of rings and corings \cite[Proposition
4.3]{Kadison:1999b}: A ring extension $i:A\to B$ is a Frobenius
extension if and only if the bimodule $_AB_A$ is endowed by a
structure of $A$-coring such that the comultiplication is a
$B$-bimodule map.

The dual result of Morita's \cite[Theorem 5.1]{Morita:1965} for
coalgebras over fields was proved in \cite[Theorem
3.5]{Castano/Gomez/Nastasescu:1999} and it states that the
corestriction functor $(-)_\varphi:\rcomod{C}\to \rcomod{D}$
associated to a morphism of coalgebras $\varphi:C\rightarrow D$ is
Frobenius if and only if $C_D$ is quasi-finite and injective and
there exists an isomorphism of bicomodules $\cohom{D}{{}_CC_D}{D}
\simeq {}_DC_C$ (here, $\cohom{D}{C}{-}$ denotes the ``cohom''
functor). Since corings generalize both rings and coalgebras over
fields, one may expect that \cite[Theorem 5.1]{Morita:1965} and
\cite[Theorem 3.5]{Castano/Gomez/Nastasescu:1999} are
specializations of a general statement on homomorphisms of corings.
In Chapter 2, we find such a result (Theorem \ref{23}) and introduce
the notion of a (right) Frobenius extension of corings (see
Definition \ref{Frobenius extension}). To prove such a result, we
study adjoint pairs and Frobenius pairs of functors between
categories of comodules over rather general corings. Precedents for
categories of modules and categories of comodules over coalgebras
over fields are contained in \cite{Morita:1965},
\cite{Takeuchi:1977}, \cite{Castano/Gomez/Nastasescu:1999} and
\cite{Caenepeel/DeGroot/Militaru:2002}. More recently,
\emph{Frobenius corings} (i.e., corings for which the functor
forgetting the coaction is Frobenius), have intensively been studied
in \cite{Brzezinski:2002,Brzezinski:2003,Brzezinski/Gomez:2003,
Brzezinski/Wisbauer:2003,Caenepeel/DeGroot/Militaru:2002}.

We think that our general approach produces results of independent
interest, beyond the aforementioned extension to the coring setting
of \cite[Theorem 5.1]{Morita:1965} and \cite[Theorem
3.5]{Castano/Gomez/Nastasescu:1999} and contributes to the
understanding of the behavior of the cotensor product functor for
corings. In fact, our general results, although they are sometimes
rather technical, have  other applications and will probably find
more. For instance, we have used them to prove new results on
equivalences of comodule categories over corings that are expounded
in Chapter 3. Moreover, when applied to corings stemming from
different algebraic theories of current interest, they boil down to
new (more concrete) results. As an illustration, we consider
entwined modules over an entwining structure in Section
\ref{entwined} and graded modules over $G$-sets in Section
\ref{graded}.

\medskip
Comatrix corings were introduced in \cite{ElKaoutit/Gomez:2003} in
order to give the structure of corings whose category of right
comodules has a finitely generated projective generator, and the
structure of cosemisimple corings which were defined and studied in
\cite{ElKaoutit/Gomez/Lobillo:2001unp}. In
\cite{Brzezinski/Gomez:2003} the authors have given an equivalent
definition of them and studied when are cosplit, coseparable, and
Frobenius.

In Section \ref{comatrixgeneralized}, we extend comatrix coring to
the case of quasi-finite comodules and we study when is cosplit,
coseparable, and Frobenius. Suppose that $_A\coring{C}$ and
$_B\coring{D}$ are flat. Let $X\in\bcomod{\coring{C}}{\coring{D}}$
and $\Lambda\in\bcomod{\coring{D}}{\coring{C}}$ and suppose that
$-\square _{\coring{C}}X$ is a left adjoint to
$-\square_{\coring{D}}\Lambda$. If $\coring{C}_A$ and $\coring{D}_B$
are flat and $_{\coring{C}}X$, $_{\coring{D}}\Lambda$ are coflat, or
$A$ and $B$ are von Neumann or if $\coring{C}$ and $\coring{D}$ are
coseparable, then the $B$-bimodule $\Lambda\cotensor{\coring{C}}X$
is endowed with a structure of $B$-coring (Theorem
\ref{Comatrixcoring}). The coproduct is
 \[
\xymatrix{\Delta:\Lambda\cotensor{\coring{C}}X \ar[r]^-\simeq &
\Lambda\cotensor{\coring{C}}(
\coring{C}\cotensor{\coring{C}}X)\ar[rr]^-{\Lambda\cotensor{\coring{C}}
(\eta_\coring{C}\cotensor{\coring{C}}X)} & &
\Lambda\cotensor{\coring{C}}((X\cotensor{\coring{D}}\Lambda)
\cotensor{\coring{C}}X)\ar[r]^-\simeq &
(\Lambda\cotensor{\coring{C}}X)\cotensor{\coring{D}}(
\Lambda\cotensor{\coring{C}}X) \ar@{^{(}->}[d] \\ & & & &
(\Lambda\cotensor{\coring{C}}X)\tensor{B}(\Lambda
\cotensor{\coring{C}}X),}
 \]
 where $\eta$ is the counit of this adjunction, and the counit is
$\xymatrix{\epsilon=\epsilon_{\coring{D}}\circ\omega
:\Lambda\square_{\coring{C}}X\ar[r] & B}$.

The most important property of this coring is the following: Let
$\Lambda\in \bcomod{\coring{D}}{\coring{C}}$ be a bicomodule,
quasi-finite as a right $\coring{C}$-comodule, such that
$_A\coring{C}$ and $_B\coring{D}$ are flat. Set
$X=\cohom{\coring{C}}{\Lambda}{\coring{C}}\in
\bcomod{\coring{D}}{\coring{C}}$.
 If
\begin{enumerate}[(a)]
\item $\coring{C}_A$ and $\coring{D}_B$ are flat,
 the cohom functor $\cohom{\coring{C}}{\Lambda}{-}$ is exact and
$_{\coring{D}}{\Lambda}$ is coflat, or
\item $A$ and $B$ are von Neumann regular rings and the cohom
functor $\cohom{\coring{C}}{\Lambda}{-}$ is exact or
\item $\coring{C}$ and
$\coring{D}$ are coseparable corings,
\end{enumerate}
then the canonical isomorphism
$\xymatrix{\delta_{\Lambda}:\e{\coring{C}}{\Lambda}\ar[r] &
\Lambda\cotensor{\coring{C}}X}$ is an isomorphism of $B$-corings
(see Proposition \ref{comatrixproperty1}). We think that, under
these conditions, this coring isomorphism gives a concrete
description of the coendomorphism coring even in the case of
coalgebras over fields. Indeed, let $C$ be a coalgebra over a field
$k$ and $\Lambda\in\rcomod{C}$ a quasi-finite and injective
comodule. Then $\e{C}{\Lambda}\simeq
\Lambda\cotensor{C}\cohom{C}{\Lambda}{C}$ as coalgebras (the case
(b)) (see Corollary \ref{coendcoalgebra}). We also think that the
last isomorphism gives a more concrete characterization of
equivalence between categories of comodules over coalgebras over a
field, than that of Takeuchi \cite{Takeuchi:1977}. Of course, this
characterization of the coendomorphism coring is also useful in the
characterization of equivalence between categories of comodules over
corings. Indeed, in Section \ref{equivalences}, we study
equivalences between categories of comodules over rather general
corings. We generalize and improve (using results we give in Chapter
2) the main results concerning equivalences between categories of
comodules given in \cite{Takeuchi:1977}, \cite{AlTakhman:2002} and
\cite{Brzezinski/Wisbauer:2003}. We also give new characterizations
of equivalences between categories of comodules over coseparable
corings and corings over QF rings. We apply our results to the
particular case of the adjoint pair of functors associated to a
morphism of corings over different base rings. Finally, when
applied to corings associated to entwining structures and that
associated to a $G$-graded algebra and a right $G$-set, we obtain
new results concerning entwined modules and graded modules. We
think that our result, Theorem \ref{graded2}, is more simple than
Del R\'{\i}o's \cite[Theorem 2.3]{DelRio:1992}.

\medskip
In \cite{Rosenberg/Zelinsky:1961} the Picard group of an Azumaya
algebra is defined and the connections between this group and the
algebra automorphisms are studied. H. Bass generalized in
\cite{Bass:1967,Bass:1968} the Picard group and these connections
for arbitrary algebras. Further study of the Picard group of
algebras is given in \cite{Frohlich:1973}. Morita theory for rings
with local units was developed in \cite{Abrams:1983} or
\cite{Anh/Marki:1987}. In \cite{Beattie/Delrio:1996} the authors
introduced the Picard group of a ring with local units and gave the
versions of the corresponding connections for this ring, in order to
study the Picard group of the category $R-gr$, where $R$ is a
$G$-graded ring. In \cite{Torrecillas/Zhang:1996}, the authors
defined the Picard group of a coalgebra over a field and studied the
corresponding connections for it.

The purpose of Chapter 4 is to introduce and study the right Picard
group of corings. The motivation is the fact that there is an
isomorphism of groups between the Picard group of the category
$\rcomod{\coring{C}}$, for a certain coring $\coring{C}$ and the
right Picard group of $\coring{C}$, which is defined as the group of
the isomorphism classes of right invertible $\coring{C}$-bicomodules
(=$T(\coring{C})$ for some $k$-autoequivalence of
$\rcomod{\coring{C}}$, $T$) with the composition law induced by the
cotensor product (see Proposition \ref{picard-cat-crg}). It is clear
that the Picard group of corings generalizes that of algebras and
coalgebras over fields. We recall two theorems, the first is
\cite[Proposition II (5.2)(3)]{Bass:1968} (see also \cite[Theorems
55.9, 55.11]{Curtis/Reiner:1987}) and the second is the particular
case of \cite[Theorem 2.7]{Torrecillas/Zhang:1996} where $R=k$.

\begin{theorem}\label{Picard1}
For a $k$-algebra $A$, there is an exact sequence
$$\xymatrix{1\ar[r] & \inn{A}\ar[r]& \aut{A}\ar[r] & \pic{A}},$$
where $\inn{A}:=\{a\mapsto bab^{-1}\mid b \;\textrm{invertible in}
\;A\},$ the \emph{group of inner automorphisms} of $A$.
\end{theorem}

\begin{theorem}\label{Picard2}
For a coalgebra $C$ over a field $k$, there is an exact sequence
$$\xymatrix{1\ar[r] & \inn{C}\ar[r]& \aut{C}\ar[r] &
\pic{C}},$$ where $\inn{C}$ (the \emph{group of inner automorphisms}
of $C$) is the set of $\varphi\in \aut{C}$ such that there is $p\in
C^*$ invertible with $\varphi(c)=\sum
p(c_{(1)})c_{(2)}p^{-1}(c_{(3)}),$ for every $c\in C.$
\end{theorem}

Our Theorem \ref{coringexseq} generalizes both Theorem \ref{Picard1}
and Theorem \ref{Picard2}. We extend the Aut-Pic property from
algebras \cite{Bolla:1984} and coalgebras over fields
\cite{Cuadra/Garcia/Torrecillas:2000} to corings, namely, a coring
$\coring{C}$, which is flat on both sides over its base ring $A$,
has the right Aut-Pic property, if the morphism $\aut{\coring{C}}\to
\rpic{\coring{C}}$ is an epimorphism. We also extend a result which
is useful to show that a given coring $\coring{C}$ has the right
Aut-Pic property (Proposition 2.7). Of course, all of the examples
of algebras and coalgebras having the Aut-Pic property that are
given in \cite{Bolla:1984,Cuadra/Garcia/Torrecillas:2000} are
corings having this property. In this chapter, we give some new
examples of corings having the Aut-Pic property. We also simplify
the computation of the right Picard group of several interesting
corings (see Proposition \ref{examples}). Finally, in Section 4, we
give the corresponding exact sequences for the category of entwined
modules over an entwining structure, the category of
Doi-Koppinen-Hopf modules over a Doi-Koppinen structure and the
category of graded modules by a $G$-set, where $G$ is a group. Our
result, regarding the Picard group of the category $gr-(A,X,G)$, is
more natural than the Beattie-Del R\'io result given in \cite[\S
2]{Beattie/Delrio:1996}.

\medskip
Now, we move to a more detailed outline of the contents of this
work.

\medskip
\textbf{Chapter 1}. The aim of this chapter is to collect some
results that will be useful in the next chapters. Several of them
are new. In Section 1, we state certain results of category and ring
theory. Sometimes, we also give some proofs in order to give the
reader clear tools that are necessary to understand the rest of this
work. Section 2 is devoted to the study of the category of comodules
over a coring. In Section 3, we study the cotensor product over
corings. Section 4 deals with the cohom functor and the
coendomorphism coring. We give complete and rigorous proofs of the
most results in Section 2 and of all the stated results in Sections
3 and 4. Finally, in Section 5, we recall the definitions and some
results regarding induction functors, entwining structures and
graded rings.

For an introduction to the theory of abelian and Grothendieck
categories, we refer to the very famous and monumental paper by A.
Grothendieck (known as the ``Tohoku paper'')
\cite{Grothendieck:1957}. For a complete treatment, we refer to
\cite{Gabriel:1962,Freyd:1964,Mitchell:1965,Popescu:1973,MacLane:1998,
Stenstrom:1975} and the very recent book by M. Kashiwara and P.
Schapira \cite{Kashiwara/Schapira:2006}. For the foundation of the
theory of categories, using Grothendieck's universes, see
\cite{Grothendieck/Verdier:1963/64}. For a detailed discussion
of corings, we refer to \cite{Brzezinski/Wisbauer:2003}. We refer to
\cite{Brzezinski/Wisbauer:2003,Caenepeel/Militaru/Zhu:2002} for a
study of entwining structures and to
\cite{Nastasescu/VanOystaeyen:1982,Nastasescu/Raianu/VanOystaeyen:1990,
Nastasescu/VanOystaeyen:2004} for graded ring theory.

\medskip
\textbf{Chapter 2}. Section \ref{Frobeniusgeneral} deals with
adjoint and Frobenius pairs on categories of comodules. Some
refinements of results from \cite{Gomez:2002} and \cite[\S
23]{Brzezinski/Wisbauer:2003} on the representation as cotensor
product functors of certain functors between comodule categories are
needed and are, thus, included in Section \ref{Frobeniusgeneral}.
From our general discussion on adjoint pairs of cotensor product
functors, we will derive our main general result on Frobenius pairs
between comodule categories (Theorem \ref{20a}) that extends the
known characterizations in the setting of modules over rings and of
comodules over coalgebras. In the first case, the key property to
derive the result on modules from Theorem \ref{20a} is the
separability of the trivial corings (see Remark \ref{damodulos}). In
the case of coalgebras, the fundamental additional property is the
duality between finite left and right comodules. We already consider
a much more general situation in Section \ref{condualidad}, where we
introduce the class of so called \emph{corings having a duality} for
which we prove characterizations of Frobenius pairs that are similar
to the coalgebra case.

Section \ref{Frobmor} is one of the principal motivations of this
study. After the technical development of sections
\ref{Frobeniusgeneral} and \ref{condualidad}, our main results
follow without difficulty. We prove, in particular, that the
induction functor $-\otimes_{A}B:\rcomod{\coring{C}}
\to\rcomod{\coring{D}}$ associated to a homomorphism
$(\varphi,\rho):\coring{C}\to\coring{D}$ of corings $\coring{C}$ and
$\coring{D}$ flat over their respective base rings $A$ and $B$ is
Frobenius if and only if the $(\coring{C} ,\coring{D})$-bicomodule
$\coring{C}\otimes_{A}B$ is quasi-finite and injector as a right
$\coring{D}$-comodule and there exists an isomorphism of
$(\coring{D},\coring{C})$-bicomodules
$\cohom{\coring{D}}{\coring{C}\otimes_{A}B}{\coring{D}}\simeq
B\otimes _{A}\coring{C}$ (Theorem \ref{23}). We show as well how
this theorem unifies previous results for ring homomorphisms
\cite[Theorem 5.1]{Morita:1965}, coalgebra maps \cite[Theorem
3.5]{Castano/Gomez/Nastasescu:1999} and Frobenius corings
\cite[27.10, 28.8]{Brzezinski/Wisbauer:2003}.

In Section \ref{entwined}, we specialize one of the general results
on corings to entwining structures.

In Section \ref{graded}, we particularize our results in the
previous sections to the coring associated to a $G$-graded algebra
and a $G$-set, where $G$ ia a group. Then, we obtain a series of new
results for graded modules by $G$-sets.

\medskip
\textbf{Chapter 3.} In Section \ref{comatrixgeneralized}, we
generalize the comatrix coring introduced in
\cite{ElKaoutit/Gomez:2003} and \cite{Brzezinski/Gomez:2003}. We
also generalize some of its interesting properties given in
\cite{Brzezinski/Gomez:2003}. Section \ref{equivalences} is devoted
to the study of equivalences of comodule categories over corings.
Our results given in \cite{Zarouali:2004} and Section
\ref{comatrixgeneralized} will allow us to generalize and improve
the main results in both \cite{AlTakhman:2002} and
\cite{Brzezinski/Wisbauer:2003} (see Propositions
\ref{equivalence1}, \ref{equivalence2}, Theorems \ref{equivalence3},
\ref{equivalence4}) and also to give new results concerning
equivalences of comodule categories over coseparable corings (see
Propositions \ref{equivalence1}, \ref{equivalence2}, Theorems
\ref{equivalence3}, \ref{equivalence5}) and over corings over QF
rings (Theorem \ref{equivalence6}). Obviously, our last theorem
generalizes \cite[Theorem 3.5]{Takeuchi:1977} and \cite[Corollary
7.6]{AlTakhman:2002}. In Section \ref{induction}, we deal with the
application of some of our results given in Section
\ref{equivalences} to the induction functor. In Section
\ref{entwined}, we apply some of our results given in the previous
sections to the corings associated to entwining structures, in
particular, those associated to a $G$-graded algebra and a right
$G$-set, where $G$ is a group.

\medskip
\textbf{Chapter 4.} The aim of Section \ref{Picard} is to introduce
and study the Picard group of corings. Section \ref{Aut-Pic} is
devoted to define the Aut-Pic property for corings and give some new
examples of corings having this property. In the final section
\ref{Picard graded}, we apply our result, Theorem \ref{coringexseq},
to the corings coming from entwining structures, specially from
graded ring theory.

\medskip
\section*{Acknowledgements}

I should be deeply grateful and thankful to my advisor, Dr Jos\'e
G\'omez-Torrecillas, for his assistance, time and advice.

\medskip
I would also like to express my deep gratitude to (in alphabetical
order):\begin{itemize}

\item Dr Tomasz Brzezi\'nski for his useful comments on Section 6 of
\cite{Zarouali:2004} and his recommendation letter for a Ph.D.
studentship which contains a goodwill to help me and an opinion on
my work on coring theory. I am very proud of this opinion.

\item Dr Manuel Bullejos Lorenzo for his assistance.

\item Dr Edgar Enochs for communicating the example of a commutative
self-injective ring which is not coherent given in Subsection
\ref{category theory}.

\item The members of the jury who agreed to attend to my defense of
thesis.

\item Dr Blas Torrecillas Jover for the helpful discussion on the
subject of \cite{Zarouali:2006}.

\item Dr Robert Wisbauer and the referee for some helpful
suggestions on \cite{Zarouali:2004}.
\end{itemize}

I do not have sufficient words to thank my brother Youssef for his
encouragement and financial support during my stay in Granada for
more than four years. I owe my carrier to him. It is only fair to
dedicate him this modest work.

\medskip
Special thanks go to my mother and sisters for their support.

\medskip
I thank my friend Youssef and my friends of Granada for their help,
as well as my friend Majid for some questions of English language.

\medskip
Finally, I would like to thank the people who are not mentioned
 above and who contributed to the realization of this work.

\newpage

\section*{Some general notations}

\begin{enumerate}[(i)]
\item Throughout this work, unless otherwise stated, $k$ denotes an
associative and commutative ring with unit, $A,$ $A',$ $A'',$ $A_1,$
$A_2,$ $B,$ $B_1$ and $B_2$ denote associative and unitary algebras
over $k$ and $\coring{C},$ $\coring{C}',$ $\coring{C}'',$
$\coring{C}_1,$ $\coring{C}_2,$ $\coring{D},$ $\coring{D}_1$ and
$\coring{D}_2$ denote corings over $A,$ $A',$ $A'',$ $A_1,$ $A_2,$
$B,$ $B_1$ and $B_2$, respectively.
\item $\textbf{Set}$ will denote the category of sets, $\textbf{Ab}$
will denote the category of abelian groups, and $\rmod{A}$,
$\lmod{A}$ and $\bmod{A}{B}$ will denote the categories of right
$A$-modules, left $A$-modules and $(A,B)$-modules, respectively.
\item The notation $\otimes$ will stand for the tensor product over
$k$.
\item We denote the dual algebra of an algebra $A$ by $A^\circ$. We
denote the dual category of a category $\cat{C}$ by $\cat{C}^\circ$.
\end{enumerate}

\chapter{Generalities}

\section{Some category and ring theory}\label{category theory}

We refer to Grothendieck's paper \cite{Grothendieck:1957} for
the definitions of the basic notions of abelian categories. We
recall only Grothendieck's following axioms for a category
$\cat{C}$:

\begin{enumerate}[{AB} 1)]
\item Every morphism in $\cat{C}$ has a kernel and cokernel.
\item For every morphism $u:A\to B$ in $\cat{C}$, the induced
morphism $\overline u:\operatorname{Coim}(u)\to
\operatorname{Im}(u)$ defined to be the unique morphism making
commutative the diagram
$$\xymatrix{\operatorname{Ker}(u)\ar[r] & A\ar[r]^u \ar[d] & B\ar[r] &
\operatorname{Coker}(u)\\
& \operatorname{Coim}(u)\ar[r]^{\overline u} &
\operatorname{Im}(u)\ar[u] &},$$ is an isomorphism.
\item The direct sum of every family of objects $(A_i)_{i\in I}$ of
$\cat{C}$ exists.
\item The axiom AB 3) holds, and the direct sum of a family of
monomorphisms is a monomorphism.
\item The axiom AB 3) holds, and if $(A_i)_{i\in I}$ is a directed
family of subobjects of an object $A\in\cat{C}$, and $B$ is an
arbitrary subobject of $A$, then
$$(\sum_iA_i)\cap B=\sum_i(A_i\cap B).$$
\item The axiom AB 3) holds, and for every object $A\in\cat{C}$ and
every family $(B^j)_{j\in J}$, where $B^j=(B^j_i)_{i\in I_j}$ is a
directed family of subobjects of $A$ for every $j\in J$,
$$ \bigcap_{j\in J}\Big(\sum_{i\in I_j}B^j_i\Big)=
\sum_{(i_j)\in\prod I_j}\Big(\bigcap_{j\in J}B^j_{i_j}\Big).$$ (This
axiom includes implicitly the existence of the greatest lower bound
of every family of subobjects of $A$.)
\end{enumerate}

Obviously the axioms AB 1) and AB 2) coincide with their duals. We
define also the dual axioms AB 3*), AB 4*), AB 5*), and AB 6*).

A category is called \emph{additive} if it is preadditive, has a
zero object and finite products (or coproducts). A category is
called \emph{abelian} if it is additive and satisfies AB 1) and
AB~2).

\medskip
\textbf{$k$-categories, $k$-functors.} A category $\cat{C}$ is said
to be $k$-\emph{category} if for every $C,C'\in\cat{C}$,
$\hom{\cat{C}}{C}{C'}$ is a $k$-module, and the composition is
$k$-bilinear. In particular, for every $C\in\cat{C}$,
$\operatorname{End}_\cat{C}(C)$ is a $k$-algebra. For a $k$-algebra
$A$, $\rmod{A}$ is a $k$-category.

A functor between $k$-categories is said to be $k$-\emph{functor} or
$k$-\emph{linear functor} if it is $k$-linear on the $k$-modules of
morphisms.

Of course, if $k=\mathbb{Z}$, we find the very well known notions of
a preadditive category and additive functor.

A $k$-category is called a $k$-\emph{abelian category} if it is an
abelian category.

\medskip
Let $\cat{C}$ be a $k$-category. Following Mac Lane (\cite[p.
194]{MacLane:1998}), a \emph{biproduct diagram} for the objects
$A,B\in \cat{C}$ is a diagram
\begin{equation}\label{biproduct}
\xymatrix{A\ar@<2pt>[r]^{i_1} &
C\ar@<2pt>[l]^{p_1}\ar@<2pt>[r]^{p_2}& B\ar@<2pt>[l]^{i_2}}
\end{equation}
such that the morphisms $p_1,p_2,i_1,i_2$ satisfy the conditions
\begin{equation}\label{biproduct-cond}
p_1i_1=1_A,\quad p_2i_2=1_B,\quad i_1p_1+i_2p_2=1_C.
\end{equation}
In particular, $p_1,p_2$ are epimorphisms, $i_1,i_2$ are
monomorphisms, $p_2i_1=0$ and $p_1i_2=0$.

If such biproduct exists, it is both a product and coproduct.

\begin{proposition}\label{$k$-functor}
Let $\cat{C}$ and $\cat{D}$ be two $k$-categories such that all two
objects have a biproduct diagram. Then a functor
$T:\cat{C}\to\cat{D}$ is $k$-linear if and only if the following
hold
\begin{enumerate}[(a)]
\item $T$ carries every biproduct diagram in $\cat{C}$ to a biproduct
diagram in $\cat{D}$,
\item $T(h1_A)=h1_{T(A)}$ for every $A\in \cat{C}$ and $h\in k$.
\end{enumerate}
\end{proposition}

\begin{proof}
First, by Grothendieck's \cite[Proposition VIII.4]{MacLane:1998},
$T$ is additive if and only if the condition (a) holds. On the other
hand, for every morphism $f:A\to B$ in $\cat{C}$ and $h\in k$,
$hf=(h1_A)f$. Hence the claimed result follows.
\end{proof}

A functor $T:\cat{C}\to\cat{D}$ between abelian categories is said
to be \emph{middle-exact} if for every short exact sequence $0\to
A\to B\to C\to 0$, the sequence $T(A)\to T(B)\to T(C)$ is exact.
Obviously every left (right) exact functor is middle-exact.

\begin{proposition}\label{middle-exact}
Every middle-exact functor is additive.
\end{proposition}

\begin{proof}
Let $A,B\in \cat{C}$. Let $p_1,p_2,i_1,i_2$ be morphisms satisfying
\eqref{biproduct}. It is easy to verify that \eqref{biproduct-cond}
is equivalent to
\begin{enumerate}[(a)]
\item $p_1i_1=1_A$ and $p_2i_2=1_B$, and
\item the sequences $\xymatrix{0\ar[r] & A\ar[r]^{i_1} &
C\ar[r]^{p_2} & B\ar[r] & 0}$ and $\xymatrix{0\ar[r] & B\ar[r]^{i_2}
& C\ar[r]^{p_1} & A\ar[r] & 0}$ are exact.
\end{enumerate}
Finally, the use of Proposition \ref{$k$-functor} achieves the
proof.
\end{proof}

\textbf{Abelian subcategory.} Let $\cat{C}$ be an abelian category.
A nonempty full subcategory of $\cat{C}$, $\cat{B}$, is called
\emph{an abelian subcategory} if $\cat{B}$ is abelian and the
injection functor $ \cat{B}\to \cat{C}$ is exact. The proof of the
following lemma is straightforward.

\begin{lemma}\label{abeliansubcategory}
Let $\cat{C}$ be an abelian category, and $\cat{B}$ be a nonempty
full subcategory of $\cat{C}$. Then
\begin{enumerate}[(1)]
\item $\cat{B}$ is additive if and only if $0\in \cat{B}$ and for
every $B_1,B_2\in \cat{B}$, $B_1\oplus B_2\in \cat{B}$.
\item $\cat{B}$ is an abelian subcategory of $\cat{C}$ if and only
if $\cat{B}$ is additive, and for every morphism $f:B\to B'$ in
$\cat{B}$, $\oper{Ker}{f}, \oper{Coker}{f}\in \cat{B}$
\emph{(\cite[Theorem 3.41]{Freyd:1964})}.
\end{enumerate}
\end{lemma}

\medskip
\textbf{Diagram categories.} We propose here a detailed exposition
of diagram categories.  This paragraph is essentially due to
Grothendieck, see \cite{Grothendieck:1957}.

A \emph{diagram scheme} is a triple $S=(I,\Phi,d)$ consisting of two
sets $I$ and $\Phi$ and a map $d:\Phi\to I\times I$. The elements of
$I$ are called the \emph{vertices}, the elements of $\Phi$ are
called the \emph{arrows} of the scheme. If $\varphi$ is an arrow of
the diagram, $d(\varphi)=(i,j)$ is called its \emph{direction},
characterized by the origin and the extremity of the arrow: $i$ and
$j$. In this situation we write $\varphi:i\to j$.

A \emph{composite arrow} $\varphi$ of origin $i$ and extremity $j$
is defined to be a nonempty finite sequence of arrows
$\varphi=(\varphi_1,\dots,\varphi_n)$ such that the origin of
$\varphi_1$ is $i$ and the extremity of $\varphi_n$ is $j$. We say
that $n$ is the \emph{length} of $\varphi$.

Let $\cat{C}$ be a category and $S$ a diagram scheme. A
\emph{diagram in $\cat{C}$ of scheme} $S$ is a map $D$ which assigns
to every $i\in I$ an object $D(i)\in\cat{C}$ and assigns to every
arrow $\varphi:i\to j$, a morphism $D(\varphi):D(i)\to D(j)$ in
$\cat{C}$. A morphism $\upsilon$ from a diagram $D$ to another $D'$
is a family of morphisms $\upsilon_i:D(i)\to D'(i)$ such that for
every arrow $\varphi:i\to j$, the diagram in $\cat{C}$
$$\xymatrix{D(i)\ar[r]^{\upsilon_i}\ar[d]_{D(\varphi)} &
D'(i) \ar[d]^{D'(\varphi)} \\
D(j)\ar[r]^{\upsilon_j} & D'(j)}$$ is commutative. Let
$\upsilon:D\to D'$ and $\upsilon': D'\to D''$ be morphisms of
diagrams. We define two morphisms of diagrams
$\upsilon'\upsilon:D\to D''$ and $1_D:D\to D$ by
$(\upsilon'\upsilon)_i=\upsilon'_i\upsilon_i$ and $(1_D)_i=
1_{D(i)}$ for every $i\in I$.

The diagrams in $\cat{C}$ of scheme $S$ with their morphisms is a
category which we denote by $\cat{C}^S$.

A diagram scheme $S=(I,\Phi,d)$ and a diagram $D\in \cat{C}^S$ are
called \emph{finite} (resp. \emph{discrete}) if $I$ is a finite set
(resp. $\Phi=\emptyset$).

Now let $D$ be a diagram $\cat{C}$ of scheme $S$ and $\varphi=
(\varphi_1,\dots,\varphi_n)$ be a composite arrow of origin $i$ and
extremity $j$. Define the morphism $D(\varphi)=D(\varphi_n)\dots
D(\varphi_1):D(i)\to D(j)$ in $\cat{C}$. We say that $D$ is a
\emph{commutative diagram} if $D(\varphi)=D(\varphi')$ for every
composite arrows $\varphi$ and $\varphi'$ which have the same origin
and the same extremity.

\begin{lemma}\label{cube} \emph{[The Cube Lemma]}\\
Consider the following diagram in a fixed category.
$$
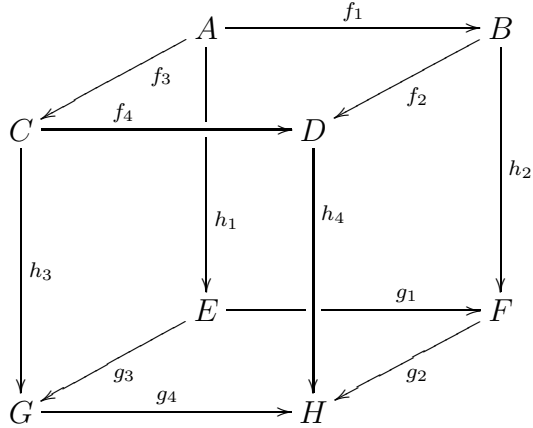
\xymatrix{
 && A \ar[rrr]^{f_1}\ar[dll]^(.3){f_3}\ar'[d][ddd]^{h_1}&&& B
\ar[dll]^{f_2}\ar[ddd]^{h_2}
\\
C \ar[rrr]^(.35){f_4}\ar[ddd]^{h_3} &&& D \ar[ddd]^(.3){h_4}
\\
&&&&&&&&&
\\
 && E \ar'[r][rrr]^{g_1} \ar[dll]^{g_3} &&& F\ar[dll]^{g_2}
\\
G \ar[rrr]^{g_4}  &&& H }.$$ If $h_4$ is a monomorphism, and all
squares, save possibly the top one, are commutative, then the cube
is also commutative.
\end{lemma}

\begin{proof}
We have
 $$h_4f_2f_1=g_2h_2f_1=g_2g_1h_1=g_4g_3h_1=g_4h_3f_3=h_4f_4f_3.$$
Since $h_4$ is a monomorphism, we obtain $f_2f_1=f_4f_3$.
\end{proof}

Let $S=(I,\Phi,d)$ be a diagram scheme. We define the \emph{dual
scheme} of $S$ to be $S^\circ=(I,\Phi,d^\circ)$ where
$d^\circ:\Phi\to I\times I$ is the map defined by
$d^\circ(\varphi)=(j,i)$ for every $\varphi\in \Phi$ such that
$d(\varphi)=(i,j)$. Consider the contravariant functor
\begin{equation}\label{diagram functor}
(-)^\circ:\cat{C}^S\to (\cat{C}^\circ)^{S^\circ}
\end{equation}
 which assigns to a diagram $D\in \cat{C}^S$ the diagram
 $D^\circ\in (\cat{C}^\circ)^{S^\circ}$ defined by $D^\circ(i)=D(i)$
 and $D^\circ(\varphi^\circ)=D(\varphi)^\circ$ for every $i\in I$
 and $\varphi\in \Phi$. It also assigns to a morphism of diagrams
 $\upsilon:D\to D'$ the morphism of diagrams $\upsilon^\circ:
 D'^\circ\to D^\circ$ such that $(\upsilon^\circ)_i=
 (\upsilon_i)^\circ$ for every $i\in I$. This functor leads to an
 isomorphism of categories
\begin{equation}\label{diagram duality}
(\cat{C}^S)^\circ \simeq (\cat{C}^\circ)^{S^\circ}.
\end{equation}

An object $A$ of a category $\cat{C}$ is called an \emph{initial
object} if $\hom{\cat{C}}{A}{B}$ is a singleton for every
$B\in\cat{C}$. An initial object is determined up to an isomorphism.
Dually we define a final object. An object is called a \emph{zero
object} if it is at the same time an initial and final object. In a
preadditive category, an object is an initial object if and only it
is a final object, if and only if it is a zero object.

Now we will see that the category $\cat{C}^S$ inherits interesting
properties from the category $\cat{C}$.

\begin{lemma}\label{DC1}
Let $\cat{C}$ be a category and $S=(I,\Phi,d)$ a diagram scheme.
Then
\begin{enumerate}[(i)]
\item If $\cat{C}$ is a $k$-category then the
same holds for $\cat{C}^S$.
\item If $\cat{C}$ has an initial (resp. a final, resp. a zero)
object then the same holds for $\cat{C}^S$.
\item If $\cat{C}$ has products (finite products)
(resp. direct sums (finite direct sums)) then the same holds for
$\cat{C}^S$.
\item If $\cat{C}$ has kernels (resp. cokernels) then the
same holds for $\cat{C}^S$.
\item Let $\cat{C}$ be a $k$-abelian category. Then so is
$\cat{C}^S$. A subobject of a diagram $D\in \cat{C}^S$ is a family
of objects $(D'(i))_{i\in I}$ such that $D'(i)$ is a subobject of
$D(i)$ for every $i\in I$, and for every arrow $\varphi:i\to j$ of
$S$, $D(\varphi).D'(i)\subset D'(j)$ as subobjects of $D(j)$. The
morphism $D'(\varphi)$ is then defined to be the unique morphism
making commutative the diagram
$$\xymatrix{0\ar[r] & D'(i)\ar[r]\ar[d]_{D'(\varphi)} &
D(i) \ar[d]^{D(\varphi)} \\
0\ar[r] & D'(j)\ar[r] & D(j).}$$ If moreover $\cat{C}$ satisfy one
of the statements AB 4), AB 5), AB 6), or one of their dual
statements AB 4*), AB 5*), AB 6*), then the same holds for
$\cat{C}^S$.
\end{enumerate}
\end{lemma}

\begin{proof}
(i) Let $\cat{C}$ be a $k$-category and $\upsilon,\omega:D\to D'$
two morphisms of diagrams in $\cat{C}^S$. Set
$$(\upsilon+\omega)_i=\upsilon_i+\omega_i,\;(h\upsilon)_i=h\upsilon_i$$
for all $i\in I$ and $h\in k$. It is easy to verify that
$\hom{\cat{C}^S}{D}{D'}$ is a $k$-module with the morphism of
diagrams $0_i=0:D(i)\to D'(i)$, $i\in I$, as neutral element for the
addition, and the composition is $k$-bilinear.

(ii) Let $C$ be an initial object in $\cat{C}$. Let $D\in\cat{C}^S$
be the diagram defined by: $D(i)=C$ for every $i\in I$, and
$D(\varphi)$ be the unique morphism $C\to C$ in $\cat{C}$ for every
$\varphi\in \Phi$. Obviously $D$ is an initial object in
$\cat{C}^S$. Analogously we prove the second statement. The third
one is obvious from the first and the second ones.

(iii) Suppose that $\cat{C}$ has products. Let
$(D_\lambda)_{\lambda\in \Lambda}$ be a nonempty family of diagrams
in $\cat{C}^S$. For every $i\in I$ let $(\prod_{\lambda\in
\Lambda}D_\lambda(i),(p_{\lambda,i}))$ be a product of
$D_\lambda(i)$, $\lambda\in \Lambda$, where
$$p_{\lambda,i}:\prod_{\mu\in \Lambda}D_\mu(i)\to D_\lambda(i),$$
$\lambda\in \Lambda$, are the canonical projections. Let
$\varphi:i\to j$ be an arrow of $S$. Consider the morphism
$$\prod_{\lambda\in \Lambda}D_\lambda(\varphi):\prod_\lambda
D_\lambda(i)\to\prod_\lambda D_\lambda(j)$$ in $\cat{C}$, i.e. the
unique morphism making commutative the diagram
$$\xymatrix{\prod_\mu D_\mu(i)\ar[rr]^{p_{\lambda,i}}
\ar[d]_{\prod_\mu D_\mu(\varphi)} &&
D_\lambda(i) \ar[d]^{D_\lambda(\varphi)} \\
\prod_\mu D_\mu(j)\ar[rr]^{p_{\lambda,j}} && D_\lambda(j)}$$ for
every $\lambda\in\Lambda$. Then we obtain the diagram defined by
$$(\prod_{\lambda\in \Lambda}D_\lambda)(i)=\prod_{\lambda\in
\Lambda}D_\lambda(i),\;(\prod_{\lambda\in
\Lambda}D_\lambda)(\varphi)=\prod_{\lambda\in
\Lambda}D_\lambda(\varphi),$$ and for every $\lambda\in\Lambda$,
$p_\lambda=(p_{\lambda,i})_{i\in I}: \prod_\mu D_\mu\to D_\lambda$
is a morphism of diagrams.

Finally, let $\upsilon_\lambda:D\to D_\lambda$, $\lambda\in
\Lambda$, be a family of morphisms of diagrams. We should prove that
there is a unique morphism of diagrams $\upsilon:D\to \prod_\mu
D_\mu$ making commutative the diagram
$$\xymatrix{D\ar[rr]^{\upsilon_\lambda} \ar[drr]^\upsilon && D_\lambda
\\ && \prod_\mu D_\mu \ar[u]_{p_\lambda}}$$
for every $\lambda\in\Lambda$.

For this let $i\in I$. Then there is a unique morphism
$\upsilon_i:D(i)\to \prod_\mu D_\mu(i)$ in $\cat{C}$ making
commutative the diagram
$$\xymatrix{D(i)\ar[rr]^{(\upsilon_\lambda)_i} \ar[drr]^{\upsilon_i}
&& D_\lambda(i)
\\ && \prod_\mu D_\mu(i) \ar[u]_{p_{\lambda,i}}}$$ for every
$\lambda\in\Lambda$. It remains to verify that $\upsilon$ is a
morphism of diagrams. Let $\varphi:i\to j$ be an arrow of $S$ and
$\lambda\in\Lambda$. Consider the diagram
$$\xymatrix{D(i) \ar[rrr]^{D(\varphi)} \ar[dd]^{\upsilon_i}
\ar[dr]^{(\upsilon_\lambda)_i} &&&
D(j)\ar'[d][dd]^{\upsilon_j}\ar[dr]^{(\upsilon_\lambda)_j} &&\\
& D_\lambda(i) \ar[rrr]^(.3){D_\lambda(\varphi)} &&&
D_\lambda(j)\\
\prod_\mu D_\mu(i)\ar[rrr]^{\prod_\mu D_\mu(\varphi)}
\ar[ur]^{p_{\lambda,i}} &&& \prod_\mu
D_\mu(j)\ar[ur]^{p_{\lambda,j}}&& .}$$ It follows that
$$p_{\lambda,j}\upsilon_j D(\varphi)=
p_{\lambda,j}(\prod_\mu D_\mu(\varphi))\upsilon_i$$ for every
$\lambda\in\Lambda$. Hence $\upsilon_j D(\varphi)= (\prod_\mu
D_\mu(\varphi))\upsilon_i.$

Analogously we prove the second statement.

Now suppose that $\cat{C}$ has direct sums. Let
$(D_\lambda)_{\lambda\in \Lambda}$ be a nonempty family of diagrams
in $\cat{C}^S$.  For every $i\in I$ let $(\oplus_{\lambda\in
\Lambda}D_\lambda(i),(u_{\lambda,i}))$ be a direct sum of
$D_\lambda(i)$, $\lambda\in \Lambda$, where
$$u_{\lambda,i}:D_\lambda(i)\to \oplus_{\mu\in \Lambda}D_\mu(i),$$
$\lambda\in \Lambda$, are the canonical injections. Define as above
the diagram in $\cat{C}^S$
$$(\oplus_{\lambda\in \Lambda}D_\lambda)(i)=\oplus_{\lambda\in
\Lambda}D_\lambda(i),\;(\oplus_{\lambda\in
\Lambda}D_\lambda)(\varphi)=\oplus_{\lambda\in
\Lambda}D_\lambda(\varphi).$$ Moreover, we have that
$u_\lambda=(u_{\lambda,i})_i:D_\lambda\to \oplus_{\mu\in
\Lambda}D_\mu$ is a morphism in $\cat{C}^S$ for every
$\lambda\in\Lambda$.

Let $(-)^\circ:\cat{C}^S\to (\cat{C}^\circ)^{S^\circ}$ be the
contravariant functor \eqref{diagram functor}. We have
$$(\oplus_{\lambda\in \Lambda}D_\lambda)^\circ(i)=\oplus_{\lambda\in
\Lambda}D_\lambda(i)=\oplus_{\lambda\in
\Lambda}D_\lambda^\circ(i)=(\prod_{\lambda\in
\Lambda}D_\lambda^\circ)(i)$$ and
$$(\oplus_{\lambda\in
\Lambda}D_\lambda)^\circ(\varphi^\circ)=\big[(\oplus_{\lambda\in
\Lambda}D_\lambda)(\varphi)\big]^\circ=\big(\oplus_{\lambda\in
\Lambda}D_\lambda(\varphi)\big)^\circ=\prod_{\lambda\in
\Lambda}D_\lambda^\circ(\varphi^\circ)=(\prod_{\lambda\in
\Lambda}D_\lambda^\circ)(\varphi^\circ),$$ for every $i\in I$ and
$\varphi\in \Phi$. Hence
\begin{equation}
(\oplus_{\lambda\in \Lambda}D_\lambda)^\circ=\prod_{\lambda\in
\Lambda}D_\lambda^\circ.
\end{equation}
Finally, $(\oplus_{\lambda\in \Lambda}D_\lambda,(u_\lambda))$ is
then a direct sum of the family $(D_\lambda)_\lambda$ in
$\cat{C}^S$.

Analogously we prove the last statement.

(iv) Suppose that $\cat{C}$ is additive and has kernels. Let
$\upsilon:D\to D'$ be a morphism in $\cat{C}^S$. Set
$K(i)=\operatorname{Ker}(\upsilon_i)$ for all $i\in I$. Let
$\varphi:i\to j$ be an arrow of $S$. Then there is a unique morphism
$K(\varphi):K(i)\to K(j)$ in $\cat{C}$ making commutative the
diagram with exact rows
\begin{equation}\label{KDC1}
\xymatrix{
 0 \ar[r] & K(i) \ar[rr]^{\operatorname{ker}(\upsilon_i)}
 \ar[d]^{K(\varphi)}&& D(i)
\ar[d]^{D(\varphi)} \ar[rr]^{\upsilon_i} &&
D'(i)\ar[d]^{D'(\varphi)} \\
0 \ar[r] &  K(j)\ar[rr]^{\operatorname{ker}(\upsilon_i)}&&
D(j)\ar[rr]^{\upsilon_j} &&  D'(j).}
\end{equation}
We obtain that $K\in\cat{C}^S$ and
$(\operatorname{ker}(\upsilon_i))_{i\in I}:K\to D$ is a morphism in
$\cat{C}^S$.

Now, let $\xi:D''\to D$ be a morphism in $\cat{C}^S$ such that
$\upsilon\xi=0$. Then, for all $i\in I$, there is a unique morphism
$\omega_i:D''(i)\to k(i)$ in $\cat{C}$ making commutative the
diagram
$$\xymatrix{&& D''(i) \ar[dl]_{\omega_i}\ar[dr]^{\xi_i}&&
\\0 \ar[r] & K(i) \ar[rr]^{\operatorname{ker}(\upsilon_i)} && D(i)
\ar[r]^{\upsilon_i} & D'(i).}$$ To complete the proof of the first
statement it is enough to verify that $\omega=(\omega_i)_i:D''\to K$
is a morphism in $\cat{C}^S$. For this, let $\varphi:i\to j$ be an
arrow of $S$ and consider the diagram
$$\xymatrix{D''(i) \ar[rrr]^{D''(\varphi)} \ar[dd]^{\omega_i}
\ar[dr]^{\xi_i} &&&
D''(j)\ar'[d][dd]^{\omega_j}\ar[dr]^{\xi_j} &&\\
& D(i) \ar[rrr]^(.3){D(\varphi)} &&&
D(j)\\
K(i)\ar[rrr]^{K(\varphi)}\ar[ur]_{\operatorname{ker}(\upsilon_i)}
&&& K(j)\ar[ur]_{\operatorname{ker}(\upsilon_j)}&& .}$$ Since
$\operatorname{ker}(\upsilon_j)$ is a monomorphism in $\cat{C}$, the
rectangle is commutative. Hence $K$ is the kernel of $\upsilon$ in
$\cat{C}^S$.

Suppose that $\cat{C}$ is additive and has cokernels. Let
$\upsilon:D\to D'$ be a morphism in $\cat{C}^S$. Set
$K'(i)=\operatorname{Coker}(\upsilon_i)$ for all $i\in I$. Let
$\varphi:i\to j$ be an arrow of $S$. Then there is a unique morphism
$K'(\varphi):K'(i)\to K'(j)$ in $\cat{C}$ making commutative the
diagram with exact rows
\begin{equation*}
\xymatrix{
  D(i)\ar[r]^{\upsilon_i}\ar[d]^{D(\varphi)}&
  D'(i)\ar[rr]^{\operatorname{coker}(\upsilon_i)}
  \ar[d]^{D'(\varphi)}
&& K'(i)\ar[d]^{K'(\varphi)} \ar[r] & 0 \\
D(j)\ar[r]^{\upsilon_j} &
D(j)\ar[rr]^{\operatorname{coker}(\upsilon_j)}&& K'(j)\ar[r] & 0.}
\end{equation*}
Then $(K'^\circ,(\operatorname{coker}(\upsilon_i)^\circ)_i)$ is the
kernel of $\upsilon^\circ:D'^\circ\to D^\circ$ in
$(\cat{C}^\circ)^{S^\circ}$. Hence
$(K',(\operatorname{coker}(\upsilon_i))_i)$ is the cokernel of
$\upsilon:D\to D'$ in $\cat{C}^S$.

Finally, (v) is obvious from the statements (i) to (iv), their
proofs, and the fact that a morphism $\upsilon:D\to D'$ in
$\cat{C}^S$ is an isomorphism if and only if $\upsilon_i:D(i)\to
D'(i)$ is an isomorphism in $\cat{C}$ for every $i\in I$.
\end{proof}

\begin{remark}
Let $\cat{C}$ be a category and $S=(I,\Phi,d)$ a diagram scheme.
\begin{enumerate} [(a)]
\item Suppose that $\cat{C}$ has products. Let
$(D_\lambda)_{\lambda\in \Lambda}$ be a nonempty family of diagrams
in $\cat{C}^S$. For a composite arrow $\varphi$, $(\prod_\lambda
D_\lambda)(\varphi)=\prod_\lambda D_\lambda(\varphi).$
\item Suppose that $\cat{C}$ is additive and has kernels. Let
$\upsilon:D\to D'$ be a morphism in $\cat{C}^S$. Let $K$ be the
diagram in $\cat{C}^S$ defined as in the proof of Lemma \ref{DC1}
(iv). Let $\varphi$ be a composite arrow of $S$ of origin $i$ and
extremity $j$. Then $K(\varphi):K(i)\to K(j)$ is the unique morphism
making commutative the diagram \eqref{KDC1}.
\end{enumerate}
\end{remark}

Now we consider very interesting full subcategories of $\cat{C}^S$.
Let $S$ be a diagram scheme and $\cat{C}$ an arbitrary category.
\begin{enumerate}[(i)]\label{Sigma}
\item Let $R$ be a set of couples $(\varphi,\varphi')$ where $\varphi$
and $\varphi'$ are two composite arrows having the same origin and
the same extremity, and of composite arrows $\varphi$ whose origin
and extremity coincide. Consider $\cat{C}^{(S,R)}$ the full
subcategory of $\cat{C}^S$ consisting of the diagrams
$D\in\cat{C}^S$ such that $D(\varphi)=D(\varphi')$ for every
$(\varphi,\varphi')\in R$ and $D(\varphi)=1_{D(i)}$
 for every $\varphi\in R$ of origin and extremity $i$.
\item Now assume that $\cat{C}$ is a $k$-category. Let
$(i,j)\in I\times I$. Let $R_{ij}$ be a set of formal linear
combinations of composite arrows of origin $i$ and extremity $j$
with coefficients in $k$, and if $i=j$ we add some formal linear
combinations of composite arrows of origin and extremity $i$, and a
fixed element $e_i$, with coefficients in $k$. Let $D$ be a diagram
in $\cat{C}^S$. We define for every $L\in R_{ij}$, the morphism
$D(L):D(i)\to D(j)$ by replacing in the expression of $L$, a
composite arrow $\varphi$ by $D(\varphi)$, and $e_i$ by $1_{D(i)}$.
Set $R=\cup_{(i,j)\in I\times I}R_{ij}$.

A diagram $D$ in $\cat{C}^S$ is called $R$-\emph{commutative} if
 $D(L)=0$ for every $L\in R$. The couple $(S,R)=\Sigma$ is
 called a \emph{diagram scheme with commutation relations}. The
 equalities $L=0$, $L\in R$, are called the \emph{commutation
 relations}. Consider $\cat{C}^\Sigma$ the full subcategory of
 $\cat{C}^S$ consisting of the $R$-commutative diagrams.

 Let $\varphi=(\varphi_1,\dots,\varphi_n)$ be a composite arrow
 of $S$ of origin $i$ and extremity $j$. Define the composite
 arrow $\varphi^\circ=(\varphi_n^\circ,\dots,\varphi_1^\circ)$ of
 $S^\circ$ of origin $j$ and extremity $i$. Let $(-)^\circ:
 \cat{C}^S\to (\cat{C}^\circ)^{S^\circ}$ be the contravariant
 functor \eqref{diagram functor}. Then $D^\circ(\varphi^\circ)=
 D(\varphi)^\circ$. Let $(i,j)\in I\times I$ and $L\in R_{ij}$.
 Let $L^\circ$ be the linear combination obtained by replacing
 in the expression of $L$ each composite arrow $\varphi$ by
 $\varphi^\circ$, and if $i=j$, $e_i$ by a fixed element
 $e_i^\circ$. Let $R^\circ_{ji}$ be the set of $L^\circ$ where
 $L\in R_{ij}$. Set $R^\circ=\cup_{(i,j)\in I\times I}R^\circ_{ji}$.
 If $L\in R_{ij}$, then $D^\circ(L^\circ)=D(L)^\circ$. Therefore,
 $D$ is $R$-commutative if and only if $D^\circ$ is
 $R^\circ$-commutative. Hence we obtain an isomorphism of categories
\begin{equation}\label{diagram duality2}
(\cat{C}^\Sigma)^\circ \simeq (\cat{C}^\circ)^{\Sigma^\circ}.
\end{equation}
\end{enumerate}

\begin{proposition}\label{DC2}
Let $\cat{C}$ a $k$-category and $(S,R)=\Sigma$ a diagram scheme
with commutation relations. Then
\begin{enumerate}[(1)]
\item $\cat{C}^\Sigma$ is a $k$-category and the injection functor
$\cat{C}^\Sigma\to\cat{C}^S$ is $k$-linear. Moreover, if $\cat{C}$
has a zero object, then the zero object of $\cat{C}^S$ belongs to
$\cat{C}^\Sigma$.
\item Suppose that $\cat{C}$ has products (resp. finite products). Let
$(D_\lambda)_{\lambda\in \Lambda}$ be a nonempty family (resp.
finite family) of diagrams in $\cat{C}^\Sigma$. Then
$\prod_{\lambda\in \Lambda}D_\lambda$ belongs to $\cat{C}^\Sigma$.
The same holds for direct sums.
\item Suppose that $\cat{C}$ has kernels. Let $\upsilon:D\to D'$ be
a morphism in $\cat{C}^\Sigma$ and $K$ its kernel in $\cat{C}^S$.
Then $K$ belongs to $\cat{C}^\Sigma$. The same holds for cokernels.
\item If $\cat{C}$ is an abelian category, then $\cat{C}^\Sigma$ is
an abelian subcategory of $\cat{C}^S$.
\end{enumerate}
\end{proposition}

\begin{proof}
(1) Obvious.

(2) Let $(i,j)\in I\times I$ and $L=\sum_{l=1}^m h_l\varphi^l$,
where $\varphi^l$ is a composite arrow of origin $i$ and extremity
$j$ and $h_l\in k$ for every $l$. Then
\begin{eqnarray*}
(\prod_\lambda D_\lambda)(L)&=& \sum_{l=1}^m h_l(\prod_\lambda
D_\lambda)(\varphi^l)\\
&=& \sum_{l=1}^m h_l\prod_\lambda D_\lambda(\varphi^l)\\
&=& \prod_\lambda(\sum_{l=1}^m h_l D_\lambda(\varphi^l))\\
&=& \prod_\lambda D_\lambda(L).
\end{eqnarray*}

Hence, if $L\in R_{ij}$ then $(\prod_\lambda D_\lambda)(L)=0$.

On the other hand,
\begin{eqnarray*}
(\prod_\lambda D_\lambda)(e_i)&=& 1_{(\prod_\lambda
D_\lambda)(i)}\\
&=& 1_{\prod_\lambda D_\lambda(i)}\\
&=&\prod_\lambda1_{D_\lambda(i)}\\
&=& \prod_\lambda D_\lambda(e_i).
\end{eqnarray*}

Let $h\in k$. Then
\begin{eqnarray*}
(\prod_\lambda D_\lambda)(L+he_i)&=&(\prod_\lambda
D_\lambda)(L)+h(\prod_\lambda D_\lambda)(e_i)\\
&=& \prod_\lambda D_\lambda(L)+h\prod_\lambda D_\lambda(e_i)\\
&=& \prod_\lambda D_\lambda(L)+h D_\lambda(e_i)\\
&=& \prod_\lambda D_\lambda(L+he_i).
\end{eqnarray*}
Hence, if $L+he_i\in R_{ii}$ then $(\prod_\lambda
D_\lambda)(L+he_i)=0.$

Analogously we prove the second statement.

The last statement follows from the first one and the isomorphism of
categories \eqref{diagram duality2}.

(3) Let $(i,j)\in I\times I$ and $L=\sum_{l=1}^m h_l\varphi^l$,
where $\varphi^l$ is a composite arrow of origin $i$ and extremity
$j$ and $h_l\in k$ for every $l$. Then
$$D(L)\operatorname{ker}(\upsilon_i)=
\operatorname{ker}(\upsilon_j)K(L).$$ Hence, if $L\in R_{ij}$ then
$K(L)=0$.

Let $h\in k$. Then
$$D(L+he_i)\operatorname{ker}(\upsilon_i)=
\operatorname{ker}(\upsilon_j)K(L+he_i).$$ Hence, if $L+he_i\in
R_{ii}$ then $K(L+he_i)=0.$

The last statement follows from the first one and the isomorphism of
categories \eqref{diagram duality2}.

(4) Obvious using Lemma \ref{abeliansubcategory} and the statements
(1), (2) and (3).
\end{proof}

Now let $F:\cat{C}\to\cat{C'}$ be a functor and $S=(I,\Phi,d)$ a
diagram scheme. Define the functor $F^S:\cat{C}^S\to\cat{C'}^S$ such
that $$F^S(D)(i)=F(D(i)),\; F^S(D)(\varphi)=F(D(\varphi))$$ for
every $D\in\cat{C}^S$, $i\in I$, $\varphi\in\Phi$, and
$$F^S(\upsilon)_i=F(\upsilon_i):F(D(i))\to F(D'(i))$$ for every
morphism $\upsilon:D\to D'$ in $\cat{C}^S$ and $i\in I$.

If $F':\cat{C}\to\cat{C'}$ is another functor and $\eta:F\to F'$ is
a natural transformation, then $\eta^S:F^S\to F'^S$ defined by
$$(\eta^S_D)_i=\eta_{D(i)}:F(D(i))\to F'(D(i)),$$ for
every $D\in\cat{C}^S$ and $i\in I$, is a natural transformation.

\begin{lemma}
\begin{enumerate}[(a)]
\item $F\mapsto F^S,\eta\mapsto\eta^S$ acts as a functor.
\item Let $G:\cat{C'}\to\cat{C''}$ be another functor. Then
$(GF)^S=G^SF^S$.
\item If $\cat{C}$ and $\cat{C'}$ are abelian categories and $F$ is
left (resp. right) exact, then $F^S$ is so.
\end{enumerate}
\end{lemma}

\begin{proof}
Easy verifications.
\end{proof}

Now assume that $\cat{C}$ and $\cat{C'}$ are $k$-categories and let
$\Sigma=(S,R)$ be a diagram scheme with commutation relations and
$F:\cat{C}\to\cat{C'}$ a $k$-linear functor. The functor $F^S$ is
$k$-linear and induces a $k$-linear functor
$F^\Sigma:\cat{C}^\Sigma\to\cat{C'}^\Sigma$. The functors $F^S$ and
$F^\Sigma$ are called \emph{canonical prolongations of $F$ to
diagrams}. Sometimes and if there is no risk of confusion we denote
the functors $F^S$ and $F^\Sigma$ by $F$.

If $F':\cat{C}\to\cat{C'}$ is another $k$-linear functor and $\eta:
F\to F'$ is a natural transformation, then $\eta^S$ induces a natural
transformation $\eta^\Sigma:F^\Sigma\to F'^\Sigma$.

\begin{proposition}\label{can ext fun diag}
\begin{enumerate}[(a)]
\item $F\mapsto F^\Sigma,\eta\mapsto\eta^\Sigma$ acts as a functor.
\item Let $G:\cat{C'}\to\cat{C''}$ be another $k$-linear functor.
Then $(GF)^\Sigma=G^\Sigma F^\Sigma$.
\item If $\cat{C}$ and $\cat{C'}$ are abelian categories and $F$ is
left (resp. right) exact, then $F^\Sigma$ is so.
\end{enumerate}
\end{proposition}

\begin{proof}
Obvious from the last lemma.
\end{proof}

\begin{examples}\label{DC examples}
Let $\cat{C}$ be an arbitrary category.
\begin{enumerate}[(1)]
\item Assume that $\cat{C}$ is a $k$-category. Let $I$ be a
singleton and $\Phi=\emptyset$. The commutation relations are of the
form $he=0$. Let $h_1,\dots,h_n$ be elements of $k$. Consider the
commutation relations $h_le=0$, $l=1,\dots,n$. Then $\cat{C}^\Sigma$
is the full subcategory of $\cat{C}$ consisting of the objects
vanished by $h_l$ for all $l$.
\item Now we take $I$ is a set which has exactly two elements $i$ and
$j$, and $\Phi=\{\varphi\}$ is a singleton with $\varphi$ of origin
$i$ and extremity $j$. Then $\cat{C}^S$ is the category of morphisms
$u:A\to B$ between objects of $\cat{C}$. If we add the commutation
relation $h\varphi=0$ where $h\in k$, we restrict ourselves to $A$
and $B$. $u$ is vanished by $h$. If $h=1$, then
$\cat{C}^\Sigma\simeq \cat{C}\times\cat{C}$.
\item \textbf{Functor categories.} Let $\textbf{I}$ be a small
category. Let $S=(I,\Phi,d)$ be the diagram scheme such that $I$ and
$\Phi$ are respectively the sets of objects and morphisms of the
category $\cat{I}$. The directions of the arrows is defined by a
natural way. Let $R$ be the set of $(gf,(f,g))$ where $f$ and $g$
are morphisms in $\textbf{I}$ such that $(f,g)$ is a composite arrow
of length 2, and the identity morphisms. Then the category of
functors $\textbf{Fun}(\textbf{I},\cat{C})$ is the category
$\cat{C}^{(S,R)}$.

Assume moreover that $\textbf{I}$ and $\cat{C}$ are $k$-categories.
$\textbf{Hom}_k(\textbf{I},\cat{C})$ denotes the full subcategory of
$\textbf{Fun}(\textbf{I},\cat{C})$ consisting of the $k$-functors.
If we add the commutation relations $(f+g)-f-g=0$ and $(hf)-h.f=0$,
where $f$ and $g$ are morphisms in $\textbf{I}$ and $h\in k$, we
obtain that $\textbf{Hom}_k(\textbf{I},\cat{C})$ can be viewed as a
category $\cat{C}^\Sigma$.
\item \textbf{Complexes with values in $\cat{C}$.} $I=\mathbb{Z}$
(the set of integers), $\Phi=\{d_n\mid n\in\mathbb{Z}\}$. For every
$n$, $d_n$ of origin $n$ and extremity $n+1$. The commutation
relations are $d_{n+1}d_n=0$. In order to consider complexes with
positive degrees we take $I=\mathbb{N}$ (the set of positive
integers).
\end{enumerate}
\end{examples}

\medskip
\textbf{Faithful functors.} Let $T:\cat{C}\to \cat{D}$ be a functor.
$T$ is called \emph{faithful} (\emph{full}, \emph{fully faithful})
if the maps
$$T:\hom{\cat{C}}{C}{C'}\to \hom{\cat{D}}{T(C)}{T(C')}$$
are injective (surjective, bijective), for every $C,C'\in \cat{C}$.

$T$ is called \emph{conservative} if it reflects isomorphisms (for
every morphism in $f$ in $\cat{C}$, $T(f)$ is an isomorphism
implies $f$ is also an isomorphism). Obviously a fully faithful
functor is conservative.

\begin{proposition}\emph{[Freyd]}
\emph{(\cite[Theorem II.7.1]{Mitchell:1965})}\label{faithful}
\\ Let $T:\cat{C}\to \cat{D}$ be a faithful functor. Then
\begin{enumerate}[(1)]
\item $T$ reflects monomorphisms, epimorphisms, and commutative
diagrams.
\item If $\cat{C}$ and $\cat{D}$ are preadditive categories with zero
and $T$ is additive, then $T$ reflects zero objects.
\item If $\cat{C}$ and $\cat{D}$ are abelian categories and $T$
is zero preserving, then reflects exact sequences.
\end{enumerate}
\end{proposition}

\begin{proposition}\label{faithful-exact}
Let $T:\cat{C}\to \cat{D}$ be a functor between abelian
categories. Consider the below statements
\begin{enumerate}[(1)]
\item $T$ is faithful.
\item $T$ is conservative.
\item $T(C)\simeq 0\Longrightarrow C\simeq 0$ for every
$C\in\cat{C}$.
\end{enumerate}
Then:
\begin{enumerate}[(i)]
\item (1) implies (2).
\item If $T$ is additive, then (2) implies (3).
\item If $T$ is exact (and then additive by Proposition
\ref{middle-exact}), then (1), (2) and (3) are equivalent.
\end{enumerate}
\end{proposition}

\begin{proof}
First $(1)\Longrightarrow(2)$ Follows from Proposition
\ref{faithful} (1).

(ii) Assume that $T$ is additive and (2) holds. Let
$C\in\cat{C}$ such that $T(C)\simeq 0$. Then $T(0)=0:
T(C)\to T(C)$ is an isomorphism. Thus, $0:C\to C$ is
also an isomorphism. Hence $\hom{\cat{C}}{C}{C}$ is a
singleton.

(iii) Assume that $T$ is exact. It is enough to show
$(3)\Longrightarrow(1)$. For this, let $f$ be a morphism
in $\cat{C}$ such that $T(f)=0$, that is $\oper{Im}{T(f)}
=0$. We have $T(\oper{Im}{f})=\oper{Im}{T(f)}$. Hence
$\oper{Im}{f}=0$, i.e. $f=0$.
\end{proof}

\begin{remark}
We refer to \cite[Exercise 8.26]{Kashiwara/Schapira:2006} for an
example which shows that (2) does not imply (1) in general, and to
\cite[Exercise 8.27]{Kashiwara/Schapira:2006} for an example which
shows that (3) does not imply (2) (and then does not imply (1)) in
general.
\end{remark}

\textbf{Completely faithful modules}. Let $R$ be a ring and $U_R$
and $_RM$ modules. The \emph{annihilator in $M$ of $U$} is
$$\operatorname{Ann}_M(U)=\{m\in M\mid u\otimes m=0 \; \textrm{in}
\; U\tensor{R}M \;\textrm{for all}\; u\in U\}.$$ We say that $U_R$
is $_RM$-\emph{faithful} if $\operatorname{Ann}_M(U)=0$. Observe
that $U_R$ is $_RR$-faithful if and only if it is faithful.

\begin{lemma}\label{tensor}
Let $f:M\to M'$ be an epimorphism in $\rmod{R}$ and
$g:N\to N'$ an epimorphism in $\lmod{R}$. Then
$h=f\tensor{R}g:M\tensor{R}N\to M'\tensor{R}N'$ is an epimorphism
and $\operatorname{Ker}(h)$ is the set $K$ of $\sum m_i\otimes n_i$
such that $m_i\in \operatorname{Ker}(f)$ or $n_i\in
\operatorname{Ker}(g)$ for all $i$.
\end{lemma}

\begin{proof}
First $h$ is surjective since it is a composition of two surjective
maps. It is obvious that $K\subset\operatorname{Ker}(h)$. Let $Q$ be
the quotient group $(M\tensor{R}N)/K$ and $p:M\tensor{R}N\to Q$ is
the canonical epimorphism. Then $h$ induces a morphism
$h^*¨:Q\to M'\tensor{R}N'$ such that $h=h^*p$.

Let $m'\in M'$ and $n'\in N'$ and suppose that $m'=f(m_1)=f(m_2)$
and $n'=g(n_1)=g(n_2)$ with $m_1,m_2\in M$ and $n_1,n_2\in N$. Then
$m=m_2-m_1\in\operatorname{Ker}(f)$ and $n=n_2-n_1\in
\operatorname{Ker}(g)$. Therefore,
$$m_2\otimes n_2-m_1\otimes n_1=m\otimes n_1+m_1\otimes n+
m\otimes n \in K.$$ Hence, $p(m_1\otimes n_1)=p(m_2\otimes n_2)$.

Now, we can define a map $j:M'\times N'\to Q$ as follows. For
$m'\in M'$ and $n'\in N'$, $j(m'\otimes n')=p(m\otimes n)$ where
$m\in M$ and $n\in N$ such that $m'=f(m)$ and $n'=g(n)$. It is
easy to verify that $j$ is $R$-balanced. By the definition of
the tensor product of modules (see \cite[\S19]{Anderson/Fuller:1992}),
there is a unique abelian group morphism $j^*:M'\tensor{R}N'\to Q$
such that $j^*(m\otimes n)=j(m,n)$ for every $m\in M,n\in N$.

Finally, we have that $j^*h^*p=j^*h=p$. Then $j^*h^*=1_Q$ and
$h^*$ is injective. Hence $\operatorname{Ker}(h)=K$.
\end{proof}

In the following lemma we collect some basic properties of the
annihilator that we need to prove Lemma \ref{completely faithful}.
We refer to \cite[Exercises 19.18, 19.20]{Anderson/Fuller:1992} for
more properties.

\begin{lemma}\label{annihilator}
Let $R$ be a ring and $U_R$ and $_RM$ modules.
\begin{enumerate}[(a)]
\item $\operatorname{Ann}_M(U)$ is the (unique) smallest submodule
$K$ of $M$ such that $U$ is $M/K$-faithful.
\item Let $(U_i)_{i\in I}$ be a family of right $R$-modules. Then
$$\operatorname{Ann}_M(\oplus_{i\in I}U_i)=
\cap_{i\in I}\operatorname{Ann}_M(U_i).$$
\item If $U_R$ generates $V_R$,  then
$\operatorname{Ann}_M(U)\leq\operatorname{Ann}_M(V)$.
\item $U$ is $M$-faithful if and only if for every morphism
$f:N\to M$ in $\lmod{R}$, $$U\tensor{R}f=0\Longrightarrow f=0.$$
\end{enumerate}
\end{lemma}

\begin{proof}
(a) is straightforward using Lemma \ref{tensor}.

(b) is clear using the canonical isomorphism
$$(\oplus_{i\in I}U_i)\tensor{R}M\to
\oplus_{i\in I}(U_i\tensor{R}M),\;
(u_i)_i\otimes m\mapsto (u_i\otimes m)_i.$$

(c) Assume that $U_R$ generates $V_R$. Then there is an
epimorphism $U^{(I)}\to V$ in $\rmod{R}$, where $I$ is a
nonempty set. Hence,
$\operatorname{Ann}_M(U)=\operatorname{Ann}_M(U^{(I)})\leq
\operatorname{Ann}_M(V)$ (by (b)).

(d) ($\Rightarrow$) Obvious. ($\Leftarrow$) Let $m\in
\operatorname{Ann}_M(U)$. To obtain $m=0$, it is enough to
take $N=Rm$ and $f:Rm\to M$ the canonical injection.
\end{proof}

Following \cite[Exercise 19.19]{Anderson/Fuller:1992}, a
module $W_R$ is said to be \emph{completely faithful} if
it is $_RM$-faithful for every left $R$-module $M$.

\begin{lemma}\label{completely faithful}
Let $R$ be a ring and $W_R$ a module.
\begin{enumerate}[(1)]
\item $R_R$ is completely faithful.
\item Every generator in $\rmod{R}$ is completely faithful.
\item $W_R$ is completely faithful if and only if the
functor $W\tensor{R}-:\lmod{R}\to \textbf{Ab}$ is faithful.
\end{enumerate}
\end{lemma}

\begin{proof}
(1) is obvious. (2) is obvious from Lemma \ref{annihilator} (c)
and (1). (3) is obvious from Lemma \ref{annihilator} (d).
\end{proof}

\medskip
\textbf{Generators and cogenerators.} Let $\cat{C}$ be a category,
and $(U_i)_{i\in I}$ a family of objects of $\cat{C}$. Following
Grothendieck \cite[p. 134]{Grothendieck:1957}, we say that
$(U_i)_{i\in I}$ is a \emph{family of generators} of $\cat{C}$ if
for all $A\in \cat{C}$ and all subobject $B\neq A$, there exist
$i\in I$ and a morphism $u:U_i\to A$ that does not factor
through $B\to A$. We say that $U$ is a \emph{generator} if the
family $\{U\}$ is a family of generators.

A category $\cat{C}$ is called \emph{locally small} if for every
object $A\in \cat{C}$, the class of subobjects of $A$ is a set. We
refer to \cite{Mitchell:1965} for the definitions and the basic
properties of pullbacks and equalizers and their duals.

\begin{proposition} \emph{[Grothendieck]}
\\ Every category $\cat{C}$ that has a family of generators and
pullbacks is locally small.
\end{proposition}

\begin{proof}
Let $(U_i)_{i\in I}$ be a family of generators of $\cat{C}$ and
$A\in \cat{C}$. Set $E=\bigcup_{i\in I}\hom{\cat{C}}{U_i}{A}$. Let
$\textsf{L}(A)$ be the class of subobjects of $A$ and
$\mathcal{P}(E)$ the set of subsets of $E$. Consider the
correspondence $\phi:\textsf{L}(A)\to\mathcal{P}(E)$ defined by: for
a subobject $B$ of $A$, $\phi(B)$ is the set of morphisms $f:U_i\to
A$ that factors through $B\to A$ for some $i\in I$. Let
$B,B'\in\textsf{L}(A)$ such that $B\neq B'$ and a pullback
$$\xymatrix{P\ar[r]^{\beta_2}\ar[d]_{\beta_1} & B'\ar[d]^{\alpha_2}
\\ B \ar[r]_{\alpha_1} & A}$$
associated to the canonical injections $\alpha_1$ and $\alpha_2$.
By \cite[proposition I.7.1]{Mitchell:1965}, $\beta_1$ and $\beta_2$
are monomorphisms. We have that $P\neq B$ or $P\neq B'$ since
$B\neq B'$. Assume for example that $P\neq B$. Obviously $\phi(P)
\subset \phi(B)$. This inclusion is strict from the definition
of a family of generators. Now suppose that $\phi(B)=\phi(B')$
and let $u\in \phi(B)$. Therefore, there is an index $i$ such
that $u:U_i\to A$ is a morphism with $u=\alpha_1\xi_1=
\alpha_2\xi_2$ for some morphisms $\xi_1:U_i\to B$ and $\xi_2:
U_i\to B'$. By the definition of a pullback, there is a morphism
$\gamma:U_i\to P$ making commutative the following diagram
 $$\xymatrix{U_i\ar@/_/[ddr]_{\xi_1} \ar@/^/[drr]^{\xi_2}
 \ar[dr]|-\gamma
 \\
 & P \ar[r]^{\beta_2}\ar[d]^{\beta_1} & B'\ar[d]^{\alpha_2}
 \\
 & B \ar[r]_{\alpha_1} & A.}$$
Then $u=\alpha_1\beta_1\gamma$. Thus $\phi(B)=\phi(P)$. We obtain
then a contradiction. Hence the correspondence $\phi$ is injective.
It follows that $\textsf{L}(A)$ is a set and
$\operatorname{Card}(\textsf{L}(A))\leq
\operatorname{Card}(\mathcal{P}(E))$.
\end{proof}

\begin{proposition}\emph{\cite[Proposition 1.9.1]{Grothendieck:1957}}
\\ Let $\cat{C}$ be an abelian category that has arbitrary direct
sums. Let $(U_i)_{i\in I}$ be a family of objects of $\cat{C}$. Set
$U=\oplus_{i\in I}U_i$. Then the below statements are equivalent
\begin{enumerate}[(1)]
\item $(U_i)_{i\in I}$ is a family of generators of $\cat{C}$;
\item $U$ is a generator of $\cat{C}$;
\item every $A\in \cat{C}$ is isomorphic to a quotient object of
$U^{(J)}$ where $J$ is a set.
\end{enumerate}
\end{proposition}

\begin{proposition}\emph{\cite[Proposition 2.3.4]{Popescu:1973}}
\\ Let $\cat{C}$ be an abelian category and $(U_i)_{i\in I}$ a
family of objects of $\cat{C}$. Then the below statements are
equivalent
\begin{enumerate}[(1)]
\item $(U_i)_{i\in I}$ is a family of generators of $\cat{C}$;
\item for every $A,C\in \cat{C}$ and every distinct morphisms
$f,g\in \hom{\cat{C}}{A}{C}$, there are an index $i_0$ and a
morphism $h\in \hom{\cat{C}}{U_{i_0}}{A}$ such that $fh\neq gh$;
\item for every $A,C\in \cat{C}$ and every nonzero morphism
$f\in \hom{\cat{C}}{A}{C}$, there are an index $i_0$ and a morphism
$h\in \hom{\cat{C}}{U_{i_0}}{A}$ such that $fh\neq 0$.
\end{enumerate}
Hence, $U\in \cat{C}$ is a generator if and only if the covariant
functor $\hom{\cat{C}}{U}{-}$ is faithful.
\end{proposition}

From the last result and Proposition \ref{faithful-exact}, we obtain

\begin{proposition}\emph{\cite[Proposition V.6.3]{Stenstrom:1975}}
\label{proj-gen}
\\ A projective object $P$ of an abelian category is a generator if
and only if there is a nonzero morphism $P\to C$ for every nonzero
object $C$.
\end{proposition}

Dually we define a family of cogenerators of $\cat{C}$.

\medskip
An abelian category is called a \emph{Grothendieck category} if it
has a generator and satisfying the AB 5) condition. The following
result is due to Grothendieck and it extends Baer's criterion for
injective modules.

\begin{proposition}\emph{\cite[Lemme 1, p. 136]{Grothendieck:1957}}
\label{Baer}
 \\Let $\cat{C}$ be a Grothendieck category with a generator $U$.
 Then $M\in \cat{C}$ is injective if and only for every subobject
 $V$ of $U$, and every morphism $v:V\to M$, there is a morphism
 $U\to M$ extending $v$.
\end{proposition}

Notice that Proposition \ref{Baer} remains true for an arbitrary
family of generators.

An abelian category is said to have \emph{enough injectives} if
every object is a subobject of an injective object. Dually we define
a category having enough projectives.

\begin{proposition}\emph{\cite[Th\'eor\`eme 1.10.1]{Grothendieck:1957}}
 \\ Every Grothendieck category has enough injectives.
\end{proposition}

\begin{proposition}\emph{[Grothendieck]}
\emph{(\cite[Lemma 3.7.12]{Popescu:1973})}\label{inj.cogen.obj}
\\Every abelian category which has arbitrary direct sums, a generator,
and enough injectives (e.g. Grothendieck categories), has an
injective cogenerator object.
\end{proposition}

\medskip
\textbf{Adjoint and Frobenius pairs of functors.} Adjoint functors
were introduced and studied in \cite{Kan:1958}.

Let $S:\cat{C}\to \cat{D}$ and $T:\cat{D}\to \cat{C}$ be functors,
and $$\eta:\hom{\cat{D}}{S(-)}{-}\to\hom{\cat{C}}{-}{T(-)}$$ a
natural transformation of functors $\cat{C}^\circ\times\cat{D}\to
\textbf{Set}$. Let $$\zeta_C:=\eta_{C,S(C)}(1_{S(C)}):C\to TS(C),$$
for every $C\in \cat{C}.$

The following five results are due to W. Shih \cite{Shih:1958} (see
also \cite{MacLane:1998}).

\begin{lemma}
$\zeta$ is a natural transformation and the correspondence
$$\nat{\hom{\cat{D}}{S(-)}{-}}{\hom{\cat{C}}{-}{T(-)}}\to
\nat{1_\cat{C}}{TS},\quad \eta\mapsto\zeta,$$ is bijective with
inverse map $\zeta\mapsto \eta^\zeta$, where
$$\eta^\xi_{C,D}:\hom{\cat{D}}{S(C)}{D}\to\hom{\cat{C}}{C}{T(D)},
\quad \beta\mapsto T(\beta)\circ \zeta_C,$$ for every $C\in \cat{C},
D\in \cat{D}.$
\end{lemma}

Now let $$\eta':\hom{\cat{C}}{-}{T(-)}\to \hom{\cat{D}}{S(-)}{-}$$
be a natural transformation, and
$$\xi_D:=\eta'_{T(D),D}(1_{T(D)}):ST(D)\to D,$$ for every $D\in
\cat{D}.$ Dually we obtain

\begin{lemma}
$\xi$ is a natural transformation and the correspondence
$$\nat{\hom{\cat{C}}{-}{T(-)}}{\hom{\cat{D}}{S(-)}{-}}\to
\nat{ST}{1_\cat{D}},\quad \eta'\mapsto\xi,$$ is bijective with
inverse map $\xi\mapsto \eta'^\xi$, where
$$\eta'^\xi_{C,D}:\hom{\cat{C}}{C}{T(D)}\to
\hom{\cat{D}}{S(C)}{D},\quad \alpha\mapsto \xi_D\circ S(\alpha),$$
for every $C\in \cat{C}, D\in \cat{D}.$
\end{lemma}

Now let $\eta$ and $\eta'$ be natural transformations as above, and
$\zeta$ ad $\xi$ the associated natural transformations.

\begin{proposition}
$\eta'\eta$ is the identity transformation of the functor
$\hom{\cat{D}}{S(-)}{-}$ if and only if the composition of natural
transformations
$$\xymatrix{S\ar[r]^-{S\zeta} & STS\ar[r]^-{\xi S} & S}$$
is the identity transformation of the functor $S$, i.e.
\begin{equation}\label{adj-fun-id-trans-2}
\xi_{S(C)}S(\zeta_C)=1_{S(C)}\textrm{ for all } C\in \cat{C}.
\end{equation}
\end{proposition}

Dually we obtain

\begin{proposition}
$\eta\eta'$ is the identity transformation of the functor
$\hom{\cat{C}}{-}{T(-)}$ if and only if the composition of natural
transformations
$$\xymatrix{T\ar[r]^-{\zeta T} & TST\ar[r]^-{T\xi} & T}$$
is the identity transformation of the functor $T$, i.e.
\begin{equation}\label{adj-fun-id-trans-1}
T(\xi_D)\zeta_{T(D)}=1_{T(D)}\textrm{ for all } D\in \cat{D}.
\end{equation}
\end{proposition}

\begin{corollary}
Let $S:\cat{C}\to \cat{D}$  and $T:\cat{D}\to \cat{C}$ be functors.
To give a natural equivalence
\begin{equation}\label{adj-fun-nat-eq}
\eta:\hom{\cat{D}}{S(-)}{-}\to\hom{\cat{C}}{-}{T(-)}
\end{equation}
is the same as to give two natural transformations
$\zeta:1_\cat{C}\to ST$ and $\xi:ST\to 1_\cat{D}$ such that the
composition of natural transformations
$$\xymatrix{S\ar[r]^-{S\zeta} & STS\ar[r]^-{\xi S} & S}$$
is the identity transformation of the functor $S$, and the
composition of natural transformations
$$\xymatrix{T\ar[r]^-{\zeta T} & TST\ar[r]^-{T\xi} & T}$$
is the identity transformation of the functor $T$, i.e.
\begin{equation}\label{adj-fun-nat-eq-equation}
\xi_{S(C)}S(\zeta_C)=1_{S(C)},\quad
T(\xi_D)\zeta_{T(D)}=1_{T(D)}\quad\textrm{for all } C\in
\cat{C},D\in \cat{D}.
\end{equation}
\end{corollary}

Let $S:\cat{C}\to \cat{D}$  and $T:\cat{D}\to \cat{C}$ be functors.
We say that $T$ is a \emph{right adjoint} to $S$ (and symmetrically,
$S$ is a \emph{left adjoint} to $T$), or $(S,T)$ is a \emph{pair of
adjoint functors}, if one of the equivalent conditions of the last
corollary holds. Sometimes and following \cite{Kan:1958} we write
$\eta:S\dashv T$.

The natural equivalence \eqref{adj-fun-nat-eq} is called the
\emph{adjunction isomorphism}. The natural transformations
$\zeta:1_\cat{C}\to ST$ and $\xi:ST\to 1_\cat{D}$ are called the
\emph{unit} and the \emph{counit} of the adjunction while and the
equalities \eqref{adj-fun-nat-eq-equation} are called \emph{the
adjunction equalities}.

\begin{proposition}\emph{\cite[Theorem 3.2]{Kan:1958}}\label{Kan 3.2}
\\Let $S,S':\cat{C}\to \cat{D}$  and $T,T':\cat{D}\to \cat{C}$ be
functors and let $\eta:S\dashv T$ and $\eta':S'\dashv T'$. Assume
that $\sigma:S'\to S$ is a natural transformation. Then there is a
unique natural transformation $\tau:T\to T'$ such that the following
diagram is commutative
$$\xymatrix{\hom{\cat{D}}{S(C)}{D}\ar[r]^\eta
\ar[d]_{\hom{\cat{D}}{\sigma_C}{D}} &
\hom{\cat{C}}{C}{T(D)}\ar[d]^{\hom{\cat{C}}{C}{\tau_D}}
\\
\hom{\cat{D}}{S'(C)}{D}\ar[r]^{\eta'} & \hom{\cat{C}}{C}{T'(D)}}$$
for every $C\in \cat{C},D\in \cat{D}$.
\end{proposition}

\begin{corollary}\label{Kan 3.2'}
\begin{enumerate}[(1)]
\item Let $S:\cat{C}\to \cat{D}$  and $T:\cat{D}\to \cat{C}$ be
functors and let $\eta:S\dashv T$. Then the natural transformation
defined as in Prop. \ref{Kan 3.2} associated to the natural
transformation $1_S:S\to S$ is $1_T:T\to T$.
\item Let $S,S',S'':\cat{C}\to \cat{D}$ and $T,T',T'':\cat{D}\to
\cat{C}$ be functors and let $\eta:S\dashv T$, $\eta':S'\dashv T'$
and $\eta'':S''\dashv T''$. Assume that $\sigma:S'\to S$ and
$\sigma':S''\to S'$ are natural transformations. Let $\tau:T\to T'$
and $\tau':T'\to T''$ be the natural transformations defined as in
Prop. \ref{Kan 3.2} associated to $\sigma$ and $\sigma'$
respectively. Then the natural transformations defined as in Prop.
\ref{Kan 3.2} associated to $\sigma\sigma'$ is $\tau'\tau$.
\item Let $S,S':\cat{C}\to \cat{D}$  and $T,T':\cat{D}\to \cat{C}$ be
functors and let $\eta:S\dashv T$ and $\eta':S'\dashv T'$. Assume
that $\sigma:S'\to S$ is a natural equivalence. Then the natural
transformations defined as in Prop. \ref{Kan 3.2} associated to
$\sigma$ is so.
\end{enumerate}
\end{corollary}

Now, we give the duality theorem for adjoint functors. Let
$S:\cat{C}\to \cat{D}$  and $T:\cat{D}\to \cat{C}$ be functors and
$S^\circ:\cat{C}^\circ\to \cat{D}^\circ$ and
$T^\circ:\cat{D}^\circ\to \cat{C}^\circ$ be their duals
($S^\circ(C^\circ)=S(C)^\circ$ and $S^\circ(f^\circ)=S(f)^\circ$ for
every object $C\in \cat{C}$ and every morphism $f$ in $\cat{C}$).

\begin{proposition}\emph{\cite[Theorem 3.4]{Kan:1958}}\label{Kan 3.4}
\\Assume that $\eta:S\dashv T$. Define for every $C\in \cat{C}$ and
$D\in \cat{D}$ a map
$$\eta^\#_{D^\circ,C^\circ}:\hom{\cat{C}^\circ}{T^\circ(D^\circ)}{C^\circ}
\to \hom{\cat{D}^\circ}{D^\circ}{S^\circ(C^\circ)}$$ by
$$\eta^\#_{D^\circ,C^\circ}=\eta^{-1}_{C,D}.$$
Then $\eta^\#$ is a natural equivalence, and hence
$\eta^\#:T^\circ\dashv S^\circ$. Furthermore, $\eta^{\#\#}=\eta$.
\end{proposition}

Notice that by this duality theorem we obtain the dual result of
Proposition \ref{Kan 3.2}.

The following result whose proof is obvious will be useful.
\begin{proposition}\label{compoadjfun}
Let $\cat{C},\cat{C}',\cat{C}''$ be categories and let
$$\xymatrix{\cat{C}\ar@<2pt>[r]^S & \cat{C}'\ar@<2pt>[l]^T
\ar@<2pt>[r]^{S'}& \cat{C}''\ar@<2pt>[l]^{T'}}$$ be functors. If
$(S,T)$ and $(S',T')$ are pairs of adjoint functors, then
$(S'S,TT')$ is a pair of adjoint functors.
\end{proposition}

\begin{proposition}\label{M IV.3}
Let $S:\cat{C}\to \cat{D}$ be a functor between $k$-categories and
$T:\cat{D}\to \cat{C}$ a right adjoint of $T$. If $S$ (resp. $T$) is
$k$-linear, then $T$ (resp. $S$) is $k$-linear and the adjunction
bijections
$$\eta_{C,D}:\hom{\cat{D}}{S(C)}{D}\simeq\hom{\cat{C}}{C}{T(D)}$$
are isomorphisms of $k$-modules ($C\in\cat{C},D\in\cat{D}$).
Conversely, if the bijections $\eta_{C,D}$ are isomorphisms of
$k$-modules, then $S$ and $T$ are $k$-linear.
\end{proposition}

\begin{proof}
Suppose that $S$ is $k$-linear. By \cite[Proposition
V.1.4]{Mitchell:1965}, $T$ is additive and the bijections
$\eta_{C,D}$ are isomorphisms of abelian groups. Moreover, for every
morphism $\alpha:C\to T(D)$ in $\cat{C}$ and $h\in k$,
\begin{eqnarray*}
\eta^{-1}_{C,D}(h\alpha)&=& \xi_D\circ S(h\alpha)\\
 &=&\xi_D\circ (hS(\alpha))\\
 &=& h\xi_D\circ S(\alpha)
 \quad\textrm{(since $\cat{D}$ is a $k$-category)}
 \\
 &=& h\eta^{-1}_{C,D}(\alpha).
\end{eqnarray*}
Hence $\eta_{C,D}$ are isomorphisms of $k$-modules. On the other
hand, for every morphism $\beta:D\to D'$ in $\cat{D}$ and $h\in k$,
\begin{eqnarray*}
\eta^{-1}_{T(D),D'}(T(h\beta))&=& \xi_{D'}\circ ST(h\beta)\\
 &=&(h\beta)\circ\xi_D
 \quad\textrm{(since $\xi$ is a natural transformation)}\\
 &=& h\beta\circ\xi_D,
\end{eqnarray*}
and
\begin{eqnarray*}
\eta^{-1}_{T(D),D'}(hT(\beta))&=& \xi_{D'}\circ S(hT(\beta))\\
 &=&\xi_{D'}\circ (hST(\beta)) \\
 &=& h\xi_{D'}\circ ST(\beta)\\
 &=& h\beta\circ\xi_D
 \quad\textrm{(since $\xi$ is a natural transformation)}.
\end{eqnarray*}
Then $T(h\beta)=hT(\beta)$. Hence $T$ is $k$-linear.

Finally, using Proposition \ref{Kan 3.4}, we prove the second
statement. The converse is proved by a similar way.
\end{proof}

\begin{proposition}\label{right-adjoint-inj-obj}
Let $\cat{C}$ and $\cat{D}$ be two abelian categories, and let
$S:\cat{C}\rightarrow \cat{D}$ be a functor. Let
$T:\cat{D}\rightarrow \cat{C}$ be a right adjoint of $S$.
\begin{enumerate}[(1)]
\item If $N\in\cat{D}$ is injective and $S$ is exact, then
$T(N)\in\cat{C}$ is injective \emph{(\cite[Lemma
3.2.7]{Popescu:1973})}.
\item If moreover $\cat{D}$ is a category with sufficiently many
injectives (e.g. if $\cat{D}$ is a Grothendieck category), then $S$
is exact if and only if $T$ preserves injective objects
\emph{(\cite[Theorem 3.2.8]{Popescu:1973})}.
\end{enumerate}
\end{proposition}

Let $S:\cat{C}\to \cat{D}$ be a functor. If $T:\cat{D}\to \cat{C}$
is both a right and a left adjoint functor to $S$, we say that
$(S,T)$ is a \emph{Frobenius pair of functors}. In this case, $S$
and $T$ are called \emph{Frobenius functors}.

\medskip
\textbf{Adjoint functors in several variables.} To give a covariant
functor $S:\cat{X}\times\cat{Y}\to\cat{Z}$ is the same as to give

(a) a covariant functor $S(-,Y):\cat{X}\to\cat{Z}$ for every
$Y\in\cat{Y}$, and

(b) a natural transformation $S(-,y):S(-,Y)\to S(-,Y')$ for every
morphism $y:Y\to Y'$ in $\cat{Y}$ such that

(i) $S(-,y'y)=S(-,Y)S(-,Y)$ for every morphisms $y:Y\to Y'$ and
$y':Y'\to Y''$ in $\cat{Y}$, and

(ii) $S(-,1_Y)=1_{S(-,Y)}$ for every $Y\in\cat{Y}$.

Now assume that for every $Y\in\cat{Y}$ there are a functor
$T_Y:\cat{Z}\to\cat{X}$ and a natural equivalence
$$\eta_Y:\hom{\cat{Z}}{S(\cat{X},Y)}{\cat{Z}}\to
\hom{\cat{X}}{\cat{X}}{T_Y(\cat{Z})},$$ i.e. $\eta_Y:S(-,Y)\dashv
T_Y$. From Proposition \ref{Kan 3.2}, for every morphism $y:Y\to Y'$
in $\cat{Y}$ there is a unique natural transformation $T_y:T_{Y'}\to
T_Y$ making commutative the diagram
$$\xymatrix{\hom{\cat{Z}}{S(\cat{X},Y)}{\cat{Z}}\ar[r]^{\eta_Y}
 & \hom{\cat{X}}{\cat{X}}{T_Y(\cat{Z})}
\\
\hom{\cat{Z}}{S(\cat{X},Y')}{\cat{Z}}
\ar[u]^{\hom{\cat{Z}}{S(\cat{X},y)}{\cat{Z}}}\ar[r]^{\eta_{Y'}} &
\hom{\cat{X}}{\cat{X}}{T_Y(\cat{Z})}
\ar[u]_{\hom{\cat{X}}{\cat{X}}{(T_y)_{\cat{Z}}}}.}$$ Let $y':Y'\to
Y''$ be a morphism $y':Y'\to Y''$ in $\cat{Y}$. By Corollary
\ref{Kan 3.2'} (1),(2), we obtain $T_{y'y}=T_yT_{y'}$ and
$T_{1_Y}=1_{T_Y}$. We have then a functor
$T:\cat{Y}\times\cat{Z}\to\cat{X}$ contravariant in $\cat{Y}$ and
covariant in $\cat{Z}$ defined by
$$T(Y,Z)=T_Y(Z),\; T(Y,-)=T_Y,\; T(y,-)=T_y$$ for every $Y\in\cat{Y},
Z\in\cat{Z}$ and $y:Y\to Y'$ in $\cat{Y}$. Hence
$$\eta:\hom{\cat{Z}}{S(\cat{X},\cat{Y})}{\cat{Z}}\to
\hom{\cat{X}}{\cat{X}}{T(\cat{Y},\cat{Z})}$$ defined by
$\eta_{(X,Y,Z)}=(\eta_Y)_{(X,Z)}$ for every
$X\in\cat{X},Y\in\cat{Y}, Z\in\cat{Z}$, is a natural equivalence.
Thus we get:

\begin{proposition}\label{Kan 4.1}
Let $S:\cat{X}\times\cat{Y}\to\cat{Z}$ be a covariant functor.
Assume that for every $Y\in\cat{Y}$ there are a functor
$T_Y:\cat{Z}\to\cat{X}$ and a natural equivalence
$$\eta_Y:\hom{\cat{Z}}{S(\cat{X},Y)}{\cat{Z}}\to
\hom{\cat{X}}{\cat{X}}{T_Y(\cat{Z})},$$ i.e. $\eta_Y:S(-,Y)\dashv
T_Y$. Then there is a unique functor
$T:\cat{Y}\times\cat{Z}\to\cat{X}$ contravariant in $\cat{Y}$ and
covariant in $\cat{Z}$ and a unique natural equivalence
$$\eta:\hom{\cat{Z}}{S(\cat{X},\cat{Y})}{\cat{Z}}\to
\hom{\cat{X}}{\cat{X}}{T(\cat{Y},\cat{Z})}$$ such that for every
$X\in\cat{X},Y\in\cat{Y}, Z\in\cat{Z}$
$$T(Y,Z)=T_Y(Z),\; T(Y,-)=T_Y,\;\eta_{(X,Y,Z)}=(\eta_Y)_{(X,Z)}.$$
\end{proposition}

Let $S:\cat{X}\times\cat{Y}\to\cat{Z}$ be a covariant functor and
$T:\cat{Y}\times\cat{Z}\to\cat{X}$ a functor contravariant in
$\cat{Y}$ and covariant in $\cat{Z}$. If there is a natural
equivalence
$$\eta:\hom{\cat{Z}}{S(\cat{X},\cat{Y})}{\cat{Z}}\to
\hom{\cat{X}}{\cat{X}}{T(\cat{Y},\cat{Z})},$$ we say that $S$ is a
\emph{left adjoint of} $T$, and symmetrically, $T$ is a \emph{right
adjoint of} $S$. Sometimes we will use the notation $\eta:S\dashv
T$.

Now using Proposition \ref{Kan 3.2} and Lemma \ref{cube} we get

\begin{proposition}\label{Kan 4.4}
Let $S,S':\cat{X}\times\cat{Y}\to\cat{Z}$ be covariant functors and
$T,T':\cat{Y}\times\cat{Z}\to\cat{X}$ functors contravariant in
$\cat{Y}$ and covariant in $\cat{Z}$ such that $\eta:S\dashv T$ and
$\eta':S'\dashv T'$. Let $\sigma:S'\to S$ be a natural
transformation. Then there is a unique natural transformation
$\tau:T\to T'$ such that the following diagram is commutative
$$\xymatrix{\hom{\cat{Z}}{S(\cat{X},\cat{Y})}{\cat{Z}}\ar[r]^\eta
\ar[d]_{\hom{\cat{Z}}{\sigma_{(\cat{X},\cat{Y})}}{\cat{Z}}} &
\hom{\cat{X}}{\cat{X}}{T(\cat{Y},\cat{Z})}
\ar[d]^{\hom{\cat{X}}{\cat{X}}{\tau_{(\cat{Y},\cat{Z})}}} \\
\hom{\cat{Z}}{S'(\cat{X},\cat{Y})}{\cat{Z}}\ar[r]^{\eta'} &
\hom{\cat{X}}{\cat{X}}{T'(\cat{Y},\cat{Z})}.}$$
\end{proposition}

Note that Corollary \ref{Kan 3.2'} is also true for adjoint functors
in two variables.

Let $S:\cat{X}\times\cat{Y}\to\cat{Z}$ be a covariant functor.
Consider the functor
$S^\#:\cat{Y}\times\cat{X}^\circ\to\cat{Z}^\circ$ contravariant in
$\cat{Y}$ and covariant in $\cat{X}^\circ$ defined by
$$S^\#(Y,X^\circ)=S(X,Y)^\circ,\quad S^\#(y,x^\circ)=S(x,y)^\circ,$$
for every objects $Y\in\cat{Y}$, $X\in\cat{X}$ and morphisms $y$ in
$\cat{Y}$ and $x$ in $\cat{X}$.

Let $T:\cat{Y}\times\cat{Z}\to\cat{X}$ be a functor contravariant in
$\cat{Y}$ and covariant in $\cat{Z}$. Consider the covariant functor
$T^\#:\cat{Z}^\circ\times\cat{Y}\to\cat{X}^\circ$ defined by
$T^{\#\#}=T$.

We have
$$\hom{\cat{X}^\circ}{T^\#(Z^\circ,Y)}{X^\circ}=
\hom{\cat{X}}{X}{T(Y,Z)},$$
$$\hom{\cat{Z}^\circ}{Z^\circ}{S^\#( Y,X^\circ)}=
\hom{\cat{Z}}{S(X,Y)}{Z}$$ for every
$X\in\cat{X},Y\in\cat{Y},Z\in\cat{Z}$.

It is easy to prove the following duality theorem:

\begin{proposition}\label{Kan 4.5}
Let $S:\cat{X}\times\cat{Y}\to\cat{Z}$ be a covariant functor and
$T:\cat{Y}\times\cat{Z}\to\cat{X}$ a functor contravariant in
$\cat{Y}$ and covariant in $\cat{Z}$ such that $\eta:S\dashv T$.
Define for every $X\in\cat{X},Y\in\cat{Y},Z\in\cat{Z}$ the map
$$\eta^\#_{(Z^\circ,Y,X^\circ)}=\eta^{-1}_{(X,Y,Z)}:
\hom{\cat{X}^\circ}{T^\#(Z^\circ,Y)}{X^\circ}\to
\hom{\cat{Z}^\circ}{Z^\circ}{S^\#(Y,X^\circ)}.$$ Then $\eta^\#$ is a
natural equivalence, i.e. $\eta^\#:T^\#\dashv S^\#$. Moreover,
$\eta^{\#\#}=\eta$.
\end{proposition}

Note that Propositions \ref{Kan 3.4} and \ref{Kan 4.5} give the dual
of Proposition \ref{Kan 4.1}. Also Proposition \ref{Kan 4.5} gives
the dual of Proposition \ref{Kan 4.4}.

Finally we consider functors in more than two variables. Let
$S:\cat{X}\times\cat{Y}_1\times\dots\times\cat{Y}_n\to\cat{Z}$ be a
covariant functor and
$T:\cat{Y}_1\times\dots\times\cat{Y}_n\times\cat{Z}\to\cat{X}$ a
functor contravariant in $\cat{Y}_1,\dots,\cat{Y}_n$ and covariant
in $\cat{Z}$. Let $Y=\prod_i\cat{Y}_i$ be the cartesian product
category. $S$ can be viewed as a covariant functor
$\cat{Y}\times\cat{Z}\to\cat{X}$ and $T$ can be viewed as a functor
$\cat{Y}\times\cat{Z}\to\cat{X}$ contravariant in $\cat{Y}$ and
covariant in $\cat{Z}$. We return then to the case of functors in
two variables.

\medskip
\textbf{Equivalences of categories.} A functor $S:\cat{C}\to
\cat{D}$ is called an equivalence, if there exists a functor
$T:\cat{D}\to \cat{C}$ with natural equivalences $\eta:1\rightarrow
TS$ and $\varepsilon:ST\rightarrow 1$. A functor $S:\cat{C}\to
\cat{D}$ is called \emph{essentially surjective} if for every $D\in
\cat{D}$, there is $C\in C$ such that $S(C)\simeq D$.

\medskip
The following result is very well known.

\begin{theorem}
A functor $S:\cat{C}\to \cat{D}$ is an equivalence if and only if it
is faithful, full and essentially surjective.
\end{theorem}

A functor between $k$-categories is said to be a
$k$-\emph{equivalence} if it is $k$-linear and an equivalence.

\medskip
\textbf{Split exact and $Z$-pure sequences.} A monomorphism $f:A\to
B$ in a category $\cat{C}$ is called a \emph{section} (or a
\emph{coretraction}) if there is a morphism $g:B\to A$ such that
$gf=1_A$. A section is an epimorphism if and only if it is an
isomorphism. Dually we define a \emph{retraction}.

\begin{proposition}\emph{\cite[Proposition 8.3.14]{Kashiwara/Schapira:2006}}
\\Let $\xymatrix{0\ar[r] & B\ar[r]^{f} & A\ar[r]^{g} & C\ar[r] & 0}$
be a short exact sequence in an abelian category $\cat{C}$. Then the
following are equivalent
\begin{enumerate}[(1)]
\item there is $h:C\to A$ such that $gh=1_C$ (i.e. $g$ is a retraction);
\item there is $k:A\to B$ such that $kf=1_B$ (i.e. $f$ is a section);
\item there is $h:C\to A$ and $k:A\to B$ such that $1_A=fk+hg;$
\item there is $h:C\to A$ and $k:A\to B$ such that $gh=1_C$, $kf=1_B$,
$kh=0$, and $1_A=fk+hg;$
\item there are $\varphi=(k,g):A\to B\oplus C$ and
$\psi=f\oplus h:B\oplus C\to A$ are isomorphisms inverse to each
other.

\end{enumerate}
\end{proposition}

In such a case we say that the sequence is \emph{split exact}.

\begin{proposition}\label{AF 5.3}
Let $\xymatrix{0\ar[r] & A_1\ar[r]^{f_1} & A\ar[r]^{g_2} & A_2\ar[r]
& 0}$ be a short sequence in an abelian category $\cat{C}$. Then the
following are equivalent
\begin{enumerate}[(1)]
\item the sequence is split exact;
\item there is a sequence $\xymatrix{0\ar[r] & A_2\ar[r]^{f_2} &
A\ar[r]^{g_1} & A_1\ar[r] & 0}$ (necessarily split exact) such that
for every $i,j\in \{1,2\}$,
$$g_if_j=\delta_{i,j}1_{A_i}\quad \textrm{and}\quad
f_1g_1+f_2g_2=1_A;$$
\item there is an isomorphism $h:A_1\oplus A_2\to A$ such that the
following diagram is commutative
$$\xymatrix{
& &  A_1\oplus A_2\ar[dd]^h\ar[dr]^{p_2} &
\\0 \ar[r] & A_1\ar[ur]^{\iota_1}\ar[dr]^{f_1} && A_2 \ar[r] & 0
\\&& A\ar[ur]^{g_2} & .}$$
where $\iota_i:A_i\to A_1\oplus A_2$ and $p_i:A_1\oplus A_2\to A_i$,
$i\in \{1,2\}$, are such that $p_j\iota_i=\delta_{i,j}$ and
$\iota_1p_1+\iota_2p_2=1_{A_1\oplus A_2}.$
\end{enumerate}
\end{proposition}

\begin{proof}
The proof of \cite[Proposition 5.3]{Anderson/Fuller:1992} works
here.
\end{proof}

\begin{corollary}\label{split exactness}
Every additive functor between abelian categories preserves split
exactness.
\end{corollary}

More generally, and following \cite[40.5]{Brzezinski/Wisbauer:2003},
we say that a sequence $\xymatrix{B\ar[r]^{f} & A\ar[r]^{g} & C}$ of
morphisms in an abelian category is \emph{split exact} if it is
exact (i.e. $\oper{Ker}{g}=\oper{Im}{f}$) and the canonical
monomorphism $\oper{Im}{g}\to C$ is a section.

Let $n\in \{2,\ldots\}\subset\mathbb{N}$. We say that an exact
sequence $$\xymatrix{A_1\ar[r]^{f_1} & A_2\ar[r]^{f_2} &
A_3\ar[r]^-{f_3} & \ldots \ar[r]^{f_{n-1}} & A_n\ar[r]^-{f_n}&
A_{n+1}}$$ in an abelian category is \emph{split exact} if the
canonical monomorphism $\oper{Im}{f_i}\to A_{i+1}$ is a section for
every $i\in \{2,\ldots,n\}$.

Let $n\in \{2,\ldots\}\subset\mathbb{N}$ and $Z\in \lmod{A}$.
Following \cite[40.13]{Brzezinski/Wisbauer:2003} we say that an
exact sequence
\begin{equation}\label{pure}
\xymatrix{M_1\ar[r]^{f_1} & M_2\ar[r]^{f_2} & M_3\ar[r]^-{f_3} &
\ldots \ar[r]^{f_{n-1}} & M_n\ar[r]^-{f_n}& M_{n+1}}
\end{equation}
in $\rmod{A}$, is $Z$-\emph{pure} if it still exact under the
functor $-\tensor{A}Z$. If this is the case for every $Z\in
\lmod{A}$, we say simply that the sequence \eqref{pure} is
\emph{pure}. Of course if $Z$ is flat, then every exact sequence in
$\rmod{A}$ is $Z$-pure.

Let $f:M\to N$ be a monomorphism in $\rmod{A}$. $f$ is called
$Z$-\emph{pure} and $M$ is called a $Z$-\emph{pure submodule} of $N$
if the exact sequence $\xymatrix{0\ar[r] & M\ar[r]^{f} & N}$ is
$Z$-pure, that is if $f\tensor{A}Z$ is a monomorphism.

\begin{lemma}
The exact sequence \eqref{pure} is $Z$-pure if and only if the
canonical monomorphism $\oper{Im}{f_i}\to M_{i+1}$ is $Z$-pure for
every $i\in \{2,\ldots,n\}$.

In particular if the sequence \eqref{pure} is split exact then it is
pure.
\end{lemma}

\begin{proof}
We begin by factorizing the map $f_i$ through its image for all
$i\in \{2,\ldots,n\}$ as follows
$$\xymatrix{M_1\ar[r]^{f_1} & M_2\ar[r]^{f_2}
\ar[dr]^{p_2}& M_3\ar[r]^-{f_3} &
\ldots M_n\ar[r]^{f_n}\ar[dr]^{p_n} & M_{n+1}&
\\& & \oper{Im}{f_2}\ar@{^{(}->}[u]^{\iota_2}\ar[dr] \ldots &&
 \oper{Im}{f_n}\ar@{^{(}->}[u]^{\iota_n}\ar[dr]&
 \\ &&& 0\ldots && 0.}$$
 Since $\iota_2$ is a monomorphism,
 $\oper{Im}{f_1}=\oper{Ker}{f_2}=\oper{Ker}{p_2}$. Then the sequence
 \\$\xymatrix{M_1\ar[r]^{f_1} & M_2 \ar[r]^-{p_2}& \oper{Im}{f_2}\ar[r] &
 0}$ is exact. Since $-\tensor{A}Z$ is right exact,
 \\$\xymatrix{M_1\tensor{A}Z\ar[r]^{f_1\tensor{A}Z} & M_2\tensor{A}Z
 \ar[r]^-{p_2\tensor{A}Z}& \oper{Im}{f_2}\tensor{A}Z\ar[r] &
 0}$ is also exact. The same thing holds for every $i\in
 \{3,\ldots,n\}$. Then the proof follows easily.
\end{proof}

\begin{corollary} Let $Z\in \lmod{A}$ and $f: M
\to N$ be a morphism in $\rmod{A}$. Then the following are
equivalent
\begin{enumerate}[(1)]
\item the sequence
$$\xymatrix{0\ar[r] & \oper{Ker}{f}\tensor{A}Z \ar[r] & M\tensor{A}Z
 \ar[r]^{f\tensor{A}Z}& N\tensor{A}Z}$$
 is exact;
 \item the canonical monomorphisms $\oper{Ker}{f}\to M$ and
 $\oper{Im}{f}\to N$ are $Z$-pure.
\end{enumerate}
\end{corollary}

In such a case we say that $f$ is $Z$-\emph{pure}. If this the case
for every $Z\in \lmod{A}$, we say simply that $f$ is \emph{pure}.

\medskip
\textbf{The Yoneda--Grothendieck Lemma and representable functors.}
Let $\cat{C}$ be a category and $X\in\cat{C}$. We have then a
covariant functor
$h_X=\hom{\cat{C}}{-}{X}:\cat{C}^\circ\to\textbf{Set}$. Let
$X',Y\in\cat{C}$ and let $u:X\to X'$ be a morphism. Then
$$h_u=\hom{\cat{C}}{-}{u}:h_X\to h_{X'}$$ is a natural transformation.
Hence we have a ``functor''
$$h:\cat{C}\to\textbf{Fun}(\cat{C}^\circ,\textbf{Set}).$$

Let $F:\cat{C}^\circ\to\textbf{Set}$ be a covariant functor. Define
$$\alpha:\nat{h_X}{F}\to F(X),\; \eta\mapsto \eta_X(1_X)\in F(X)$$
and $$\beta:F(X)\to\nat{h_X}{F},\;\xi\mapsto \eta$$ where $\eta$ is
the natural transformation defined by
$$\eta_Y:\hom{\cat{C}}{Y}{X}\to F(X),\; \upsilon\mapsto
F(\upsilon)(\xi)$$ for every $Y\in\cat{C}$.

\begin{proposition}\emph{[Yoneda--Grothendieck]}\label{Yoneda}
\\ The maps $\alpha$ and $\beta$ are bijections inverse to each
other. (This includes that $\nat{h_X}{F}$ is a set.)
\end{proposition}

If we take $F=h_{X'}$ for $X'\in\cat{C}$, the map
$$\beta:\hom{\cat{C}}{X}{X'}\to\nat{h_X}{h_{X'}}$$ is exactly the map
$u\mapsto h_u$. Thus:

\begin{corollary}\label{Yoneda'}
The canonical ``functor''
$h:\cat{C}\to\textbf{Fun}(\cat{C}^\circ,\textbf{Set})$ is fully
faithful.
\end{corollary}

We say that a contravariant functor $F:\cat{C}\to \textbf{Set}$ is
\emph{representable} if there is $X\in \cat{C}$ such that $F\simeq
h_X$. Since $h$ is conservative then $X$ is determined up to an
isomorphism. Then $F$ is representable if and only if there exist
$X\in \cat{C},\xi\in F(X)$ such that $$\hom{\cat{C}}{Y}{X}\to
F(Y),\; \upsilon\mapsto F(\upsilon)(\xi)$$ is an isomorphism, for
every $Y\in\cat{C}$. We say that the couple $(X,\xi)$
\emph{represents} $F$. Usually we omit $\xi$ and we say only that
$X$ \emph{represents} $F$. It follows from Corollary \ref{Yoneda'}
that every category $\cat{C}$ is equivalent to the category of
all representable contravariant functors $\cat{C}\to\textbf{Set}$.

For the details we refer to \cite{Grothendieck:1961} or
\cite{Grothendieck/Dieudonne:1971}.

Using the axiom of choice and the Yoneda--Grothendieck Lemma, we
can easily prove the following useful result:

\begin{proposition}\label{Gabriel I.10}
Let $S:\cat{C}\to \cat{D}$ be a functor. Then the following are
equivalent
\begin{enumerate}[(1)]
\item the functor $S$ has a right adjoint;
\item for every object $D\in \cat{D}$, the contravariant functor
$C\mapsto\hom{\cat{D}}{S(C)}{D}$ is representable.
\end{enumerate}
\end{proposition}

We also define for $X\in \cat{C}$ a covariant functor
$h'_X=\hom{\cat{C}}{X}{-}:\cat{C}\to \textbf{Set}$. Dually, we
define representable covariant functors.

Using also the Yoneda--Grothendieck Lemma we prove the result below:

\begin{proposition}\emph{\cite[Theorem IV.3.1]{MacLane:1998}}
\label{Maclane-adj-fun}
\\Let $S:\cat{C}\to \cat{D}$ be a functor. Suppose that there is an
adjunction $(S,T)$ with unit $\zeta$ and counit $\xi$.
\begin{enumerate}[(i)]
\item $T$ is faithful if and only if $\xi_D$ is an epimorphism for
every $D\in \cat{D}$.
\item $T$ is full if and only if $\xi_D$ is a section for
every $D\in \cat{D}$.
\end{enumerate}
Hence, $T$ is faithful and full if and only if $\xi_D$ is an
isomorphism for every $D\in \cat{D}$.
\end{proposition}

\medskip
\textbf{Limits (Projective limits).} Let $\cat{C}$ be a category,
$S=(I,\Phi,d)$ a diagram scheme, and $D$ a diagram in $\cat{C}$
of scheme $S$. A family of morphisms $\{u_i:X\to D(i)\mid i\in I\}$,
is called \emph{compatible} for $D$ if $D(\varphi)u_i=u_j$ for every
arrow $\varphi:i\to j$ of $S$. Assume that $(u_i)$ is compatible and
let $v_i:Y\to D(i)$ be another compatible family. A \emph{morphism}
from $(v_i)$ to $(u_i)$ is a morphism $v:Y\to X$ in $\cat{C}$ such
that $v_i=u_iv$ for every $i\in I$. Then we obtain the category of
all compatible families for $D$. A final object for this category
is called a \emph{limit} for $D$. Hence a limit is determined up to
an isomorphism. We denote it by
$\{\underleftarrow{\lim}D\to D(i)\mid i\in I\}$.

Obviously pullbacks, products and equalizers are special cases of
limits. We say that $\cat{C}$ is $S$-\emph{complete} if every
diagram in $\cat{C}^S$ has a limit. We say that $\cat{C}$ is
\emph{complete} (resp. \emph{finitely complete}) if it is $\cat{C}$
is $S$-complete for every (resp. finite) diagram scheme $S$. Every
finitely complete category has pullbacks, equalizers and finite
products, and every complete category has products.

Dually we define colimits (inductive limits).

Assume that $\cat{C}$ is $S$-complete. Let $\upsilon:D\to D'$ be a
morphism in $\cat{C}^S$. By the definition of limit there exists a
unique morphism
$\underleftarrow{\lim}\upsilon:\underleftarrow{\lim}D\to
\underleftarrow{\lim}D'$ making the following diagram commutative
$$\xymatrix{\underleftarrow{\lim}D
\ar[r]^{\underleftarrow{\lim}\upsilon}\ar[d] &
\underleftarrow{\lim}D'\ar[d]\\
D(i)\ar[r]^{\upsilon_i}& D'(i)}$$ for every $i\in I$. Then we obtain
a functor $$\underleftarrow{\lim}:\cat{C}^S\to\cat{C}.$$ Now let
$A\in\cat{C}$. Define $D^A\in\cat{C}^S$ by $D^A(i)=A$ and
$D^A(\varphi)=1_A$ for every $i\in I$ and $\varphi\in\Phi$. A
morphism $D^A\to D$ in $\cat{C}^S$ can be identified with a
compatible family $\{A\to D(i)\mid i\in I\}$ for $D$. Let $u:A\to B$
be a morphism in $\cat{C}$. Consider the morphism $\upsilon^u:D^A\to
D^B$ in $\cat{C}^S$ defined by $\upsilon^u_i=u$ for every $i\in I$.
We obtain then a functor $I:\cat{C}\to\cat{C}^S$. Let $X\in\cat{C}$
and $\xi\in\hom{\cat{C}^S}{D^X}{D}$. Then $(X,\xi)$ is a limit for
$D$ if and only if $(X,\xi)$ represents the contravariant functor
$h_DI=\hom{\cat{C}^S}{-}{D}I:\cat{C}\to \textbf{Set}$. Hence, by
Proposition \ref{Gabriel I.10}, $\cat{C}$ is $S$-complete if and
only if the functor $I$ has a right adjoint. In such a case, $I$ is
a left adjoint to $\underleftarrow{\lim}$.

Let $F:\cat{C}\to\cat{C}'$ be a functor and $S$ a diagram scheme. We
know that $F$ induces a functor $F^S:\cat{C}^S\to\cat{C}'^S$. We say
that $F$ \emph{preserves limits} if $(F(u_i))_i$ is the limits for
$F(D)$ whenever $(u_i)$ is the limits for $D$. Of course, if $F$
preserves limits then it preserves pullbacks, equalizers and
products. Dually, we say that a functor $F:\cat{C}\to\cat{C}'$
\emph{preserves colimits} if its dual
$F^\circ:\cat{C}^\circ\to\cat{C}'^\circ$ preserves limits.

Let $G:\cat{C}\to\cat{C}'$ be another functor and $\eta:F\to G$ a
natural transformation. If $\{v_i:Y\to F(D(i))\mid i\in I\}$ is a
compatible family for the diagram $F(D)\in \cat{C}'^S$, then the
$(\eta_{D(i)}v_i)_i$ is a compatible family for $G(D)$. Assume that
$\eta$ is a natural equivalence and let $u_i:X\to D(i)$, $i\in I$,
be a family of morphism in $\cat{C}$. Then $(F(u_i))_i$ is a limit
for $F(D)$ if and only if $(G(u_i))_i$ is a limit for $G(D)$. Hence
$F$ preserves limits if and only if $G$ preserves limits.

It is also easy to verify the following two facts. A morphism $u$ in
$\cat{C}$ is a monomorphism if and only if $h'_A(u)$ is an injective
map for every $A\in \cat{C}$. Let $\{u_i:X\to D(i)\mid i\in I\}$ be
a family of morphism in $\cat{C}$. Then $(u_i)_i$ is a limit for $D$
if and only if $(h'_A(u_i))_i$ is a limit for $h'_A(D)$ for every
$A\in \cat{C}$. Using the fact that
$h'_{A^\circ}=h_A:\cat{C}^\circ\to \textbf{Set}$, we obtain the dual
of these considerations.

Obviously, an object $X\in\cat{C}$ is a final object if and only if
$h'_A(X)$ is a final object of the category $\textbf{Set}$ for every
$A\in \cat{C}$ (Since the final objects of $\textbf{Set}$ are the
sets with one element). Let $\cat{C}$ and $\cat{C}'$ be
$k$-categories and $F:\cat{C}\to\cat{C}'$ a functor. We have $F$ is
$k$-linear if and only if $h'_XF$ is $k$-linear for every
$X\in\cat{C}$ (of course the functors $h'_X$ are $k$-linear).
Moreover, $F$ is $k$-linear if and only if its dual $F^\circ$ is
$k$-linear.

A functor $F$ is called a \emph{monofunctor} if it preserves
monomorphisms. Dually, a functor $F$ is an \emph{epifunctor} if its
dual $F^\circ$ is a monofunctor, or equivalently, if it preserves
epimorphisms.

From the dual of Proposition \ref{Gabriel I.10}, a functor
$G:\cat{C}'\to\cat{C}$ has a left adjoint if and only if $h'_AG$ is
representable for every $A\in \cat{C}$.

Let $S:\cat{X}\times\cat{Y}\to\cat{Z}$ be a covariant functor and
$T:\cat{Y}\times\cat{Z}\to\cat{X}$ a functor contravariant in
$\cat{Y}$ and covariant in $\cat{Z}$ such that $\eta:S\dashv T$.
Consider the covariant functors $S_X=S(X,-):\cat{Y}\to\cat{Z}$ and
$T_Z=T(-,Z):\cat{Y}^\circ\to \cat{X}$. We have then a natural
equivalence
$$h'_{Z^\circ}(S_X)^\circ\simeq h'_XT_Z:\cat{Y}^\circ\to \textbf{Set}$$
for every $X\in\cat{X}, Z\in \cat{Z}$.

Then from the above considerations we prove the following result:

\begin{proposition}\label{adj fun-lim}
\begin{enumerate}[(1)]
\item Let $F:\cat{C}\to \cat{C}'$ be a functor and
$G:\cat{C}'\to\cat{C}$ a right adjoint functor to $F$. Then
\begin{enumerate}[(i)]
\item $G$ is a monofunctor and
preserves limits. $G$ preserves also final objects.
\item Dually, $F$ is an
epifunctor and preserves colimits. $F$ preserves also initial
objects.
\end{enumerate}
\item Let $S:\cat{X}\times\cat{Y}\to\cat{Z}$ be a covariant functor and
$T:\cat{Y}\times\cat{Z}\to\cat{X}$ a functor contravariant in
$\cat{Y}$ and covariant in $\cat{Z}$ such that $\eta:S\dashv T$.
Consider the covariant functors $S_X=S(X,-):\cat{Y}\to\cat{Z}$ and
$T_Z=T(-,Z):\cat{Y}^\circ\to\cat{X}$. Then
\begin{enumerate}[(a)]
\item $S_X$ preserves colimits for every $X\in\cat{X}$ if and only
if $T_Z$ preserves limits (as a covariant functor) for every $Z\in
\cat{Z}$.
\item Assume moreover that the categories $\cat{X}$, $\cat{Y}$ and
$\cat{Z}$ are $k$-categories. Then $S_X$ is $k$-linear for every
$X\in\cat{X}$ if and only if $T_Z$ is $k$-linear (as a covariant
functor) for every $Z\in \cat{Z}$.
\end{enumerate}
\end{enumerate}
\end{proposition}

\begin{corollary}
Let $\cat{C}$ be a category and $S$ a diagram scheme. If $\cat{C}$
is $S$-complete, then the functor
$\underleftarrow{\lim}:\cat{C}^S\to\cat{C}$ is a monofunctor and
preserves limits.
\end{corollary}

\begin{proposition}\emph{[Freyd]}
\emph{(see \cite[Theorem II.7.1]{Mitchell:1965})}\label{faithful'}
\\ Let $F:\cat{C}\to \cat{C}'$ be a faithful functor. Then
\begin{enumerate}[(1)]
\item If $\cat{C}$ is an abelian category and $\cat{C}'$ is additive,
then $F$ reflects limits and colimits of finite diagrams.
\item If moreover $F$ is full, then $F$ reflects limits and
colimits.
\end{enumerate}
\end{proposition}

For many more interesting properties of limits, see
\cite{Mitchell:1965}.

\medskip
\textbf{Separable functors.} Separable functors were introduced and
studied in \cite{Nastasescu/VanDerBergh/VanOystaeyen:1989}. Further
study of these functors is given in \cite{Rafael:1990}.

Let $T:\cat{C}\to \cat{D}$ be a functor between arbitrary
categories. $T$ is called a \emph{separable functor} if for every
$C,C'\in \cat{C}$, there exists a map
$$\Phi=\Phi_{C,C'}:\hom{\cat{D}}{T(C)}{T(C')}\to \hom{\cat{C}}{C}{C'}$$
such that
\begin{enumerate}[(SF1)]
\item For every morphism $f:C\to C'$ in $\cat{C}$, $\Phi(T(f))=f.$
\item For every morphisms $f:C_1\to C_1'$ and $g:C_2\to C_2'$ in
$\cat{C}$, and every commutative diagram in $\cat{D}$
$$\xymatrix{T(C_1)\ar[rr]^-h \ar[d]_{T(f)} && T(C_2) \ar[d]^{T(g)}
\\ T(C_1') \ar[rr]^-{h'} && T(C_2')},$$
the following diagram
$$\xymatrix{C_1\ar[rr]^-{\Phi(h)} \ar[d]_f && C_2 \ar[d]^g
\\ C_1' \ar[rr]^-{\Phi(h')} && C_2'}$$
is also commutative.
\end{enumerate}

Obviously, from (SF1), every separable functor is faithful. Let
$T:\cat{C}\to \cat{D}$ be a separable functor and
$T^\circ:\cat{C}^\circ\to \cat{D}^\circ$ its dual. Then, for every
objects $C,C'\in\cat{C}$,
$\hom{\cat{D}^\circ}{T^\circ(C^\circ)}{T^\circ(C'^\circ)}=
\hom{\cat{D}}{T(C')}{T(C)}$ and
$\hom{\cat{C}^\circ}{C^\circ}{C'^\circ}=\hom{\cat{C}}{C'}{C}$.
Define for every objects $C,C'\in\cat{C}$ a map
$$\Phi^\#_{C^\circ,C'^\circ}:
\hom{\cat{D}^\circ}{T^\circ(C^\circ)}{T^\circ(C'^\circ)}
\to\hom{\cat{C}^\circ}{C^\circ}{C'^\circ}$$ by
$$\Phi^\#_{C^\circ,C'^\circ}=\Phi_{C',C}.$$ We have, for every
morphism $h:T(C')\to T(C)$ in $\cat{C}$,
$$\Phi^\#_{C^\circ,C'^\circ}(h^\circ)=\Phi_{C',C}(h)^\circ.$$
The functor $T^\circ:\cat{C}^\circ\to \cat{D}^\circ$ is clearly
separable. We use this fact to prove the dual results of separable
functors.

\begin{examples} \cite[Lemma 1.1]{Nastasescu/VanDerBergh/VanOystaeyen:1989}
\begin{enumerate}[(1)]
\item An equivalence of categories is separable.
\item Let $F:\cat{C}\to\cat{D}$ and $G:\cat{D}\to\cat{E}$ be
functors. Then
\begin{enumerate}[(i)]
\item If $F$ and $G$ are separable then $GF$ is separable.
\item If $GF$ is separable then $F$ is separable.
\end{enumerate}
\end{enumerate}
\end{examples}

\begin{proposition}\emph{\cite[Proposition
1.2]{Nastasescu/VanDerBergh/VanOystaeyen:1989}}\label{separable}
\\Let $F:\cat{C}\to\cat{D}$ be a functor and $C\in\cat{C}$. Then
\begin{enumerate}[(1)]
\item $F$ reflects sections and retractions.

Now, assume that $\cat{C}$ and $\cat{D}$ are abelian categories.
Then
\item If $F$ is a monofunctor, then $F$
reflects quasi-simple objects (=every subject splits off).
\item If $F$ is a monofunctor, then $F$ reflects injective objects.
Dually, if $F$ is an epifunctor, then $F$ reflects projective
objects.
\end{enumerate}
\end{proposition}

\begin{proposition}\emph{[Rafael] (\cite[Theorem 1.2]{Rafael:1990})}
\label{Rafael}
\\ Let $S:\cat{C}\to\cat{D}$ be a functor and $T:
\cat{D}\to\cat{C}$ a right adjoint to $S$. Let $\zeta:1_\cat{C}\to
TS$ and $\xi:ST\to 1_\cat{D}$ be the unit and the counit of the
adjunction, respectively.
\begin{enumerate}[(1)]
\item $S$ is separable if and only if $\zeta$ splits
(i.e. there is a natural transformation $\zeta':TS\to 1_\cat{C}$
such that $\zeta'_C\circ\zeta_C=1_C$ for all $C\in \cat{C}$).
\item Dually, $T$ is separable if and only if $\xi$ cosplits
(i.e. there is a natural transformation $\xi':1_\cat{D}\to ST$ such
that $\xi_D\circ\xi'_D=1_D$ for all $D\in \cat{D}$).
\end{enumerate}
\end{proposition}

\medskip
\textbf{Locally finitely generated and locally noetherian
categories.} Let $\cat{C}$  be a Grothendieck category. An object
$C$ of $\cat{C}$ is said to be \emph{finitely generated} \cite[p.
121]{Stenstrom:1975} if the (upper continuous) lattice of all
subobjects $\textsf{L}(C)$ is compact, that is, if $C = \sum_i C_i$
for a direct family of subobjects $C_i$ of $C$, then $C = C_{i_0}$
for some index $i_0$. Notice that a subobject $B$ of $C$ is finitely
generated as object if and only if $B$ is compact in the lattice
$\textsf{L}(C)$ (Since there is a canonical order isomorphism (and
then a lattice isomorphism) $[0,B]\to\textsf{L}(B)$).

An object $C\in \cat{C}$ is called \emph{finitely presented} if it
is finitely generated and every epimorphism $B\to C$ where $B$ is
finitely generated  has a finitely generated kernel. An object $C\in
\cat{C}$ is called \emph{noetherian} if the lattice (of all
subobjects) $\textsf{L}(C)$ is noetherian, that is, every infinite
ascending chain in it is stationary.

\begin{lemma}\emph{\cite[Lemma V.3.1]{Stenstrom:1975}}\label{V.3.1}
\\Let $0\to C'\to C\to C''\to 0$ be a short exact sequence in $\cat{C}$.
Then
\begin{enumerate}[(a)]
\item If $C$ is finitely generated, then so is $C''$.
\item If $C'$ and $C''$ are finitely generated, then so is $C$.
\end{enumerate}
\end{lemma}

\begin{proposition}\emph{\cite[Proposition V.3.2]{Stenstrom:1975}}
\label{V.3.2}
\\The following are equivalent
\begin{enumerate}[(1)]
\item $C$ is finitely generated;
\item for every directed family of subobjects $(D_i)_{i\in I}$ of a
certain object $D$ of $\cat{C}$, and every monomorphism $\alpha:C\to
\sum_ID_i$ in $\cat{C}$, there exists some $i_0\in I$ and a morphism
$\beta:C\to D_{i_0}$ making commutative the following diagram
$$\xymatrix{C\ar[r]^-\alpha \ar[dr]^\beta & \sum D_i
\\ & D_{i_0}\ar[u]^{};}$$
\item the canonical monomorphism of abelian groups
$$\underset{\underset{I}{\longrightarrow}}\lim\hom{\cat{C}}{C}{D_i}
\to \hom{\cat{C}}{C}{\sum D_i}$$ is an isomorphism. In such a case,
we say that the functor $\hom{\cat{C}}{C}{-}$ \emph{preserves direct
unions}.
\end{enumerate}
\end{proposition}

\begin{proof}
$(2)\Longrightarrow(1)$ Obvious. $(1)\Longrightarrow(2)$ From Lemma
\ref{V.3.1}, $\oper{Im}{\alpha}$ is finitely generated. We have
$(D_i)_{i\in I}$ is a directed family of subobjects of $\sum D_i$.
Hence $\oper{Im}{\alpha}\subset D_{i_0}$ for some $i_0\in I$.
Finally, (3) is a reformulation of (2).
\end{proof}

\begin{lemma}\label{lefadjfunf.g.obj}
Let $S:\cat{C}\to \cat{D}$ be a functor between Grothendieck
categories, and $T:\cat{D}\to \cat{C}$ a right adjoint functor to
$S$, that is exact and preserves coproducts. Then $S$ preserves
finitely generated objects.

In particular, a Frobenius functor between Grothendieck categories
preserves finitely generated objects.
\end{lemma}

\begin{proof}
Let $C$ be a finitely generated object in $\cat{C}$, $(D_i)_{i\in
I}$ be a directed family of subobjects of a certain object
$D\in\cat{C}$, and $$\eta:\hom{\cat{C}}{-}{T(-)}\to
\hom{\cat{D}}{S(-)}{-}$$ be the natural equivalence.

Since $\eta_{C,-}$ is natural in $D_i\to \sum D_i$, we have the
commutative diagram
\begin{equation}\label{dlefadjfunf.g.obj}
\xymatrix{\hom{\cat{C}}{C}{T(D_i)}\ar[rr]^{\eta_{C,D_i}} \ar[d] &&
\hom{\cat{D}}{S(C)}{D_i}\ar[d]\\
\hom{\cat{C}}{C}{T(\sum D_i)}\ar[rr]^{\eta_{C,\sum D_i}} &&
\hom{\cat{D}}{S(C)}{\sum D_i}.}
\end{equation}
Since $T$ is exact, $(T(D_i))_{i\in I}$ is a directed family of
subobjects of $T(D)$. Let $\iota_i:D_i\to D$ be the canonical
injection for every $i\in I$. Since $T$ is exact and preserves
coproducts, we have
$$\sum_I T(D_i)=\oper{Im}{T(\oplus_I\iota_i)}=
T\big(\oper{Im}{\oplus_I\iota_i}\big)=T\Big(\sum_I D_i\Big).$$ (To
obtain the last statement we can also use that $T$ preserves
inductive limits.) Now, from the commutativity of the diagram
\eqref{dlefadjfunf.g.obj} and Proposition \ref{V.3.2}, $S(C)$ is
finitely generated.
\end{proof}

\begin{proposition}\emph{\cite[Proposition V.3.4]{Stenstrom:1975}}
\label{V.3.4}
\\Let $\cat{C}$ be a locally finitely generated category. An object
$C\in \cat{C}$ is finitely presented if and only if the functor
$\hom{\cat{C}}{C}{-}$ preserves direct limits.
\end{proposition}

Part of the following result is \cite[Proposition
V.4.1]{Stenstrom:1975}.

\begin{proposition}\label{V.4.1}
Let $C\in\cat{C}$. The following are equivalent
\begin{enumerate}[(1)]
\item $C$ is noetherian;
\item every subobject of $C$ is finitely generated;
\item every nonempty set of subobjects of $C$ has a maximal element.
\end{enumerate}
\end{proposition}

Since the above result is stated in \cite{Stenstrom:1975} without
proof, we propose a proof of it.

\begin{proof}
$(1)\Longrightarrow(3)$ Let $\textsf{L}$ be a nonempty set of
subobjects of $C$. Suppose that $\textsf{L}$ does not have a maximal
element, that is, for every $A\in\textsf{L}$, the set
$\{A'\in\textsf{L}\mid A<A'\}$ is not empty. Hence, by the axiom of
choice, there is an infinite strictly ascending chain of subobjects
of $C$.

$(3)\Longrightarrow(2)$ Let $B$ be a subobject of $C$, and $B=\sum
B_i$, where $(B_i)_i$ is a directed family of subobjects of $B$.
This family has a maximal element, say $B_{i_m}$. Hence $B=B_{i_m}$.

$(2)\Longrightarrow(1)$ Suppose that every subobject of $C$ is
finitely generated. Let $$C_1\subset C_2\subset\ldots \subset
C_n\subset\ldots$$ be an ascending chain of subobjects of $C$. Set
$B=\sum_{\mathbb{N}-\{0\}} C_n$. Since $B$ is finitely generated,
there is $n$ such that $B=C_n$. Hence $B=C_{n+i}$ for every $i\in
\mathbb{N}$.
\end{proof}

\begin{remark}
The equivalence $(1)\Leftrightarrow(3)$ of Proposition \ref{V.4.1}
is true in a more general context. Indeed, a lattice $L$ is
noetherian if and only if every nonempty subset of $L$ has a maximal
element (see \cite[Proposition III.3.4]{Stenstrom:1975}).
\end{remark}

\begin{proposition}\emph{\cite[Proposition V.4.2]{Stenstrom:1975}}
\label{V.4.2}
\\ Let $0\to C'\to C\to C''\to 0$ be a short exact sequence in
$\cat{C}$. Then $C$ is noetherian if and only if both $C'$ and
$C''$ are noetherian.
\end{proposition}

\begin{proposition}\label{fgfpnoeth}
\begin{enumerate}[(1)]
\item Let $C=C_1\oplus\ldots\oplus C_n\in \cat{C}$. Then $C$ is
finitely generated if and only each $C_i$ ($i=1,\ldots,n$) is
finitely generated.
\item If $C$ is the sum of many finitely generated subobjects,
then $C$ is finitely generated.
\item Every finitely generated projective object is finitely
presented.
\item Let $C=C_1\oplus\ldots\oplus C_n\in \cat{C}$. Then $C$ is
noetherian if and only each $C_i$ ($i=1,\ldots,n$) is finitely
noetherian.
\item If $C$ is the sum of many noetherian subobjects, then
$C$ is noetherian.
\item The following statements are equivalent
\begin{enumerate}[(a)]
\item every finitely generated object is noetherian;
\item every finitely generated object is finitely presented.
\end{enumerate}
In particular, a ring $R$ is right noetherian if and only if
every finitely generated right $R$-module is finitely
presented.
\end{enumerate}
\end{proposition}

\begin{proof}
(1) It is enough to prove it for $n=2$. In such a case, there are
morphisms $\iota_i:C_i\to C$ and $\pi_i:C\to C_i$ ($i=1,2$) such
that $\pi_i\iota_j=\delta_{ij}$ and $\iota_1\pi_1 +\iota_2\pi_2=1$.
Hence we have two short exact sequences
$$\xymatrix{0\ar[r] & C_1\ar[r]^{\iota_1} &
C \ar[r]^{\pi_2} & C_2\ar[r] & 0}$$
and
$$\xymatrix{0\ar[r] & C_2\ar[r]^{\iota_2} &
C \ar[r]^{\pi_1} & C_1\ar[r] & 0.}$$ By Lemma \ref{V.3.1}, (1)
follows.

(2) follows directly from (1) and Lemma \ref{V.3.1}(i). Notice that
we can prove it directly using \cite[Lemma III.5.2]{Stenstrom:1975}.

(3) Let $C$ be a finitely generated projective object, and
$\alpha:B\to C$ be an epimorphism, where $B$ is finitely generated.
Then we obtain the short exact sequence
$$\xymatrix{0\ar[r] & \oper{Ker}{\alpha}\ar[r] &
B \ar[r]^{\alpha} & C\ar[r] & 0.}$$ Since $C$ is projective, this
sequence splits and $B\simeq \oper{Ker}{\alpha}\oplus C$. Hence, by
(1), $\oper{Ker}{\alpha}$ is finitely generated.

The proofs of (4) and (5) are analogous to that of (1) and (2),
respectively, using Proposition \ref{V.4.2}.

(6) Follows directly from Propositions \ref{V.4.2}, \ref{V.4.1}, and
Lemma \ref{V.3.1}(i).
\end{proof}

The category $\cat{C}$ is \emph{locally finitely generated} if the
lattice $\textsf{L}(C)$ is compactly generated for all $C\in
\cat{C}$, i.e. every object of $\cat{C}$ is a sum of a (directed)
family of finitely generated subobjects (see \cite[p.
73]{Stenstrom:1975}).

\begin{proposition}\label{l.f.g.}
Let $\cat{C}$ be a Grothendieck category. The following are
equivalent
\begin{enumerate}[(1)]
\item $\cat{C}$ is locally finitely generated;
\item $\cat{C}$ has a family of finitely generated generators.
\end{enumerate}
\end{proposition}

\begin{proof}
$(1)\Longrightarrow (2)$ Let $U$ be a generator of $\cat{C}$. Then
there is a family of finitely generated subobjects of $C$,
$(U_i)_{i\in I}$, such that there is an epimorphism $\bigoplus_{i\in
I}U_i\to U$. Hence $(U_i)_{i\in I}$ is a family of finitely
generated generators of $\cat{C}$.

$(2)\Longrightarrow (1)$ Let $(U_i)_{i\in I}$ be a family of
finitely generated generators of $\cat{C}$. Set $U=\bigoplus_{i\in
I}U_i$. Let $C\in \cat{C}$. Then there are a set $J$ and an
epimorphism $U^{(J)}\to C$. Consider the family
$(U_{(i,j)})_{(i,j)\in I\times J}$ such that $U_{(i,j)}=U_i$ for
every $(i,j)\in I\times J$. We have
$U^{(J)}\simeq\bigoplus_{(i,j)\in I\times J}U_{(i,j)}$. It follows
that there is an epimorphism $f:\bigoplus_{(i,j)\in I\times
J}U_{(i,j)}\to C$. Let $f_{(i,j)}:U_{(i,j)}\to C$, $(i,j)\in I\times
J$, be the morphisms defining $f$. Let
$f_{(i,j)}=\oper{im}{f_{(i,j)}}\circ p_{(i,j)}$ be the canonical
factorization, for every $(i,j)\in I\times J$. Then we obtain the
commutative diagram
$$\xymatrix{\bigoplus_{(i,j)}U_{(i,j)}
\ar[rr]^-f \ar[dr]_{\bigoplus_{(i,j)}p_{(i,j)}}&& C
\\ & \bigoplus_{(i,j)}\oper{Im}{f_{(i,j)}} \ar[ur] & .}$$
Since $f$ is an epimorphism, then the morphism $\bigoplus_{(i,j)\in
I\times J}\oper{Im}{f_{(i,j)}}\to C$ is also an epimorphism.
Finally, by Lemma \ref{V.3.1} (a), each $\oper{Im}{f_{(i,j)}}$ is
finitely generated.
\end{proof}

The category $\cat{C}$ is \emph{locally noetherian} \cite[p.
123]{Stenstrom:1975} if it has a family of noetherian generators.

\begin{proposition}\label{l.n.}
Let $\cat{C}$ be a Grothendieck category. The following are
equivalent
\begin{enumerate}[(1)]
\item $\cat{C}$ is locally noetherian;
\item every object of $\cat{C}$ is a sum of a (directed)
family of noetherian subobjects;
\item $\cat{C}$ is locally finitely generated and every finitely
generated object of $\cat{C}$ is noetherian;
\item $\cat{C}$ is locally finitely generated and every finitely
generated object of $\cat{C}$ is finitely presented.
\end{enumerate}
\end{proposition}

\begin{proof}
First observe that if an object of $\cat{C}$ is a sum of a family
of noetherian subobjects, then it is a sum of a directed family
of noetherian subobjects (by Proposition \ref{fgfpnoeth} (5)).

$(1)\Longleftrightarrow (2)$ Analogous to that of Proposition
\ref{l.f.g.} using Proposition \ref{V.4.2}.

$(2)\Longrightarrow(3)$ Obvious from Proposition \ref{V.4.1}.

$(3)\Longrightarrow(1)$ Obvious from Proposition \ref{l.f.g.}.

$(3)\Longleftrightarrow (4)$ Obvious from Proposition
\ref{fgfpnoeth} (6).
\end{proof}

\medskip
\textbf{Small objects.} Let $C$ be an additive category and $X$ an
object in $\cat{C}$.

\begin{lemma}\emph{\cite[Lemma II.16.1]{Mitchell:1965}}
\\Let $f:X\to \bigoplus_{i\in I}X_i$ be a morphism in $\cat{C}$ and
$F$ is a nonempty finite subset of $I$. The following are equivalent
\begin{enumerate}[(a)]
\item there exists a factorization of $f$
\begin{equation}\label{small diag}
\xymatrix{X \ar[rr]^f \ar[dr]^{f'}&& \bigoplus_{i\in I}X_i
\\ & \bigoplus_{j\in F}X_j \ar[ur]^{u_F}&},
\end{equation}
where $u_F$ is the unique morphism satisfying $u_Fu'_j=u_j$ for
every $j\in F$, with $u_i$, $i\in I$, are the canonical injections
of the direct sum $\bigoplus_{i\in I}X_i$ and $u'_j$, $j\in J$, are
these of $\bigoplus_{j\in F}X_j$;
\item $f=\sum_{j\in J}u_jp_jf$, where $u_j$, $j\in J$, are the
canonical injections and $p_j$, $j\in J$, are the associated
projections.
\end{enumerate}
\end{lemma}

\begin{proposition}\emph{\cite[Proposition II.16.2]{Mitchell:1965}}
\\The following are equivalent
\begin{enumerate}[(1)]
\item  Every morphism $f:X\to \bigoplus_{i\in I}X_i$ in
$\cat{C}$ has the factorization \eqref{small diag} for some nonempty
finite subset $F$ of $I$;
\item the functor $h'_X=\hom{\cat{C}}{X}{-}:\cat{C}\to\textbf{Ab}$
preserves direct sums.
\end{enumerate}
\end{proposition}

$X$ is called a \emph{small object} if the conditions of the above
proposition hold.

\begin{proposition}\emph{\cite[Exercise V.12, p.
134]{Stenstrom:1975}}\label{small}
\\Let $C$ be a Grothendieck category. The following are equivalent
\begin{enumerate}[(1)]
\item $X$ is small;
\item the functor $h'_X=\hom{\cat{C}}{X}{-}$
preserves denumerable direct sums;
\item the functor $h'_X=\hom{\cat{C}}{X}{-}$
preserves denumerable direct unions;
\item if $X=\sum X_i$ where $X_1\subset X_2\subset\dots$ is a
denumerable ascending chain of subobjects of $X$, then $X=X_n$ for
some $n$.
\end{enumerate}
\end{proposition}

In \cite{Stenstrom:1975} the terminology ``$\sum$-generated
object'' was used for small object.

Let $C$ be a Grothendieck category. From Proposition \ref{small},
every finitely generated of $\cat{C}$ is small. Following \cite[p.
76]{Nastasescu/VanOystaeyen:2004}, we say that $\cat{C}$ is a
\emph{steady category} if every small object of $\cat{C}$ is
finitely generated. For instance, every locally noetherian category
is a steady category (see \cite[Exercice 2.12.22, p.
76]{Nastasescu/VanOystaeyen:2004}).

\medskip
\textbf{Ring action on a $k$-category.} To give a $k$-category with
only one object is the same as to give a $k$-algebra. Let $A$ be a
$k$-algebra and $\cat{C}$ a $k$-category. Consider $A$ as the
category with only one object $*$. To give a $k$-functor $T:A\to
\cat{C}$ is the same as to give an object $T(*)=C$ and a $k$-algebra
morphism
$$\rho:A=\hom{A}{*}{*}\to \hom{\cat{C}}{T(*)}{T(*)}=\hom{\cat{C}}{C}{C}.$$
The functor $T:A\to \cat{C}$ is called a left $A$-\emph{object} in
$\cat{C}$ and is denoted by $(C,\rho)$. The category of left
$A$-objects in $\cat{C}$ is the category of $k$-functors
$\textbf{Hom}_k(A,\cat{C})$ and it is denoted by $_A\cat{C}$. Thus,
a morphism $(C,\rho)\to (C',\rho')$ in $_A\cat{C}$ is a morphism $
f:C\to C'$ in $\cat{C}$ such that, for all $a\in A$, the following
diagram is commutative
$$\xymatrix{C\ar[r]^f\ar[d]_{\rho(a)} & C'\ar[d]^{\rho'(a)}
\\ C\ar[r]^f & C'.}$$
From Proposition \ref{DC2} and Example \ref{DC examples} (3), if
$\cat{C}$ is a $k$-category (resp. $k$-abelian category), then so is
$_A\cat{C}$.

Dually, a right $A$-\emph{object} is a $k$-functor $A^\circ\to
\cat{C}$. The category of right $A$-objects in $\cat{C}$ is the
category of $k$-functors $\textbf{Hom}_k(A^\circ,\cat{C})$ and is
denoted by $\cat{C}_A$. Obviously, $\cat{C}_A\simeq
{}_{A^\circ}\cat{C}$.

If $\cat{C}=\rmod{k}$, then $_A\cat{C}=\lmod{A}$ and
$\cat{C}_A=\rmod{A}.$

Let $A$ be $k$-algebra and $F:\cat{C}\to \cat{D}$ a $k$-functor.
If $C\in{}_A\cat{C}$, then $F(C)\in{}_A\cat{D}$. Hence, we obtain a
functor $_AF=F:{}_A\cat{C}\to {}_A\cat{D}$. Notice that $_AF$ is
nothing else than the functor $F^\Sigma$ defined in ``Diagram
categories''. From Proposition \ref{can ext fun diag}(c), if $F$
is left (resp. right) exact, then the functor $_AF$ is so. If $F$
is a contravariant functor, we obtain a functor
$F:{}_A\cat{C}\to \cat{D}_A$. For a $k$-category $\cat{C}$ and a
$k$-algebra morphism $A\to B$, we obtain a functor
$_A(-):{}_B\cat{C}\to {}_A\cat{C}$ (the restriction functor).

\begin{lemma}\label{$A$-objects}
Let $F,G:\cat{C}\to \cat{D}$ be $k$-functors and $\eta:F\to G$ a
natural transformation. Then $_A\eta=\eta:{}_AF\to {}_AG$ is a
natural transformation. In particular, if $C\in{}_A\cat{C}$ then
$\eta_C:F(C)\to G(C)$ is a morphism in $_A\cat{D}$.
\end{lemma}

\begin{proof}
Obvious. (This lemma is also an immediate consequence from the
definition of the natural transformation $\eta^\Sigma$, see
``Diagram categories''.)
\end{proof}

Now, let $T:\cat{C}\times\cat{D}\to\cat{E}$ be a $k$-linear
bifunctor (i.e. $k$-linear in each variable) and $(C,\rho)$ is a
left $A$-object in $\cat{C}$, then $T(C,D)$ is a left $A$-object in
$\cat{E}$ for all $D\in \cat{D}$. Moreover, it is easy to verify
that, for a morphism $f:C\to C'$ in $\cat{C}$ and a morphism $g:D\to
D'$ in $\cat{D}$, $T(f,g)$ is a morphism in $_A\cat{E}$. Hence, we
obtain a functor $_A\cat{C}\times\cat{D} \to{} _A\cat{E}$. If $T$ is
contravariant in $\cat{C}$ and covariant in $\cat{D}$, then we
obtain a bifunctor $_A\cat{C}\times\cat{D} \to \cat{E}_A$. A special
case is the following: Let $\cat{C}$ be a $k$-category. The
$k$-linear bifunctor $\hom{\cat{C}}{-}{-}:\cat{C}\times\cat{C}\to
\rmod{k}$ induces a bifunctor
$\hom{\cat{C}}{-}{-}:{}_A\cat{C}\times\cat{C}\to \rmod{A}$.

\begin{proposition}\label{KS 8.5.5}
Let $\cat{C}$ be a $k$-category and $A$ a $k$-algebra.
\begin{enumerate}[(1)]
\item Assume that $\cat{C}$ has direct sums and cokernels. Then
\begin{enumerate}[(i)]
\item Let $X\in {}_A\cat{C}$ and $N\in \rmod{A}$. The functor
$\cat{C}\to \rmod{k}$ defined by $$Y\mapsto
\hom{A}{N}{\hom{\cat{C}}{X}{Y}}$$ is representable. Set
$N\tensor{A}X$ its representative.
\item The functor contravariant
in the first variable and covariant in the second one
$$\hom{\cat{C}}{-}{-}:{}_A\cat{C}\times\cat{C}\to\rmod{A}$$ has a
left adjoint functor $$-\tensor{A}-:\rmod{A}\times{}_A\cat{C}\to
\cat{C}$$ which assigns to each $(N,X)\in
\rmod{A}\times{}_A\cat{C}$, $N\tensor{A}X$. In particular, the
functor $-\tensor{A}-$ is $k$-linear and preserves limits in each
variable.
\end{enumerate}
\item Assume that $\cat{C}$ has products and kernels. Then
\begin{enumerate}[(i)]
\item let $X\in {}_A\cat{C}$ and $M\in \lmod{A}$. the functor
$\cat{C}\to \rmod{k}$ defined by $$Y\mapsto
\hom{A}{M}{\hom{\cat{C}}{Y}{X}}$$ is representable. Set
$\hom{A}{M}{X}$ its representative.
\item There is a functor
$$\hom{A}{-}{-}:(\lmod{A})^\circ\times{}_A\cat{C}\to \cat{C}$$
which is $k$-linear and preserves colimits in each variable, and
assigns to each $(M,X)\in (\lmod{A})^\circ\times{}_A\cat{C}$,
$\hom{A}{M}{X}$, such that the isomorphism
$$\hom{A}{M}{\hom{\cat{C}}{Y}{X}}\simeq
\hom{\cat{C}}{Y}{\hom{A}{M}{X}}$$ is natural in $M$, $Y$ and $X$.
\end{enumerate}
\end{enumerate}
\end{proposition}

\begin{proof}
(1) (i) First we consider the particular case where $N=A^{(I)}$ with
$I$ a set. Then there is an isomorphism which is natural in
$Y\in\cat{C}$:
$$\hom{A}{A^{(I)}}{\hom{\cat{C}}{X}{Y}}\simeq\hom{\cat{C}}{X}{Y}^I
\simeq\hom{\cat{C}}{X^{(I)}}{Y}.$$ For the general case, let
$Y\in\cat{C}$ and let
$\xymatrix{A^{(J)}\ar[r]^f & A^{(I)}\ar[r] & N\ar[r] & 0}$ be an
exact sequence with $I$ and $J$ sets. By the Yoneda--Grothendieck
Lemma, there is a morphism $g:X^{(J)}\to X^{(I)}$ in $\cat{C}$ and a
unique isomorphism of $k$-modules
$$\eta_Y:\hom{A}{N}{\hom{\cat{C}}{X}{Y}}\to
\hom{\cat{C}}{\operatorname{Coker}(g)}{Y}$$ making commutative the
diagram
$$\xymatrix{0\ar[r] &
\hom{A}{N}{\hom{\cat{C}}{X}{Y}}\ar[r]\ar[d]^{\eta_Y} &
\hom{A}{A^{(I)}}{\hom{\cat{C}}{X}{Y}}\ar[r]\ar[d]^\simeq &
\hom{A}{A^{(J)}}{\hom{\cat{C}}{X}{Y}}\ar[d]^\simeq\\
0\ar[r] & \hom{\cat{C}}{\operatorname{Coker}(g)}{Y}\ar[r] &
\hom{\cat{C}}{X^{(I)}}{Y}\ar[r]^{-\circ g} &
\hom{\cat{C}}{X^{(J)}}{Y}.}$$ By Lemma \ref{cube}, $\eta_Y$ is
natural in $Y\in\cat{C}$.

(ii) Follows immediately using Proposition \ref{Gabriel I.10}, the
dual of Proposition \ref{Kan 4.1}, and Propositions \ref{adj
fun-lim}, \ref{M IV.3}.

(2) Apply (1) to $\cat{C}^\circ$ and $A^\circ$.
\end{proof}

\begin{remark}
If $\cat{C}$ has finite direct sums and cokernels (resp. $\cat{C}$
has finite direct sums and kernels) and $N$ (resp. $M$) is finitely
presented, then, by a slightly different proof, we obtain that
$N\tensor{A}X$ (resp. $\hom{A}{M}{X}$) is well defined.
\end{remark}

Now, we will recall a well known result of Gabriel--Popescu
(\cite[Theorem 3.7.9]{Popescu:1973}). Let $\cat{C}$ be a
Grothendieck category and $U\in\cat{C}$. Set
$A=\operatorname{End}_\cat{C}(U)$. Then $U\in {}_A\cat{C}$. Consider
the functors
$$T=\hom{\cat{C}}{U}{-}:\cat{C}\to \rmod{A},
\quad S=-\tensor{A}U:\rmod{A}\to \cat{C}.$$ We have that $(S,T)$ is
an adjoint pair.

\begin{theorem}\emph{[Gabriel--Popescu]}\label{Gabriel-Popescu}
\\ Let $\cat{C}$ be a Grothendieck category and $U\in\cat{C}$. The
assertions below are equivalent
\begin{enumerate}[(1)]
\item $U$ is a generator of $\cat{C}$;
\item $T$ is fully faithful and $S$ is
exact.
\end{enumerate}
\end{theorem}

For an alternative proof of it see that of \cite[Theorem
8.5.8]{Kashiwara/Schapira:2006}.

\medskip
The following result is contained in \cite[Proposition
3]{Harada:1968}.

\begin{theorem}\emph{[Harada]}\label{Harada}
\\ Let $\cat{C}$ be a Grothendieck category and $U\in\cat{C}$.
The assertions below are equivalent
\begin{enumerate}[(1)]
\item $S$ is an equivalence of categories;
\item $U$ is a generator, projective and small in $\cat{C}$.
\end{enumerate}
\end{theorem}

\textbf{Some Hom-tensor relations.} Let $R$ and $S$ be rings.
\begin{lemma}\label{hom-tens-relations}
\begin{enumerate}[(1)]
\item Let $\cat{C}$ and $\cat{D}$ be abelian categories, $F$ and $G$
be covariant (contravariant) additive functors from $\cat{C}$
to $\cat{D}$, and let $\eta:F\to G$ be a natural transformation.
\begin{enumerate}[(a)]
\item If $$0\to C'\to C\to C''\to 0$$ is split exact in $\cat{C}$,
then $\eta_C$ is a monomorphism (resp. an epimorphism, resp.
an isomorphism) if and only if $\eta_{C'}$ and $\eta_{C''}$ are
monomorphisms (resp. epimorphisms, resp. isomorphisms).
\item If $C_1,\ldots,C_n\in\cat{C}$, then $\eta_{\oplus_{i=1}^nC_i}$
is a monomorphism (resp. an epimorphism, resp. an isomorphism) if and
only if $\eta_{C_1},\ldots,\eta_{C_n}$ are monomorphisms (resp.
epimorphisms, resp. isomorphisms).

In particular, if $\cat{C}$ and $\cat{D}$ are the categories of left
or right modules over $R$ and $S$, respectively, and if $\eta_R$ is
a monomorphism (resp. epimorphism, resp. isomorphism), then so is
$\eta_P$ for every finitely generated projective $R$-module $P$.
\end{enumerate}
\item Let $F$ and $G$ be left exact contravariant (additive)
functors from $\lmod{R}$ to $\textbf{Ab}$, and let $\eta:F\to G$ be
a natural transformation. If $\eta_R:F(R)\to G(R)$ is a monomorphism
(resp. isomorphism), then $\eta_M:F(M)\to G(M)$ is also a
monomorphism (resp. isomorphism) for every finitely generated (resp.
finitely presented) module $_RM$.

The version of this statement for right exact covariant
(additive) functors is obtained by replacing monomorphism
by epimorphism.
\item In the situation $(_SL,{}_SU_R,{}_RQ)$, there is a
homomorphism natural in $L,U$ and $Q$,
$$\eta:\hom{S}{L}{U}\tensor{R}Q\to \hom{S}{L}{U\tensor{R}Q}$$
defined by
$$\eta(\gamma\otimes q):l\mapsto \gamma(l)\otimes q.$$
\begin{enumerate}[(i)]
\item If $_SL$ is finitely generated and $_RQ$ is flat,
then $\eta$ is a monomorphism.
\item If
\begin{enumerate}[(a)]
\item $_SL$ is finitely generated and projective, or
\item $_SL$ is finitely presented and $_RQ$ is flat, or
\item $_RQ$ is finitely generated and projective,
\end{enumerate}
then $\eta$ is an isomorphism.
\end{enumerate}
\item In the situation $(L_S,{}_RU_S,Q_R)$, there is a
homomorphism natural in $L,U$ and $Q$,
$$\eta:Q\tensor{R}\hom{S}{L}{U}\to \hom{S}{L}{Q\tensor{R}U}$$
defined by
$$\eta(q\otimes \gamma):l\mapsto q\otimes\gamma(l).$$
\begin{enumerate}[(i)]
\item If $L_S$ is finitely generated and $Q_R$ is flat,
then $\eta$ is a monomorphism.
\item If
\begin{enumerate}[(a)]
\item $L_S$ is finitely generated and projective, or
\item $L_S$ is finitely presented and $Q_R$ is flat, or
\item $Q_R$ is finitely generated and projective,
\end{enumerate}
then $\eta$ is an isomorphism.
\end{enumerate}
\item Let $M$ and $N$ be left $A$-modules. If
\begin{enumerate}[(a)]
\item $_AN$ is finitely generated and projective, or
\item $_AN$ is finitely presented and $_AM$ is flat, or
\item $_AM$ is finitely generated and projective,
\end{enumerate}
then the map
$$\theta_{N,M}:\rdual{N}\tensor{A}M\to \hom{A}{N}{M}$$
defined by
$$\theta_{N,M}(\gamma\otimes m):n\mapsto \gamma(n)m,$$
is an isomorphism which is natural in $M$ and $N$.
\item Let $M$ and $N$ be right $A$-modules. If
\begin{enumerate}[(a)]
\item $N_A$ is finitely generated and projective, or
\item $N_A$ is finitely presented and $M_A$ is flat, or
\item $M_A$ is finitely generated and projective,
\end{enumerate}
then the map
$$\theta_{N,M}:M\tensor{A}\rdual{N}\to \hom{A}{N}{M}$$
defined by
$$\theta_{N,M}(m\otimes\gamma):n\mapsto m\gamma(n),$$
is an isomorphism which is natural in $M$ and $N$.
\item In the situation $(L_S,{}_RU_S,{}_RQ)$, there is a
homomorphism natural in $L,U$ and $Q$,
$$\nu:\hom{S}{U}{L}\tensor{R}Q\to
\hom{S}{\hom{R}{Q}{U}}{L}$$
defined by
$$\nu(\gamma\otimes q):\alpha\mapsto \gamma(\alpha(q)).$$
\begin{enumerate}[(i)]
\item If $_RQ$ is finitely generated and $L_S$ is injective,
then $\nu$ is an epimorphism.
\item If
\begin{enumerate}[(a)]
\item $_RQ$ is finitely generated and projective, or
\item $_RQ$ is finitely presented and $L_S$ is injective,
\end{enumerate}
then $\nu$ is an isomorphism.
\end{enumerate}
\item In the situation $(_SL,{}_SU_R,Q_R)$, there is a
homomorphism natural in $L,U$ and $Q$,
$$\nu:Q\tensor{R}\hom{S}{U}{L}\to
\hom{S}{\hom{R}{Q}{U}}{L}$$
defined by
$$\nu(q\otimes \gamma):\alpha\mapsto \gamma(\alpha(q)).$$
\begin{enumerate}[(i)]
\item If $Q_R$ is finitely generated and $_SL$ is injective,
then $\nu$ is an epimorphism.
\item If
\begin{enumerate}[(a)]
\item $Q_R$ is finitely generated and projective, or
\item $Q_R$ is finitely presented and $_SL$ is injective,
\end{enumerate}
then $\nu$ is an isomorphism.
\end{enumerate}
\end{enumerate}
 \end{lemma}

\begin{proof}
(1) (a) From Proposition \ref{AF 5.3}, there are two split
exact short sequences
$$0\to C'\to C\to C''\to 0$$ and $$0\to C''\to C\to C'\to 0.$$
By Corollary \ref{split exactness}, we obtain then, in both cases
(the ``covariant'' case and the ``contravariant'' one), two
diagrams with split exact rows
$$\xymatrix{0\ar[r]& F(C')\ar[r]\ar[d]^{\eta_{C'}}
& F(C) \ar[r]\ar[d]^{\eta_C}& F(C'')\ar[r]\ar[d]^{\eta_{C''}}& 0
\\ 0\ar[r]& G(C')\ar[r]
& G(C) \ar[r]& G(C'')\ar[r]& 0}$$ and
$$\xymatrix{0\ar[r]& F(C'')\ar[r]\ar[d]^{\eta_{C''}}
& F(C) \ar[r]\ar[d]^{\eta_C}& F(C')\ar[r]\ar[d]^{\eta_{C'}}& 0
\\ 0\ar[r]& G(C'')\ar[r]
& G(C) \ar[r]& G(C')\ar[r]& 0.}$$

The short five lemma achieves then the proof.

(b) It follows by induction using (a).

(2) Let $_RM$ be a finitely generated (resp. finitely presented)
module. Then there is an exact sequence $R^{(I)}\to R^{n}\to M\to
0$, where $I$ is a (resp. a finite) nonempty index set and $n\in
\mathbb{N}-\{0\}$. Consider the commutative diagram with exact rows
$$\xymatrix{0\ar[r]& F(M)\ar[r]\ar[d]^{\eta_M}
& F(R^n) \ar[r]\ar[d]^{\eta_{R^n}}&
F(R^{(I)})\ar[d]^{\eta_{R^{(I)}}}
\\ 0\ar[r]& G(M)\ar[r]
& G(R^n) \ar[r]& G(R^{(I)}).}$$ By (1), $\eta_{R^n}$ is a
monomorphism (resp. $\eta_{R^n}$ and $\eta_{R^{(I)}}$ are
isomorphisms). Hence $\eta_M$ is also a monomorphism (resp. an
isomorphism).

Analogously we prove the second statement.

(3), (4), (7) and (8) are straightforward from (1) and (2). (5) is
obvious from (3) while (6) is obvious from (4).
\end{proof}

\medskip
\textbf{Separable bimodules.} Separable bimodules are introduced by
Sugano in \cite{Sugano:1971}. Let $_AM_B$ be a bimodule. Define the
evaluation map $ev:M\tensor{B}{}^*M\to A$ by $m\otimes f\mapsto
f(m)$ which is an $A$-bimodule morphism. $M$ is \emph{separable} or
$A$ \emph{is} $M$-\emph{separable over} $B$ if $ev$ is a retraction
in $\bmod{A}{A}$.

A \emph{ring extension} $A/S$ is a ring morphism $\iota:S\to A$. A
ring extension $A/S$ is called a \emph{separable extention} if
$$\mu_S:A\tensor{S}A\to A,\;a\otimes b\mapsto ab,$$
is a retraction in $\bmod{A}{A}$.

Let $_AM_B$ be a bimodule. We have, $E=\operatorname{End}_B(M)$
(resp. $E'=\operatorname{End}_A(M)$) is a ring extension over $A$
(resp. $B$) via left (resp. right) multiplication.

\begin{proposition}\emph{\cite[Theorem 1, Proposition 2]{Sugano:1971}
(see also \cite[Theorem 3.1]{Kadison:1999})}\label{Sugano}
\\Let $_AM_B$ be a bimodule.
\begin{enumerate}[1)]
\item If $M$ is separable then $E/A$ is a separable extension.
\item If $M_B$ is finitely generated projective and $E/A$ is a separable
extension, then $M$ is separable.
\item If $_AM$ is finitely generated projective and $M$ is separable, then
$E'/B$ is a separable extension.
\item If $_AM$ is a generator and $E'/B$ is a separable extension,
then $M$ is separable.
\end{enumerate}
\end{proposition}

\medskip
\textbf{Locally projective modules.} Let $R$ be a ring.
Following Zimmermann-Huisgen \cite{Zimmermann-Huisgen:1976},
a right $R$-module $M$ is called \emph{locally projective} if
for every diagram in $\rmod{A}$ with exact rows:
$$\xymatrix{0\ar[r] & F \ar[r]^i & M \ar[d]^g \ar@{-->}[dl]_h &
\\ & L \ar[r]^f & N \ar[r] & 0,}$$
where $F$ is finitely generated, there exists $h:M\to L$ with
$gi=fhi$. Obviously, every projective module is locally projective.

For every $M\in\rmod{R}$ and $N\in \lmod{R}$, let
$$\alpha_{M,N}:M\tensor{R}N\to \hom{R}{\rdual{M}}{N}$$
be the map defined by $\alpha_{M,N}(m\otimes n)(f)=f(m)n$,
for every $m\in M,n\in N,f\in \rdual{M}.$

Following Bass, $M$ is called \emph{torsionless} if $\alpha_{M,R}$
is injective and \emph{reflexive} if $\alpha_{M,R}$ is an isomorphism.

Following Garfinkel \cite{Garfinkel:1976}, $M$ is called
\emph{universally torsionless} (UTL) if $\alpha_{M,N}$ is
injective for all $N\in\lmod{R}$ and \emph{universally ring
torsionless} (URTL) if $M\tensor{R}S$ is $S$-torsionless for all
ring extensions $S/R$.

\begin{proposition}\emph{\cite[Theorem 2.2]{Garfinkel:1976}}
\\Let $M$ be a module over a commutative ring $R$. The following
are equivalent
\begin{enumerate}[(1)]
\item $M$ is UTL;
\item $M$ is URTL;
\item $M\tensor{R}S$ is $S$-torsionless for all
commutative ring extensions $S/R$.
\end{enumerate}
\end{proposition}

\begin{proposition}\emph{\cite[Proposition 2.4]{Garfinkel:1976}}
\\Let $M$ be a right $R$-module. The following are
equivalent
\begin{enumerate}[(1)]
\item $M$ is UTL;
\item $M$ is URTL and flat.
\end{enumerate}
\end{proposition}

Let $\beta:G\to M$ be a map in $\rmod{R}$ and $n\in\mathbb{N}-\{0\}$.
Let us consider the property $P(n)$: for all $x_1,\dots,x_n\in G$,
there is a map $\varphi:M\to G$ in $\rmod{R}$ such that
$\beta\varphi\beta(x_i)=\beta(x_i)$ for $i=1,\dots,n$.

By Lemma \cite[Lemma 3.1]{Garfinkel:1976}, we have that, for every $n$,
$P(1)$ and $P(n)$ are equivalent.

A submodule $K$ of $M_R$ is called \emph{ideal pure} if for every
left ideal $I$ of $R$, $KI=K\cap MI$. Following
\cite[p. 124]{Garfinkel:1976}, a map $\beta:G\to M$ in $\rmod{R}$
is said to be \emph{split} if there a map $\varphi:M\to G$ in $\rmod{R}$
such that $\beta=\beta\varphi\beta$, \emph{finitely split}
if it satisfies $P(1)$, and \emph{(ideal) pure} if its image is (ideal)
pure. Note that, if the image of $\beta$ is finitely generated, then
it is split if and only if it is finitely split.

From \cite[Theorem 2.6]{Garfinkel:1976} and
\cite[Theorem 2.1]{Zimmermann-Huisgen:1976} we obtain the following
proposition. \cite[Theorem 2.1]{Zimmermann-Huisgen:1976} provides more
equivalent conditions to the (1) to (11) ones.

\begin{proposition}
Let $M$ be a right $R$-module. The following are
equivalent
\begin{enumerate}[(1)]
\item $M$ is locally projective;
\item $M$ is UTL;
\item $M$ is a (ideal) pure submodule of a right UTL module $P$;
\item $\alpha_{M,N}$ is injective for every cyclic left $R$-module $N$;
\item every $m\in M$ belongs to $M.\rdual M(m)$ with
$\rdual M(m)=\{f(m)\mid f\in \rdual M\}$;
\item for every $m\in M$, there are $m_1,\dots,m_n\in M$ and
$f_1,\dots,f_n\in \rdual M$ such that $m=\sum m_if_i(m)$;
\item every ideal pure map $\beta:G\to M$ finitely splits;
\item every epimorphism $\beta:G\to M$ finitely splits;
\item for every finitely generated submodule $K$ of $M$ there is a
finitely generated free right module $F$ and maps $\alpha:F\to M$ and
$\gamma:M\to F$ such that $\alpha\gamma$ is the identity on $K$;
\item for every finitely generated submodule $K$ of $M$ there are
$m_1,\dots,m_n\in M$ and $f_1,\dots,f_n\in \rdual M$ such that
$x=\sum m_if_i(x)$ for $x\in K$;
\item for every epimorphism $\beta:G\to M$ the induced map
$$\hom{R}{C}{\beta}:\hom{R}{C}{G}\to\hom{R}{C}{M}$$ is an epimorphism
for every finitely generated module $C$ such that for every $m\in M$
there is an $\alpha\in \operatorname{End}_R(M)$ with $\alpha(m)=m$
and the image of $\alpha$ is contained in a finitely generated
submodule of $M$.
\end{enumerate}
\end{proposition}

\begin{corollary}\emph{\cite[p. 126]{Garfinkel:1976}}
\begin{enumerate}[(a)]
\item Let $M_R$ be UTL and $S/R$ a ring extension. Then the right $S$-module
$M\tensor{R}S$ is also UTL.
\item Each projective module is UTL.
\item Any finitely generated UTL module is projective.
\item Any finitely generated pure submodule of a UTL module is a
projective direct summand.
\item Any countably generated submodule of a UTL module $M$ is
contained in a countably generated pure projective submodule.
\end{enumerate}
\end{corollary}

We refer to \cite{Garfinkel:1976,Zimmermann-Huisgen:1976}
for more details.

\medskip
The following example has been communicated to me by Edgard Enochs.

\medskip
\textbf{An example of a self-injective commutative ring which is not
coherent.} Let $(R_i)_{i\in I}$ be any family of rings indexed by
the set $I$ and let $R=\prod_{i\in I} R_i$. Then the following are
easy to prove. If $I$ is a finitely generated ideal of $R$ such that
$I$ is generated by $n\geq 1$ elements, then for each $i\in I$ there
are ideals $I_i \in R_i$ such that $I_i$ is generated by $n$
elements and such that $I=\prod I_i$. Conversely, given an $n\geq 1$
and such ideal $I_i$ for each $i\in I$ we get that $I=\prod I_i$ is
a finitely generated ideal of $R$ and moreover that $I$ is generated
by $n$ elements.

  The above can be generalized to finitely generated submodules of $R^m$
for any $m\geq 1$.

  Now using the above facts it is not hard to see that if $R$ is coherent
 the following two conditions hold:
\begin{enumerate}[(a)]\item each $R_i$ is coherent,
\item for each $i$ and $n\geq 1$ if $I_i$ is an ideal of $R_i$
generated by $n$ elements and if $R_i^n\rightarrow I_i$ is
surjective and linear then there is some $m\geq 1$ such that each of
these kernels is generated by $m$ elements.
\end{enumerate}

Now for an example, let $I=P$ ($P$ the set of positive natural
numbers) and let $R_n=k[[x_1,\dots ,x_n]]$ modulo the ideal
generated by $x_1^2,\dots ,x_n^2$ Then each $R_n$ is artinian and so
coherent. But now look at \cite[Proposition 21.5 (p. 530), Corollary
21.19]{Eisenbud:1995}. Then we see that each $R_n$ has a unique
minimal ideal $I_n$. But it is not hard to see that if
$R_n\rightarrow I_n$ is surjective then the kernel (which is the
ideal generated by the cosets of $x_1,\dots , x_n$) cannot be
generated by fewer than $n$ elements. So in this case, the product
ring $R$ is not coherent.

By the results in \cite{Eisenbud:1995} mentioned above, we have that
each $R_n$ is self-injective (in his language, Gorenstein of
dimension 0). Finally, from \cite[Corollary 3.11B]{Lam:1999}, the
product of self-injective rings is also self-injective. So this is
an example of a commutative ring $R$ that is self-injective but
which is not coherent.

\section{The category of comodules}

We recall from \cite{Sweedler:1975} the definition of a coring.
\begin{definition}
An $A$-\emph{coring} is an $A$-bimodule $\coring{C}$ with two
$A$-bimodule maps $\Delta : \coring{C} \rightarrow \coring{C}
\tensor{A} \coring{C}$ (\emph{coproduct} or \emph{comultiplication})
and $\epsilon:\coring{C} \rightarrow A$ (\emph{counit}) such that
$(\coring{C}\tensor{A} \Delta) \circ \Delta = (\Delta \tensor{A}
\coring{C})\circ \Delta$ (\emph{coassociative property}) and
$(\epsilon \tensor{A}\coring{C}) \circ \Delta = ( \coring{C}
\tensor{A} \epsilon) \circ\Delta = 1_\coring{C}$ (\emph{counit
property}), that is the following diagrams are commutative
$$\xymatrix{\coring{C} \ar[rr]^\Delta \ar[d]^\Delta &&
\coring{C}\tensor{A}\coring{C}\ar[d]^{\coring{C} \tensor{A} \Delta} \\
\coring{C}\tensor{A}\coring{C}\ar[rr]^{\Delta\tensor{A}\coring{C}}
&& \coring{C}\tensor{A}\coring{C}\tensor{A}\coring{C}}
\xymatrix{\coring{C} \ar[r]^-\Delta \ar[dr]^\simeq &
\coring{C}\tensor{A}\coring{C}\ar[d]^{\epsilon\tensor{A}\coring{C}}
\\ & A\tensor{A}\coring{C}}
\xymatrix{\coring{C} \ar[r]^-\Delta \ar[dr]^\simeq &
\coring{C}\tensor{A}\coring{C}\ar[d]^{\coring{C}\tensor{A}\epsilon}
\\ & \coring{C}\tensor{A}A}.$$
\end{definition}

As for coalgebras, we will use the Sweedler's sigma notation, for
every $c\in \coring{C}$, we write
$$\Delta(c)=\sum c_{(1)}\tensor{A}c_{(2)}.$$

Then, as for coalgebras, the coassociative property can be expressed
by $$\sum\Delta(c_{(1)})\tensor{A}c_{(2)}=\sum
c_{(1)}\tensor{A}\Delta(c_{(2)})=\sum
c_{(1)}\tensor{A}c_{(2)}\tensor{A}c_{(3)},$$ and the counit property
can be expressed by
$$\sum\epsilon(c_{(1)})c_{(2)}=c=\sum c_{(1)}\epsilon(c_{(2)}).$$

Sometimes the symbol $\sum$ will be omitted.

\begin{examples}\label{examples corings}
\begin{enumerate}[(1)]
\item If we take $A=k$ we find the notion of a $k$-coalgebra.
\item $\coring{C}=A$ endowed with the obvious structure maps is an
$A$-coring. This coring is called \emph{the trivial coring}.
\item Let $\rho:A\to B$ be a morphism of $k$-algebras, and
$\coring{C}$ be an $A$-coring.
$B\coring{C}B=B\tensor{A}\coring{C}\tensor{A}B$ is a $B$-coring with
coproduct and counit:
$$\Delta_{B\coring{C}B}:=B\tensor{A}\upsilon\tensor{A}B:
B\tensor{A}\coring{C}\tensor{A}B\to
B\tensor{A}\coring{C}\tensor{A}B\tensor{A}\coring{C}\tensor{A}B,$$
$$\epsilon_{B\coring{C}B}:B\tensor{A}\coring{C}\tensor{A}B\to B,
\quad b\tensor{A}c\tensor{A}b'\mapsto b\rho(\epsilon(c))b',$$ where
$$\upsilon:\coring{C}\to \coring{C}\tensor{A}B\tensor{A}\coring{C},\quad
c\mapsto \sum c_{(1)}\tensor{A}1_B\tensor{A}c_{(2)}.$$ This coring
is called \emph{a base ring extension} of $\coring{C}$.
\item Let $X$ be a nonempty set. $AX=A^{(X)}$ is an $A$-coring with
coproduct and counit:
$$\Delta(x)=x\tensor{A}x,\quad\epsilon(x)=1_A,$$ for every $x\in X$.
This coring is called the \emph{grouplike coring on $X$.} If we take
$A=k$ we find the definition of a grouplike coalgebra
(cocommutative).
\item Let $M$ be a $B-A$-bimodule with $M_A$ is finitely generated
projective. Let $\{e_i,e_i^*\}_{i\in I}$ be a finite dual basis of
$M$. Define $\tau:B\to M\tensor{A}M^*, \quad b\mapsto\sum_{i\in
I}be_i\otimes e_i^*= \sum_{i\in I}e_i\otimes e_i^*b,$ and
$\epsilon_{M^*\tensor{B}M}:M^*\tensor{B}M\to B, \quad \varphi\otimes
m\mapsto\varphi(m).$

Then we have
$(A,B,{}_AM^*_B,{}_BM_A,\epsilon_{M^*\tensor{B}M},\tau)$ is a
comatrix coring context, i.e. $\epsilon_{M^*\tensor{B}M}$ and $\tau$
are bimodule maps, and the following diagrams are commutative
$$\xymatrix{
M^*\tensor{B}M\tensor{A}M^*\ar[rrr]^-{\epsilon_{M^*\tensor{B}M}\tensor{A}M^*}
&&&  A\tensor{A}M^*\ar[d]^{\simeq}\\
 M^*\tensor{B}B\ar[u]^{M^*\tensor{B}B}\ar[rrr]^-{\simeq}&&& N}
\xymatrix{
M\tensor{A}M^*\tensor{B}M\ar[rrr]^-{M\tensor{A}\epsilon_{M^*\tensor{B}M}}
&&&  M\tensor{A}A\ar[d]^{\simeq}\\
 B\tensor{B}M\ar[u]^{\tau\tensor{B}M}\ar[rrr]^-{\simeq}&&& M}.$$
 By \cite[Proposition 2.1]{ElKaoutit/Gomez:2003} or
 \cite[Theorem 2.4(2)]{Brzezinski/Gomez:2003},
 $M^*\tensor{B}M$ is an $A$-coring with coproduct
 $$\Delta_{M^*\tensor{B}M}:M^*\tensor{B}M\to M^*\tensor{B}M\tensor{A}
 M^*\tensor{B}M, \quad \varphi\otimes m\mapsto \sum_{i\in I}\varphi
 \otimes e_i\otimes e_i^*\otimes m,$$ and counit $\epsilon_{M^*\tensor{B}M}.$
 This coring is called the \emph{comatrix coring} associated to the
 bimodule $M$. We propose a generalization of this coring in
 Chapter 3. In fact, we consider, under certain suitable conditions, a coring
 associated to a quasi-finite comodule over a coring.
\item We take in $(5)$ $M=A^n$ for $n\in \mathbb{N}-\{0\}$. $M^*\tensor{A}M$
 can be identified with the ring of all $n\times n$ matrices with
 entries in $A$, $M_n(A)$. Let $\{e_{i,j}\}_{1\leq i,j\leq n}$ be the
 canonical $A$-basis of $M_n(A)$. Then $M_n(A)$ is an $A$-coring with
 coproduct and counit defined as follows
 $$\Delta(e_{i,j})=\sum_ke_{i,k}\tensor{A}e_{k,j}, \quad \epsilon(e_{i,j})=
 \delta_{i,j}.$$ This coring is called the $(n,n)$-\emph{matrix coring} over
 $A$ and it is denoted by $M_n^c(A)$.
\end{enumerate}
\end{examples}

Following \cite{Sweedler:1975}, for an $A$-coring, $\coring{C}$, we
define:
$$\coring{C}^*=\hom{-A}{\coring{C}}{A},
\quad^*\coring{C}=\hom{A-}{\coring{C}}{A},\quad
^*\coring{C}^*=\hom{(A,A)}{\coring{C}}{A}=\coring{C}^*\cap{}^*\coring{C}.$$

The following result that gives the elementary properties of these
$k$-modules is essentially due to Sweedler, see \cite[Proposition
3.2]{Sweedler:1975} (see also
\cite[17.8]{Brzezinski/Wisbauer:2003}).

\begin{proposition}\label{BW 17.8}
\begin{enumerate}[(1)]
\item $\coring{C}^*$ is a $k$-algebra with product
$f*^rg=g(f\tensor{A}\coring{C})\Delta$, and unit $\epsilon$. More
explicitly, $f*^rg(c)=\sum g(f(c_{(1)})c_{(2)}),$ for every $c\in
\coring{C}$. Moreover the map $i_R:A\to \coring{C}^*,\quad a\mapsto
\epsilon(a-),$ is an anti-morphism of $k$-algebras.
\item $^*\coring{C}$ is a $k$-algebra with product
$f*^lg=f(\coring{C}\tensor{A}g)\Delta$, and unit $\epsilon$. More
explicitly, $f*^lg(c)=\sum f(c_{(1)}g(c_{(2)})),$ for every $c\in
\coring{C}$. Moreover the map $i_L:A\to {}^*\coring{C},\quad
a\mapsto \epsilon(-a),$ is an anti-morphism of $k$-algebras.
\item $^*\coring{C}^*$ is a $k$-algebra with product
$f*g(c)=\sum f(c_{(1)})g(c_{(2)}),$ for every $f,g\in
{}^*\coring{C}^*,c\in \coring{C}$, and unit $\epsilon$. Moreover the
map $i:\oper{Z}{A}\to \oper{Z}{^*\coring{C}^*},\quad a\mapsto
\epsilon(a-)=\epsilon(-a),$ is a morphism of $k$-algebras.
($\oper{Z}{A}$ and $\oper{Z}{^*\coring{C}^*}$ denote the center of
$A$ and $^*\coring{C}^*$ respectively.)
\item $\oper{Z}{\coring{C}^*}\subset \oper{Z}{^*\coring{C}^*}$
and $\oper{Z}{^*\coring{C}}\subset\oper{Z}{^*\coring{C}^*}$.
 \item If $A=k$ ($\coring{C}$ is a $k$-coalgebra), then $\coring{C}^*=
 {}^*\coring{C}={}^*\coring{C}^*.$
\item If $\coring{C}$ is the trivial $A$-coring, then
$\coring{C}^*\simeq A^\circ$ (the opposite algebra of $A$) and
$^*\coring{C}\simeq A$.
\end{enumerate}
\end{proposition}

\begin{definitions}
A right $\coring{C}$-\emph{comodule} is a pair $(M,\rho_M)$
consisting of a right $A$-module $M$ and an $A$-linear map $\rho_M:
M \rightarrow M \tensor{A}\coring{C}$ (right coaction of
$\coring{C}$ on $M$) satisfying $(M \tensor{A} \Delta) \circ \rho_M
= (\rho_M \tensor{A} \coring{C}) \circ \rho_M$ and $(M \tensor{A}
\epsilon) \circ \rho_M = 1_M$. That means that the following
diagrams are commutative
\begin{equation}\label{comod}
\xymatrix{M \ar[rr]^{\rho_M} \ar[d]^{\rho_M} &&
M\tensor{A}\coring{C}\ar[d]^{M \tensor{A}\Delta} \\
M\tensor{A}\coring{C}\ar[rr]^{\rho_M\tensor{A}\coring{C}} &&
M\tensor{A}\coring{C}\tensor{A}\coring{C}}
\xymatrix{M\ar[r]^-{\rho_M} \ar[dr]^\simeq &
M\tensor{A}\coring{C}\ar[d]^{M\tensor{A}\epsilon}
\\ & M\tensor{A}A}.
\end{equation}

A \emph{comodule morphism} $f: M \to N$ of right
$\coring{C}$-comodules $(M,\rho_M)$ and $(N,\rho_N)$ is a right
$A$-linear map such that $(f \tensor{A} \coring{C}) \circ \rho_M =
\rho_N \circ f$. That means that the following diagram is
commutative
\begin{equation}\label{comodmor}
\xymatrix{M\ar[rr]^f \ar[d]^{\rho_M} && N \ar[d]^{\rho_N}
\\ M\tensor{A}\coring{C}\ar[rr]^{f\tensor{A}\coring{C}}
&& N\tensor{A}\coring{C}}.
\end{equation}
The set of all such morphisms is a $k$-submodule of $\hom{A}{M}{N}$,
and it will be denoted by $\hom{\coring{C}}{M}{N}$.

Analogously we define a left \emph{comodule} $M$, and a
\emph{comodule morphism} of left comodules. The left coaction of
$\coring{C}$ on $M$, $M\to \coring{C}\tensor{A}M$, is denoted by
$\lambda_M.$
\end{definitions}
Since $\rmod{A}$ is a $k$-category, right comodules over
$\coring{C}$ and their morphisms form a $k$-category, and it is
denoted by $\rcomod{\coring{C}}$. When $\coring{C}$ is the trivial
$A$-coring, $\rcomod{A}$ is the category of right $A$-modules
$\rmod{A}$. Analogously left comodules over $\coring{C}$ and their
morphisms form a $k$-category, and it is denoted by
$\lcomod{\coring{C}}$.

Coproducts and cokernels (and then inductive limits) in
$\rcomod{\coring{C}}$ exist and they coincide respectively with
coproducts and cokernels in the category of right $A$-modules
$\rmod{A}$ (see Proposition \ref{ker}). If ${}_A\coring{C}$ is flat,
then $\rcomod{\coring{C}}$ is a $k$-abelian category (see
Proposition \ref{ker}). Moreover it is a Grothendieck category (see
also Lemma \ref{comodgenerators}).

As for coalgebras, we will use the Sweedler's sigma notation, for
every $m\in M$, we write
$$\rho_M(m)=\sum m_{(0)}\tensor{A}m_{(1)}.$$

Again as for coalgebras, the commutativity of the diagrams
\eqref{comod} can be expressed by
$$\sum\rho_M(m_{(0)})\tensor{A}m_{(1)}=\sum
m_{(0)}\tensor{A}\Delta(m_{(1)})=\sum m_{(0)}
\tensor{A}m_{(1)}\tensor{A}m_{(2)},$$ and
$$\sum m_{(0)}\epsilon(m_{(1)})=m.$$
Also the commutativity of the diagram \eqref{comodmor} can be
expressed by $$\sum f(m)_{(0)}\tensor{A}f(m)_{(1)}=\sum
f(m_{(0)})\tensor{A}m_{(1)}.$$

Sometimes the symbol $\sum$ will be omitted.

\begin{proposition}
\begin{enumerate}[(1)]
\item Let $T$ be a $k$-algebra and $M\in {}_T\rcomod{\coring{C}}$.
If $X\in \rmod{T}$, then $X\tensor{T}M$ is a right
$\coring{C}$-comodule. If $f$ is a morphism in $\rmod{T}$, then
$f\tensor{T}M$ is a morphism in $\rcomod{\coring{C}}$. This yields a
functor $-\tensor{T}M:\rmod{T}\to \rcomod{\coring{C}}$. In
particular, we obtain a functor $-\tensor{A}\coring{C}:\rmod{A}\to
\rcomod{\coring{C}}$.

Right $\coring{C}$-comodules of the type $X\tensor{A}\coring{C}$
where $X\in \rmod{A}$, are called \emph{standard comodules}.
\item The forgetful functor (forgets the coaction)
$U:\rcomod{\coring{C}}\to \rmod{A}$ is a left adjoint to the functor
$-\tensor{A}\coring{C}:\rmod{A}\to \rcomod{\coring{C}}$. In
particular, the functor $-\tensor{A}\coring{C}$ preserves projective
limits.
\item Let $T$ be a $k$-algebra and $M\in {}_T\rcomod{\coring{C}}$.
The functor $-\tensor{T}M:\rmod{T}\to \rcomod{\coring{C}}$ is a left
adjoint to the functor
$\hom{\coring{C}}{M}{-}:\rcomod{\coring{C}}\to\rmod{T}$. In
particular, the functor $-\tensor{T}M$ preserves inductive limits.

Let $I$ be a nonempty index set.
$A^{(I)}\tensor{A}\coring{C}\simeq\coring{C}^{(I)}$ as comodules,
and for $M\in\rcomod{\coring{C}}$, there is a surjective
$\coring{C}$-comodule $\coring{C}^{(I')}\to M\tensor{A}\coring{C}$,
where $I'$ is a set.
\item Let $M\in\rcomod{\coring{C}}$. $M$ is a
direct summand of the right $A$-module $M\tensor{A}\coring{C}$. The
coaction $\rho_M:M\to M\tensor{A}\coring{C}$ is a morphism of right
$\coring{C}$-comodules, and then $M$ is a subcomodule of a standard
comodule.
\end{enumerate}
\end{proposition}

\begin{proof}
(1) Easy verifications.

(2) For every $M\in \rcomod{\coring{C}}$ and $X\in \rmod{A}$, the
$k$-linear map
$$\varphi:\hom{\coring{C}}{M}{X\tensor{A}\coring{C}}\to
\hom{A}{M}{X},\, f\mapsto (X\tensor{A}\epsilon)f,$$ is bijective,
with inverse map $g\mapsto (g\tensor{A}\coring{C})\rho_M$. Moreover,
$\varphi$ is natural in $M$ and $X$.

(3) From the Hom-tensor relations, we obtain that for every $X\in
\rmod{T}$ and $N\in \rcomod{\coring{C}}$, the $k$-linear map
$$\psi:\hom{T}{X}{\hom{A}{M}{N}}\to \hom{A}{X\tensor{T}M}{N},$$
$\psi(f)(x\otimes m)=f(x)(m),$ is bijective, with inverse map
$\psi^{-1}(g)(x)=g(x\otimes-)$. Moreover, $\psi$ is natural in both
$X$ and $N$. It is easy to verify that $\psi$ induces a bijection
$$\psi':\hom{T}{X}{\hom{\coring{C}}{M}{N}}\to
\hom{\coring{C}}{X\tensor{T}M}{N}.$$ Since $\psi$ is natural in $X$
and $N$, $\psi'$ is also natural in $X$ and $N$. This achieves the
proof of the first statement.

The first part of the last statement follows directly from the first
statement. To show the second part, take a surjective $A$-linear map
$f:A^{(M)}\to M$. Then,
$f\tensor{A}\coring{C}:A^{(M)}\tensor{A}\coring{C}\to
M\tensor{A}\coring{C}$ is a surjective comodule morphism, and
$A^{(M)}\tensor{A}\coring{C}\simeq \coring{C}^{(M)}$ as comodules.

(4)  By the counit property, $\rho_M$ is a section in $\rmod{A}$
($(M\tensor{A}\epsilon)\rho_M=M$). In particular, $\rho_M$ is an
injective $A$-linear map. On the other hand, from the
coassociativity property, $\rho_M$ is a comodule morphism.
\end{proof}

\begin{definitions}
An $A'-A$-bimodule $M$ is a
$\coring{C}'-\coring{C}$-\emph{bicomodule} if
\begin{enumerate}[(1)]
\item $M$ is at the same time a right and
a left comodule over $\coring{C}$ and $\coring{C}'$ respectively.
\item $\rho_M$ is $A'$-linear and $\lambda_M$ is $A$-linear.
\item one of the following equivalent statements holds
\begin{enumerate}[(a)]
\item $\rho_M : M \to M \tensor{A}\coring{C}$ is a morphism
of left $\coring{C}'$-comodules;
\item $\lambda_M : M \rightarrow \coring{C}' \tensor{A'} M$
is a morphism of right $\coring{C}$-comodules;
\item the following diagram is commutative
$$\xymatrix{
M\ar[rr]^{\rho_M} \ar[d]^{\lambda_M}&&
M\tensor{A}\coring{C}\ar[d]^{\lambda_M\tensor{A}\coring{C}}
\\\coring{C}'\tensor{A'} M \ar[rr]^-{\coring{C}'\tensor{A'}\rho_M}
&& \coring{C}'\tensor{A'}M\tensor{A}\coring{C}.}$$
\end{enumerate}
\end{enumerate}
A \emph{morphism of bicomodules} is a morphism of right and left
comodules.
\end{definitions}
$\coring{C}'-\coring{C}$-bicomodules with their morphisms form a
$k$-category $\bcomod{\coring{C}'}{\coring{C}}$. Coproducts and
cokernels (and then inductive limits) in this category exist and
they coincide respectively with coproducts and cokernels in the
category of right $A'-A$-modules, $\bmod{A'}{A}$ (see Proposition
\ref{ker}). If $_A\coring{C}$ and $\coring{C'}_{A'}$ are flat,
$\bcomod{\coring{C}'}{\coring{C}}$ is a $k$-abelian category (see
Proposition \ref{ker}). Moreover it is a Grothendieck category (see
also Corollary \ref{tensflat}, Proposition \ref{bcomod}, and Lemma
\ref{comodgenerators}). If $\coring{C}$ and $\coring{C}'$ are the
trivial $A$-coring and $A'$-coring respectively, then
$\bcomod{\coring{C'}}{\coring{C}}$ is exactly the category
$\bmod{A'}{A}$.

\medskip
Let $R$, $S$ and $T$ be $k$-algebras. In the situation
$(M_{R-S},{}_RN_T)$, $M\tensor{R}N$ is a right $S\otimes T$-module
via $(m\otimes n)(s\otimes t)=ms\otimes nt$. The situation $_{R-S}M$
is equivalent to that of $_{R\otimes S}M$ via $(r\otimes
s)m=s(rm)=r(sm)$.

\begin{lemma}
Let $R$, $S$ and $T$ be $k$-algebras. In the situation
$(M_{R-S},{}_RN_T,{}_{S-T}L)$, there is a unique morphism
$$\alpha:(M\tensor{R}N)\tensor{S\otimes T}L\to
M\tensor{R\otimes S}(N\tensor{T}L)$$ such that $\alpha((m\otimes
n)\otimes l)=m\otimes (n\otimes l).$ Moreover, this morphism is an
isomorphism and yields a natural isomorphism of functors.
\end{lemma}

\begin{proof}
Easy verifications.
\end{proof}

\begin{corollary}\label{tensflat}
Let $S$ and $T$ be $k$-algebras. In the situation $(M_S,{}N_T)$, if
$M_S$ and $N_T$ are flat, then $(M\otimes N)_{S\otimes T}$ is also
flat.
\end{corollary}

\begin{lemma}
In the situation $(_AM_A,{}_AN_A,{}_BL_B,{}_BP_B)$, there is a
unique morphism
$$\beta:(M\tensor{A}N)\otimes(L\tensor{B}P)\to
(M\otimes L)\tensor{A\otimes B}(N\otimes P)$$ such that
$\beta((m\otimes n)\otimes (l\otimes p))=(m\otimes l)\otimes
(n\otimes p)$. Moreover, this morphism is an isomorphism.
\end{lemma}

\begin{proof}
Easy verifications.
\end{proof}

\begin{proposition}\emph{\cite[Proposition 1.5]{Gomez/Louly:2003}}
\\Let $\coring{C}$ be an $A$-coring and $\coring{D}$ be a $B$-coring.
Then $\coring{C}\otimes\coring{D}$ is an $A\otimes B$-coring with
coproduct
$$\xymatrix{\coring{C}\otimes\coring{D}
\ar[rr]^-{\Delta_\coring{C}\otimes\Delta_\coring{D}}
&&(\coring{C}\tensor{A}\coring{C})\otimes(\coring{D}\tensor{B}\coring{D})
\ar[r]^-\simeq & (\coring{C}\otimes\coring{D})\tensor{A\otimes B}
(\coring{C}\otimes\coring{D}),}$$ and counit
$$\xymatrix{\coring{C}\otimes\coring{D}
\ar[rr]^{\epsilon_\coring{C}\otimes\epsilon_\coring{D}}&& A\otimes
B.}$$
\end{proposition}
For example, if $\coring{C}$ is the trivial $A$-coring and
$\coring{D}$ is the trivial $B$-coring, then the $A\otimes B$-coring
$\coring{C}\otimes\coring{D}$ is noting else that the trivial
$A\otimes B$-coring.

\medskip
Every $(A,B)$-bimodule $M$ yields a $(B^\circ,A^\circ)$-bimodule
$M^\circ$. Let $\coring{C}$ be an $A$-coring. Then
$\coring{C}^\circ$ is an $A^\circ$-bimodule. The \emph{opposite
coring} of $\coring{C}$, that is denoted by $\coring{C}^\circ$, is
the $A^\circ$-coring with coproduct
$$\xymatrix{\Delta^\circ:\coring{C}^\circ\ar[r]^\Delta &
\coring{C}\tensor{A}\coring{C} \ar[r]^-\tau &
\coring{C}^\circ\tensor{A^\circ}\coring{C}^\circ},$$ where $\tau$ is
the isomorphism of $k$-modules defined by $\tau(c\otimes
c')=c'\otimes c$, and counit
$\epsilon^\circ=\epsilon:\coring{C}^\circ\to A^\circ$.

As examples, the opposite coring of the trivial $A$-coring is the
trivial $A^\circ$-coring. The opposite coring of a $k$-coalgebra $C$
is the opposite coalgebra of $C$.

\medskip
The following result generalizes both the algebra case and the
coalgebra case.

\begin{proposition}\emph{\cite[Proposition 1.8]{Gomez/Louly:2003}}\label{bcomod}
\\Let $\coring{C}$ be an $A$-coring and $\coring{D}$ be a $B$-coring.
Then there are isomorphisms of categories
$$\bcomod{\coring{C}}{\coring{D}}\simeq \lcomod{\coring{C}
\otimes\coring{D}^\circ}\simeq \rcomod{\coring{C}^\circ
\otimes\coring{D}}.$$
\end{proposition}

\begin{lemma}\label{mor}
\begin{enumerate}[(1)]
\item Let $M,N,K\in \rcomod{\coring{C}}$, $p:M\to N$ be a morphism in
$\rcomod{\coring{C}}$ which is a surjective map, and $f:N\to K$ be a
map such that $fp$ is a morphism in $\rcomod{\coring{C}}$. Then $f$
is a morphism in $\rcomod{\coring{C}}$.
\item  Let $M,N,K\in \rcomod{\coring{C}}$, $f:M\to N$ be a map, and
$i:N\to K$ be a morphism in $\rcomod{\coring{C}}$, such that $i$ and
$i\tensor{A}\coring{C}$ are injective maps, and $if$ is a morphism
in $\rcomod{\coring{C}}$. Then $f$ is a morphism in
$\rcomod{\coring{C}}$.
\end{enumerate}
\end{lemma}

\begin{proof}
The proof is easy. We prove for example (2). It is obvious that $f$
is $A$-linear. We have moreover,
$(i\tensor{A}\coring{C})(f\tensor{A}\coring{C})\rho_M=
(if\tensor{A}\coring{C})\rho_M=
\rho_Kif=(i\tensor{A}\coring{C})\rho_Nf$. Since
$i\tensor{A}\coring{C}$ is injective, then
$(f\tensor{A}\coring{C})\rho_M=\rho_Nf$, that is, $f$ is a morphism
in $\rcomod{\coring{C}}$.
\end{proof}

\begin{proposition}\label{ker}
\begin{enumerate}[(1)]
\item Let $\{M_i\}_I$ be a family of $\bcomod{\coring{C'}}{\coring{C}}$.
Put $M=\bigoplus_I M_i$ in $\bmod{A'}{A}$. Then there is a unique
morphism $\rho_M:M\to M \tensor{A}\coring{C}$ in $\bmod{A'}{A}$
making the following diagrams commutative $(i\in I)$
$$\xymatrix{
M_i \ar[rr]^-{\rho_{M_i}} \ar[d]^{\iota_i}&& M_i\tensor{A}\coring{C}
\ar[d]^{\iota_i\tensor{A}\coring{C}}
\\ M\ar[rr]^-{\rho_M} && M \tensor{A}\coring{C}
}.$$ $(M,\rho_M)$ is a right $\coring{C}$-comodule and
$\iota_i:M_i\to M$ is a morphism of right $\coring{C}$-comodules.
Moreover $\{M,\{\iota_i\}_I\}$ is the coproduct of $\{M_i\}_I$ in
$\rcomod{\coring{C}}$.

Analogously we define $\lambda_M:M\to \coring{C}'\tensor{A'}M$ to be
the unique morphism in $\bmod{A'}{A}$ making commutative the
diagrams $(i\in I)$
$$\xymatrix{
M_i \ar[rr]^-{\lambda_{M_i}} \ar[d]^{\iota_i}&&
\coring{C}'\tensor{A'}M_i \ar[d]^{\coring{C}'\tensor{A'}\iota_i}
\\ M\ar[rr]^-{\lambda_M} && \coring{C}'\tensor{A'}M
}.$$ We have $(M,\rho_M,\lambda_M)$ is a
$\coring{C'}-{\coring{C}}$-bicomodule and $\iota_i:M_i\to M$ is a
morphism of $\coring{C'}-{\coring{C}}$-bicomodules. Moreover
$\{M,\{\iota_i\}_I\}$ is the coproduct of $\{M_i\}_I$ in
$\bcomod{\coring{C'}}{\coring{C}}$.
\item The category $\bcomod{\coring{C'}}{\coring{C}}$ has cokernels,
which coincide with those in $\bmod{A'}{A}$.
\item Let $f:M\to N$ be a morphism in $\bcomod{\coring{C}'}{\coring{C}}$,
and let $i:\oper{Ker}{f}\to M$ be its kernel in $\bmod{A'}{A}$. If
$f$ is $\coring{C}'_{A'}$-pure and $_A\coring{C}$-pure, and the
following
\begin{equation*}
i\otimes_A\coring{C} \otimes_A\coring{C}, \quad
\coring{C}'\otimes_{A'}\coring{C}'\otimes_{A'}i \quad \textrm{and}
\quad \coring{C}'\otimes_{A'}i \otimes_A\coring{C}
\end{equation*}
are injective maps, then $\oper{Ker}{f}$ is the kernel of $f$ in
$\bcomod{\coring{C}'}{\coring{C}}$. This is the case if $f$ is
$(\coring{C}'\otimes_{A'}\coring{C}')_{A'}$-pure, $_A(\coring{C}
\otimes_A\coring{C})$-pure, and $\coring{C}'\otimes_{A'}f$ is
$_A\coring{C}$-pure (e.g. if $\coring{C}'_{A'}$ and $_A\coring{C}$
are flat).
\end{enumerate}
\end{proposition}

\begin{proof}
(1) The proof of the first part is well known and we leave it to the
reader. For the second part consider, for all $i\in I$, the diagram

$$\xymatrix{
 && M_i \ar[rrr]^{\rho_{M_i}}\ar[dll]^(.3){\lambda_{M_i}}
 \ar'[d][ddd]^{\iota_i}&&&
 M_i \tensor{A}\coring{C} \ar[dll]^{\lambda_{M_i}\tensor{A}\coring{C}}
 \ar[ddd]^{\iota_i\tensor{A}\coring{C}}
\\
 \coring{C}'\tensor{A'}M_i \ar[rrr]^(.35){\coring{C}'\tensor{A'}\rho_{M_i}}
 \ar[ddd]^{\coring{C}'\tensor{A'}\iota_i}
 &&& \coring{C}'\tensor{A'}M_i\tensor{A}\coring{C}
 \ar[ddd]^(.3){\coring{C}'\tensor{A'}\iota_i \tensor{A}\coring{C}}
\\
&&&&&&&&&
\\
 && M \ar'[r][rrr]^{\rho_M} \ar[dll]^{\lambda_M} &&&
 M \tensor{A}\coring{C}\ar[dll]^{\lambda_M\tensor{A}\coring{C}}
\\
\coring{C}'\tensor{A'}M \ar[rrr]^{\coring{C}'\tensor{A'}\rho_M} &&&
\coring{C}'\tensor{A'}M \tensor{A}\coring{C}}.$$

Then for all $i\in I$, all the squares, save possibly the bottom
one, are commutative. Hence the bottom square is also commutative.

(2) The proof is well known and we leave it to the reader.

(3) First put $M'=\oper{Ker}{f}$. From the assumptions there is a
unique morphism $\rho_{M'}:M'\to M'\tensor{A}\coring{C}$ in
$\bmod{A'}{A}$ making commutative in $\bmod{A'}{A}$ the following
diagram
$$\xymatrix{0 \ar[r] & M'\ar[r]^i \ar[d]^{\rho_{M'}}&
M \ar[r]^f \ar[d]^{\rho_M}
& N\ar[d]^{\rho_N}
\\ 0 \ar[r] & M'\tensor{A}\coring{C}\ar[r]^{i\tensor{A}\coring{C}} &
M\tensor{A}\coring{C}\ar[r]^{f\tensor{A}\coring{C}} &
N\tensor{A}\coring{C}}.$$
Now consider the diagram
\begin{equation}\label{coass}
\xymatrix{M'\ar[rr]^{\rho_{M'}} \ar[d]^{\rho_{M'}}&&
M'\tensor{A}\coring{C} \ar[rr]^{i\tensor{A}\coring{C}}
\ar[d]^{M'\tensor{A}\Delta_\coring{C}}&& M\tensor{A}\coring{C}
\ar[d]^{M\tensor{A}\Delta_\coring{C}}
\\ M'\tensor{A}\coring{C}\ar[rr]^-{\rho_{M'}\tensor{A}\coring{C}} &&
M'\tensor{A}\coring{C}\tensor{A}\coring{C}\ar[rr]^
{i\tensor{A}\coring{C}\tensor{A}\coring{C}}
 & & M\tensor{A}\coring{C}\tensor{A}\coring{C}}.
 \end{equation}
 $i\tensor{A}\coring{C}$ is a morphism in $\rcomod{\coring{C}}$ means that
 the second rectangle of \eqref{coass} is commutative. Since
 \begin{eqnarray*}
 (M\tensor{A}\Delta_\coring{C})(i\tensor{A}\coring{C})\rho_{M'}&=
 &(M\tensor{A}\Delta_\coring{C})\rho_Mi\\
 &=& (\rho_M\tensor{A}\coring{C})\rho_Mi\\
 &=& (\rho_M\tensor{A}\coring{C})(i\tensor{A}\coring{C})\rho_{M'}\\
 &=& (\rho_Mi\tensor{A}\coring{C})\rho_{M'}\\
 &=&
 (i\tensor{A}\coring{C}\tensor{A}\coring{C})
 (\rho_{M'}\tensor{A}\coring{C})\rho_{M'}.
\end{eqnarray*}
Then the diagram \eqref{coass} is commutative. Hence the
coassociative property of $\rho_{M'}$ holds.

The counit property of $\rho_{M'}$ follows from the commutativity of
the following diagram and the counit property of $\rho_M$
$$\xymatrix{M'\ar[r]^-{\rho_{M'}} \ar[d]^i & M'\tensor{A}\coring{C}
\ar[rr]^{M'\tensor{A}\epsilon_\coring{C}}
\ar[d]^{i\tensor{A}\coring{C}}&& M'\tensor{A}A\ar[d]^{i\tensor{A}A}
\\ M\ar[r]^-{\rho_M} &
M\tensor{A}\coring{C}\ar[rr]^{M\tensor{A}\epsilon_\coring{C}} &&
M\tensor{A}A}.$$ Hence $(M',\rho_{M'})$ is a right
$\coring{C}$-comodule.

By the same way we define $\lambda_{M'}:M'\to
\coring{C'}\tensor{A'}M'$ which endow $M'$ with a structure of left
$\coring{C'}$-comodule.

Now consider the diagram
$$\xymatrix{
 && M' \ar[rrr]^{\rho_{M'}}\ar[dll]^(.3){\lambda_{M'}} \ar'[d][ddd]^i &&&
 M'\otimes_A\coring{C} \ar[dll]^{\lambda_{M'}\otimes_A\coring{C}}
 \ar[ddd]^{i\otimes_A\coring{C}}
\\
\coring{C}'\otimes_{A'}M'
\ar[rrr]^(.35){\coring{C}'\otimes_{A'}\rho_{M'}}
\ar[ddd]^{\coring{C}'\otimes_{A'}i} &&&
\coring{C}'\otimes_{A'}M'\otimes_A\coring{C}
\ar[ddd]^(.3){\coring{C}'\otimes_{A'}i\otimes_A\coring{C}}
\\
&&&&&&&&&
\\
&& M \ar'[r][rrr]^{\rho_M} \ar[dll]^{\lambda_M} &&&
M\otimes_A\coring{C}\ar[dll]^{\lambda_M\otimes_A\coring{C}}
\\
\coring{C}'\otimes_{A'}M \ar[rrr]^{\coring{C}'\otimes_{A'}\rho_M}
&&& \coring{C}'\otimes_{A'}M\otimes_A\coring{C}}.$$ It follows from
Lemma \ref{cube} ($\coring{C}'\otimes_{A'}i\otimes_A\coring{C}$ is
an injective map), that $\rho_M$ is a morphism of left
$\coring{C}'$-comodules (or $\lambda_M$ is a morphism of right
$\coring{C}$-comodules). Hence
$M'\in\bcomod{\coring{C'}}{\coring{C}}$ and $i$ is a monomorphism in
$\bcomod{\coring{C'}}{\coring{C}}$.

Finally, let $\xi:X\to M$ be a monomorphism in
$\bcomod{\coring{C'}}{\coring{C}}$ such that $f\xi=0$. Then there is
a unique morphism $\gamma:X\to M'$ in $\bmod{A'}{A}$ such that
$\xi=i\gamma$. By Lemma \ref{mor} ($i\otimes_A\coring{C}$ and
$\coring{C}'\otimes_{A'}i$ are injective), $\gamma$ is a morphism in
$\bcomod{\coring{C'}}{\coring{C}}$. Therefore $(M',i)$ is the kernel
of $f$ in $\bcomod{\coring{C'}}{\coring{C}}$.
\end{proof}

\begin{corollary}\emph{\cite[Proposition 1.1]{Guzman:1989}}
\\Let $f:M\to N$ be a morphism in $\bcomod{\coring{C}'}{\coring{C}}$,
and let $i:\oper{Ker}{f}\to M$ be its kernel in $\bmod{A'}{A}$. If
$i$ is section in $\bmod{A'}{A}$ and $\oper{Im}{f}$ is a
$\coring{C}'_{A'}$-pure submodule and a $_A\coring{C}$-pure
submodule of $N$, then $\oper{Ker}{f}$ is the kernel of $f$ in
$\bcomod{\coring{C}'}{\coring{C}}$.
\end{corollary}

\begin{corollary}\emph{\cite[18.14]{Brzezinski/Wisbauer:2003}}
\\ For an $A$-coring $\coring{C}$, the following are equivalent
\begin{enumerate}[(1)]
\item $_A\coring{C}$ is flat;
\item the forgetful functor $\rcomod{\coring{C}}\to\rmod{A}$ is a
monofunctor;
\item every monomorphism $U\to \coring{C}$ in $\rcomod{\coring{C}}$
is injective (in $\rmod{A}$).
\end{enumerate}
\end{corollary}

\begin{proof}
$(1)\Longrightarrow (2)$ Clear. $(2)\Longrightarrow (3)$ Trivial.

$(3)\Longrightarrow (1)$ By the Flat Test Lemma
\cite[19.17]{Anderson/Fuller:1992}, it is enough to show that
$\coring{C}$ is $A_A$-flat, that is, for every right ideal $I$ of
$A$, the map $I\tensor{A}\coring{C}\to A\tensor{A}\coring{C}$ is
injective. Since $A\tensor{A}\coring{C}\simeq \coring{C}$ in
$\rcomod{\coring{C}}$ and the functor $-\tensor{A}\coring{C}$
preserves monomorphisms, the claimed statement holds.
\end{proof}

\begin{lemma}\label{comodgenerators}
If $_A\coring{C}$ is flat, then the subcomodules of $\coring{C}^n$,
$n\in \mathbb{N}$, form a generating family for the category
$\rcomod{\coring{C}}$.
\end{lemma}

\begin{proof}
Let $M\in \rcomod{\coring{C}}$. $M$ is a subcomodule of a standard
comodule $M'$. Then there is an epimorphism
$f:\coring{C}^{(\Lambda)}\to M'$ for some index set $\Lambda$. Let
$K$ be a subcomodule of $M$ such that $K\neq M$, and let $m\in M-K$.
There are $n\in \mathbb{N}$ and a monomorphism $i_n:\coring{C}^{(n)}
\to\coring{C}^{(\Lambda)}$ such that $m\in \operatorname{Im}(fi_n)$.
Then we obtain the following diagram
$$\xymatrix{
\coring{C}^{(n)}\ar[r]^{i_n} & \coring{C}^{(\Lambda)}\ar[r]^{f} & M'
\\
(fi_n)^{-1}(M) \ar[rr]^{g}\ar[u] & & M \ar@{^{(}->}[u]}.$$ Then, by
Lemma \ref{mor}, $g$ is a morphism in $\rcomod{\coring{C}}$. $g$
cannot be factorized through a morphism $(fi_n)^{-1}(M)\to K$. Hence
the mentioned family is a family of generators of
$\rcomod{\coring{C}}$.
\end{proof}

Now, we end this section by recalling some corings of particular
interest.

\medskip
\textbf{Coseparable corings.} Following \cite{Guzman:1989}, a coring
$\coring{C}$ is said to be \emph{coseparable} if the
comultiplication map $\Delta_{\coring{C}}$ is a section in the
category $\bcomod{\coring{C}}{\coring{C}}$. Obviously the trivial
coring is coseparable.

\begin{definition}
A comodule $M\in \rcomod{\coring{C}}$ is called $A$-relative
injective comodule or $(\coring{C},A)$-injective comodule if for
every morphism in $\rcomod{\coring{C}}$, $i:N\to L$, that is a
section in $\rmod{A}$, every morphism in $\rcomod{\coring{C}}$,
$f:N\to M$, there is a morphism in $\rcomod{\coring{C}}$, $g:L\to
M$, such that $f=gi$.
\end{definition}

\begin{proposition}\emph{\cite[18.18]{Brzezinski/Wisbauer:2003}}
\label{relatinj}
\begin{enumerate}[(1)]
\item Let $M\in \rcomod{\coring{C}}$. The following are equivalent
\begin{enumerate}[(a)]
\item $M$ is $(\coring{C},A)$-injective;
\item every morphism in $\rcomod{\coring{C}}$, $i:M\to L$, that is
a section in $\rmod{A}$, is also a section in $\rcomod{\coring{C}}$;
\item the coaction $\rho_M:M\to M\tensor{A}\coring{C}$ is a section
in $\rcomod{\coring{C}}$.
\end{enumerate}
\item For every $X\in \rmod{A}$, $X\tensor{A}\coring{C}$ is
$(\coring{C},A)$-injective.
\item For every $M\in \rcomod{\coring{C}}$ that is
$(\coring{C},A)$-injective,
and every $L\in \rcomod{\coring{C}}$, the canonical sequence
$$\xymatrix{
0\ar[r] & \hom{\coring{C}}{L}{M}\ar[r]^i &
\hom{A}{L}{M}\ar[r]^-\gamma & \hom{A}{L}{M\tensor{A}\coring{C}}}$$
is split exact in $\rmod{B}$, where
$B=\operatorname{End}_\coring{C}(L)$ and
$\gamma(f)=\rho_Mf-(f\tensor{A}\coring{C})\rho_L$.
\end{enumerate}
\end{proposition}
For instance the comodule $\coring{C}\in \rcomod{\coring{C}}$  is
$(\coring{C},A)$-injective.

Let $\coring{C}$ be an $A$-coring such that $_A\coring{C}$ is flat.
An exact sequence in $\rcomod{\coring{C}}$ is called
$(\coring{C},A)$-\emph{exact} if it is split exact in $\rmod{A}$. A
functor on $\rcomod{\coring{C}}$ is called right (left)
$(\coring{C},A)$-\emph{exact} if it is right (left) exact on short
$(\coring{C},A)$-exact sequences. For instance $M$ is
$(\coring{C},A)$-injective if and only if $\hom{\coring{C}}{-}{M}$
is $(\coring{C},A)$-exact.

\begin{proposition}\emph{\cite[26.1]{Brzezinski/Wisbauer:2003}}
\label{cosepcoring}
\\Let $\coring{C}$ be an $A$-coring. Then the following are
equivalent
\begin{enumerate}[(1)]
\item $\coring{C}$ is coseparable;
\item there is an $A$-bimodule map $\delta:\coring{C}\tensor{A}
\coring{C}\to A$ satisfying
\begin{equation}\label{cointegral-eq}
\delta\Delta=\epsilon \quad \textrm{and}\quad
(\coring{C}\tensor{A}\delta)(\Delta\tensor{A}\coring{C})=
(\delta\tensor{A}\coring{C})(\coring{C}\tensor{A}\Delta);
\end{equation}
\item the forgetful functor $(-)_A:\rcomod{\coring{C}}\to \rmod{A}$
is separable;
\item the forgetful functor $_A(-):\lcomod{\coring{C}}\to \lmod{A}$
is separable;
\item the forgetful functor
$_A(-)_A:\bcomod{\coring{C}}{\coring{C}} \to \bmod{A}{A}$ is
separable;
\item $\coring{C}$ is $(A,A)$-relative semisimple as a
$\coring{C}$-bicomodule, that is any monomorphism in
$\bcomod{\coring{C}}{\coring{C}}$ that splits as an $A$-bimodule map
also splits in $\bcomod{\coring{C}}{\coring{C}}$;
\item $\coring{C}$ is $(A,A)$-relative injective as a
$\coring{C}$-bicomodule;
\item $\coring{C}$ is $(A,\coring{C})$-relative injective as a
$\coring{C}$-bicomodule;
\item $\coring{C}$ is $(\coring{C},A)$-relative injective as a
$\coring{C}$-bicomodule.
\end{enumerate}
In such a case, $\coring{C}$ is right and left
$(\coring{C},A)$-semisimple, that is, all comodules in
$\rcomod{\coring{C}}$ and $\lcomod{\coring{C}}$ are
$(\coring{C},A)$-injective.
\end{proposition}

\begin{definition}\label{cointegral}
An $A$-bimodule map $\delta:\coring{C}\tensor{A} \coring{C}\to A$
satisfying conditions \eqref{cointegral-eq} is called a
\emph{cointegral} in $\coring{C}$.
\end{definition}

\textbf{Cosplit corings.} A coring is called a \emph{cosplit coring}
if it satisfy the conditions of the proposition below.

\begin{proposition}
\label{cosplitcoring} Let $\coring{C}$ be an $A$-coring. Then the
following are equivalent
\begin{enumerate}[(1)]
\item the functor $-\tensor{A}\coring{C}:\rmod{A}\to\rcomod{\coring{C}}$ is
separable;
\item the functor $\coring{C}\tensor{A}-:\lmod{A}\to\lcomod{\coring{C}}$
is separable;
\item $\epsilon$ is a retraction in $\bmod{A}{A}$.
\end{enumerate}
\end{proposition}

\begin{proof}
By Rafael's Proposition \ref{Rafael}, $-\tensor{A}\coring{C}$ is
separable if and only if there exists a natural transformation
$\xi': 1_{\rmod{A}}\to -\tensor{A}\coring{C}$ such that
$\xi_X\circ\xi'_X=1_X$ for every $X\in\rmod{A}$, where
$\xi_X=X\tensor{A}\epsilon$, $X\in\rmod{A}$, is the counit of the
adjunction $(U_r,-\tensor{A}\coring{C})$. ($U_r$ is the forgetful
functor.) By Lemma \ref{$A$-objects} (or by Lemma \ref{18}), a
natural transformation $\xi': 1_{\rmod{A}}\to -\tensor{A}\coring{C}$
is entirely determined by the data of an $A$-bimodule map
$\xi':A\to\coring{C}$ such that $\xi'_X=X\tensor{A}\xi'$ for every
$X\in\rmod{A}$. Hence the equivalence $(1)\Leftrightarrow(3)$
follows. $(2)\Leftrightarrow(3)$ follows by symmetry.
\end{proof}

\medskip
\textbf{Cosemisimple corings.} Let $\cat{C}$ be an abelian category.
$\cat{C}$ is called a \emph{spectral category} if every short exact
sequence in $\cat{C}$ splits, or equivalently, if every object is
injective (resp. projective). As for modules, an object $S$ of
$\cat{C}$ is \emph{simple} if $S\neq 0$ has exactly two object 0 and
$S$, and it is \emph{semisimple} if it a sum of simple subjects. 0
is, by definition, a semisimple object. A spectral category is
called \emph{discrete} if every object is semisimple.

\begin{proposition}\emph{\cite[Propsition 6.7]{Stenstrom:1975}}\label{spectral category}
\\Let $\cat{C}$ be a Grothendieck category. Then the following are
equivalent
\begin{enumerate}[(1)]
\item $\cat{C}$ is a discrete spectral category;
\item $\cat{C}$ is a locally finitely generated spectral category;
\item every object of $\cat{C}$ is semisimple;
\item $\cat{C}$ has a family of simple generators;
\item $\cat{C}$ is equivalent to a product category
$\prod_I\rmod{K_i}$, where $K_i$, $i\in I$, is a family of division
rings.
\end{enumerate}
\end{proposition}

For more details we refer to \cite[Chap.V, \S6]{Stenstrom:1975}.

\begin{proposition}\emph{\cite[Theorem 3.1]{ElKaoutit/Gomez/Lobillo:2001unp}}\label{cosemisimple}
\\Let $\coring{C}$ be an $A$-coring. Then the following are
equivalent
\begin{enumerate}[(1)]
\item $\rcomod{\coring{C}}$ is an abelian and a discrete spectral category;
\item $\lcomod{\coring{C}}$ is an abelian and a discrete spectral category;
\item $\coring{C}$ is semisimple in $\rcomod{\coring{C}}$ and
$_A\coring{C}$ is flat;
\item $\coring{C}$ is semisimple in $\lcomod{\coring{C}}$ and
$\coring{C}_A$ is flat;
\item $\coring{C}$ is a semisimple right $\rdual{\coring{C}}$-module
and $\coring{C}_A$ is projective;
\item $\coring{C}$ is a semisimple left $\ldual{\coring{C}}$-module
and $_A\coring{C}$ is projective.
\end{enumerate}
\end{proposition}

If the above conditions hold, we say that $\coring{C}$ is a
\emph{cosemisimple coring}.

\medskip
\textbf{Semiperfect corings.} In \cite{Bass:1960}, H. Bass defined
and studied perfect and semiperfect rings. In \cite{Lin:1977}, L.-P.
Lin introduced and studied semiperfect coalgebras over fields. In
\cite[\S3]{Harada:1973}, M. Harada defined perfect and semiperfect
Grothendieck categories and characterized them with respect to a
generating set. Let $\cat{C}$ be a Grothendieck category. $\cat{C}$
is called \emph{perfect} (resp. \emph{semiperfect}) if every object
of $\cat{C}$ has a projective cover (which is defined as for
modules).

Let $\coring{C}$ be an $A$-coring such that $_A\coring{C}$ is flat.
$\coring{C}$ is called right \emph{semiperfect} if
$\rcomod{\coring{C}}$ is a semiperfect category.

\begin{proposition}\emph{\cite[Theorem 3.1]{Caenepeel/Iovanov:2005}}\label{semiperfect}
\\Let $\coring{C}$ be an $A$-coring such that $A$ is a right
artinian ring and $_A\coring{C}$ is projective. Then the following
are equivalent
\begin{enumerate}[(1)]
\item $\coring{C}$ is right semiperfect;
\item every simple object of $\rcomod{\coring{C}}$ has a projective
cover;
\item every finitely generated object of $\rcomod{\coring{C}}$ has a
finitely generated projective cover;
\item every simple object of $\rcomod{\coring{C}}$ has a
finitely generated projective cover;
\item the category $\rcomod{\coring{C}}$ has enough projectives;
\item the category $\rcomod{\coring{C}}$ has a projective generator.
\end{enumerate}
\end{proposition}

For a study of semiperfect corings over QF rings, see
\cite{ElKaoutit/Gomez:2002unp} or
\cite[\S4]{Caenepeel/Iovanov:2005}.

\section{Cotensor product over corings}

First we recall that the cotensor product over coalgebras over
fields was introduced by Milnor and Moore in
\cite{Milnor/Moore:1965}. A detailed study of the cotensor product
over corings is given in \cite{Brzezinski/Wisbauer:2003}.

Let
 $M\in{}\bcomod{\coring{C}'}{\coring{C}}$ and
$N\in{}\bcomod{\coring{C}}{\coring{C}''}$. The map
\[
\omega_{M,N}=\rho_M\otimes_AN-M\otimes_A\lambda_N:M\otimes_AN
\rightarrow M\otimes_A\coring{C}\otimes_AN
\]
is a $\coring{C}'-\coring{C}''$-bicomodule map.

Notice that $W\tensor{A'}\omega_{M,N}=\omega_{W\tensor{A'}M,N}$ and
$\omega_{M,N}\tensor{A''}V=\omega_{M,N\tensor{A''}V}$ for every
$W\in \rmod{A'}$ and $V\in \lmod{A''}$.

\begin{definition}
The kernel of $\omega_{M,N}$ in $\bmod{A'}{A''}$ is the
\emph{cotensor product} of $M$ and $N$, and it is denoted by
$M\cotensor{\coring{C}}N$ . We denote the canonical injection
$M\cotensor{\coring{C}}N\to M\tensor{A}N$ by $\iota_{M,N}$.
\end{definition}

\begin{proposition}\label{bcomodcotensor}
If $\omega_{M,N}$ is $\coring{C}'_{A'}$- and
$_{A''}\coring{C''}$-pure, and the following
\begin{equation*}
\iota_{M,N}\tensor{A''}\coring{C}''\tensor{A''}\coring{C}'', \quad
\coring{C}'\otimes_{A'}\coring{C}'\otimes_{A'}\iota_{M,N} \quad
\textrm{and} \quad \coring{C}'\tensor{A'}\iota_{M,N}
\tensor{A''}\coring{C}''
\end{equation*}
are injective maps, then $M\cotensor{\coring{C}}N$ is the kernel of
$\omega_{M,N}$ in $\bcomod{\coring{C}'}{\coring{C}''}$.

 These conditions are fulfilled if $\omega_{M,N}$ is
$(\coring{C}'\otimes_{A'}\coring{C}')_{A'}$- and
$_{A''}(\coring{C}'' \otimes_{A''}\coring{C}'')$-pure, and
$\coring{C}'\otimes_{A'}\omega_{M,N}$ is $_{A''}\coring{C}''$-pure
(e.g. if $\coring{C}'_{A'}$ and $_{A''}\coring{C}''$ are flat).
\end{proposition}

\begin{proof}
Follows immediately from Proposition \ref{ker}.
\end{proof}

\begin{definition}
Let $f:M\to M'$ be a morphism in $\bcomod{\coring{C}'}{\coring{C}}$,
and $g:N\to N'$ a morphism in $\bcomod{\coring{C}}{\coring{C}''}$.
There is a unique morphism
$f\cotensor{\coring{C}}g:M\cotensor{\coring{C}}N\to
M'\cotensor{\coring{C}}N'$ in $\bmod{A'}{A''}$ making commutative
the diagram
$$\xymatrix{
0\ar[r] & M\cotensor{\coring{C}}N
\ar[r]\ar[d]^{f\cotensor{\coring{C}}g} & M\tensor{A}N
\ar[rr]^-{\omega_{M,N}}\ar[d]^{f\tensor{A}g} &&
M\tensor{A}\coring{C}\tensor{A}N
\ar[d]^{f\tensor{A}\coring{C}\tensor{A}g}
\\ 0 \ar[r] & M'\cotensor{\coring{C}}N'\ar[r] & M'\tensor{A}N'
\ar[rr]^-{\omega_{M',N'}} && M'\tensor{A}\coring{C}\tensor{A}N'.}
$$
$f\cotensor{\coring{C}}g$ is called the \emph{cotensor product} of
$f$ and $g$.
\end{definition}

\begin{proposition}\label{cotensorbifunctor}
If for every $M\in{}\bcomod{\coring{C}'}{\coring{C}}$ and
$N\in{}\bcomod{\coring{C}}{\coring{C}''}$, $\omega_{M,N}$ is
$\coring{C}'_{A'}$- and $_{A''}\coring{C}''$-pure, then we have a
$k$-linear bifunctor
\begin{equation}\label{cotenbifun}
\xymatrix{-\cotensor{\coring{C}}-:{}^{\coring{C}'}\mathcal{M}^{\coring{C}}
\times{}^{\coring{C}}\mathcal{M}^{\coring{C}''}\ar[r]&
{}^{\coring{C}'}\mathcal{M}^{\coring{C}''}}.
\end{equation}
In particular, if $\coring{C}'_{A'}$ and $_{A''}\coring{C}'' $ are
flat, then the bifunctor \eqref{cotenbifun} is well defined.
\end{proposition}

\begin{proof}
Follows immediately from Proposition \ref{bcomodcotensor} and Lemma
\ref{mor}.
\end{proof}

\begin{lemma}\label{assotenscotens}
\begin{enumerate}[(1)]
\item Let $M\in{}\bcomod{\coring{C}'}{\coring{C}}$ and
$N\in{}\bcomod{\coring{C}}{\coring{C}''}$ and $W\in \rmod{A'}$. We
define $\psi_{W,M,N}: W\tensor{A'}(M\cotensor{\coring{C}}N)\to
(W\tensor{A'}M)\cotensor{\coring{C}}N$ to be the unique morphism
making commutative the following diagram
$$\xymatrix{ &
W\tensor{A'}(M\cotensor{\coring{C}}N) \ar[r]\ar[d]^{\psi_{W,M,N}} &
W\tensor{A'}(M\tensor{A}N)
\ar[rr]^-{W\tensor{A'}\omega_{M,N}}\ar[d]^{\simeq} &&
W\tensor{A'}(M\tensor{A}\coring{C}\tensor{A}N) \ar[d]^{\simeq}
\\ 0 \ar[r] & (W\tensor{A'}M)\cotensor{\coring{C}}N\ar[r] &
(W\tensor{A'}M)\tensor{A}N \ar[rr]^-{\omega_{W\tensor{A'}M,N}} &&
(W\tensor{A'}M)\tensor{A}\coring{C}\tensor{A}N.}
$$ ($\psi_{W,M,N}\big(\sum w_i\otimes (m_j\otimes n_k)\big)=
\sum (w_i\otimes m_j)\otimes n_k$.)
Then
\begin{enumerate}[(a)]
\item $\psi_{W,M,N}$ is natural in $W$, $M$ and
$N$.
\item $\psi_{W,M,N}$ is an isomorphism if and only if $\omega_{M,N}$
is $W_{A'}$-pure. In particular this the case if $W_{A'}$ is flat.
\end{enumerate}
\item We have analogous statements for $W\in \lmod{A''}$.
\end{enumerate}
\end{lemma}

\begin{proof}
(1) $(a)$ For example we verify that $\psi_{W,M,N}$ is natural in
$N$. For this, let $f:N\to N'$ be a morphism in
$\bcomod{\coring{C}}{\coring{C}''}$, and consider the diagram
$$\xymatrix@C 2pt{
 && W\tensor{A'}(M\cotensor{\coring{C}}N)
  \ar[rrr]^{W\tensor{A'}(M\cotensor{\coring{C}}f)}
  \ar[dll]^(.3){\psi_{W,M,N}}
  \ar'[d][ddd]&&& W\tensor{A'}(M\cotensor{\coring{C}}N')
  \ar[dll]^{\psi_{W,M,N'}}\ar[ddd]
\\
(W\tensor{A'}M)\cotensor{\coring{C}}N
\ar[rrr]^(.35){(W\tensor{A'}M)\cotensor{\coring{C}}f}\ar[ddd] &&&
(W\tensor{A'}M)\cotensor{\coring{C}}N' \ar[ddd]
\\
&&&&&&&&&
\\
 && W\tensor{A'}(M\tensor{A}N) \ar'[r][rrr]^{W\tensor{A'}(M\tensor{A}f)}
 \ar[dll]^{\simeq} &&& W\tensor{A'}(M\tensor{A}N')\ar[dll]^{\simeq}
\\
(W\tensor{A'}M)\tensor{A}N \ar[rrr]^{(W\tensor{A'}M)\tensor{A}f} &&&
(W\tensor{A'}M)\tensor{A}N'.}$$ Then the claimed statement follows
from Lemma \ref{cube}.

$(b)$ Obvious.

(2) Analogous to (1).
\end{proof}

\begin{lemma}\label{relatinjsplit}
Let $A$ and $T$ be $k$-algebras. Let $M\in \rcomod{\coring{C}}$ and
$N\in \lcomod{\coring{C}}_T$.
\begin{enumerate}[(1)]
\item If the functor $M\cotensor{\coring{C}}-$ is right exact, then
$\omega_{M,N}$ is a pure morphism in $\rmod{T}$.
\item If $M$ is $(\coring{C},A)$-injective, then the sequence
$$\xymatrix{0\ar[r] & M\cotensor{\coring{C}}N\ar[r] & M\tensor{A}N
\ar[r]^-{\omega_{M,N}} & M\tensor{A}\coring{C}\tensor{A}N}$$ is
split exact in $\rmod{T}$. In particular it is pure in $\rmod{T}$.
\end{enumerate}
Analogous statements are true for $M\in {}_S\rcomod{\coring{C}}$ and
$N\in \lcomod{\coring{C}}$, where $S$ is a $k$-algebra.
\end{lemma}

\begin{proof}
(1) Let $X\in \lmod{T}$. Then there is an exact sequence in
$\lmod{T}$, $F'\to F\to X\to 0$, where $F$ and $F'$ are free. We
obtain the commutative diagram with exact rows
$$\xymatrix{(M\cotensor{\coring{C}}N)\tensor{T}F'
\ar[r]\ar[d]^{\psi'_{M,N,F'}}_\simeq &
(M\cotensor{\coring{C}}N)\tensor{T}F
\ar[r]\ar[d]^{\psi'_{M,N,F}}_\simeq &
(M\cotensor{\coring{C}}N)\tensor{T}X \ar[r]\ar[d]^{\psi'_{M,N,X}} &
0
\\ M\cotensor{\coring{C}}(N\tensor{T}F') \ar[r] &
M\cotensor{\coring{C}}(N\tensor{T}F)\ar[r] &
M\cotensor{\coring{C}}(N\tensor{T}X) \ar[r] & 0. }$$ $\psi'_{M,N,F}$
and $\psi'_{M,N,F'}$ are isomorphisms since $F$ and $F'$ are flat.
Hence $\psi'_{M,N,X}$ is also an isomorphism. The statement follows
then from Lemma \ref{assotenscotens}.

(2) It enough to show that the canonical injections in $\rmod{T}$,
$M\cotensor{\coring{C}}N\to M\tensor{A}N$ and
$\oper{Im}{\omega_{M,N}} \to M\tensor{A}\coring{C}\tensor{A}N$ are
sections. Let $\lambda:M\tensor{A}\coring{C}\to M$ be a morphism in
$\rcomod{\coring{C}}$ such that $\lambda\circ\rho_M=1_M$. Consider
$\alpha=(\lambda\tensor{A}N)\circ(M\tensor{A}\lambda_N):M\tensor{A}N\to
M\tensor{A}N$ which is a morphism in $\rmod{T}$.

First we have $\oper{Im}{\alpha}\subset M\cotensor{\coring{C}}N$
since
\begin{eqnarray*}
 (\rho_M\tensor{A}N)\circ \alpha&=&
 (\rho_M\circ\lambda\tensor{A}N)\circ(M\tensor{A}\lambda_N)\\
 &=& (\lambda\tensor{A}\coring{C}\tensor{A}N)
 \circ(M\tensor{A}\Delta_\coring{C}\tensor{A}N)
 \circ(M\tensor{A}\lambda_N) \\
 && (\textrm{Since}\, \lambda \,\textrm{is a morphism in}
\,\rcomod{\coring{C}}\,\textrm{i.e.} \, \rho_M\circ\lambda=
(\lambda\tensor{A}\coring{C})\circ(M\tensor{A}\Delta_\coring{C}))\\
&=& (\lambda\tensor{A}\coring{C}\tensor{A}N)
 \circ\Big(M\tensor{A}\big((\Delta_\coring{C}\tensor{A}N)\lambda_N\big)\Big),
\end{eqnarray*}
and
 \begin{eqnarray*}
 (M\tensor{A}\lambda_N)\circ\alpha&=& (\lambda\tensor{A}\coring{C}\tensor{A}N)
 \circ(M\tensor{A}\coring{C}\tensor{A}\lambda_N)\circ(M\tensor{A}\lambda_N)\\
 &=& (\lambda\tensor{A}\coring{C}\tensor{A}N)
 \circ\Big(M\tensor{A}\big((\Delta_\coring{C}\tensor{A}N)\lambda_N\big)\Big)\\
&& (\textrm{Since}\, \lambda_N \,\textrm{is a morphism in}
\,\lcomod{\coring{C}}\,\textrm{i.e.} \,
(\coring{C}\tensor{A}\lambda_N)\circ\lambda_N=
(\Delta_\coring{C}\tensor{A}N)\circ\lambda_N).
\end{eqnarray*}
On the other hand, for every $x\in M\cotensor{\coring{C}}N$,
$\alpha(x)=(\lambda\tensor{A}N)\circ(\rho_M\tensor{A}N)(x)=x$. Hence
$M\cotensor{\coring{C}}N\to M\tensor{A}N$ is a section.

Now consider the canonical diagram
$$\xymatrix{M\tensor{A}N
\ar[r]^-{\omega_{M,N}}\ar[d]^\pi & M\tensor{A}\coring{C}\tensor{A}N
\\ (M\tensor{A}N)/(M\cotensor{\coring{C}}N)
\ar[r]^-{\bar\omega_{M,N}}_-\simeq &
\oper{Im}{\omega_{M,N}}\ar@{^{(}->}[u]^\iota},$$ and the map
$\gamma=\bar\omega_{M,N}\circ\pi\circ(\lambda\tensor{A}N):
M\tensor{A}\coring{C}\tensor{A}N\to
\oper{Im}{\omega_{M,N}}$. We have
$(\lambda\tensor{A}N)\circ\omega_{M,N}=1_{M\tensor{A}N}-\alpha$.
Then $\pi\circ(\lambda\tensor{A}N)\circ\omega_{M,N}=\pi$. Hence
$\gamma\circ\iota=1_{\oper{Im}{\omega_{M,N}}}$.

The last assertion is proved by a similar method.
\end{proof}

\begin{corollary}\label{}
Let
 $M\in{}\bcomod{\coring{C}'}{\coring{C}}$ and
$N\in{}\bcomod{\coring{C}}{\coring{C}''}$. If $\coring{C}$ is a
coseparable $A$-coring, then $M\cotensor{\coring{C}}N$ is the kernel
of $\omega_{M,N}$ in ${}\bcomod{\coring{C}'}{\coring{C}''}$.
\end{corollary}

\begin{proof}
Follows immediately from Propositions \ref{bcomodcotensor},
\ref{cosepcoring}, and \ref{relatinjsplit}.
\end{proof}

\begin{corollary}
If $\coring{C}$ is a coseparable $A$-coring, then the bifunctor
\eqref{cotenbifun} is well defined. In the special case when
$\coring{C}$ is the trivial $A$-coring,
$-\cotensor{\coring{C}}-=-\tensor{A}-$.
\end{corollary}

\begin{proof}
Follows immediately from Propositions \ref{cotensorbifunctor},
\ref{cosepcoring}, and \ref{relatinjsplit}.
\end{proof}

\begin{proposition}\label{cotensdirlim}
For every $M\in \rcomod{\coring{C}}$, the functor
$M\cotensor{\coring{C}}-$ preserves direct limits.
\end{proposition}

\begin{proof}
Let $I$ be a direct set, $\textbf{I}$ be the associated category,
and let $(N_i,u_{ij})$ be a directed system in $\lcomod{\coring{C}}$
indexed by $I$. Let $(N,u_i)$ be a direct limit in $\lmod{A}$ of
this direct system. It is clear that $N$ is a left
$\coring{C}$-comodule via the unique morphism $\lambda_N$ making
commutative the following diagram
$$\xymatrix{N_i\ar[rr]^{u_i} \ar[d]^{\lambda_{N_i}} &&
N\ar[d]^{\lambda_N}
\\ \coring{C}\tensor{A}N_i \ar[rr]^{\coring{C}\tensor{A}u_i} &&
\coring{C}\tensor{A}N. }$$ We have moreover $u_i$ is a morphism in
$\lcomod{\coring{C}}$. Let $\upsilon_i:N_i\to L$, $i\in I$, be a
compatible family with $(N_i,u_{ij})$ in $\lcomod{\coring{C}}$. Let
$\upsilon:N\to L$ be the unique morphism in $\lmod{A}$ making
commutative the diagram
$$\xymatrix{N_i\ar[r]^{u_i} \ar[dr]^{\upsilon_i}& N\ar[d]^{\upsilon}
\\ & L.
}$$ Therefore,
$(\coring{C}\tensor{A}\upsilon)\lambda_Nu_i=(\coring{C}\tensor{A}\upsilon
u_i)\lambda_{N_i}=(\coring{C}\tensor{A}\upsilon_i)\lambda_{N_i}$ and
$\lambda_L\upsilon u_i=
\lambda_L\upsilon_i=(\coring{C}\tensor{A}\upsilon_i)\lambda_{N_i}$.
Hence $(\coring{C}\tensor{A}\upsilon)\lambda_Nu_i=\lambda_L\upsilon
u_i$ for every $i\in I$, and by the definition of direct limit,
$(\coring{C}\tensor{A}\upsilon)\lambda_N=\lambda_L\upsilon$, i.e.
$\upsilon$ is a morphism in $\lcomod{\coring{C}}$. It follows that
$(N,u_i)$ is the direct limit of $(N,u_i)$ in $\lcomod{\coring{C}}$.

We have moreover, for every $i\in I$, a commutative diagram in
$\rmod{k}$ with exact rows
$$\xymatrix{
0 \ar[r] & M\cotensor{\coring{C}}N_i \ar[r]
\ar[d]^{M\cotensor{\coring{C}}u_i}& M\tensor{A} N_i
\ar[rr]^-{\omega_{M,N_i}}\ar[d]^{M\tensor{A}u_i} &&
M\tensor{A}\coring{C}\tensor{A}N_i
\ar[d]^{M\tensor{A}\coring{C}\tensor{A}u_i}
\\ 0 \ar[r] & M\cotensor{\coring{C}}N \ar[r] & M\tensor{A} N
\ar[rr]^-{\omega_{M,N}} && M\tensor{A}\coring{C}\tensor{A}N. }$$ The
first rows of the the above diagrams define an exact sequence in
$\textbf{Fun}(\textbf{I},\rmod{k})$. Since direct limits are exact
in $\rmod{k}$, $(M\cotensor{\coring{C}}N,M\cotensor{\coring{C}}u_i)$
is the direct limit of the direct system
$(M\cotensor{\coring{C}}N_i,M\cotensor{\coring{C}}u_{ij})$ in
$\rmod{k}$.
\end{proof}

It is well known that the cotensor functor, for coalgebras over
fields, is left exact (see \cite{Takeuchi:1977}). But for corings,
the cotensor functor is not left exact nor right exact in general.
For a counterexample, see \cite[Counterexample 2.5, p.
338]{Grunenfelder/Pare:1987} which is an example of a
$\mathbb{Z}$-coalgebra $C$ such that the cotensor product
$\cotensor{C}$ is not associative. We will give in the next result
sufficient conditions to have left exactness of the cotensor
functor.

\begin{proposition}
Let $_A\coring{C}$ be flat and $M\in {}\lcomod{\coring{C}}$. Let
$$\xymatrix{ 0 \ar[r] & N_1 \ar[r]^f & N_2 \ar[r]^g & N_3\ar[r] & 0}$$
be an exact sequence in $\rcomod{\coring{C}}$. Then the resulting
sequence
$$\xymatrix{ 0 \ar[r] & N_1\cotensor{\coring{C}}M
\ar[rr]^{f\cotensor{\coring{C}}M} && N_2\cotensor{\coring{C}}M
\ar[rr]^{g\cotensor{\coring{C}}M} && N_3\cotensor{\coring{C}}M}$$ is
exact if $f\tensor{A}M$ and $f\tensor{A}\coring{C}\tensor{A}M$ are
injective maps.

In particular, the functor
$-\cotensor{\coring{C}}M:\rcomod{\coring{C}}\to \rmod{k}$ is left
$(\coring{C},A)$-exact, and it is left exact if $_AM$ is flat.
\end{proposition}

\begin{proof}
Consider the following commutative diagram with exact columns
$$\xymatrix{& 0\ar[d] && 0\ar[d] && 0\ar[d]
\\0 \ar[r] & N_1\cotensor{\coring{C}}M
\ar[rr]^{f\cotensor{\coring{C}}M} \ar[d] &&
N_2\cotensor{\coring{C}}M \ar[rr]^{g\cotensor{\coring{C}}M}\ar[d] &&
N_3\cotensor{\coring{C}}M\ar[d]
\\ 0 \ar[r] & N_1\tensor{A}M
\ar[rr]^{f\tensor{A}M} \ar[d]^{\omega_{N_1,M}} && N_2\tensor{A}M
\ar[rr]^{g\tensor{A}M} \ar[d]^{\omega_{N_2,M}} &&
N_3\tensor{A}M\ar[d]^{\omega_{N_3,M}}\ar[r] & 0
\\ 0 \ar[r] & N_1\tensor{A}\coring{C}\tensor{A}M
\ar[rr]^{f\tensor{A}\coring{C}\tensor{A}M} &&
N_2\tensor{A}\coring{C}\tensor{A}M
\ar[rr]^{g\tensor{A}\coring{C}\tensor{A}M}&&
N_3\tensor{A}\coring{C}\tensor{A}M\ar[r] & 0. }$$ By assumptions,
the second and the third rows are exact. From Snake Lemma
(\cite[Lemmas III.3.2, III.3.3]{Cartan/Eilenberg:1956}) the top row
is also exact.
\end{proof}

\begin{remark}\label{leftexactcotenfun}
As a corollary of the above result we state: Let $_A\coring{C}$ is
flat and $M\in {}\lcomod{\coring{C}}$. Let
$$\xymatrix{ 0 \ar[r] & N_1 \ar[r]^f & N_2 \ar[r]^g & N_3}$$
be an exact sequence in $\rcomod{\coring{C}}$ which is $_AM$- and
$_A(\coring{C}\tensor{A}M)$-pure. Then the resulting sequence
$$\xymatrix{ 0 \ar[r] & N_1\cotensor{\coring{C}}M
\ar[rr]^{f\cotensor{\coring{C}}M} && N_2\cotensor{\coring{C}}M
\ar[rr]^{g\cotensor{\coring{C}}M} && N_3\cotensor{\coring{C}}M}$$ is
exact.
\end{remark}

\begin{lemma}\label{caniso}
Let $M\in \bcomod{\coring{D}}{\coring{C}}$ and $V\in \lmod{A}$. Then
there is a canonical isomorphism
$$M\cotensor{\coring{C}}(\coring{C}\tensor{A}V)\simeq M\tensor{A}V.$$
If moreover $\omega_{M,\coring{C} \tensor{A} V}$ is
$(\coring{D}\tensor{B}\coring{D})_{B}$-pure for every $V\in
\lmod{A}$, then
$$M\cotensor{\coring{C}}(\coring{C}\tensor{A}-)\simeq M\tensor{A}-:
\lmod{A}\to \lcomod{\coring{D}}.$$
\end{lemma}

\begin{proof}
First consider the diagram
\[
\xymatrix{ 0 \ar[r] & M \cotensor{\coring{C}} (\coring{C} \tensor{A}
V) \ar^{\iota}[rr] && M \tensor{A} \coring{C} \tensor{A} V
\ar^-{\omega_{M,\coring{C} \tensor{A} V}}[rr] \ar@<2pt>[dll]^{M
\tensor{A}\epsilon_\coring{C}\tensor{A} V}& & M \tensor{A}
\coring{C} \tensor{A} \coring{C} \tensor{A} V \\
 & M \tensor{A} V \ar^{\psi_M}@{-->}[u] \ar^-{\rho_M \tensor{A} V}[urr]
 &&&&},
\]
By the universal property of kernel, there is a unique map $\psi_M$
such that $\iota\psi_M=\rho_M \tensor{A} V.$ Moreover we have
$(M\tensor{A}\epsilon_\coring{C}\tensor{A} V)(\rho_M \tensor{A}
V)=M\tensor{A} V$ and
\begin{eqnarray*}
(\rho_M \tensor{A} V)(M\tensor{A}\epsilon_\coring{C}\tensor{A}
V)\iota&=&(M\tensor{A}\coring{C}\tensor{A}\epsilon_\coring{C}
\tensor{A}V)(\rho_M\tensor{A}\coring{C}\tensor{A} V)\iota
\\&=& (M\tensor{A}\coring{C}\tensor{A}\epsilon_\coring{C}
\tensor{A}N)(M\tensor{A}\Delta_\coring{C}\tensor{A} V)\iota
\\ &=& (M\tensor{A}\coring{C}\tensor{A} V)\iota
\\ &=& \iota.
\end{eqnarray*}
Then $\psi_M$ is bijective with inverse map
$(M\tensor{A}\epsilon_\coring{C}\tensor{A}V)\iota$.

Now let $f:V\to V'$ be a morphism in $\lmod{A}$. Consider the
diagram
$$\xymatrix{
 && M\tensor{A}V \ar[rrr]^\simeq\ar[dll]^{M\tensor{A}f}
 \ar[ddrrr]^(.3){\rho_M\tensor{A}V}&&&
 M\cotensor{\coring{C}}(\coring{C}\tensor{A}V)
 \ar[dll]^{M\cotensor{\coring{C}}(\coring{C}\tensor{A}f)}\ar[dd]
\\
M\tensor{A}V' \ar[rrr]^\simeq \ar[ddrrr]^{\rho_M\tensor{A}V'}&&&
M\cotensor{\coring{C}}(\coring{C}\tensor{A}V') \ar[dd]
\\
 &&&&&  M \tensor{A}\coring{C}\tensor{A}V
 \ar[dll]^{M\tensor{A}\coring{C}\tensor{A}f}
\\
  &&& M \tensor{A}\coring{C}\tensor{A}V'.}$$
All sides, except perhaps the top one, are commutative, then the top
side is also commutative.
\end{proof}

\begin{definition}
Let $\coring{C}_A$ be flat. A comodule $M\in \rcomod{\coring{C}}$ is
called coflat (resp. faithfully coflat) if the functor
$M\cotensor{\coring{C}}-:\lcomod{\coring{C}}\to \rmod{k}$ is exact
(resp. faithfully exact).
\end{definition}

\begin{corollary}\label{coflat flat}
\begin{enumerate}[(1)]
\item Let $M\in \bcomod{\coring{D}}{\coring{C}}$. If
$\oper{Im}{\omega_{M,\coring{C}}}$ is a $_A\coring{C}$-pure
submodule of $M\tensor{A}\coring{C}\tensor{A}\coring{C}$, then
$$M\cotensor{\coring{C}}\coring{C}\simeq M$$ as
$(\coring{D},\coring{C})$-bicomodules.
\item Let $M\in \rcomod{\coring{C}}$. If $\coring{C}_A$ is flat
and $M$ is coflat in $\rcomod{\coring{C}}$ then $M_A$ is flat.
\end{enumerate}
\end{corollary}

\begin{proof}
(1) Follows from Proposition \ref{bcomodcotensor}, Lemmas
\ref{caniso}, \ref{relatinjsplit} ($\coring{C}$ is
$(\coring{C},A)$-injective), and Lemma \ref{mor}. (2) Obvious from
Lemma \ref{caniso}.
\end{proof}

\begin{proposition}\label{assocotens}
Suppose that $_A\coring{C}$ is flat and for every
$X\in{}\bcomod{\coring{C}'}{\coring{C}}$,
$Y\in{}\bcomod{\coring{C}}{\coring{C}''}$ and
$Z\in{}\bcomod{\coring{C}'''}{\coring{C}'}$, we have $\omega_{X,Y}$
is $\coring{C}'_{A'}$- and $_{A''}\coring{C}''$-pure and
$\omega_{Z,X}$ is $\coring{C}'''_{A'''}$-pure.

Let
 $M\in{}\bcomod{\coring{C}'}{\coring{C}}$,
$N\in{}\bcomod{\coring{C}}{\coring{C}''}$ and
$L\in{}\bcomod{\coring{C}'''}{\coring{C}'}$. If the following hold
\begin{enumerate}[(a)]
\item $\omega_{M,N}$ is $L_{A'}$-,
$(L\tensor{A'}\coring{C}')_{A'}$-pure,
\item $\omega_{L,M}$ is
$_AN$-, $_A(\coring{C}\tensor{A}N)$-pure,
\item $\omega_{L\tensor{A'}M,N}$ is $_{A''}\coring{C}''$-,
$\coring{C}'''_{A'''}$-pure,
\item $\omega_{L,M\cotensor{\coring{C}}N}$ is $_{A''}\coring{C}''$-,
$\coring{C}'''_{A'''}$-pure,
\end{enumerate}
then there is a canonical isomorphism of
$\coring{C}'''-\coring{C}'$-bicomodules
$$L\cotensor{\coring{C}'}(M\cotensor{\coring{C}}N)\simeq
(L\cotensor{\coring{C}'}M)\cotensor{\coring{C}}N.$$
\end{proposition}

\begin{proof}
Consider the following commutative diagram
$$\xymatrix{0 \ar[r] &
L\cotensor{\coring{C}'}(M\cotensor{\coring{C}}N)\ar[r]
\ar@{-->}[d]^\psi & L\tensor{A'}(M\cotensor{\coring{C}}N)
\ar[rr]^-{\omega_{L,M\cotensor{\coring{C}}N}}\ar[d]^{\psi_{L,M,N}}
&& L\tensor{A'}\coring{C}'\tensor{A'}(M\cotensor{\coring{C}}N)
\ar[d]^{\psi_{L\tensor{A'}\coring{C}',M,N}}
\\0 \ar[r] & (L\cotensor{\coring{C}'}M)\cotensor{\coring{C}}N \ar[r]
& (L\tensor{A'}M)\cotensor{\coring{C}}N
\ar[rr]^-{\omega_{L,M}\cotensor{\coring{C}}N}&&
(L\tensor{A'}\coring{C}'\tensor{A'}M)\cotensor{\coring{C}}N.}$$

The top row is exact from the definition of the cotensor product and
the second one is exact from Remark \ref{leftexactcotenfun} (see
Condition (b)). Then there is a unique map
$\psi:L\cotensor{\coring{C}'}(M\cotensor{\coring{C}}N)\to
(L\cotensor{\coring{C}'}M)\cotensor{\coring{C}}N$ making the above
diagram commutative. From Lemma \ref{assotenscotens} (see Condition
(a)), $\psi_{L,M,N}$ and $\psi_{L\tensor{A'}\coring{C}',M,N}$ are
isomorphisms. Hence $\psi$ is also an isomorphism. From Lemmas
\ref{mor} and \ref{assotenscotens} (see Condition (c)),
$\psi_{L,M,N}$ is an isomorphism of
$\coring{C}'''-\coring{C}'$-bicomodules. Again from Lemma \ref{mor}
(see Condition (d)), we obtain that $\psi$ is an isomorphism of
$\coring{C}'''-\coring{C}'$-bicomodules.
\end{proof}

\begin{proposition}\label{assocotens1}
All the assumed conditions in Proposition \ref{assocotens} are
fulfilled if at least one of the following holds
\begin{enumerate}[(1)]
\item $_A\coring{C}$, $\coring{C}'_{A'}$, $_{A''}\coring{C}''$,
$\coring{C}'''_{A'''}$, $L_{A'}$ and $_AN$ are flat;
\item $_A\coring{C}$, $\coring{C}'_{A'}$, $_{A''}\coring{C}''$,
$\coring{C}'''_{A'''}$ are flat, $N$ is coflat in
$\lcomod{\coring{C}}$;
\item $_A\coring{C}$, $\coring{C}'_{A'}$, $_{A''}\coring{C}''$,
$\coring{C}'''_{A'''}$ are flat, $L$ is coflat in
$\rcomod{\coring{C}'}$;
\item $_A\coring{C}$, $\coring{C}'_{A'}$, $_{A''}\coring{C}''$,
$\coring{C}'''_{A'''}$ are flat, $N$ is $(\coring{C},A)$-injective
(e.g. if $\coring{C}$ is coseparable), and $L$ is
$(\coring{C}',A')$-injective (e.g. if $\coring{C}'$ is coseparable);
\item If $_A\coring{C}$ is flat and $\coring{C}$ and
$\coring{C}'$ are coseparable.
\end{enumerate}
\end{proposition}

\begin{proof}
Immediate from the aforementioned results.
\end{proof}

\begin{proposition}\label{assocotens2}
Let
 $M\in{}\bcomod{\coring{C}'}{\coring{C}}$,
$N\in{}\bcomod{\coring{C}}{\coring{C}''}$ and
$L\in{}\bcomod{\coring{C}'''}{\coring{C}'}$. Then the canonical
isomorphism of $\coring{C}'''-\coring{C}'$-bicomodules
$$L\cotensor{\coring{C}'}(M\cotensor{\coring{C}}N)\simeq
(L\cotensor{\coring{C}'}M)\cotensor{\coring{C}}N$$ is natural in
$L$, $M$ and $N$ provided that it is well defined.
\end{proposition}

\begin{proof}
For example let $f:L\to L'$ be a morphism in
$\bcomod{\coring{C}'''}{\coring{C}'}$. Consider the diagram
$$\xymatrix@C 2pt{
 && L\cotensor{\coring{C}'}(M\cotensor{\coring{C}}N)
  \ar[rrr]^\simeq \ar[dll]_{f\cotensor{\coring{C}'}
  (M\cotensor{\coring{C}}N)}\ar'[d][ddd]&&&
  (L\cotensor{\coring{C}'}M)\cotensor{\coring{C}}N
   \ar[dll]^{(f\cotensor{\coring{C}'}M)\cotensor{\coring{C}}N}\ar[ddd]
\\
L'\cotensor{\coring{C}'}(M\cotensor{\coring{C}}N)\ar[rrr]^(.35)\simeq
\ar[ddd] &&&
(L'\cotensor{\coring{C}'}M)\cotensor{\coring{C}}N\ar[ddd]
\\
&&&&&&&&&
\\
 && L\tensor{A'}(M\cotensor{\coring{C}}N)
 \ar'[r][rrr]^\simeq  \ar[dll]^{f\tensor{A'}(M\cotensor{\coring{C}}N)}
  &&& (L\tensor{A'}M)\cotensor{\coring{C}}N
 \ar[dll]^{(f\tensor{A'}M)\cotensor{\coring{C}}N}
\\
L'\tensor{A'}(M\cotensor{\coring{C}}N) \ar[rrr]^\simeq  &&&
(L'\tensor{A'}M)\cotensor{\coring{C}}N.}$$ By Lemmas \ref{cube} and
\ref{assotenscotens}, the top square is commutative.
\end{proof}

\section{Cohom functors and coendomorphism corings}\label{Cohom functors}

\begin{definition}
A bicomodule $N \in \bcomod{\coring{C}}{\coring{D}}$ is said to be
\emph{quasi-finite} as a right $\coring{D}$-comodule or
$(A,\coring{D})$-\emph{quasi-finite}, if the functor $-\tensor{A}N:
\rmod{A} \rightarrow \rcomod{\coring{D}}$ has a left adjoint
$\cohom{\coring{D}}{N}{-} : \rcomod{\coring{D}} \rightarrow
\rmod{A}$, and we call it the \emph{cohom functor} associated to
$N$.
\end{definition}

\begin{example}\label{quasifinite}
If $\coring{D}$ is the trivial $B$-coring, then a
$(\coring{C},B)$-bicomodule $N$ is $(A,B)$-quasi-finite if and only
if $_AN$ is finitely generated and projective. The cohom functor is
$-\tensor{B}\ldual{N}:\rmod{B} \rightarrow \rmod{A}$, where
$\ldual{N}=\hom{A}{N}{A}$. In Chapter 2 we will give a
generalization of this fact for bigraded modules.
\end{example}

\begin{remark}
Let $N \in \bcomod{\coring{C}}{\coring{D}}$. It is clear that, if
$_B\coring{D}$ is flat and $N$ is $(A,\coring{D})$-quasi-finite,
then $_AN$ is flat.
\end{remark}

Now, let $N \in \bcomod{\coring{C}}{\coring{D}}$ be a
$(A,\coring{D})$-quasi-finite comodule, and let
\begin{equation}\label{cohomiso}
\Phi_{X,Y}:\hom{A}{\cohom{\coring{D}}{N}{X}}{Y}\to
\hom{\coring{D}}{X}{Y\tensor{A}N},
\end{equation}
where $Y\in \rmod{A}$ and $X\in \rcomod{\coring{D}}$, be the natural
isomorphism giving the adjunction
$(-\tensor{A}N,\cohom{\coring{D}}{N}{-})$. Let
$\theta:1_{\rcomod{\coring{D}}}\to \cohom{\coring{D}}{N}{-}
\tensor{A}N$ be the unit of this adjunction. Then
$\Phi_{X,Y}(f)=(f\tensor{A}N)\theta_X$ for every $f\in
\hom{A}{\cohom{\coring{D}}{N}{X}}{Y}$.

Let $X\in \rcomod{\coring{D}}$. Then there is a unique $A$-linear
map $$\rho_{\cohom{\coring{D}}{N}{X}}: \cohom{\coring{D}}{N}{X}\to
\cohom{\coring{D}}{N}{X} \tensor{A}\coring{C}$$ such that
$$(\cohom{\coring{D}}{N}{X}\tensor{A}\lambda_N)\theta_X=
(\rho_{\cohom{\coring{D}}{N}{X}}\tensor{A}N)\theta_X.$$

\begin{proposition}
We have that
$(\cohom{\coring{D}}{N}{X},\rho_{\cohom{\coring{D}}{N}{X}})$ is a
right $\coring{C}$-comodule, and we obtain a functor
$$\cohom{\coring{D}}{N}{-} : \rcomod{\coring{D}}
\rightarrow \rcomod{\coring{C}}.$$
\end{proposition}

\begin{proof}
We have
\begin{eqnarray*}
 (\rho_{\cohom{\coring{D}}{N}{X}}\tensor{A}\coring{C}\tensor{A}N)
 (\rho_{\cohom{\coring{D}}{N}{X}}\tensor{A}N)\theta_X
 &=& (\rho_{\cohom{\coring{D}}{N}{X}}\tensor{A}\coring{C}\tensor{A}N)
 (\cohom{\coring{D}}{N}{X}\tensor{A}\lambda_N)\theta_X\\
 &=& (\rho_{\cohom{\coring{D}}{N}{X}}\tensor{A}\lambda_N)\theta_X,
\end{eqnarray*}
and
\begin{eqnarray*}
(\cohom{\coring{D}}{N}{X}\tensor{A}\Delta_\coring{C}\tensor{A}N)
(\rho_{\cohom{\coring{D}}{N}{X}}\tensor{A}N)\theta_X &=&
(\cohom{\coring{D}}{N}{X}\tensor{A}\Delta_\coring{C}\tensor{A}N)
(\cohom{\coring{D}}{N}{X}\tensor{A}\lambda_N)\theta_X\\
&=& \big(\cohom{\coring{D}}{N}{X}\tensor{A}
(\coring{C}\tensor{A}\lambda_N)\lambda_N\big)\theta_X
\\ &=& (\cohom{\coring{D}}{N}{X}\tensor{A}\coring{C}\tensor{A}\lambda_N)
(\rho_{\cohom{\coring{D}}{N}{X}}\tensor{A}N)\theta_X\\
&=& (\rho_{\cohom{\coring{D}}{N}{X}}\tensor{A}\lambda_N)\theta_X.
\end{eqnarray*}
Hence $(\rho_{\cohom{\coring{D}}{N}{X}}\tensor{A}\coring{C})
\rho_{\cohom{\coring{D}}{N}{X}}=(\cohom{\coring{D}}{N}{X}
\tensor{A}\Delta_\coring{C})\rho_{\cohom{\coring{D}}{N}{X}}.$

On the other hand,
\begin{eqnarray*}
(\cohom{\coring{D}}{N}{X}\tensor{A}\epsilon_\coring{C}\tensor{A}N)
(\rho_{\cohom{\coring{D}}{N}{X}}\tensor{A}N)\theta_X &=&
(\cohom{\coring{D}}{N}{X}\tensor{A}\epsilon_\coring{C}\tensor{A}N)
(\cohom{\coring{D}}{N}{X}\tensor{A}\lambda_N)\theta_X\\
&=& \theta_X.
\end{eqnarray*}
Hence $(\cohom{\coring{D}}{N}{X}\tensor{A}\epsilon_\coring{C})
\rho_{\cohom{\coring{D}}{N}{X}}=\cohom{\coring{D}}{N}{X}$.

Now let $f:X\to X'$ be a morphism in $\rcomod{\coring{D}}$. Now, we
will verify that $\cohom{\coring{D}}{N}{f}$ is
$\coring{C}$-colinear. First, since $\theta$ is a natural
transformation, $\theta_{X'}f=
(\cohom{\coring{D}}{N}{f}\tensor{A}N)\theta_X$. Then
\begin{eqnarray*}
(\rho_{\cohom{\coring{D}}{N}{X'}}\tensor{A}N)
(\cohom{\coring{D}}{N}{f}\tensor{A}N)\theta_X & =&
(\rho_{\cohom{\coring{D}}{N}{X'}}\tensor{A}N)\theta_{X'}f\\
&=& (\cohom{\coring{D}}{N}{X'}\tensor{A}\lambda_N)\theta_{X'}f\\
&=& (\cohom{\coring{D}}{N}{X'}\tensor{A}\lambda_N)
(\cohom{\coring{D}}{N}{f}\tensor{A}N)\theta_X \\
&=& (\cohom{\coring{D}}{N}{f}\tensor{A}\lambda_N)\theta_X.
\end{eqnarray*}
On the other hand,
\begin{eqnarray*}
(\cohom{\coring{D}}{N}{f}\tensor{A}\coring{C}\tensor{A}N)
(\rho_{\cohom{\coring{D}}{N}{X}}\tensor{A}N)\theta_X &=&
(\cohom{\coring{D}}{N}{f}\tensor{A}\coring{C}\tensor{A}N)
(\cohom{\coring{D}}{N}{X}\tensor{A}\lambda_N)\theta_X\\
&=& (\cohom{\coring{D}}{N}{f}\tensor{A}\lambda_N)\theta_X.
\end{eqnarray*}
Hence $\rho_{\cohom{\coring{D}}{N}{X'}}\cohom{\coring{D}}{N}{f}=
(\cohom{\coring{D}}{N}{f}\tensor{A}\coring{C})
\rho_{\cohom{\coring{D}}{N}{X}}.$
\end{proof}

\begin{lemma}\label{cohomind}
Let $X\in \rcomod{\coring{D}}$, $Y\in \rcomod{\coring{C}}$, and
$f\in \hom{A}{\cohom{\coring{D}}{N}{X}}{Y}$. Then $f$ is
$\coring{C}$-colinear if and only if
$\oper{Im}{(f\tensor{A}N)\theta_X}\subset Y\cotensor{\coring{C}}N.$
\end{lemma}

\begin{proof}
We have,
\begin{eqnarray*}
(f\tensor{A}\coring{C}\tensor{A}N)
(\rho_{\cohom{\coring{D}}{N}{X}}\tensor{A}N)\theta_X &=&
(f\tensor{A}\coring{C}\tensor{A}N)(\cohom{\coring{D}}{N}{f}
\tensor{A}\lambda_N)\theta_X \\
&=& (Y\tensor{A}\lambda_N)(f\tensor{A}N)\theta_X.
\end{eqnarray*}
Hence,
\begin{eqnarray*}
f \textrm{ is } \coring{C}\textrm{-colinear} &\Longleftrightarrow &
(\rho_Y\tensor{A}N)(f\tensor{A}N)\theta_X=
(f\tensor{A}\coring{C}\tensor{A}N)
(\rho_{\cohom{\coring{D}}{N}{X}}\tensor{A}N)\theta_X\\
&\Longleftrightarrow & (\rho_Y\tensor{A}N)(f\tensor{A}N)\theta_X=
(Y\tensor{A}\lambda_N)(f\tensor{A}N)\theta_X\\
&\Longleftrightarrow& \oper{Im}{(f\tensor{A}N)\theta_X}\subset
Y\cotensor{\coring{C}}N.
\end{eqnarray*}
\end{proof}

\begin{proposition}\label{adjointcohomcotens}
Let $N\in\bcomod{\coring{C}}{\coring{D}}$. Assume that
$\omega_{Y,N}:=\rho_Y\otimes_AN-Y\otimes_A\lambda_N$ is
$_B(\coring{D} \tensor{B}\coring{D})$-pure for every right
$\coring{C}$-comodule $Y$ (e.g., $_B\coring{D}$ is flat or
$\coring{C}$ is coseparable).
\begin{enumerate}[(1)]
\item If $N$ is $(A,\coring{D})$-quasi-finite then
$-\cotensor{\coring{C}} N:\rcomod{\coring{C}}\to
\rcomod{\coring{D}}$ has a left adjoint, which is
$\cohom{\coring{D}}{N}{-}$.
\item If $-\cotensor{\coring{C}} N:\rcomod{\coring{C}}\to
\rcomod{\coring{D}}$ has a left adjoint, then $N$ is
$(A,\coring{D})$-quasi-finite.
\end{enumerate}
\end{proposition}

\begin{proof}
(1) From Lemmas \ref{cohomind} and \ref{mor}, the isomorphism
\eqref{cohomiso} induces an isomorphism
$\Phi_{X,Y}:\hom{\coring{C}}{\cohom{\coring{D}}{N}{X}}{Y}\to
\hom{\coring{D}}{X}{Y\cotensor{\coring{C}}N}.$ Hence (1) follows.

(2) From Lemma \ref{caniso}, $(-\tensor{A}\coring{C})
\cotensor{\coring{C}}N\simeq -\tensor{A}N:\rmod{A}\to
\rcomod{\coring{D}}$. We know that $-\tensor{A}\coring{C}$ is a
right adjoint to the forgetful functor $U:\rcomod{\coring{C}}\to
\rmod{A}$. Hence $-\tensor{A}N$ has a left adjoint (by Proposition
\ref{compoadjfun}).
\end{proof}

\begin{definition}
A bicomodule $N \in \bcomod{\coring{C}}{\coring{D}}$ is said to be
\emph{injector} as a right $\coring{D}$-comodule or
$(A,\coring{D})$-\emph{injector} if the functor $-\tensor{A}N:
\rmod{A} \to \rcomod{\coring{D}}$ preserves injective objects.
\end{definition}

By Proposition \ref{right-adjoint-inj-obj}, if $_A\coring{C}$ is
flat, then $\coring{C}$ is $(A,\coring{C})$-injector.

\begin{proposition}\label{injector}
Let $N\in\bcomod{\coring{C}}{\coring{D}}$ be a bicomodule with
$_A\coring{C}$ and $_B\coring{D}$ are flat. Suppose that $N$ is
$(A,\coring{D})$-quasi-finite. Then the following are equivalent
\begin{enumerate}[(1)]
\item $N$ is $(A,\coring{D})$-injector;
\item the cohom functor $h_{\coring{D}}(N,-):
\rcomod{\coring{D}}\to \rmod{A}$ is exact;
\item the cohom functor $h_{\coring{D}}(N,-):
\rcomod{\coring{D}}\to \rcomod{\coring{C}}$ is exact;
\item the functor $-\cotensor{\coring{C}}N:
\rcomod{\coring{C}}\to \rcomod{\coring{D}}$ preserves injective
objects.
\end{enumerate}
\end{proposition}

\begin{proof}
Follows directly from Propositions \ref{right-adjoint-inj-obj},
\ref{adjointcohomcotens}.
\end{proof}

\medskip
Now, let $N \in {}_A\rcomod{\coring{D}}$ be
$(A,\coring{D})$-quasi-finite. Set
$\e{\coring{D}}{N}=\cohom{\coring{D}}{N}{N}$. From Subsection
\ref{Frobeniusgeneral}, $\e{\coring{D}}{N}$ is an $A$-bimodule.
From Lemma \ref{$T$-objects}, $\theta_N$ is a morphism
in $_A\rcomod{\coring{D}}$. Consider the map in $\rcomod{\coring{D}}$
$$\xymatrix{N\ar[r]^-{\theta_N} & \e{\coring{D}}{N}\tensor{A}N
\ar[rr]^-{\e{\coring{D}}{N}\tensor{A}\theta_N} &&
\e{\coring{D}}{N}\tensor{A}\e{\coring{D}}{N}\tensor{A}N.}$$ Then
there is a unique right $A$-linear map
$$\Delta_e:\e{\coring{D}}{N}\to
\e{\coring{D}}{N}\tensor{A}\e{\coring{D}}{N}$$ such that
$$\Phi_{N,\e{\coring{D}}{N}\tensor{A}\e{\coring{D}}{N}}(\Delta_e)=
(\Delta_e\tensor{A}N)\theta_N=
(\e{\coring{D}}{N}\tensor{A}\theta_N)\theta_N.$$

Let $\rho^N:A\to\operatorname{End}_{\coring{D}}(N)$ and
$\rho:A\to\operatorname{End}_{\coring{C}}(\e{\coring{D}}{N})$ be the
maps defining the left $A$-action on $N$ and $\e{\coring{D}}{N}$,
respectively.

We want to verify that $\Delta_e$ is left $A$-linear, that is, for
every $a\in A$, the diagram in $\rmod{A}$
$$\xymatrix{\e{\coring{D}}{N}\ar[rr]^-{\Delta_e}\ar[d]_{\rho(a)} &&
\e{\coring{D}}{N}\tensor{A}\e{\coring{D}}{N}
\ar[d]^{\rho(a)\tensor{A}\e{\coring{D}}{N}} \\
\e{\coring{D}}{N}\ar[rr]^-{\Delta_e} &&
\e{\coring{D}}{N}\tensor{A}\e{\coring{D}}{N}}$$ commutes.

Let $a\in A$. We have the commutativity of the following diagram
\begin{equation}\label{e}
\xymatrix{N\ar[rr]^-{\theta_N}\ar[d]_{\rho^N(a)} &&
\e{\coring{D}}{N}\tensor{A}N
\ar[d]^{\rho(a)\tensor{A}N} \\
N\ar[rr]^-{\theta_N} && \e{\coring{D}}{N}\tensor{A}N.}
\end{equation}
Then
\begin{eqnarray*}
(\rho(a)\tensor{A}\e{\coring{D}}{N}\tensor{A}N)
(\Delta_e\tensor{A}N)\theta_N &=&
(\rho(a)\tensor{A}\theta_N)\theta_N
\\
&=&
(\e{\coring{D}}{N}\tensor{A}\theta_N)
(\rho(a)\tensor{A}N)\theta_N
\\
&=& (\e{\coring{D}}{N}\tensor{A}\theta_N)
\theta_N\rho^N(a),
\end{eqnarray*}
and
\begin{eqnarray*}
(\Delta_e\tensor{A}N)(\rho(a)\tensor{A}N)\theta_N &=&
(\Delta_e\tensor{A}N)\theta_N\rho^N(a)
\\ &=&
(\e{\coring{D}}{N}\tensor{A}\theta_N) \theta_N\rho^N(a).
\end{eqnarray*}
Hence $\Delta_e$ is left $A$-linear.

There is also a unique right $A$-linear map
$$\epsilon_e:\e{\coring{D}}{N}\to A$$
such that
$$\Phi_{N,A}(\epsilon_e)=(\epsilon_e\tensor{A}N)\theta_N=N.$$

Now we want to verify that $\epsilon_e$ is left $A$-linear, that is,
for every $a\in A$, the diagram in $\rmod{A}$
$$\xymatrix{\e{\coring{D}}{N}\ar[rr]^-{\epsilon_e}\ar[d]_{\rho(a)} &&
A\ar[d]^{a-} \\
\e{\coring{D}}{N}\ar[rr]^-{\epsilon_e} && A}$$ commutes.

Indeed, for every $a\in A$, we have
\begin{eqnarray*}
(\epsilon_e\tensor{A}N)(\rho(a)\tensor{A}N)\theta_N &=&
(\epsilon_e\tensor{A}N)\theta_N\rho^N(a)
\\ &=&
\rho^N(a),
\end{eqnarray*}
and $({a-}\tensor{A}N)(\epsilon_e\tensor{A}N)\theta_N=\rho^N(a).$

\begin{proposition}\label{BW 23.8}
Let $N \in {}_A\rcomod{\coring{D}}$ be
$(A,\coring{D})$-quasi-finite. Set
$\e{\coring{D}}{N}=\cohom{\coring{D}}{N}{N}$. Then
$(\e{\coring{D}}{N},\Delta_e,\epsilon_e)$ is an $A$-coring, where
$\Delta_e$ and $\epsilon_e$ defined above. Moreover, $N$ is a
$(\e{\coring{D}}{N},\coring{D})$-bicomodule by
$\theta_N:N\to\e{\coring{D}}{N}\tensor{A}N$, and there is an algebra
anti-isomorphism
$$\xymatrix{\rdual{\e{\coring{D}}{N}}=\hom{A}{\e{\coring{D}}{N}}{A}
\ar[r]^-{\Phi_{N,A}}& \operatorname{End}_\coring{D}(N)}.$$
\end{proposition}

\begin{proof}
First set $\coring{C}=\e{\coring{D}}{N}$. The coassociativity of
$\Delta_e$ follows from
\begin{eqnarray*}
(\coring{C}\tensor{A}\Delta_e\tensor{A}N)
(\Delta_e\tensor{A}N)\theta_N
&=&
(\coring{C}\tensor{A}\Delta_e\tensor{A}N)
(\coring{C}\tensor{A}\theta_N)\theta_N
\\ &=&
(\coring{C}\tensor{A}\coring{C}\tensor{A}\theta_N)
(\coring{C}\tensor{A}\theta_N)\theta_N,
\end{eqnarray*}
and
\begin{eqnarray*}
(\Delta_e\tensor{A}\coring{C}\tensor{A}N)
(\Delta_e\tensor{A}N)\theta_N
&=&
(\Delta_e\tensor{A}\coring{C}\tensor{A}N)
(\coring{C}\tensor{A}\theta_N)\theta_N
\\ &=&
(\Delta_e\tensor{A}\theta_N)\theta_N
\\ &=&
(\coring{C}\tensor{A}\coring{C}\tensor{A}\theta_N)
(\coring{C}\tensor{A}\theta_N)\theta_N.
\end{eqnarray*}
The counit property follows from
\begin{eqnarray*}
(\epsilon_e\tensor{A}\coring{C}\tensor{A}N)
(\Delta_e\tensor{A}N)\theta_N
&=&(\epsilon_e\tensor{A}\coring{C}\tensor{A}N)
(\coring{C}\tensor{A}\theta_N)\theta_N
\\&=&
(A\tensor{A}\theta_N)(\epsilon_e\tensor{A}N)\theta_N
\\&=&
(A\tensor{A}\theta_N)N
\\&=& \theta_N,
\end{eqnarray*}
and
\begin{eqnarray*}
(\coring{C}\tensor{A}\epsilon_e\tensor{A}N)
(\Delta_e\tensor{A}N)\theta_N &=&
(\coring{C}\tensor{A}\epsilon_e\tensor{A}N)
(\coring{C}\tensor{A}\theta_N)\theta_N
\\&=&
\big(\coring{C}\tensor{A}
(\epsilon_e\tensor{A}N)\theta_N\big)\theta_N
\\&=& \theta_N.
\end{eqnarray*}
Finally, we know that $\theta_N:N\to\coring{C}\tensor{A}N$ is a
morphism in $_A\rcomod{\coring{D}}$. Then, by the definitions of
$\Delta_e$ and $\epsilon_e$, $N$ is a
$(\coring{C},\coring{D})$-bicomodule.

Let $f,g\in\rdual{\coring{C}}$. Since $\theta_N$ is left $A$-linear,
$(f\tensor{A}\coring{C}\tensor{A}N)(\coring{C}\tensor{A}\theta_N)=
\theta_N(f\tensor{A}N).$ Hence
\begin{eqnarray*}
\Phi_{N,A}(f*^rg)
&=&
(g\tensor{A}N)(f\tensor{A}\coring{C}\tensor{A}N)
(\Delta_e\tensor{A}N)\theta_N
\\&=&
(g\tensor{A}N)(f\tensor{A}\coring{C}\tensor{A}N)
(\coring{C}\tensor{A}\theta_N)\theta_N
\\&=&
(g\tensor{A}N)\theta_N(f\tensor{A}N)\theta_N
\\&=&
\Phi_{N,A}(g)\Phi_{N,A}(f).
\end{eqnarray*}
\end{proof}

Now, we recall from \cite{Sweedler:1975} the definition of a
morphism of corings.

\begin{definition}
Let $\coring{C}$ and $\coring{C}'$ be two $A$-corings. An
$A$-bilinear map $f:\coring{C}\to\coring{C}'$ is called a
\emph{morphism of corings} if the diagrams
$$\xymatrix{\coring{C}\ar[rr]^-f\ar[d]_{\Delta_\coring{C}} &&
\coring{C}'\ar[d]^{\Delta_{\coring{C}'}} \\
\coring{C}\tensor{A}\coring{C}\ar[rr]^-{f\tensor{A}f} &&
\coring{C}'\tensor{A'}\coring{C}'}
\qquad\xymatrix{\coring{C}\ar[rr]^-f\ar[drr]^{\epsilon_\coring{C}}
&& \coring{C}' \ar[d]^{\epsilon_{\coring{C}'}}
\\&& A}$$
are commutative.
\end{definition}

\begin{proposition}\label{KG 5.2}
Let $N \in \bcomod{\coring{C}}{\coring{D}}$ be
$(A,\coring{D})$-quasi-finite. Let
$f:\e{\coring{D}}{N}\to\coring{C}$ be the unique morphism in
$\rmod{A}$ such that
$$(f\tensor{A}N)\theta_N=\lambda_N.$$
Then $f$ is a morphism of $A$-corings.
\end{proposition}

\begin{proof}
First we must verify that $f$ is left $A$-linear, that is, for every
$a\in A$, the diagram in $\rmod{A}$
$$\xymatrix{\e{\coring{D}}{N}\ar[rr]^-f\ar[d]_{\rho(a)} &&
\coring{C}\ar[d]^{a-} \\
\e{\coring{D}}{N}\ar[rr]^-f && \coring{C}}$$ commutes, where
$\rho:A\to\operatorname{End}_{\coring{C}}(\e{\coring{D}}{N})$ is the
map defining the left $A$-action on $\e{\coring{D}}{N}$. We have for
every $a\in A$,
\begin{eqnarray*}
(f\tensor{A}N)(\rho(a)\tensor{A}N)\theta_N &=&
(f\tensor{A}N)\theta_N\rho^N(a) \textrm{ (see Diagram \ref{e})}
\\ &=&
\lambda_N\rho^N(a),
\end{eqnarray*}
and
$({a-}\tensor{A}N)(f\tensor{A}N)\theta_N=({a-}\tensor{A}N)\lambda_N.$
Hence, since $\lambda_N$ is left $A$-linear, $f$ is also left
$A$-linear.

Moreover,
\begin{eqnarray*}
(f\tensor{A}f\tensor{A}N)(\Delta_e\tensor{A}N)\theta_N
&=&
(f\tensor{A}f\tensor{A}N)
(\e{\coring{D}}{N}\tensor{A}\theta_N)\theta_N
\\&=&
(f\tensor{A}\lambda_N)\theta_N
\\&=&
(\coring{C}\tensor{A}\lambda_N)\lambda_N,
\end{eqnarray*}
and
\begin{eqnarray*}
(\Delta_\coring{C}\tensor{A}N)(f\tensor{A}N)\theta_N
&=&
(\Delta_\coring{C}\tensor{A}N)\lambda_N
\\&=&
(\coring{C}\tensor{A}\lambda_N)\lambda_N.
\end{eqnarray*}
Hence, $(f\tensor{A}f)\Delta_e=\Delta_\coring{C}f.$

Finally, we have
\begin{eqnarray*}
(\epsilon_{\coring{C}}\tensor{A}N)(f\tensor{A}N)\theta_N
&=&
(\epsilon_{\coring{C}}\tensor{A}N)\lambda_N
\\&=&
N,
\end{eqnarray*}
and $(\epsilon_e\tensor{A}N)\theta_N=N.$ Hence,
$\epsilon_\coring{C}f=\epsilon_e$.
\end{proof}

\begin{remark}\label{KG 5.2'}
By Lemma \ref{cohomind}, the map $f$ defined in Proposition \ref{KG
5.2} is right $\coring{C}$-colinear. Hence, if
$\omega_{Y,N}:=\rho_Y\otimes_AN-Y\otimes_A\lambda_N$ is
$_B(\coring{D} \tensor{B}\coring{D})$-pure for every right
$\coring{C}$-comodule $Y$ (e.g., $_B\coring{D}$ is flat or
$\coring{C}$ is coseparable), then
$f=\chi_\coring{C}\circ\cohom{\coring{D}}{N}{\lambda_N}$, where
$\chi:\cohom{\coring{D}}{N}{-\cotensor{\coring{C}}N}\to
 1_{\rcomod{\coring{C}}}$ is the counit of the adjunction
 $(\cohom{\coring{D}}{N}{-},-\cotensor{\coring{C}}N)$. Hence,
 Proposition \ref{KG 5.2} and this remark recover
 \cite[Proposition 5.2]{ElKaoutit/Gomez:2003}.
\end{remark}

\section{Induction functors, the category of entwining modules,
the category of graded modules}

\subsection{Induction functors}\label{IndunctionFunctors}

We start this subsection by recalling from \cite{Gomez:2002} and
\cite{Brzezinski/Wisbauer:2003} the following

\begin{definition}\label{coring}
A coring homomorphism from the coring $(\coring{C}:A)$ to the coring
$(\coring{D}:B)$ is a pair $(\varphi,\rho)$, where $\rho:A\to B$ is
a homomorphism of $k$-algebras and $\varphi:\coring{C}\to\coring{D}$
is a homomorphism of $A$-bimodules such that the diagrams
$$\xymatrix{\coring{C}\ar[rrr]^{\Delta_\coring{C}}\ar[dd]_\varphi
&&&
\coring{C}\tensor{A}\coring{C}\ar[dr]^{\varphi\tensor{A}\varphi}
\\ &&&& \coring{D}\tensor{A}\coring{D}
\ar[dl]^{\omega_{\coring{D},\coring{D}}}
\\\coring{D} \ar[rrr]^{\Delta_\coring{D}} &&&
\coring{D}\tensor{B}\coring{D}} \qquad\xymatrix{\coring{C}
\ar[rrr]^{\epsilon_\coring{C}}\ar[dd]_\varphi &&& A \ar[dd]^\rho
\\&&&&
\\\coring{D} \ar[rrr]^{\epsilon_\coring{D}} &&& B}$$
are commutative, where
$\omega_{\coring{D},\coring{D}}:\coring{D}\tensor{A}\coring{D}
\to\coring{D}\otimes_B\coring{D}$ is the canonical map induced by
$\rho:A\to B.$  Equivalently, the map
$$\xymatrix{\varphi':B\tensor{A}\coring{C}\tensor{A}B
\ar[rr]^{B\tensor{A}\varphi\tensor{A}B}&&
B\tensor{A}\coring{D}\tensor{A}B\ar[r]^-\omega &
B\tensor{B}\coring{D}\tensor{B}B\ar[r]^-\simeq & \coring{D}}$$
($\varphi'(b\tensor{A}c\tensor{A}b')=b\varphi(c)b'$) is a morphism
of $B$-corings, where $\omega$ is the canonical map and
$B\tensor{A}\coring{C}\tensor{A}B$ is the base ring extension of
$\coring{C}$ (See Example \ref{examples corings} (3)).
\end{definition}

\begin{example}
Let $\rho:A\to B$ be a morphism of $k$-algebras and $\coring{C}$ an
$A$-coring. Consider the map $\varphi:\coring{C}\to
B\tensor{A}\coring{C}\tensor{A}B,\; c\mapsto 1_B\otimes c\otimes
1_B$. Then $(\varphi,\rho):(\coring{C}:A)\to
(B\tensor{A}\coring{C}\tensor{A}B:B)$ is a morphism of corings.
\end{example}

Let $(\varphi,\rho):(\coring{C}:A)\to(\coring{D}:B)$ be a coring
morphism. Now we introduce the \emph{induction functor}
$$-\tensor{A}B:\rcomod{\coring{C}} \to \rcomod{\coring{D}}.$$
It assigns to a comodule $M\in \rcomod{\coring{C}}$, $M\tensor{A}B$
which is a right $\coring{D}$-comodule by the coaction
\[
\rho_{M \tensor{A} B}(m \tensor{A}b) = \sum m_{(0)} \tensor{A} 1_B
\tensor{B} \varphi(m_{(1)})b,
\]
where $\rho_M(m)=\sum m_{(0)}\tensor{A} m_{(1)}$, and assigns to a
morphism $f$ in $\rcomod{\coring{C}}$, the map $f\tensor{A}B$.
Analogously the induction functor $B\tensor{A}-:\lcomod{\coring{C}}
\to\lcomod{\coring{D}}$ is defined.

\begin{definition}
We say that $(\varphi,\rho)$ is a \emph{pure morphism of corings}
\cite[24.8]{Brzezinski/Wisbauer:2003}, if
$\omega_{Y,B\tensor{A}\coring{C}}$ is
${}_A(\coring{C}\tensor{A}\coring{C})$-pure for every right
$\coring{D}$-comodule $Y$,  where the left
$\coring{D}$-coaction on the $B-\coring{C}$-bicomodule
$B\tensor{A}\coring{C}$ is given by:
$$\lambda_{B\tensor{A}\coring{C}}:B\tensor{A}\coring{C}\to
\coring{D}\tensor{B}B\tensor{A}\coring{C}\simeq
\coring{D}\tensor{A}\coring{C},\quad b\tensor{A}c\mapsto \sum
b\varphi(c_{(1)})\tensor{A}c_{(2)},$$ with
$\Delta_{\coring{C}}(c)=\sum c_{(1)}\tensor{A}c_{(2)}$.
\end{definition}

Let $(\varphi,\rho):(\coring{C}:A)\to(\coring{D}:B)$ be a pure
morphism of corings. We also define the \emph{coinduction functor}
$$-\cotensor{\coring{D}}(B\tensor{A}\coring{C}):
\rcomod{\coring{D}}\rightarrow \rcomod{\coring{C}}.$$
\medskip
The following result is given in \cite[Proposition 5.4]{Gomez:2002}
or \cite[24.11]{Brzezinski/Wisbauer:2003}, but there is a missing
condition, in both results, that guaranties that the coinduction
functor is well defined.

\begin{proposition}\emph{}\label{BW 24.11}
If $(\varphi,\rho)$ is a pure morphism of corings, then we have an
adjoint pair of functors
$$(-\otimes_AB,-\cotensor{\coring{D}}(B\tensor{A}\coring{C})).$$
\end{proposition}

From the proof of \cite[24.11]{Brzezinski/Wisbauer:2003}, the counit
of the adjunction
$(-\tensor{A}B,-\cotensor{\coring{D}}(B\tensor{A}\coring{C}))$ is
given by
$$\psi_N:(N\cotensor{\coring{D}}(B\tensor{A}\coring{C}))\tensor{A}B\to
N, \quad (\sum_in_i\tensor{A}c_i)\tensor{A}b\mapsto
\sum_in_i\rho(\epsilon_{\coring{C}}(c_i))b$$ $(N\in
\rcomod{\coring{D}})$. We obtain a commutative diagram
\begin{equation}\label{counit-ind-fun-diag}
\xymatrix{(\coring{D}\cotensor{\coring{D}}(B\tensor{A}\coring{C}))
\tensor{A}B\ar[rr]^-{\psi_\coring{D}} \ar[dr]_-\simeq && \coring{D}
\\ & B\tensor{A}\coring{C}\tensor{A}B\ar[ur]_-{\varphi'} &,}
\end{equation}
where $\varphi'$ is the map defined in Definition \ref{coring}.

\begin{proposition}\emph{\cite[23.9]{Brzezinski/Wisbauer:2003}}\label{BW 23.9}
\\Let $\rho:A\to B$ be a morphism of $k$-algebras and $\coring{C}$ an
$A$-coring such that $_A\coring{C}$ is flat. The coendomorphism
coring of the $(B,\coring{C})$-quasi-finite comodule
$B\tensor{A}\coring{C}$ is isomorphic as coring to the $B$-coring
$B\tensor{A}\coring{C}\tensor{A}B$.
\end{proposition}

\begin{proof}
Let $\varphi:\coring{C}\to B\tensor{A}\coring{C}\tensor{A}B$ be the
map defined by $\varphi(c)=1_B\otimes c\otimes 1_B$. By Example
1.5.2, $(\varphi,\rho):(\coring{C}:A)\to
(B\tensor{A}\coring{C}\tensor{A}B:B)$ is a morphism of coring. In
this case $\varphi'$ is the identity map. From the commutativity of
the diagram \eqref{counit-ind-fun-diag},
$\psi_{B\tensor{A}\coring{C}\tensor{A}B}$ is an isomorphism. Hence,
by Remark \ref{KG 5.2'}, there is an isomorphism of corings
$f:\e{\coring{C}}{B\tensor{A}\coring{C}}\to
B\tensor{A}\coring{C}\tensor{A}B$.
\end{proof}

\subsection{The category of entwined modules}\label{CategoryEntwined}

Entwining structures were introduced in \cite{Brzezinski/Majid:1998}
with the aim of preserving in noncommutative geometry symmetry
properties of principal bundles. Entwined modules over entwining
structures were introduced in \cite{Brzezinski:1999} as a
generalization of Doi-Koppinen Hopf modules and, in particular,
graded modules by $G$-sets, Hopf modules and Yetter-Drinfeld
modules.

Concerning the subject of this subsection and for the rest of this
work, we adopt the notations of \cite{Caenepeel/Militaru/Zhu:2002}.

\medskip
We recall from \cite{Brzezinski/Wisbauer:2003} that a
\emph{right-right entwining structure} over $k$ is a triple
$(A,C,\psi)$, where $A$ is a $k$-algebra, $C$ is a $k$-coalgebra,
and $\psi:C\otimes A\to A\otimes C$ is a $k$-linear map, such that
\begin{enumerate}[(ES1)]
\item $\psi\circ(1_C\otimes m)=(m\otimes 1_C)\circ(1_A\otimes
\psi)\circ(\psi\otimes 1_A)$, or equivalently, for all $a,b\in A$,
$c\in C$, $\sum(ab)_\psi\otimes c^\psi=\sum a_\psi b_\Psi\otimes
c^{\psi\Psi}$,
\item $(1_A\otimes \Delta)\circ \psi=(\psi\otimes1_C)\circ (1_C\otimes
\psi)\circ(\Delta\otimes 1_A)$, or equivalently, for all $a\in A$,
$c\in C$, $\sum a_\psi\otimes \Delta(c^\psi)=\sum
a_{\psi\Psi}\otimes c_{(1)}^\Psi\otimes c_{(2)}^\psi$,
\item $\psi\circ(1_C\otimes\eta)=\eta\otimes1_C$, or equivalently, for
all $c\in C$, $\sum(1_A)_\psi\otimes c^\psi=1_A\otimes c$,
\item $(1_A\otimes\epsilon)\circ\psi=\epsilon\otimes1_A$, or
equivalently, for all $a\in A$, $c\in C$, $\sum
a_\psi\epsilon(c^\psi)=a\epsilon(c)$.
\end{enumerate}
where $m$ and $\eta$ are respectively the multiplication and the
unit maps of $A$, and $\psi(c\otimes a)=\sum a_\psi\otimes
c^\psi=\sum a_\Psi\otimes c^\Psi.$ (Sometimes the symbol $\sum$
will be omitted.)

We say that $A$ and $C$ are \emph{entwined} by $\psi$ while $\psi$
is called the \emph{entwining map}.

The conditions (ES1) to (ES4) mean that the following diagram is
commutative

$$
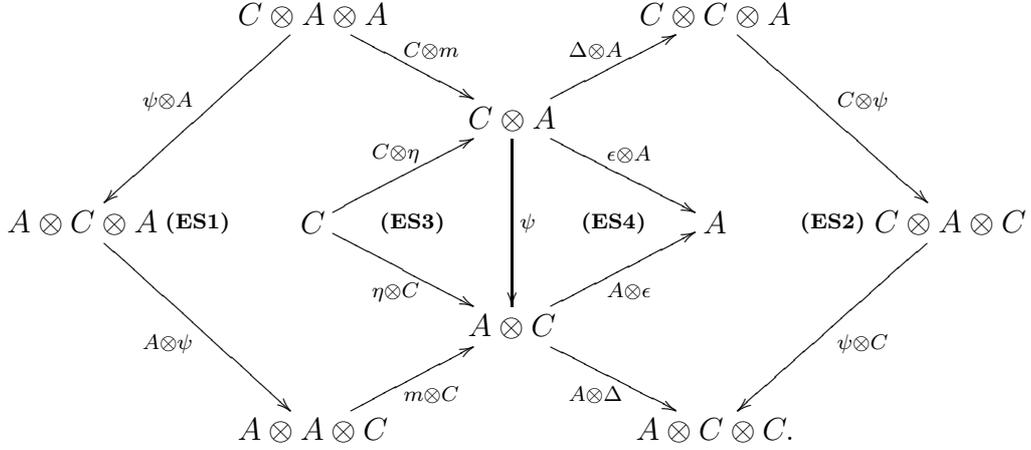
\xymatrix{
& C\otimes A\otimes A\ar[dr]^{C\otimes m} \ar[ddl]_{\psi\otimes A}
&& C\otimes C\otimes A\ar[ddr]^{C\otimes \psi}
\\
&& C\otimes A\ar[ur]^{\Delta\otimes A}\ar[dr]^{\epsilon\otimes A}
\ar[dd]^\psi &&
\\
\ar@{}[r]|{\textbf{(ES1)}} A\otimes C\otimes A \ar[ddr]_{A\otimes
\psi} & \ar@{}[r]|{\textbf{(ES3)}} C \ar[ur]^{C\otimes \eta}
\ar[dr]_{\eta\otimes C} && \ar@{}[l]|{\textbf{(ES4)}} A &
\ar@{}[l]|{\textbf{(ES2)}} C\otimes A\otimes C\ar[ddl]^{\psi\otimes
C}
\\
&& A\otimes C \ar[ur]_{A\otimes \epsilon} \ar[dr]_{A\otimes
\Delta}&&
\\
& A\otimes A\otimes C \ar[ur]_{m\otimes C} && A\otimes C\otimes C
.}$$

A \emph{morphism} $(A,C,\psi)\to(A',C',\psi')$ is a pair
$(\alpha,\gamma)$ with  $\alpha:A\to A'$ is a morphism of algebras,
and $\gamma:C\to C'$ is a morphism of coalgebras such that the
following diagram is commutative
$$\xymatrix{C\otimes A \ar[rr]^\psi\ar[d]_{\gamma\otimes\alpha} &&
A\otimes C\ar[d]^{\alpha\otimes\gamma}
\\ C'\otimes A'\ar[rr]^{\psi'} && A'\otimes C',}$$
or equivalently, for all $a\in A$, $c\in C$,
$$\sum\alpha(a_\psi)\otimes\gamma(c^\psi)=\sum\alpha(a)_{\psi'}\otimes
\gamma(c)^{\psi'}.$$ We denote this category by
$\mathbb{E}_\bullet^\bullet(k)$.

We can also define ($\tau$ and $\tau'$ are the twist maps)
\begin{itemize}
\item Left-right entwining structures over $k$: triples
$(A,C,\psi)$, where $\psi:A\otimes C\to A\otimes C$ such that
$(A^\circ,C,\psi\circ\tau)\in \mathbb{E}_\bullet^\bullet(k).$ We
denote this category by $_\bullet\mathbb{E}^\bullet(k)$.
\item Right-left entwining structures over $k$: triples
$(A,C,\psi)$, where $\psi:C\otimes A\to C\otimes A$ such that
$(A,C^\circ,\tau\circ\psi)\in \mathbb{E}_\bullet^\bullet(k).$ We
denote this category by $^\bullet\mathbb{E}_\bullet(k)$.
\item Left-left entwining structures over $k$: triples
$(A,C,\psi)$, where $\psi:A\otimes C\to C\otimes A$ such that
$(A^\circ,C^\circ,\tau'\circ\psi\circ\tau)\in
\mathbb{E}_\bullet^\bullet(k).$ We denote this category by
$_\bullet^\bullet\mathbb{E}(k)$.
\end{itemize}

\begin{proposition}
The categories $\mathbb{E}_\bullet^\bullet(k)$,
$_\bullet\mathbb{E}^\bullet(k)$, $^\bullet\mathbb{E}_\bullet(k)$,
and $_\bullet^\bullet\mathbb{E}(k)$ are isomorphic.
\end{proposition}

Again we recall from \cite{Brzezinski/Wisbauer:2003}, that a
\emph{right-right entwined modules} over a right-right entwining
structure $(A,C,\psi)$ is a $k$-module $M$ which is a right
$A$-module with multiplication $\psi_M$, and a right $C$-comodule
with comultiplication $\rho_M$ such that the following diagram is
commutative
$$\xymatrix{M\otimes A \ar[rrr]^{\rho_M\otimes A}\ar[dd]_{\psi_M}
&&& M\otimes C\otimes A\ar[dr]^{M\otimes \psi}
\\ &&&& M\otimes A\otimes C\ar[dl]^{\psi_M\otimes C}
\\M \ar[rrr]^{\rho_M} &&&
M\otimes C,}$$
or equivalently, for all $m\in M$ and $a\in A$,
$$\rho_M(ma)=\sum m_{(0)}a_\psi\otimes m_{(1)}^\psi.$$

A \emph{morphism} between entwined modules is a morphism of right
$A$-modules and right $C$-comodules at the same time. We denote this
category by $\mathcal{M}(\psi)_A^C$.

The following result shows that entwining structures and entwined
modules are very related to corings and comodules over corings
respectively.

\begin{theorem}\emph{[Takeuchi]}\label{Takeuchi}
\begin{enumerate}[(a)]
\item Let $A$ be an algebra, $C$ be a $k$-coalgebra, and let
$\psi:C\otimes A\to A\otimes C$ be a $k$-linear map. Obviously
$A\otimes C$ has a structure of left $A$-module by $b(a\otimes
c)=ba\otimes c$ for $a,b\in A$ and $c\in C$. Define the right
$A$-module action on $A\otimes C$, $$\xymatrix{\psi_{A\otimes
C}^r:A\otimes C\otimes A \ar[r]^-{A\otimes\psi} & A\otimes A\otimes
C \ar[r]^-{m\otimes C} & A\otimes C},$$ that is, $(a\otimes
c)b=a\psi(c\otimes b)$ for $a,b\in A$ and $c\in C$. Define also
$$\xymatrix{\Delta:A\otimes C\ar[r]^-{A\otimes \Delta_C} & A\otimes
C\otimes C\simeq (A\otimes C)\tensor{A}(A\otimes C)}$$ (for every
$a\in A,c\in C,$ $\Delta(a\otimes c)=\sum(a\otimes
c_{(1)})\tensor{A}(1_A\otimes c_{(2)})$, where $\Delta_C(c)=\sum
c_{(1)}\otimes c_{(2)}$, and
$$\xymatrix{\epsilon:A\otimes C\ar[r]^-{A\otimes \epsilon_C} &
A\otimes k\simeq A}$$ ($\epsilon(a\otimes c)=a\epsilon_C(c)$). Thus,
$(A\otimes C,\Delta,\epsilon)$ is an $A$-coring if and only if
$(A,C,\psi)$ is a right-right entwining structure.
\item Let $(A,C,\psi)$ be a right-right entwining
structure, $M$ be a right $A$-module with the action
$\psi^r_M:M\otimes A\to M$, and let $\rho_M:M\to M\otimes C$ be a
$k$-linear map. Define $$\xymatrix{\rho'_M:M\ar[r]^-{\rho_M} &
M\otimes C\simeq M\tensor{A}(A\otimes C)}$$ (for every $m\in M,$
$\rho'_M(m)=\sum m_{(0)}\tensor{A}(1_A\otimes m_{(1)})$, where
$\rho_M(m)=\sum m_{(0)}\otimes m_{(1)}$). Then, $(M,\rho'_M)$ is a
right $A\otimes C$-comodule if and only if $(M,\psi^r_M,\rho_M)$ is
a right-right entwined module over $(A,C,\psi)$. Hence, there is an
isomorphism of categories
$$\rcomod{A\otimes C}\simeq \mathcal{M}_A^C(\psi).$$
\end{enumerate}
\end{theorem}

\begin{proof}(Sketch)
(a) It is easy to verify that
\begin{itemize}
\item For all $a,b\in A,c\in C$, $(1\otimes
c)(ab)=((1\otimes c)a)b \Longleftrightarrow (\textrm{ES}1)$.
\item For all $c\in
C$, $(1\otimes c)1=1\otimes c \Longleftrightarrow (\textrm{ES}3)$.
\item For all $a\in A,c\in C$, $\Delta((1\otimes c)a)=\Delta(1\otimes
c)a \Longleftrightarrow (\textrm{ES}2)$. \item For all $a\in A,c\in
C$, $\epsilon((1\otimes c)a)=\epsilon(1\otimes c)a
\Longleftrightarrow (\textrm{ES}4)$.
\end{itemize}
Moreover, the coassociativity of $\Delta$ and the counit property of
$\epsilon$ follow from that of $\Delta_C$ and $\epsilon_C$. Hence
(a) follows.

(b) It is easy to verify that\begin{itemize}
\item $\rho'_M$ is $A$-linear if and only if for all $m\in M$ and $a\in
A$, $\rho_M(ma)=\sum m_{(0)}a_\psi\otimes m_{(1)}^\psi$.
\item $\rho'_M$ is coassociative if and only if $\rho_M$ is
coassociative.
\item The counit property of $\epsilon$ holds if and only if that of
$\epsilon_C$ holds.
 Hence (b) follows.
\end{itemize}
\end{proof}

The last theorem has a left-handed version:

\begin{theorem}\label{Takeuchi'}
\begin{enumerate}[(a)]
\item Let $A$ be an algebra, $C$ be a $k$-coalgebra, and let
$\varphi:A\otimes C\to C\otimes A$ be a $k$-linear map. Obviously
$C\otimes A$ has a structure of right $A$-module by $(c\otimes
a)b=c\otimes ab$ for $a,b\in A$ and $c\in C$. Define the left
$A$-module action on $C\otimes A$, $$\xymatrix{\psi_{C\otimes
A}^l:A\otimes C\otimes A \ar[r]^-{\varphi\otimes A} & C\otimes
A\otimes A \ar[r]^-{C\otimes m} & C\otimes A},$$ that is,
$b(c\otimes a)=\varphi(b\otimes c)a$ for $a,b\in A$ and $c\in C$.
Define also
$$\xymatrix{\Delta:C\otimes A\ar[r]^-{\Delta_C\otimes A} & C\otimes
C\otimes A\simeq (C\otimes A)\tensor{A}(C\otimes A)}$$ (for every
$a\in A,c\in C,$ $\Delta(c\otimes a)=\sum(c_{(1)}\otimes
1_A)\tensor{A}(c_{(2)}\otimes a)$, where $\Delta_C(c)=\sum
c_{(1)}\otimes c_{(2)}$, and
$$\xymatrix{\epsilon:C\otimes A\ar[r]^-{\epsilon_C\otimes A} &
k\otimes A\simeq A}$$ ($\epsilon(c\otimes a)=\epsilon_C(c)a$). Thus,
$(C\otimes A,\Delta,\epsilon)$ is an $A$-coring if and only if
$(A,C,\varphi)$ is a left-left entwining structure.
\item Let $(A,C,\psi)$ be a left-left entwining
structure, $M$ be a left $A$-module with the action
$\psi^l_M:A\otimes M\to M$, and let $\lambda_M:M\to C\otimes M$ be a
$k$-linear map. Define $$\xymatrix{\lambda'_M:M\ar[r]^-{\lambda_M} &
C\otimes M\simeq (C\otimes A)\tensor{A}M}$$ (for every $m\in M,$
$\lambda'_M(m)=\sum (m_{(-1)}\tensor{A}1_A)\otimes m_{(0)}$, where
$\rho_M(m)=\sum m_{(-1)}\otimes m_{(0)}$). Then, $(M,\lambda'_M)$ is
a left $C\otimes A$-comodule if and only if $(M,\psi^l_M,\lambda_M)$
is a left-left entwined module over $(A,C,\varphi)$. Hence, there is
an isomorphism of categories
$$\lcomod{C\otimes A}\simeq{}_A^C\mathcal{M}(\varphi).$$
\end{enumerate}
\end{theorem}

Let $(A,C,\psi)\in \mathbb{E}_\bullet^\bullet(k)$ be a right-right
entwining structure and $(B,D,\varphi)\in
{}_\bullet^\bullet\mathbb{E}(k)$ a left-left entwining structure. A
\emph{two-sided entwined module} \cite[pp.
68--69]{Caenepeel/Militaru/Zhu:2002} is a $k$-module $M$ that have a
structure of a left-left $(B,D,\varphi)$-module and a right-right
$(A,C,\psi)$-module such that the following extra conditions hold:
\begin{enumerate}[(1)]
\item $M$ is a $(B,A)$-bimodule;
\item $M$ is a $(D,C)$-bicomodule;
\item the right $A$-action is left $D$-colinear, that is,
$$\lambda_M(ma)=m_{(-1)}\otimes m_{(0)}a,\textrm{ for all } m\in M,a\in
A;$$
\item the left $B$-action is right $C$-colinear, that is,
$$\rho_M(bm)=bm_{(0)}\otimes m_{(1)},\textrm{ for all } m\in M,b\in B.$$
\end{enumerate}

We denote the category of two-sided entwined modules by
$_B^D\mathcal{M}(\varphi,\psi)_A^C$.

\subsection{The category of graded modules}\label{CategoryGraded}

For the rest of this work we adopt the notations of
\cite{Nastasescu/Raianu/VanOystaeyen:1990} and \cite{DelRio:1992}.

\medskip
In this subsection $G$ stand for a group, and $A$ for a $G$-graded
$k$-algebra (i.e. $A=\bigoplus_{g\in G}A_g$, where $A_g$, $g\in G$,
are $k$-submodules of $A$ such that $A_gA_\sigma\subset A_{g\sigma}$
for all $g,\sigma\in G$). We denote the unit of $G$ by $e$ and the
unit of $A$ by $1$.

Let $X$ be a right $G$-set (i.e. there is a map $X\times G\to
X,\;(x,g)\mapsto xg$ such that $(xg)\sigma=x(g\sigma)$ and $xe=x$
for all $x\in X, g,\sigma\in G$). Of course, $X$ is a left $G$-set
via $g.x=xg^{-1}$ for all $g\in G,x\in X$.

\medskip
We consider the category of $X$-graded right $A$-modules
$gr-(A,X,G)$. This category is introduced and studied in
\cite{Nastasescu/Raianu/VanOystaeyen:1990}.

A right $A$-module $M$ is called a \emph{$X$-graded right
$A$-module} if $M=\bigoplus_{x\in X}M_x$, where $M_x$, $x\in X$, is
$k$-submodules of $M$ such that $M_xA_g\subset M_{xg}$ for all $g\in
G,x\in X$.

Let $M$ and $N$ be two $X$-graded right $A$-modules. A morphism of
right $A$-modules $f:M\to N$ is called a \emph{morphism of
$X$-graded right $A$-modules} if $f(M_x)\subset N_x$ for all $x\in
X$. Notice that $\hom{gr-(A,X,G)}{M}{N}$ is a $k$-submodule of
$\hom{A}{M}{N}$.

Since $\rmod{A}$ is a $k$-category, $X$-graded right $A$-modules and
their morphisms define the $k$-category of $X$-graded right
$A$-modules. We denote it by $gr-(A,X,G)$.

Analogously, if $X$ be a left $G$-set, then we define the
$k$-category of $X$-graded left $A$-modules, and we denote it by
$(G,X,A)-gr$. For every $M\in (G,X,A)-gr$, we denote the $x$-th
component of $M$ by $_xM$, for all $x\in X$.

If $X=G$ with the obvious action, then $gr-(A,G,G)$ is the category
of graded right $A$-modules, and it is denoted by $gr-A$. If $X$ is
a singleton, then $gr-(A,X,G)=\rmod{A}$.

Let $M\in gr-A$ ($M\in A-gr$ and $X$ be a left $G$-set) and $x\in
X$. The $x$-th \emph{suspension} of $M$ is, by definition, the
$X$-graded right (left) $A$-module $M(x)$ ($(x)M$) that is equal to
$M$ as right (left) $A$-module endowed with the grading
$M(x)_y=\bigoplus\{M_g\mid g\in G, xg=y\}$
\big($_y(x)M=\bigoplus\{_gM\mid g\in G, gx=y\}$\big), $y\in X$.

By \cite[Theorem 2.8]{Nastasescu/Raianu/VanOystaeyen:1990},
$gr-(A,X,G)$ is a Grothendieck category with $\{A(x)\mid x\in X\}$
as a family of projective generators. We can derive that
$gr-(A,X,G)$ is a Grothendieck category from the isomorphism of
categories \eqref{graded entwined} and Theorem \ref{Takeuchi} (2).

\medskip
Now, let $G$ and $G'$ be two groups, $X$ a right $G$-set, $X'$ a
right $G'$-set, $A$ a $G$-graded algebra, and $A'$ a $G'$-graded
algebra.

\begin{definition}
Let $M\in gr-(A,X,G)$ and $N\in (G,X,A)-gr$. $M\widehat{\otimes}_AN$
is defined as the $k$-submodule of $M\tensor{A}N$ generated by the
elements $m\otimes n$ where $x\in X,m\in M_x, n\in {}_xN$.
\end{definition}

\begin{definition}
Let $Z$ be a $(G,G')$-set (i.e. $Z$ a right $G$-set and a right
$G'$-set such that $(g.z)g'=g.(zg')$ for all $g\in G,g'\in G',z\in
Z$). An $(A,A')$-bimodule $M$ is called a \emph{$Z$-graded
$(A,A')$-bimodule} if it has a $Z$-grading $M=\bigoplus_{z\in Z}M_z$
such that $A_gM_zA'_{g'}\subset M_{gzg'}$ for all $g\in G,g'\in
G',z\in Z$.
\end{definition}

Obviously, $X\times X'$ is a $(G,G')$-set via
$$g.(x,x')g'=(g.x,x'g'),\; g\in G,g'\in
G',x\in X,x'\in X'.$$

Let $P=\bigoplus_{(x,x')\in X\times X'}P_{(x,x')}$ be an $X\times
X'$-graded $(A,A')$-bimodule. Set $$_xP=\bigoplus_{x'\in X'}P_{(x,x')}
\quad \textrm{and}\quad P_{x'}=\bigoplus_{x\in X}P_{(x,x')}$$
for every $x\in X,x'\in X'.$

We have that, for every $x\in X$, $_xP$ is a submodule of $P_{A'}$
and $_xP\in gr-(A',X',G')$ with the $X'$-grading $(_xP)_{x'}=
P_{(x,x')}$, $x'\in X'$. Furthermore,
$P=\bigoplus_{x'\in X'}P_{x'}\in gr-(A',X',G')$, and $P$ is the
coproduct of the family $\{_xP,x\in X\}$ in $gr-(A',X',G')$.

Analogously, for every $x'\in X'$, $P_{x'}\in (G,X,A)-gr$ with
the $X$-grading $_x(P_{x'})=P_{(x,x')}$, $x\in X$, $P=
\bigoplus_{x\in X}{}_xP\in (G,X,A)-gr$, and $P$ is the coproduct
of the family $\{P_{x'},x'\in X'\}$ in $(G,X,A)-gr$.

\medskip
Now we will introduce two functors associated to $P$.

\begin{enumerate}[(1)]
\item \textbf{The functor} $-\widehat{\otimes}_AP:gr-(A,X,G)\to
gr-(A',X',G').$

Let $M\in gr-(A,X,G)$. $M\widehat{\otimes}_AP$ is a submodule of
$(M\tensor{A}P)_{A'}$. Then $M\widehat{\otimes}_AP\in gr-(A',X',G')$
with the $X'$-grading
$$(M\widehat{\otimes}_AP)_{x'}=M\widehat{\otimes}_AP_{x'},\; x'\in X'.$$
For a morphism $f:M\to N$ in $gr-(A,X,G)$, we define
$f\widehat{\otimes}_AP:m\otimes p\mapsto f(m)\otimes p.$
\item \textbf{The functor} $\h{P_{A'}}{-}:gr-(A',X',G')\to gr-(A,X,G).$

Let $M\in gr-(A',X',G')$. We have that $\hom{A'}{P}{M'}$ is a right
$A$-module via $(fa)(p)=f(ap),\;f\in \hom{A'}{P}{M'},a\in A, p\in
P.$ Consider the $A$-submodule of $\hom{A'}{P}{M'}$
$$\h{P_{A'}}{M'}=\{g\in \hom{gr-(A',X',G')}{P}{M'}\mid f(_xP)=0
\textrm{ for almost all } x\in X\}.$$ $\h{P_{A'}}{M'}$ can be viewed
as a $X$-graded right $A$-module via the grading
$$(\h{P_{A'}}{M'})_x=\big\{g\in \hom{gr-(A',X',G')}{P}{M'}\mid f(_yP)=0
\textrm{ for all } y\in X-\{x\}\big\},\;x\in X.$$ Let $f:M'\to N'$
be a morphism in $gr-(A',X',G')$. $\h{P_{A'}}{f}:\h{P_{A'}}{M'}\to
\h{P_{A'}}{N'}$ is defined by
$\h{P_{A'}}{f}(\varphi)=f\circ\varphi,\;\varphi\in \h{P_{A'}}{M'}.$
\end{enumerate}

\begin{proposition}\emph{\cite[Proposition 1.2]{Menini:1993}}\label{Menini 1.2}
\\For every $X\times X'$-graded $(A,A')$-bimodule $P$, we have an
adjunction $$(-\widehat{\otimes}_AP,\h{P_{A'}}{-}).$$
\end{proposition}

The unit and the counit of this adjunction are given respectively by
$$\eta_M:M\to \h{P_{A'}}{M\widehat{\otimes}_AP},$$
$\eta_M(m)(p)=\sum_{x\in X}m_x\tensor{A}{_xp}$ ($m=\sum_{x\in
X}m_x\in M,p=\sum_{x\in X}{}_xp\in P$), and
$$\varepsilon_N:\h{P_{A'}}{N}\widehat{\otimes}_AP\to N,$$
$\varepsilon_N(f\tensor{A}{P})=f(p)$ ($f\in \h{P_{A'}}{N}_x,p\in
{}_xP, x\in X$).

\medskip
Now, set, for $a\in A$, $\lambda_a=a-:A\to A$. We have, for every
$a\in A_g,x\in X$, $\lambda_a:A(x)\to A(g.x)$ is a morphism in
$gr-(A,X,G)$.

\begin{definition}
Let $F:gr-(A,X,G)\to gr-(A',X',G')$ be a (covariant) functor. The
\emph{$X\times X'$-graded $(A,A')$-bimodule associated to} $F$ is
the $X'$-graded right $A'$-module $P=\bigoplus_{x\in X}F(A(x))$
endowed with the following left $A$-module structure
$$ap=\big(F(\lambda_a)\big)(p),\;
g\in G,a\in A_g, p\in F(A(x)),$$
and the $X\times X'$-grading
$$P_{(x,x')}=F(A(x))_{x'},\; x\in X,x'\in X'.$$
\end{definition}

We will denote the $X\times X$-graded $(A,A)$-bimodule associated to
the identity functor $1_{gr-(A,X,G)}$ by $\widehat A$.

\medskip
For more details we refer to \cite{Menini:1993} and
\cite{DelRio:1992}.

\medskip
Now, we will recall the definition and some well known results of
bialgebras.

\begin{proposition}
Let $H$ be a $k$-algebra with the multiplication map $m_H$ and the
unit map $\eta_H$, and a $k$-coalgebra. The following are equivalent
\begin{enumerate}[(1)]
\item $m_H$ and $\eta_H$ are morphisms of coalgebras;
\item $\Delta_H$ and $\epsilon_H$ are morphisms of algebras;
\item for all $g,h\in H$,
\begin{eqnarray}
 \Delta(gh)&=& g_{(1)}h_{(1)}\otimes g_{(2)}h_{(2)}\label{bialg1}\\
 \epsilon(gh)&=& \epsilon(g)\epsilon(h)\label{bialg2}\\
 \Delta(1_H)&=& 1_H\otimes 1_H\label{bialg3}\\
 \epsilon(1_H) &=& 1_k\label{bialg4}.
\end{eqnarray}
\end{enumerate}
In this case, we say that $H$ is a \emph{bialgebra}. A map between
bialgebras over the same ring $k$ is called a \emph{morphism of
bialgebras} if it is both a morphism of algebras and a morphism of
coalgebras.
\end{proposition}

\begin{proof}
It is a direct consequence from the following remarks

$m_H$ is a morphism of coalgebras if and only if \eqref{bialg1} and
\eqref{bialg2} hold,

$\eta_H$ is a morphism of coalgebras if and only if \eqref{bialg3}
and \eqref{bialg4} hold,

$\Delta_H$ is a morphism of algebras if and only if \eqref{bialg1}
and \eqref{bialg3} hold,

$\epsilon_H$ is a morphism of algebras if and only if \eqref{bialg2}
and \eqref{bialg4} hold.
\end{proof}

\begin{definition}
Let $H$ be a bialgebra. We say that $H$ is a \emph{Hopf algebra} if
$1_H$ is invertible in the convolution algebra $\hom{k}{H}{H}$ (the
multiplication is defined by $f*g=m_H(f\otimes g)\Delta_H$ for all
$g,h\in \hom{k}{H}{H}$, and $\eta_H\epsilon_H$ as the unit map),
that is, there is $S\in\hom{k}{H}{H}$ such that
\begin{equation}
S(h_{(1)})h_{(2)}=h_{(1)}S(h_{(2)})=\eta_H(\epsilon_H(h)), \textrm{
for all } h\in H.
\end{equation}
$S$ is called the \emph{antipode} of $H$.
\end{definition}

Let $H$ and $K$ be two Hopf algebras and $f:H\to K$ be a morphism of
bialgebras. It is very-known that $f$ preserves the antipode, i.e.
$$S_Kf=fS_H.$$ In this case, we say that $f$ is a \emph{morphism of
Hopf algebras}.

\begin{example}
Let $G$ be a group. Then $kG$ is a $k$-algebra and a $k$-coalgebra
(see Example \ref{examples corings} (4)). Moreover $kG$ is a Hopf
algebra with the antipode $S(g)=g^{-1}$ for all $g\in G$.
\end{example}

To differentiate the opposite algebra from that of coalgebra, we
will denote the opposite coalgebra of a $k$-coalgebra $C$ by
$C^{c\circ}$.

If $H$ is a bialgebra, $H^\circ,H^{c\circ}$ and $H^{\circ c\circ}$
are also bialgebras. If $H$ has an antipode $S$, then $S$ is also an
antipode of $H^{\circ c\circ}$. Then, if $H^\circ$ has an antipode
$\overline S$, that is, $\overline S$ satisfies
\begin{equation}
h_{(2)}\overline S(h_{(1)})=\overline
S(h_{(2)})h_{(1)}=\eta_H(\epsilon_H(h)), \textrm{ for all } h\in H,
\end{equation}
then $\overline S$ is also an antipode of $H^{c\circ}$. $\overline
S$ is called a \emph{twisted antipode}.

\begin{proposition}\emph{\cite[Proposition
2]{Caenepeel/Militaru/Zhu:2002}}\label{antipode}
\\Let $H$ be a Hopf algebra. Then $S:H\to H^{\circ c\circ}$ is a
morphism of bialgebras. If $S$ is bijective, then $S^{-1}$ is a
twisted antipode. If $H$ is commutative or cocommutative, then
$S^2=1_H$. In particular, $S=\overline S$.
\end{proposition}

\begin{lemma}
Let $H$ be a bialgebra and $M$ and $N$ two right $H$-modules. Then
$M\otimes N$ is also a right $H$-module via
$$(m\otimes n)h=mh_{(1)}\otimes nh_{(2)}$$
for all $m\in M,n\in N, h\in H$.
\end{lemma}

\begin{proof}
Let $\psi_M:M\otimes H\to M$ and $\psi_N:N\otimes H\to N$ be the
right $H$-actions on $M$ and $N$, respectively. Consider the
$k$-linear map
$$\psi:\xymatrix{M\otimes N\otimes H\ar[rr]^-{M\otimes N\otimes\Delta_H}
&& M\otimes N\otimes H\otimes H\ar[r]^\simeq & M\otimes H\otimes
N\otimes H \ar[rr]^-{\psi_M\otimes\psi_N} && M\otimes N.}$$ Set
$(m\otimes n)h=\psi(m\otimes n\otimes h)=mh_{(1)}\otimes nh_{(2)}$
for all $m\in M,n\in N, h\in H$.

Now, let  $m\in M,n\in N, g,h\in H$. By \eqref{bialg1}, $(m\otimes
n)(gh)=((m\otimes n)g)h$. Finally, by \eqref{bialg3}, $(m\otimes
n)1_H=m\otimes n$.
\end{proof}

\begin{definition}
Let $H$ be a bialgebra. A $k$-module $C$ which is a $k$-coalgebra
and a right $H$-module is called a right \emph{$H$-module coalgebra}
if the comultiplication and the counit maps are right $H$-linear,
that is,
\begin{equation}
\Delta_C(ch)=c_{(1)}h_{(1)}\otimes c_{(2)}h_{(2)} \textrm{ and }
\epsilon_C(ch)=\epsilon_C(c)\epsilon_H(h)
\end{equation}
for all $h\in H,c\in C$.
\end{definition}

Analogously, we define left module coalgebras. Moreover, if $C$ is a
right $H$-module coalgebra, then $C^{c\circ}$ is a left $H^{\circ
c\circ}$-module coalgebra.

\begin{example}
Let $G$ be a group and $X$ a right $G$-set. Then the $k$-coalgebra
$kX$ is a right $kG$-module coalgebra. Indeed, $kX\otimes kG$ is a
free $k$-module with the basis $\{x\otimes g\mid x\in X,g\in G\}$.
Define $\psi_{kX}:kX\otimes kG\to kX$ such that $\psi_{kX}(x\otimes
g)=xg$ for all $x\in X,g\in G$. It is obvious that $(kX,\psi_{kX})$
is a right $kG$-module. From Example \ref{examples corings} (4),
$kX$ is also a $k$-coalgebra. It is easy to verify that $kX$ is a
right $kG$-module coalgebra.
\end{example}

Let $A$ be a $k$-algebra and $X$ a nonempty set. Let $G=\{e\}$ be
the trivial group. Then $A$ is trivially a $G$-graded algebra and
$X$ is trivially a right $G$-set. Moreover, $M\in gr-(A,X,G)$ if and
only if $M$ is a right $A$-module and $M=\bigoplus_{x\in X}M_x$ with
each $M_x$ is an $A$-submodule of $M$. In the following example we
will consider the particular case where $G=\{e\}$ and $A=k$.

\begin{example}\label{graded modules}
Let $C=kX$ be a grouplike coalgebra with $X$ is a nonempty set. Let
$M$ be a $X$-graded $k$-module, that is, $M=\bigoplus_{x\in X}M_x$
such that each $M_x$ is a $k$-submodule of $M$. Then $M$ is a right
$C$-comodule via $$\rho_M(m)=m\otimes x$$ for all $x\in X, m\in
M_x$. If $M\in \rcomod{C}$, then $M$ is an $X$-graded $k$-module via
the grading $$M_x=\{m\in M\mid \rho_M(m)=m\otimes x\},\;x\in X.$$
This defines an isomorphism from $\rcomod{C}$ to the category of
$X$-graded $k$-modules $gr-(k,X,\{e\})$.
\end{example}

\begin{lemma}
Let $H$ be a bialgebra and $M$ and $N$ are right $H$-comodules. Then
$M\otimes N$ is also a right $H$-comodule via
$$\rho_{M\otimes N}(m\otimes n)=m_{(0)}\otimes n_{(0)}\otimes
m_{(1)}n_{(1)}$$ for all $m\in M,n\in N$.
\end{lemma}

\begin{proof}
Let $\rho_M:M\to M\otimes H$ and $\rho_N:N\to N\otimes H$ be the
right $H$-coactions on $M$ and $N$, respectively, and $\eta$ the
multiplication map of $H$. Consider the $k$-linear map
$$\rho:\xymatrix{M\otimes N\ar[rr]^-{\rho_M\otimes \rho_N}
&& M\otimes H\otimes N\otimes H\ar[r]^\simeq & M\otimes N\otimes
H\otimes H \ar[rr]^-{M\otimes N\otimes\eta} && M\otimes N\otimes
H.}$$ Then $\rho(m\otimes n)=m_{(0)}\otimes n_{(0)}\otimes
m_{(1)}n_{(1)}$ for all $m\in M,n\in N$.

Finally, for all $m\in M,n\in N$,
\begin{eqnarray*}
 (M\otimes N\otimes\Delta)\rho(m\otimes n)
 &=& m_{(0)}\otimes n_{(0)}\otimes\Delta(m_{(1)}n_{(1)})
 \\
 &=& m_{(0)}\otimes n_{(0)}\otimes m_{(1)(1)}n_{(1)(1)}
 \otimes m_{(1)(2)}n_{(1)(2)}
 \textrm{ (by \eqref{bialg1}}) \\
 &=& m_{(0)(0)}\otimes n_{(0)(0)}\otimes m_{(0)(1)}n_{(0)(1)}
 \otimes m_{(1)}n_{(1)}\\
 &=&(\rho\otimes H)\rho(m\otimes n),
\end{eqnarray*}
and
\begin{eqnarray*}
 (M\otimes N\otimes\epsilon)\rho(m\otimes n)
 &=& m_{(0)}\otimes n_{(0)}\epsilon(m_{(1)}n_{(1)})
 \\&=&
 m_{(0)}\otimes n_{(0)}\epsilon(m_{(1)})\epsilon(n_{(1)})
 \textrm{ (by \eqref{bialg2}})
 \\&=& m\otimes n.
\end{eqnarray*}
\end{proof}

\begin{definition}
Let $H$ be a bialgebra. Since the map $\eta_H:k\to H$ is a morphism
of coalgebras, then $k$ is a right $H$-comodule with the coaction
$\xymatrix{k\ar[r]^-\simeq & k\otimes k\ar[r]^-{k\otimes\eta_H}&
k\otimes H}$.

A $k$-module $A$ which is a $k$-algebra and a right $H$-comodule is
called a right \emph{$H$-comodule algebra} if the multiplication and
the unit maps are right $H$-colinear, that is,
\begin{equation}
\rho_A(ab)=a_{(0)}b_{(0)}\otimes a_{(1)}b_{(1)} \textrm{ and }
\rho_A(1_A)=1_A\otimes 1_H
\end{equation}
for all $a,b\in A$.
\end{definition}

Analogously, we define left comodule algebras. Moreover, if $A$ is a
right $H$-comodule algebra, then $A^\circ$ is a left $H^{\circ
c\circ}$-comodule algebra.

\begin{example}
Let $G$ be a group, and $H=kG$ the associated Hopf algebra. To give
a $kG$-comodule algebra is the same as to give a $G$-graded
$k$-algebra.
\end{example}

\medskip
A right-right \emph{Doi-Koppinen structure or simply DK structure}
over $k$ \cite{Caenepeel/Militaru/Zhu:2002} is a triple $(H,A,C)$,
where $H$ is a bialgebra, $A$ is a right $H$-comodule algebra, and
$C$ is a right $H$-module coalgebra. A \emph{morphism of DK
structures} is a triple
$(\hbar,\alpha,\gamma):(H,A,C)\to(H',A',C')$, where $\hbar:H\to H'$,
$\alpha:A\to A'$, and $\gamma:C\to C'$ are respectively a bialgebra
morphism, an algebra morphism, and a coalgebra morphism such that
$$\rho_{A'}(\alpha(a))=\alpha(a_{(0)})\otimes \hbar(a_{(1)})$$ and
$$\gamma(ch)=\gamma(c)\hbar(h),$$ for all $a\in A, c\in C, h\in H.$
This yields a category which we denote by
$\mathbb{DK}_\bullet^\bullet(k)$.

\begin{proposition}\emph{(see \cite[Proposition 17]{Caenepeel/Militaru/Zhu:2002})}\label{CMZ 17}
\\We have a faithful functor
$$G:\mathbb{DK}_\bullet^\bullet(k)\to
\mathbb{E}_\bullet^\bullet(k)$$ defined by $G((H,A,C))=(A,C,\psi)$
with $\psi:C\otimes A\to A\otimes C$ such that $\psi(c\otimes
a)=a_{(0)}\otimes ca_{(1)}$, and
$G((\hbar,\alpha,\gamma))=(\alpha,\gamma)$.
\end{proposition}

The category of right \emph{Doi-Koppinen-Hopf modules} over the
right-right DK structure $(H,A,C)$ is exactly the category of
$\mathcal{M}(\psi)_A^C$, and it is denoted by $\mathcal{M}(H)_A^C$.

If $H$ has a twisted antipode $\overline S$ (see Subsection
\ref{CategoryGraded}), then $\psi$ is bijective and
$(A,C,\psi^{-1})\in{}_\bullet^\bullet\mathbb{E}(k)$. In this case,
$_A^C\mathcal{M}(\psi^{-1})$ will denoted by $_A^C\mathcal{M}(H)$.
The objects of this category are left $A$-modules and left
$C$-comodules such that
\begin{equation}
\rho_C(am)=m_{(-1)}\overline S(a_{(1)})\otimes a_{(0)}m_{(0)}.
\end{equation}

Analogous considerations are true for left-right, right-left, and
left-left Doi-Koppinen structures.

\medskip
Let $G$ be a group, $X$ a right $G$-set, and $A$ a $G$-graded
$k$-algebra. We know that $(kG,A,kX)$ is a DK structure with $kG$
a Hopf algebra. Then, $(A,kX,\psi)$ is an entwining structure where
$\psi:kX\otimes A\rightarrow A\otimes kX$ is the map defined by
$\psi(x\otimes a_{g})=a_{g}\otimes xg$ for all $x\in X,g\in G,a_g\in
A_g$. Moreover, we have
\begin{equation}\label{graded entwined}
gr-(A,X,G)\simeq \mathcal{M}(kG)_A^{kX}.
\end{equation}
From Theorem \ref{Takeuchi}, we have an $A$-coring, $A\otimes kX$.
The comultiplication and the counit maps of the coring $A\otimes kX$
are defined by:
\begin{equation}
\Delta_{A\otimes kX}(a\otimes x)=(a\otimes x)\tensor{A}(1_A\otimes
x), \quad \epsilon_{A\otimes kX}(a\otimes x)=a \quad (a\in A,x\in
X).
\end{equation}

We will verify in Section \ref{graded} that this coring is
coseparable.

\chapter{Adjoint and Frobenius Pairs of Functors for Corings}

\section{Frobenius functors between categories of comodules}
\label{Frobeniusgeneral}

Let $T$ be a $k$-algebra and $M\in\bmod{T}{A}$. Let
$\rho:T\rightarrow \operatorname{End}_A(M)$ be the morphism of
$k$-algebras given by the left $T$-module structure of the bimodule
$_TM_A$. Now, suppose moreover that $M\in\rcomod{\coring{C}}$. Then
$\operatorname{End}_{\coring{C}}(M)$ is a subalgebra of
$\operatorname{End}_A(M).$ We have that $\rho(T)\subset
\operatorname{End}_{\coring{C}}(M)$ if and only if $\rho_M$ is
$T$-linear. Indeed, let $m\in M$ and $t\in T$, and let
$\rho_M(m)=\sum m_{(0)}\otimes m_{(1)}.$ We have,
$(\rho(t)\otimes_A\coring{C})\rho_M(m)=\sum tm_{(0)}\otimes
m_{(1)}=t\rho_M(m),$ and $\rho_M\rho(t)(m)=\rho_M(tm)$. Hence the
left $T$-module structure of a $T-\coring{C}$-bicomodule $M$ can be
described as a morphism of $k$-algebras $\rho:T\to
\operatorname{End}_{\coring{C}}(M)$, that is, as a left $T$-object
structure on $M\in \rcomod{\coring{C}}$ (see Section
\ref{category theory}). Moreover, the category
$^T\rcomod{\coring{C}}$ is noting else than the category of left
$T$-objects in $\rcomod{\coring{C}}$, $_T\rcomod{\coring{C}}$.

Given a $k$-linear functor
$F:\rcomod{\coring{C}}\to\rcomod{\coring{D}}$ and
$M\in{}_T\rcomod{\coring{C}}$. We know that
$F(M)\in{}_T\rcomod{\coring{D}}$ via $(F(M),\{F(\rho(t))\mid t\in
T\})$. We have then two $k$-linear bifunctors
\[
-\otimes_TF(-),\,F(-\otimes_T-):\mathcal{M}_T
\times{}_T\rcomod{\coring{C}}\to\rcomod{\coring{D}}.
\]
Let $\Upsilon_{T,M}$ be the unique isomorphism of
$\coring{D}$-comodules making the following diagram commutative

\begin{equation}\label{T}
\xymatrix{  T\otimes_TF(M)\ar[rd]_\simeq
\ar[rr]^{\Upsilon_{T,M}} & &  F(T\otimes_TM)\ar[ld]^\simeq \\
 & F(M) & }
\end{equation}
for every $M\in{}_T\rcomod{\coring{C}}.$ It is easy to verify that
$\Upsilon_{T,M}$ is natural in $T$.

\begin{theorem}\emph{[Mitchell]} \emph{(\cite[Theorem 3.6.5]{Popescu:1973})}
\label{Mitchell-Th}
\\Let $\cat{C}$ be an abelian category satisfying AB 5) and
$\cat{C}'$ be a full subcategory whose objects form a set of
projective generators of $\cat{C}$. Let $T,S:\cat{C}\to\cat{D}$
be additive covariant functors where $\cat{D}$ is an abelian
category which has direct sums. Suppose that $T$ preserves
inductive limits. Let $T'$ and $S'$ be the restrictions of
$T$ and $S$ to $\cat{C}'$, respectively. Then, every natural
transformation $\eta:T'\to S'$ can be uniquely extended to a
natural transformation $\overline\eta:T\to S$.
\end{theorem}

\begin{corollary}\emph{[Mitchell]} \emph{(\cite[Corollary
3.6.6]{Popescu:1973})} \label{Mitchell-Cor}
\\In the set up of Theorem \ref{Mitchell-Th}, if $S$ preserves
inductive limits and $\eta$ is a natural equivalence, then so is
$\overline\eta$.
\end{corollary}

\begin{remark}\label{Mitchell}
Theorem \ref{Mitchell-Th} holds also if we only suppose that the
target category $\cat{D}$ is preadditive and has direct sums, or if
the category $\cat{D}$ is preadditive and the functor $S$ preserves
direct sums. Note also that Corollary \ref{Mitchell-Cor} also holds
if we suppose only that the category $\cat{D}$ is preadditive.

Now we give an interesting example of a preadditive category which
is not abelian and has direct sums. A \emph{filtration} on an
abelian group $M$ is a family is an increasing sequence
$(M_n)_{n\in\mathbb{Z}}$ of abelian subgroup of $M$ such that
$M=\bigcup_{n\in\mathbb{Z}}M_n=M_\infty$. An abelian group is called
a \emph{filtered abelian group} if it is endowed with a filtration.
Let $M$ and $N$ be filtered abelian groups. A morphism of abelian
groups $f:M\to N$ is a \emph{morphism of filtered abelian groups} if
$f(M_n)\subset N_n$ for every $n\in\mathbb{Z}$. Filtered abelian
groups with their morphisms form a category. We denote it by
\textbf{FAb}. We check easily that \textbf{FAb} is an additive
category. Moreover, this category has kernels, cokernels, product
sums, and direct sums, but it does not abelian. For the details, see
\cite [3.1.1]{Schneiders:1999}.
\end{remark}

Now, by Remark \ref{Mitchell}, there exists a unique natural
transformation
\[
\Upsilon_{-,M}:-\otimes_TF(M)\to F(-\otimes_TM)
\]
extending the natural transformation $\Upsilon_{T,M}$. (The natural
transformation $\Upsilon_{-,M}$ exists even if the category
$\rcomod{\coring{D}}$ is not abelian.)

\begin{remark}\label{Upsilon'}
Let $\cat{C}$ and $\cat{D}$ be $k$-categories having direct sums and
cokernels, $F:\cat{C}\to\cat{D}$ a $k$-linear functor, and $T$ a
$k$-algebra. We have two $k$-linear bifunctors
\[
-\otimes_TF(-),\,F(-\otimes_T-):\mathcal{M}_T
\times{}_T\cat{C}\to\cat{D}.
\]
Let $M=(M,\rho)\in{}_T\cat{C}$. By the definition of tensor product
(see Proposition \ref{KS 8.5.5} and its proof),
$T\otimes_TF(M)=F(M)$ and $T\otimes_TM=M$. Define
$$\Upsilon_{T,M}=1_{F(M)}:T\otimes_TF(M)\to F(T\otimes_TM).$$
That $\Upsilon_{T,M}$ is natural in $T$ follows from the fact that
for a morphism $f:T\to T$ in $\rmod{T}$,
$f\otimes_TM=\rho(f(1)):M\to M$. By Remark \ref{Mitchell}, there
exists a unique natural transformation
$$\Upsilon_{-,M}:-\otimes_TF(M)\to F(-\otimes_TM):\rmod{T}\to
\cat{D}$$ extending the natural transformation $\Upsilon_{T,M}$.

Let $g:M\to M'$ be a morphism in $_T\cat{C}$. Consider the two
natural transformations
$$F(-\otimes_Tg)\Upsilon_{-,M}, \Upsilon_{-,M}(-\otimes_TF(g)):
-\otimes_TF(M)\to F(-\otimes_TM').$$ Since $T\otimes_Tg=g$,
$F(T\otimes_Tg)\Upsilon_{T,M}= \Upsilon_{T,M}(T\otimes_TF(g))$.
Then, by Remark \ref{Mitchell}, these natural transformations
coincide. Hence $\Upsilon_{T,M}$ is natural in $M$. Again, by Remark
\ref{Mitchell}, if moreover $F$ preserves inductive limits, then
$\Upsilon$ is a natural equivalence of bifunctors. Now we give two
interesting applications:
\begin{enumerate}[(i)]
\item Take $\cat{C}=\rmod{T}$ and $M=T$ and assume that
$F:\rmod{T}\to\cat{D}$ is $k$-linear and preserves inductive limits.
Then we obtain a natural equivalence
$$F\simeq F(-\otimes_TT)\simeq -\otimes_TF(T).$$
This result generalizes \cite[Exercise 5, p. 157]{Mitchell:1965}
which generalizes a result due independently to Watts
\cite{Watts:1960} and Eilenberg \cite{Eilenberg:1960}.
\item Take $\cat{C}=\rmod{R}$ with $R$ a $k$-algebra,
$F=-\otimes_RN:\rmod{R}\to \cat{D}$ with $N\in{}_R\cat{D}$ and let
$X\in \rmod{T}$ and $M\in\bmod{T}{R}$. Then we obtain a natural
equivalence of bifunctors:
$$\Upsilon_{X,M}:X\otimes_T(M\otimes_RN)\to (X\otimes_TM)\otimes_RN.$$
Using Remark \ref{Mitchell} this isomorphism is also natural in $N$.
This result generalizes \cite[Exercise 4, p. 157]{Mitchell:1965}.
\end{enumerate}
\end{remark}

The following lemma is a generalization of \cite[Folgerung
III.4.3]{AlTakhman:1999}.

\begin{lemma}
\begin{enumerate}[(i)]
\item If $F$ preserves direct sums then $\Upsilon_{X,M}$ is an
isomorphism for every projective right $T$-module $X$.
\item If $F$ preserves direct limits then $\Upsilon_{X,M}$ is an
isomorphism for every flat right $T$-module $X$.
\item If $F$ preserves inductive limits then $\Upsilon_{X,M}$ is an
isomorphism for every right $T$-module $X$.
\end{enumerate}
\end{lemma}

\begin{proof}
Set $\eta=\Upsilon_{-,M}$, $S=-\otimes_TF(M)$, and
$T=F(-\otimes_TM)$. First observe that $S$ preserves inductive
limits.

(i) Let $(X_i)_{i\in I}$ be a family of right $T$-modules. It is
clear that $\eta_{\oplus_{i\in I}X_i}=\oplus_{i\in I}\eta_{X_i}$.
Then $\eta_{\oplus_{i\in I}X_i}$ is an isomorphism if and only if
$\oplus_{i\in I}\eta_{X_i}$ is so for every $i\in I$. Hence (i)
follows since a projective right module is a direct summand of a
free module.

(ii) Let $X$ be a flat right $T$-module. By Lazard's structural
theorem \cite[Th\'eor\`eme 1.2]{Lazard:1969}, every flat module is a
direct limit of finitely generated free modules. Then
$X=\underset{\underset{I}{\longrightarrow}}\lim X_i$ for a directed
set $I$ and a direct system  of finitely generated free modules
$(X_i,u_{ij})$. Then,
$S(X)=\underset{\underset{I}{\longrightarrow}}\lim S(X_i)$ and
$T(X)=\underset{\underset{I}{\longrightarrow}}\lim T(X_i)$. Hence,
$\eta_{X_i}:S(X_i)\to T(X_i)$, $i\in I$, is a morphism from
$(S(X_i),S(u_{ij}))$ to $(T(X_i),T(u_{ij}))$, and
$\underset{\underset{I}{\longrightarrow}}\lim \eta_{X_i}=\eta_X$.
Finally, from (i), $\eta_X$ is an isomorphism.

(iii) Let $X$ be a right $T$-module. Consider a free presentation in
$\rmod{T}$ of $X$ $$T^{(J)}\to T^{(I)}\to X\to 0,$$ where $I$ and
$J$ are sets. This yields the commutative diagram with exact rows
($F$ is right exact)
$$\xymatrix{T^{(J)}\tensor{T}F(M)\ar[r] \ar[d]^{\eta_{T^{(J)}}}&
T^{(I)}\tensor{T}F(M)\ar[r] \ar[d]^{\eta_{T^{(I)}}} &
X\tensor{T}F(M)\ar[r] \ar[d]^{\eta_X} & 0
\\
F(T^{(J)}\tensor{T}M)\ar[r] & F(T^{(I)}\tensor{T}M)\ar[r] &
F(X\tensor{T}M)\ar[r] & 0.}$$ From (i), $\eta_{T^{(I)}}$ and
$\eta_{T^{(J)}}$ are isomorphisms. Hence so is $\eta_X$.
\end{proof}

Clearly we have

\begin{proposition}\label{G1}
$\Upsilon$ is a natural transformation of bifunctors. Moreover,
$\Upsilon_{X,-}$ is a natural isomorphism for all projective right
$T$-module $X$. If $F$ preserves direct limits, then
$\Upsilon_{X,-}$ is a natural isomorphism for all flat right
$T$-module $X$. If $F$ preserves inductive limits, then $\Upsilon$
is a natural equivalence.
\end{proposition}

\begin{example}
Let $N\in\bcomod{\coring{C}}{\coring{D}}$. Suppose that
$_B\coring{D}$ is flat or $\coring{C}$ is coseparable. Consider the
functor $F=-\cotensor{\coring{C}}N:\rcomod{\coring{C}}\to
\rcomod{\coring{D}}$ and the natural transformation associated to
$F$, $\Upsilon$. By Mitchell's Theorem, $\Upsilon_{W,M}$ is exactly
the canonical map (see Proposition \ref{assotenscotens})
$$W\tensor{T}(M\cotensor{\coring{C}}N)\to
(W\tensor{T}M)\cotensor{\coring{C}}N,$$ for every $W\in
\rmod{T},M\in {}_T\rcomod{\coring{C}}$.
\end{example}

From Lemma \ref{$A$-objects} we obtain

\begin{lemma}\label{$T$-objects}
Let $F,G:\rcomod{\coring{C}}\to \rcomod{\coring{D}}$ be $k$-linear
functors, and let $\eta:F\to G$ be a natural transformation. For
each $M\in{}_T\rcomod{\coring{C}}$, $\eta_M:F(M)\to G(M)$ is a
morphism of $(T,\coring{D})$-bicomodules.
\end{lemma}

Using Remark \ref{Mitchell}, we prove easily the two following
lemmas. The condition ``$F$ (and $G$) preserve(s) direct sums'' is
superfluous. Moreover, these lemmas remain true for $k$-categories
which have direct sums and cokernels (see Remark \ref{Upsilon'}).

\begin{lemma}\emph{\cite[Lemma 3.2]{Gomez:2002}}\label{G2}
\\Let $F,G:\rcomod{\coring{C}}\to \rcomod{\coring{D}}$ be $k$-linear
functors and $\eta:F\to G$ a natural transformation. For $X\in
\rmod{T}$ and $M\in{}_T\rcomod{\coring{C}}$, the following diagram
is commutative
$$\xymatrix{
X\tensor{T}F(M)\ar[rr]^{X\tensor{T}\eta_M}\ar[d]^{\Upsilon^F_{X,M}}
&&
X\tensor{T}G(M)\ar[d]^{\Upsilon^G_{X,M}}\\
F(X\tensor{T}M)\ar[rr]^{\eta_{X\tensor{T}M}} && G(X\tensor{T}M).}$$
\end{lemma}

\begin{lemma}\emph{\cite[Lemma 3.3]{Gomez:2002}}\label{G3}
\\Let $S,T$ be $k$-algebras and $F:\rcomod{\coring{C}}\to
\rcomod{\coring{D}}$ a $k$-linear functor. In the situation
$(X_S,{}_SY_T,{}_TM_\coring{C})$,
$$\Upsilon_{X,Y\tensor{T}M}\circ(X\tensor{S}\Upsilon_{Y,M})=
\Upsilon_{X\tensor{S}Y,M}.$$
\end{lemma}

Let $M\in \bcomod{\coring{C'}}{\coring{C}}$ be a bicomodule. We say
that a $k$-linear functor $F:\rcomod{\coring{C}}\to
\rcomod{\coring{D}}$ is $M$-\emph{compatible} if
$\Upsilon_{\coring{C'},M}$ and
$\Upsilon_{\coring{C'}\tensor{A'}\coring{C'},M}$ are isomorphisms.
For example, $F$ is $M$-compatible for every bicomodule $M$ either
if $F$ preserves inductive limits, or $\coring{C'}_{A'}$ is flat and
$F$ preserves direct limits, or $\coring{C'}_{A'}$ is projective and
$F$ preserves coproducts (since, by \cite[Proposition
II.5.3]{Cartan/Eilenberg:1956},
$(\coring{C'}\tensor{A'}\coring{C'})_{A'}$ is projective).

Now suppose that $F$ is $M$-compatible. We define
$$\lambda_{F(M)}:=\Upsilon_{\coring{C'},M}^{-1}\circ
F(\lambda_M):F(M)\to \coring{C'}\tensor{A'}F(M).$$ We have that
$\lambda_{F(M)}$ is a morphism of $A'-\coring{D}$-bicomodules.

The following lemma will be useful in the proof of the next theorem.

\begin{lemma}\label{2}
If $M\in\bcomod{\coring{C'}}{\coring{C}}$ and
$F:\rcomod{\coring{C}}\to\rcomod{\coring{D}}$ is a $M$-compatible
$k$-linear functor, then for all $X\in\mathcal{M}_{A'},$
\[
\Upsilon_{X\otimes_{A'}\coring{C}',M}\circ(X\otimes_{A'%
}\lambda_{F(M)})=F(X\otimes_{A'}\lambda_M)\circ\Upsilon_{X,M}.
\]
\end{lemma}

\begin{proof}
Let us consider the diagram
\[
\xymatrix{
 X\otimes_{A'}F(M) \ar[rrrr]^{X\otimes_{A'}\lambda
_{F(M)}} \ar[rrd]^{X\otimes_{A'}F(\lambda_M)}
\ar[dd]_{\Upsilon_{X,M}}& & & & X
\otimes_{A'}\coring{C}'\otimes_{A'}F(M)
\ar[lld]^{X\otimes_{A'}\Upsilon_{\coring{C}',M}}
\ar[dd]^{\Upsilon_{X\otimes
_{A'}\coring{C}',M}}\\
& & X\otimes_{A'}F(\coring{C} '\otimes_{A'}M)
\ar[rrd]^{\Upsilon_{X,\coring{C}'\otimes_{A'}M}} & & \\
F(X\otimes_{A'}M) \ar[rrrr]_{F(X\otimes_{A'}\lambda _M)}& &  & &
F(X\otimes_{A'}\coring{C}'\otimes_{A'}M) }
\]
The commutativity of the top triangle follows from the definition of
$\lambda_{F(M)}$, while the right triangle commutes by Lemma
\ref{G3} (we take $S=T=A'$ and $Y=\coring{C}'$), and the left
triangle is commutative since $\Upsilon_{X,-}$ is natural.
Therefore, the commutativity of the rectangle holds.
\end{proof}

\begin{proposition}\emph{\cite[Proposition 3.4]{Gomez:2002}}\label{G4}
\\Let $M\in\bcomod{\coring{C'}}{\coring{C}}$ and let $F:\rcomod{\coring{C}}
\to\rcomod{\coring{D}}$ be a $M$-compatible functor. $F(M)$ endowed
with $\rho_{F(M)}$ and $\lambda_{F(M)}$, is a
$\coring{C'}-\coring{D}$-bicomodule. Moreover, if
$F,G:\rcomod{\coring{C}}\to \rcomod{\coring{D}}$ are $M$-compatible,
and $\eta:F\to G$ is a natural transformation, then the map
$\eta_M:F(M)\to G(M)$ is a morphism of
$\coring{C'}-\coring{D}$-bicomodules.
\end{proposition}

\begin{proof}
To prove the first statement, consider the diagram
$$\xymatrix{
 && F(M) \ar[rrr]^{\lambda_{F(M)}}
 \ar[dll]^(.3){\lambda_{F(M)}}\ar@{=}'[d][ddd]&&&
\coring{C}'\tensor{A'}F(M)
\ar[dll]^{\Delta_{\coring{C}'}\tensor{A'}F(M)}
\ar[ddd]^{\Upsilon_{\coring{C'},M}}
\\
\coring{C}'\tensor{A'}F(M)
\ar[rrr]^(.35){\coring{C}'\tensor{A'}\lambda_{F(M)}}
\ar[ddd]^{\Upsilon_{\coring{C'},M}} &&&
\coring{C}'\tensor{A'}\coring{C}'\tensor{A'}F(M)
\ar[ddd]^(.3){\Upsilon_{\coring{C}'\tensor{A'}\coring{C'},M}}
\\
&&&&&&&&&
\\
 && F(M) \ar'[r][rrr]^(.35){F(\lambda_M)} \ar[dll]^{F(\lambda_M)} &&&
F(\coring{C}'\tensor{A'}M)
\ar[dll]^{F(\Delta_{\coring{C}'}\tensor{A'}M)}
\\
F(\coring{C}'\tensor{A'}M)
\ar[rrr]^{F(\coring{C}'\tensor{A'}\lambda_M)} &&&
F(\coring{C}'\tensor{A'}\coring{C}'\tensor{A'}M)}.$$ Since
$\Upsilon_{-,M}$ is a natural transformation, and from the
definition of $\lambda_{F(M)}$, the coassociative property of
$\lambda_M$, and Lemma \ref{2}, all squares, save possibly the top
one, are commutative. Hence, by Lemma \ref{cube}
($\Upsilon_{\coring{C}'\tensor{A'}\coring{C}',M}$ is injective), the
top square is also commutative.

To show the counity property, consider the diagram
$$\xymatrix{\coring{C}'\tensor{A'}F(M)
\ar[rr]^{\Upsilon_{\coring{C'},M}}
\ar[dd]_{\epsilon_{\coring{C}'}\tensor{A'}F(M)} &&
F(\coring{C}'\tensor{A'}M)
\ar[dd]^{F(\epsilon_{\coring{C}'}\tensor{A'}M)}\\
& F(M) \ar[ur]^{F(\lambda_M)}\ar[ul]^{\lambda_{F(M)}}
\ar[dr]^\simeq\ar[dl]^\simeq &\\
A'\tensor{A'}F(M) \ar[rr]^{\Upsilon_{A',M}} && F(A'\tensor{A'}M).}$$
From the definition of $\Upsilon_{A',M}$ and $\lambda_{F(M)}$, the
counity property of the left $\coring{C}'$-comodule $M$, and the
fact that $\Upsilon_{-,M}$ is natural, all the triangles, save
possibly the left one, and the rectangle are commutative. Since
$\Upsilon_{A',M}$ is injective, then the left triangle is also
commutative.

Now we will verify the last statement. For this consider the diagram
$$\xymatrix{F(M) \ar[rrr]^{\eta_M} \ar[dd]^{\lambda_{F(M)}}
\ar[dr]^{F(\lambda_M)} &&&
G(M)\ar'[d][dd]^{\lambda_{G(M)}}\ar[dr]^{G(\lambda_M)} &&\\
& F(\coring{C}'\tensor{A'}M)
\ar[rrr]^(.3){\eta_{\coring{C}'\tensor{A'}M}} &&&
G(\coring{C}'\tensor{A'}M)\\
\coring{C}'\tensor{A'}F(M)\ar[rrr]^{\coring{C}'\tensor{A'}\eta_M}
\ar[ur]^{\Upsilon^F_{\coring{C'},M}} &&&
\coring{C}'\tensor{A'}G(M)\ar[ur]^{\Upsilon^G_{\coring{C'},M}}&&
.}$$ The commutativity of the triangles follows by the definition of
$\lambda_{F(M)}$ and $\lambda_{G(M)}$. The upper parallelogram is
commutative since $\eta$ is natural while the bottom one commutes
from Lemma \ref{G2}. Since $\Upsilon^G_{\coring{C'},M}$ is an
injective map, the rectangle is also commutative.
\end{proof}

From now on, every result for comodules over coseparable corings
applies in particular for modules over rings (since the trivial
$A$-coring is coseparable).
\medskip

 Now we propose a generalization of \cite[Proposition
2.1]{Takeuchi:1977}, \cite[23.1(1)]{Brzezinski/Wisbauer:2003} and
\cite[Theorem 3.5]{Gomez:2002}.

\begin{theorem}\label{3}
Let $F:\rcomod{\coring{C}}\to\rcomod{\coring{D}}$ be a $k$-linear
functor, such that
\begin{enumerate}[(I)]
\item $_B\coring{D}$ is flat and $F$ preserves the kernel of
$\rho_N\otimes_A\coring{C}-N\otimes_A\Delta_{\coring{C}}$ for every
$N \in \rcomod{\coring{C}}$, or
\item $\coring{C}$ is a coseparable $A$-coring and the categories
$\mathcal{M}^{\coring{C}}$ and $\rcomod{\coring{D}}$ are abelian.
\end{enumerate}
Assume that at least one of the following statements holds
\begin{enumerate}
\item $\coring{C}_A$ is projective, $F$ preserves coproducts,
and $\Upsilon_{N,\coring{C}}$, $\Upsilon_{N\otimes_A\coring{C}%
,\coring{C}}$ are isomorphisms for all $N\in\rcomod{\coring{C}}$
(e.g., if $A$ is semisimple and $F$ preserves coproducts), or \item
$\coring{C}_A$ is flat, $F$ preserves direct limits, and $\Upsilon
_{N,\coring{C}}$, $\Upsilon_{N\otimes_A\coring{C},\coring{C}}$ are
isomorphisms for all $N\in\rcomod{\coring{C}}$ (e.g., if $A$ is a
von Neumann regular ring and $F$ preserves direct limits), or
\item $F$ preserves inductive limits (e.g., if $F$ has a right
adjoint).
\end{enumerate}
 Then $F$ is naturally equivalent to
$-\cotensor{\coring{C}}F(\coring{C}).$
\end{theorem}

\begin{proof}
In each case, we have $F$ is $\coring{C}$-compatible where
$\coring{C}\in{}^{\coring{C}}\mathcal{M}^{\coring{C}}.$ Therefore,
by Proposition \ref{G4}, $F(\coring{C})$ can be viewed as a
$\coring{C}-\coring{D}$-bicomodule. From Lemma \ref{2}, and since
$\Upsilon_{-,\coring{C}}$ is a natural transformation, we have, for
every $N\in \mathcal{M}^{\coring{C}}$, the commutativity of the
following diagram with
exact rows in $\mathcal{M}^{\coring{D}}$%
\[%
\xymatrix{
 0 \ar[r] & N\square_{\coring{C}}F(\coring{C})\ar[r] &
 N\otimes_AF(\coring{C})
\ar[d]^\simeq_{\Upsilon_{N,\coring{C}}}
\ar[rrr]^-{\rho_N\otimes_AF(\coring{C})
-N\otimes_A%
\lambda_{F(\coring{C})}} &  & & N\otimes
_A\coring{C}\otimes_AF(\coring{C})\ar[d]^{\simeq
}_{\Upsilon_{N\otimes_A\coring{C},\coring{C}}} \\
 0 \ar[r] &  F(N)\ar[r]^-{F(\rho
_N)} & F(N\otimes_A\coring{C})
\ar[rrr]^-{F(\rho_N\otimes_A\coring{C}-N\otimes_A%
\Delta_{\coring{C}})} & & & F(N\otimes
_A\coring{C}\otimes_A\coring{C}).}
\]
The exactness of the bottom sequence is assumed in the case (I). For
the case (II), it follows by factorizing the map
$\omega_{N,\coring{C}} =
\rho_N\otimes_A\coring{C}-N\otimes_A%
\Delta_{\coring{C}}$ through its image, and using the facts that the
sequence $\xymatrix{0 \ar[r] &  N\ar[r]^-{\rho _N} & N
\otimes_A\coring{C} \ar[r]^-{\omega_{N,\coring{C}}}  & N \otimes
_A\coring{C}\otimes_A\coring{C}}$ is split exact in
$\rcomod{\coring{C}}$, and that additive functors between abelian
categories preserve split exactness. By the universal property of
kernel, there exists a unique isomorphism $\eta_N:N\square
_{\coring{C}}F(\coring{C}) \to F(N)$ in $\rcomod{\coring{D}}$ making
commutative the above diagram.

Now, we will verify that $\eta$ is natural. For this let $f:N\to N'$
be a morphism in $\rcomod{\coring{C}}$, and consider the following
diagram
$$\xymatrix{
 &&  N\cotensor{\coring{C}}F(\coring{C}) \ar[rrr]^{\eta_N}
 \ar[dll]^(.3){f\cotensor{\coring{C}}F(\coring{C})}\ar'[d][ddd]&&& F(N)
 \ar[dll]^{F(f)}\ar[ddd]^{F(\rho_N)}
\\
N'\cotensor{\coring{C}}F(\coring{C})
\ar[rrr]^(.3){\eta_{N'}}\ar[ddd] &&& F(N')
\ar[ddd]^(.3){F(\rho_{N'})}
\\
&&&&&&&&&
\\
 && N\tensor{A}F(\coring{C})\ar'[r][rrr]^{\Upsilon_{N,\coring{C}}}
 \ar[dll]^{f\tensor{A}F(\coring{C})} &&& F(N\tensor{A}\coring{C})
 \ar[dll]^{F(f\tensor{A}\coring{C})}
\\
N'\tensor{A}F(\coring{C}) \ar[rrr]^{\Upsilon_{N',\coring{C}}}  &&&
F(N'\tensor{A}\coring{C}).}$$ From Lemma \ref{cube}, the top square
is commutative as desired. Hence
$F\simeq-\square_{\coring{C}}F(\coring{C}).$
\end{proof}

\begin{remark}
Let $\cat{C}$ be a $k$-category having direct sums and cokernels,
$A$ a $k$-algebra, $\coring{C}$ an $A$-coring, and
$M\in{}_A\cat{C}$. A left $\coring{C}$-\emph{coaction} on $M$ is a
morphism $\lambda_M: M \to\coring{C}\tensor{A}M$ in $_A\cat{C}$ such
that the following diagrams in $\cat{C}$ are commutative
\begin{equation*}
\xymatrix{M \ar[rr]^{\lambda_M} \ar[d]^{\lambda_M} &&
\coring{C}\tensor{A}M\ar[d]^{\Delta \tensor{A}M} \\
\coring{C}\tensor{A}M\ar[rr]^-{\coring{C}\tensor{A}\lambda_M} &&
\coring{C}\tensor{A}\coring{C}\tensor{A}M}
\xymatrix{M\ar[r]^-{\lambda_M} \ar@{=}'[dr] &
\coring{C}\tensor{A}M\ar[d]^{\epsilon\tensor{A}M}
\\ & A\tensor{A}M.}
\end{equation*}

Let $N\in{}_A\cat{C}$ and $\lambda_M$ and $\lambda_N$ be left
$\coring{C}$-coactions on $M$ and $N$ respectively. A morphism $f: M
\to N$ in $_A\cat{C}$ is a \emph{comodule morphism} if the following
diagram in $\cat{C}$ is commutative
\begin{equation*}\label{comodmor}
\xymatrix{M\ar[rr]^f \ar[d]^{\lambda_M} && N \ar[d]^{\lambda_N}
\\ \coring{C}\tensor{A}M\ar[rr]^{\coring{C}\tensor{A}f}
&& \coring{C}\tensor{A}N.}
\end{equation*}
The composition of morphisms is the same as in $\cat{C}$. We obtain
a category which we denote by $^\coring{C}\cat{C}$. For every
$M,N\in {}^\coring{C}\cat{C}$, $\hom{{}^\coring{C}\cat{C}}{M}{N}$ is
a $k$-submodule of $\hom{{}_A\cat{C}}{M}{N}$. Then
$^\coring{C}\cat{C}$ is a $k$-category. If $\cat{C}=\rmod{k}$, then
$^\coring{C}\cat{C}=\lcomod{\coring{C}}$. If
$\cat{C}=\rcomod{\coring{D}}$ where $\coring{D}$ is a $B$-coring,
then $^\coring{C}\cat{C}=\bcomod{\coring{C}}{\coring{D}}$.

Now assume moreover that $\cat{C}$ has kernels. Let
$M\in\rcomod{\coring{C}}$ and $N\in{}^\coring{C}\cat{C}$. Consider
the morphism in $\cat{C}$:
\[
\omega_{M,N}=\rho_M\otimes_AN-M\otimes_A\lambda_N:M\otimes_AN \to
M\otimes_A\coring{C}\otimes_AN.
\]

The kernel of $\omega_{M,N}$ in $\cat{C}$ is called the
\emph{cotensor product} of $M$ and $N$.

Let $f:M\to M'$ be a morphism in $\rcomod{\coring{C}}$ and $g:N\to
N'$ a morphism in $^\coring{C}\cat{C}$. There is a unique morphism
$f\cotensor{\coring{C}}g:M\cotensor{\coring{C}}N\to
M'\cotensor{\coring{C}}N'$ in $\cat{C}$ making commutative the
diagram
$$\xymatrix{
0\ar[r] & M\cotensor{\coring{C}}N
\ar[r]\ar[d]^{f\cotensor{\coring{C}}g} & M\tensor{A}N
\ar[rr]^-{\omega_{M,N}}\ar[d]^{f\tensor{A}g} &&
M\tensor{A}\coring{C}\tensor{A}N
\ar[d]^{f\tensor{A}\coring{C}\tensor{A}g}
\\ 0 \ar[r] & M'\cotensor{\coring{C}}N'\ar[r] & M'\tensor{A}N'
\ar[rr]^-{\omega_{M',N'}} && M'\tensor{A}\coring{C}\tensor{A}N'.}
$$
$f\cotensor{\coring{C}}g$ is called the \emph{cotensor product} of
$f$ and $g$. We obtain then a bifunctor $k$-linear in each variable
$$-\cotensor{\coring{C}}-:
\rcomod{\coring{C}}\times {}^\coring{C}\cat{C}\to\cat{C}.$$

Now let $\cat{D}$ be a $k$-category having direct sums and
cokernels. Let $M\in \bcomod{\coring{C'}}{\coring{C}}$ be a
bicomodule. We say that a $k$-linear functor
$F:\rcomod{\coring{C}}\to \cat{D}$ is $M$-\emph{compatible} if
$\Upsilon_{\coring{C'},M}$ and
$\Upsilon_{\coring{C'}\tensor{A'}\coring{C'},M}$ are isomorphisms in
$\cat{D}$. If $F$ preserves inductive limits then $F$ is
$M$-compatible for every bicomodule $M$.

Assume that $F$ is $M$-compatible. Let $\lambda_{F(M)}$ be the
unique morphism in $_{A'}\cat{D}$ making commutative the diagram
$$\xymatrix{F(M)\ar[rr]^{\lambda_{F(M)}} \ar[dr]_{F(\lambda_M)}&&
\coring{C'}\tensor{A'}F(M)\ar[dl]^{\Upsilon_{\coring{C'},M}}\\ &
F(\coring{C'}\tensor{A'}M).}$$ Clearly Lemma \ref{2} and Proposition
\ref{G4} remain true for $F$. We have that $F(M)$, endowed with
$\lambda_{F(M)}$, belongs to $^{\coring{C}'}\cat{D}$. Moreover, if
$F,G:\rcomod{\coring{C}}\to \cat{D}$ are $M$-compatible, and
$\eta:F\to G$ is a natural transformation, then the map
$\eta_M:F(M)\to G(M)$ is a morphism in $^{\coring{C}'}\cat{D}$.

Finally, let $\coring{C}$ be an $A$-coring such that
$\rcomod{\coring{C}}$ is an abelian category, $\cat{D}$ an abelian
category which has direct sums, and
$F:\rcomod{\coring{C}}\to\cat{D}$ a $k$-linear functor which
preserves inductive limits. If $\coring{C}$ is coseparable or $F$ is
left exact then
$$F\simeq -\cotensor{\coring{C}}F(\coring{C}).$$ The proof is
analogous to that of Theorem \ref{3}. This result is another
generalization of \cite[Exercise 5, p. 157]{Mitchell:1965} (see
Remark \ref{Upsilon'}).
\end{remark}

As an immediate consequence of the last theorem we have the
following generalization of Eilenberg--Watts Theorem
\cite[Proposition VI.10.1]{Stenstrom:1975}.

\begin{corollary}\label{4}
Let $F:\rcomod{\coring{C}}\rightarrow\rcomod{\coring{D}}$ be a
$k$-linear functor.
\begin{enumerate}[(1)]
\item If $_B\coring{D}$ is
flat and $A$ is a semisimple ring (resp. a von Neumann regular
ring), then the following statements are
equivalent\begin{enumerate}[(a)]
\item $F$ is left exact and
preserves coproducts (resp. left exact and preserves direct limits);
\item $F\simeq-\square_{\coring{C}}M$ for some bicomodule
$M\in{}^{\coring{C}%
}\mathcal{M}^{\coring{D}}$.
\end{enumerate}
 \item If
$_A\coring{C}$ and $_B\coring{D}$ are flat, then the following
statements are equivalent
\begin{enumerate}[(a)]
\item $F$ is exact
and preserves inductive limits;\item $F\simeq-\square_{\coring{C}}M$
for some
bicomodule $M\in{}^{\coring{C}%
}\mathcal{M}^{\coring{D}}$ which is coflat in
$^{\coring{C}}\mathcal{M}%
$.
\end{enumerate}
\item If $\coring{C}$ is a coseparable $A$-coring
and the categories $\rcomod{\coring{C}}$ and
$\mathcal{M}^{\coring{D}}$ are abelian, then the following
statements are equivalent
\begin{enumerate}[(a)]\item $F$ preserves
inductive limits;\item $F$ preserves cokernels and
$F\simeq-\square_{\coring{C}}M$ for some bicomodule
$M\in{}^{\coring{C}}\mathcal{M}^{\coring{D}}$.
\end{enumerate}\item If $\coring{C}=A$ and the category
$\mathcal{M}^{\coring{D}}$ is abelian, then the following statements
are equivalent
\begin{enumerate}[(a)]
\item $F$ has a right adjoint;\item $F$
preserves inductive limits;
\item $F\simeq-\otimes_AM$ for some bicomodule
$M\in{}_A\rcomod{\coring{D}}.$
\end{enumerate}
\end{enumerate}
\end{corollary}

\begin{proof}
Follows directly from Theorem \ref{3} and Proposition
\ref{cotensdirlim}.
\end{proof}

The particular case of the following statement when the cohom is
exact generalizes \cite[Corollary 3.12]{AlTakhman:2002}.

\begin{corollary}\label{5}
Let $N\in{}^{\coring{C}}\mathcal{M}^{\coring{D}}$ be a bicomodule,
quasi-finite as a right $\coring{D}$-comodule, such that
$_A\coring{C}$ and $_B\coring{D}$ are flat. If the cohom functor
$\cohom{\coring{D}}{N}{-}$ is exact or if $\coring{D}$ is a
coseparable $B$-coring, then we have
\[
\cohom{\coring{D}}{N}{-}
\simeq-\square_{\coring{D}}\cohom{\coring{D}}{N}{\coring{D}}:
\mathcal{M}^{\coring{D}}\rightarrow \mathcal{M}^{\coring{C}}.
\]
\end{corollary}

\begin{proof}
The functor $\cohom{\coring{D}}{N}{-}$ is $k$-linear and preserves
inductive limits, since it is a left adjoint to the $k$-linear
functor $-\cotensor{\coring{C}}N:\rcomod{\coring{C}}\to
\rcomod{\coring{D}}$ (see Propositions \ref{adjointcohomcotens},
\ref{M IV.3}). Hence Theorem \ref{3} achieves the proof.
\end{proof}

Now we will use the following generalization of \cite[Lemma
2.2]{Takeuchi:1977}.

\begin{lemma}\label{18}
Let $\Lambda$, $\Lambda'$ be bicomodules in
$\bcomod{\coring{D}}{\coring{C}}$. Suppose that
\begin{enumerate}[(a)]
\item $\oper{Im}{\omega_{\coring{D},\Lambda}}$ is a $\coring{D}_B$-pure
submodule of $\coring{D}\tensor{B}\coring{D}\tensor{B}\Lambda$, and
$\oper{Im}{\omega_{\coring{D},\Lambda'}}$ is a $\coring{D}_B$-pure
submodule of $\coring{D}\tensor{B}\coring{D}\tensor{B}\Lambda'$,
\item $\omega_{X,\Lambda}$ and $\omega_{X,\Lambda'}$ are
$_A(\coring{C}\tensor{A}\coring{C})$-pure for every $X\in
\rcomod{\coring{D}}$ (e.g. $_A\coring{C}$ is flat or $\coring{D}$ is
coseparable).
\end{enumerate}
Let
$G=-\cotensor{\coring{D}}\Lambda,G'=-\cotensor{\coring{D}}\Lambda'
:\rcomod{\coring{D}}\to \rcomod{\coring{C}}$. Then
\[
\nat{G}{G'}\simeq \hom{(\coring{D},\coring{C})}{\Lambda}{\Lambda'}.
\]
\end{lemma}

\begin{proof}
We want to show that the correspondence
\begin{equation}\label{Nat-Hom}
\hom{(\coring{D},\coring{C})}{\Lambda}{\Lambda'}\to \nat{G}{G'},
\quad \widetilde \eta\mapsto -\cotensor{\coring{D}}\widetilde \eta
\end{equation}
is bijective.

Let $\eta:G\to G'$ be a natural transformation. By Lemma
\ref{$T$-objects}, $\eta_\coring{D}$ is left $B$-linear.

Let $M\in\rmod{B}$ and $m\in M$. Define $f_m:\coring{D}\to
M\tensor{B}\coring{D}$, $f_m(d)=m\otimes d$, $m\in
M,d\in\coring{D}$. $f_m$ is clearly right $\coring{D}$-colinear.
Then, by Lemma \ref{caniso}, we obtain the commutative diagram
$$\xymatrix{G(\coring{D})\simeq\Lambda
\ar[rr]^{\eta_\coring{D}}\ar[d]_{G(f_m)} && \Lambda'\simeq
G'(\coring{D})\ar[d]^{G'(f_m)}
\\ G(M\tensor{B}\coring{D})\simeq M\tensor{B}\Lambda
\ar[rr]^{\eta_{M\tensor{B}\coring{D}}} && M\tensor{B}\Lambda'\simeq
G'(M\tensor{B}\coring{D}),}$$ and $G(f_m)(\lambda)=m\otimes\lambda$
for all $\lambda\in\Lambda$. Hence
$\eta_{M\tensor{B}\coring{D}}=M\tensor{B}\eta_\coring{D}$.

Now, let $M\in \rcomod{\coring{D}}$. Again from Lemma \ref{caniso}
we have the commutative diagram
\begin{equation}\label{Nat-Hom-Diag}
\xymatrix{M\cotensor{\coring{D}}\Lambda
\ar[rr]^{\eta_M}\ar[d]_{\rho_M\cotensor{\coring{D}}\Lambda} &&
M\cotensor{\coring{D}}\Lambda'\ar[d]^{\rho_M\cotensor{\coring{D}}\Lambda'}
\\ (M\tensor{B}\coring{D})\cotensor{\coring{D}}\Lambda
\ar[rr]^{\eta_{M\tensor{B}\coring{D}}}\ar[d]_\simeq &&
(M\tensor{B}\coring{D})\cotensor{\coring{D}}\Lambda'\ar[d]^\simeq
\\ M\tensor{B}\Lambda
\ar[rr]^{M\tensor{B}\eta_\coring{D}} && M\tensor{B}\Lambda'}
\end{equation}
and the compositions of the vertical maps are the inclusion maps.
From the commutativity of the diagram \eqref{Nat-Hom-Diag}, we get
the following commutative diagram
$$\xymatrix{\Lambda
\ar[rr]^{\eta_\coring{D}}\ar[d]_{\lambda_\Lambda} &&
\Lambda'\ar[d]^{\lambda_{\Lambda'}}
\\ \coring{D}\tensor{B}\Lambda
\ar[rr]^{\coring{D}\tensor{B}\eta_\coring{D}} &&
\coring{D}\tensor{B}\Lambda'},$$ i.e. $\eta_\coring{D}$ is left
$\coring{D}$-colinear. Therefore,
$\eta_M=M\cotensor{\coring{D}}\eta_\coring{D}$. Hence the
correspondence \eqref{Nat-Hom} is bijective. Finally, it is clearly
an isomorphism of abelian groups.
\end{proof}

The following proposition generalizes \cite[Theorem
2.1]{Morita:1965} from bimodules over rings to bicomodules over
corings.

\begin{proposition}\label{19a}
Suppose that $_A\coring{C}$, $\coring{C}_A$, $_B\coring{D}$ and
$\coring{D}_B$ are flat. Let $X\in\bcomod{\coring{C}}{\coring{D}}$
and $\Lambda\in\bcomod{\coring{D}}{\coring{C}}$. Consider the
following properties:
\begin{enumerate}[(1)]
\item \label{adj1} $-\square _{\coring{C}}X$ is left adjoint to
$-\square_{\coring{D}}\Lambda$; \item \label{adj2} $\Lambda$ is
quasi-finite as a right $\coring{C}$-comodule and $-
\square_{\coring{C}} X \simeq \cohom{\coring{C}}{\Lambda}{-}$;
\item \label{adjsep} $\Lambda$ is quasi-finite as a right $\coring{C}%
$-comodule and $X\simeq \cohom{\coring{C}}{\Lambda}{\coring{C}}$ in
$\bcomod{\coring{C}}{\coring{D}}$; \item \label{semicontext} there
exist bicolinear maps
\[
\psi:\coring{C}\rightarrow X\square_{\coring{D}}\Lambda\text{ and
}\omega:\Lambda\square_{\coring{C}}X\rightarrow\coring{D}
\]
in $^{\coring{C}}\mathcal{M}^{\coring{C}}$ and $^{\coring{D}}
\mathcal{M}^{\coring{D}}$ respectively, such that
\begin{equation}\label{unitcounit}
(\omega\square_{\coring{D}}\Lambda)\circ(\Lambda
\square_{\coring{C}}\psi)=\Lambda\text{ and }(X\square
_{\coring{D}}\omega)\circ(\psi\square_{\coring{C}}X)=X;
\end{equation}
\item \label{adj1symm} $\Lambda \cotensor{\coring{C}}-$ is left
adjoint to $X \cotensor{\coring{D}}-$.
\end{enumerate}
Then

(\ref{adj1}) and (\ref{adj2}) are equivalent, and they imply
(\ref{adjsep}).

(\ref{adjsep}) implies (\ref{adj2}) if $\coring{C}$ is a coseparable
$A$-coring. If ${}_A X$ and ${}_B \Lambda$ are flat, and
$\omega{}_{X,\Lambda}=\rho{}_X\otimes_B\Lambda-X\otimes_A\rho{}_\Lambda$
is pure as an $A$-linear map and
$\omega{}_{\Lambda,X}=\rho{}_\Lambda\otimes_AX-\Lambda\otimes_B\rho{}_X$
is pure as a $B$-linear map (e.g. if ${}_{\coring{C}}X$ and
${}_{\coring{D}}\Lambda$ are coflat (see Lemma \ref{relatinjsplit})
or $A$ and $B$ are von Neumann regular rings), or if $\coring{C}$
and $\coring{D}$ are coseparable, then (\ref{semicontext}) implies
(\ref{adj1}). The converse is true if ${}_{\coring{C}}X$ and
${}_{\coring{D}}\Lambda$ are coflat, or if $A$ and $B$ are von
Neumann regular rings, or if $\coring{C}$ and $\coring{D}$ are
coseparable. Finally, if $\coring{C}$ and $\coring{D}$ are
coseparable, or if $X$ and $\Lambda$ are coflat on both sides, or if
$A, B$ are von Neumann regular rings, then (\ref{adj1}),
(\ref{semicontext}) and (\ref{adj1symm}) are equivalent.
\end{proposition}

\begin{proof}
The equivalence between (\ref{adj1}) and (\ref{adj2}) follows from
Proposition \ref{adjointcohomcotens}. That (\ref{adj2}) implies
(\ref{adjsep}) is a consequence of Proposition \ref{G4}. If
$\coring{C}$ is coseparable and we assume (\ref{adjsep}) then, by
Corollary \ref{5}, $\cohom{\coring{C}}{\Lambda}{-} \simeq -
\cotensor{\coring{C}} \cohom{\coring{C}}{\Lambda}{\coring{C}} \simeq
- \cotensor{\coring{C}} X$. That (\ref{adj1}) implies
(\ref{semicontext}) follows from Lemma \ref{18} by evaluating the
unit and the counit of the adjunction at $\coring{C}$ and
$\coring{D}$, respectively. Conversely, if we put
$F=-\cotensor{\coring{C}}X$ and $G=-\cotensor{\coring{D}}\Lambda$,
we have $GF \simeq
-\square_{\coring{C}%
}(X\square_{\coring{D}}\Lambda)$ and $FG \simeq -\square
_{\coring{D}}(\Lambda\square_{\coring{C}}X)$ by Proposition
\ref{assocotens1}$. $ Define natural transformations
$$\xymatrix@1{\eta :1_{\mathcal{M}^{\coring{C}}}\ar[r]^\simeq &
-\square _{\coring{C}}\coring{C}\ar[r]^-{-\square_{\coring{C}}\psi}
& GF}$$
 and
$$ \xymatrix@1{\varepsilon:FG\ar[r]^{-\square_{\coring{D}}\omega} &
-\square_{\coring{D}}\coring{D}\ar[r]^\simeq &
1_{\mathcal{M}^{\coring{D}}},} $$ which become the unit and the
counit of an adjunction by \eqref{unitcounit}. This gives the
equivalence between (\ref{adj1}) and (\ref{semicontext}). The
equivalence between (\ref{semicontext}) and (\ref{adj1symm}) follows
by symmetry.
\end{proof}

\begin{proposition}\label{19}
Suppose that $_A\coring{C}$ and $_B\coring{D}$ are flat. Let
$X\in\bcomod{\coring{C}}{\coring{D}}$ and
$\Lambda\in\bcomod{\coring{D}}{\coring{C}}.$ The following
statements are equivalent
\begin{enumerate}[(i)]
\item $-\square _{\coring{C}}X$ is left adjoint to
$-\square_{\coring{D}}\Lambda$, and $- \cotensor{\coring{C}} X$ is
left exact (or ${}_AX$ is flat or ${}_{\coring{C}}X$ is coflat);
\item $\Lambda$ is quasi-finite as a right $\coring{C}$-comodule,
$- \square_{\coring{C}} X \simeq \cohom{\coring{C}}{\Lambda}{-}$,
and $- \square_{\coring{C}}X$ is left exact (or ${}_AX$ is flat or
${}_{\coring{C}}X$ is coflat); \item
 $\Lambda$ is quasi-finite and injector as a right
$\coring{C}$-comodule and $X\simeq
\cohom{\coring{C}}{\Lambda}{\coring{C}}$ in
$^{\coring{C}}\mathcal{M}^{\coring{D}}$.
\end{enumerate}
\end{proposition}

\begin{proof}
First, observe that if $_{\coring{C}}X$ is coflat, then $_AX$ is
flat (see Corollary \ref{coflat flat} (2)), and that if $_AX$ is
flat, then the functor $- \cotensor{\coring{C}} X$ is left exact.
Thus, in view of Proposition \ref{19a}, it suffices if we prove that
the version of $(ii)$ with $-\cotensor{\coring{C}} X$ left exact
implies $(iii)$, and this last implies the version of $(ii)$ with
${}_{\coring{C}}X$ coflat. Assume that $- \cotensor{\coring{C}}X
\simeq \cohom{\coring{C}}{\Lambda}{-}$ with $-
\cotensor{\coring{C}}X$ left exact. By Proposition \ref{G4}, $X
\simeq \cohom{\coring{C}}{\Lambda}{\coring{C}}$ in
$\bcomod{\coring{C}}{\coring{D}}$. Being a left adjoint,
$\cohom{\coring{C}}{\Lambda}{-}$ is right exact and, henceforth,
exact. By Proposition \ref{injector}, $\Lambda_{\coring{C}}$ is an
injector and we have proved $(iii)$. Conversely, if
$\Lambda_{\coring{C}}$ is a quasi-finite injector and $X \simeq
\cohom{\coring{C}}{\Lambda}{\coring{C}}$ as bicomodules, then $ -
\cotensor{\coring{C}} X \simeq - \cotensor{\coring{C}}
\cohom{\coring{C}}{\Lambda}{\coring{C}}$ and, by Proposition
\ref{injector}, we get that $\cohom{\coring{C}}{\Lambda}{-}$ is an
exact functor. By Corollary \ref{5}, $\cohom{\coring{C}}{\Lambda}{-}
\simeq -\cotensor{\coring{C}}
\cohom{\coring{C}}{\Lambda}{\coring{C}} \simeq -
\cotensor{\coring{C}} X$, and ${}_{\coring{C}}X$ is coflat.
\end{proof}

From the foregoing propositions, it is easy to deduce our
characterization of Frobenius functors between categories of
comodules over corings.

\begin{theorem}\label{20a}
Suppose that $_A\coring{C}$ and $_B\coring{D}$ are flat. Let
$X\in{}^{\coring{C}}\mathcal{M}%
^{\coring{D}}$ and
$\Lambda\in{}^{\coring{D}}\mathcal{M}^{\coring{C}}.$ The following
statements are equivalent
\begin{enumerate}[(i)]
\item $(-\cotensor{\coring{C}} X, - \cotensor{\coring{D}}
\Lambda)$ is a Frobenius pair; \item $ - \cotensor{\coring{C}} X$ is
a Frobenius functor, and $\cohom{\coring{D}}{X}{\coring{D}} \simeq
\Lambda$ as bicomodules; \item there is a Frobenius pair $(F,G)$ for
$\rcomod{\coring{C}}$ and $\rcomod{\coring{D}}$ such that
$F(\coring{C}) \simeq \Lambda$ and $G(\coring{D}) \simeq X$ as
bicomodules; \item $\Lambda_{\coring{C}}, X_{\coring{D}}$ are
quasi-finite injectors,  and $X \simeq
\cohom{\coring{C}}{\Lambda}{\coring{C}}$ and $\Lambda \simeq
\cohom{\coring{D}}{X}{\coring{D}}$ as bicomodules; \item
$\Lambda_{\coring{C}}, X_{\coring{D}}$ are quasi-finite, and $-
\cotensor{\coring{C}} X \simeq \cohom{\coring{C}}{\Lambda}{-}$ and $
- \cotensor{\coring{D}} \Lambda \simeq \cohom{\coring{D}}{X}{-}$.
\end{enumerate}
\end{theorem}

\begin{proof}
$(i) \Leftrightarrow (ii) \Leftrightarrow (iii)$ This is obvious,
after Theorem \ref{3} and Proposition \ref{G4}.

$(i) \Leftrightarrow (iv)$ Follows from Proposition \ref{19}.

$(iv) \Leftrightarrow (v)$ If $X_{\coring{D}}$ and
$\Lambda_{\coring{C}}$ are quasi-finite, then ${}_AX$ and
${}_B\Lambda$ are flat. Now, apply Proposition \ref{19}.
\end{proof}

From Proposition \ref{19a} and Proposition \ref{19} (or Theorem
\ref{20a}) we get the following

\begin{theorem}\label{leftrightFrobenius}
Let $X\in{}^{\coring{C}}\mathcal{M}%
^{\coring{D}}$ and
$\Lambda\in{}^{\coring{D}}\mathcal{M}^{\coring{C}}.$ Suppose that
$_A\coring{C}$, $\coring{C}_A$, $_B\coring{D}$ and $\coring{D}_B$
are flat. The following statements are equivalent
\begin{enumerate} \item $(-\cotensor{\coring{C}} X, -
\cotensor{\coring{D}}\Lambda)$ is a Frobenius pair, with
$X_{\coring{D}}$ and $\Lambda_{\coring{C}}$ coflat; \item $(\Lambda
\cotensor{\coring{C}} -, X \cotensor{\coring{D}}-)$ is a Frobenius
pair, with ${}_{\coring{C}}X$ and ${}_{\coring{D}}\Lambda$ coflat;
\item $X$ and $\Lambda$ are coflat quasi-finite injectors on both
sides, and $X\simeq \cohom{\coring{C}}{\Lambda}{\coring{C}}$ in
$^{\coring{C}}\mathcal{M}^{\coring{D}}$ and $\Lambda\simeq
\cohom{\coring{D}}{X}{\coring{D}}$ in
$^{\coring{D}}\mathcal{M}^{\coring{C}}$.
 \end{enumerate}
  If moreover $\coring{C}$ and $\coring{D}$ are coseparable (resp. $A$
and $B$ are von Neumann regular rings), then the following
statements are equivalent
\begin{enumerate}
\item $(-\cotensor{\coring{C}} X,- \cotensor{\coring{D}}\Lambda)$
is a Frobenius pair; \item $(\Lambda \cotensor{\coring{C}} -, X
\cotensor{\coring{D}}-)$ is a Frobenius pair;\item $X$ and $\Lambda$
are quasi-finite (resp. quasi-finite injectors) on both sides, and
$X\simeq \cohom{\coring{C}}{\Lambda}{\coring{C}}$ in
$^{\coring{C}}\mathcal{M}^{\coring{D}}$ and $\Lambda\simeq
\cohom{\coring{D}}{X}{\coring{D}}$ in
$^{\coring{D}}\mathcal{M}^{\coring{C}}$.
\end{enumerate}
\end{theorem}

\begin{remark}\label{damodulos}
In the case of rings (i.e., $\coring{C} = A$ and $\coring{D} = B$,
Theorem \ref{20a} and the second part of Theorem
\ref{leftrightFrobenius} give \cite[Theorem
2.1]{Castano/Gomez/Nastasescu:1999}. To see this, observe that
${}_AX_B$ is quasi-finite as a right $B$-module if and only if
$-\tensor{A} X : \rmod{A} \rightarrow \rmod{B}$ has a left adjoint,
that is, if and only if ${}_AX$ is finitely generated and
projective. In such a case, the left adjoint is
$-\tensor{B}\hom{A}{X}{A} : \rmod{B} \rightarrow \rmod{A}$. Of
course, $- \cotensor{A} X = - \tensor{A} X$.

The dual characterization in the framework of coalgebras over fields
\cite[Theorem 3.3]{Castano/Gomez/Nastasescu:1999} will be deduced in
Section \ref{condualidad}.
\end{remark}

\section{Frobenius functors between corings with a
duality}\label{condualidad}

We will look to Frobenius functors for corings closer to coalgebras
over fields, in the sense that the categories of comodules share a
fundamental duality.

First we recall a finiteness theorem:

\begin{lemma}\emph{\cite[19.12(1)]{Brzezinski/Wisbauer:2003}}\label{BW 19.12}
\\Assume that $_A\coring{C}$ is locally projective and let $M\in
\rcomod{\coring{C}}$. Then every nonempty finite subset of $M$ is a
subset of a subcomodule of $M$ which is finitely generated as a
right $A$-module. In particular, minimal
$\ldual{\coring{C}}$-submodules are finitely generated as right
$A$-modules.
\end{lemma}

\begin{proposition}\label{locally}
Let $\coring{C}$ be an $A$-coring such that $_A\coring{C}$ is flat.
\begin{enumerate}[(1)]
\item\label{fg} A comodule $M \in \rcomod{\coring{C}}$ is finitely
generated if and only if $M_A$ is finitely generated.
\item\label{fp} A comodule $M \in \rcomod{\coring{C}}$ is finitely
presented if $M_A$ is finitely presented. The converse is true
whenever $\rcomod{\coring{C}}$ is locally finitely generated.
\item\label{lfg} If $_A\coring{C}$ is locally projective, then
$\rcomod{\coring{C}}$ is locally finitely generated.
\end{enumerate}
\end{proposition}

\begin{proof}
(1) From the definition of a finitely generated object, every
comodule that is finitely generated as right $A$-module is a
finitely generated comodule. The forgetful functor $U :
\rcomod{\coring{C}} \to \rmod{A}$ has a left adjoint $ -\tensor{A}
\coring{C} : \rmod{A} \to \rcomod{\coring{C}}$ which is exact and
preserves coproducts. Thus, from Lemma \ref{lefadjfunf.g.obj}, $U$
preserves finitely generated objects.

(2) If $M \in \rcomod{\coring{C}}$ is such that $M_A$ is finitely
presented, then for every exact sequence $0 \to K \to L \to M \to 0$
in $\rcomod{\coring{C}}$ with $L$ finitely generated, we get an
exact sequence $0 \to K_A \to L_A \to M_A \to 0$ with $M_A$ finitely
presented. Thus, $K_A$ is finitely generated and, by (\ref{fg}), $K
\in \rcomod{\coring{C}}$ is finitely generated. This proves that $M$
is a finitely presented comodule. The converse follows directly from
Lemma \ref{V.3.4} and the fact that the forgetful functor preserves
inductive limits.

(3) is a consequence of (\ref{fg}) and Lemma \ref{BW 19.12}.
\end{proof}

\begin{proposition}
\begin{enumerate}[(1)]
Let $\coring{C}$ be an $A$-coring.
\item If $_A\coring{C}$ is flat, $A$ is a right noetherian algebra,
and $\coring{C}_A$ is finitely generated module, then the category
$\rcomod{\coring{C}}$ is locally finitely noetherian.
\item If $_A\coring{C}$ is flat and $A$ is a right noetherian
algebra, then every finitely generated object of
$\rcomod{\coring{C}}$ is noetherian.
\item If $_A\coring{C}$ is locally projective and $A$ is a right
noetherian algebra, then the category $\rcomod{\coring{C}}$ is
locally finitely noetherian.
\end{enumerate}
\end{proposition}

\begin{proof}
(1) Immediate from Lemma \ref{comodgenerators} and the fact that
every $M\in \rcomod{\coring{C}}$ such that $M_A$ is noetherian, is
noetherian in $\rcomod{\coring{C}}$ (by Proposition \ref{V.4.1}).

(2) Let $M\in \rcomod{\coring{C}}$ be a finitely generated comodule.
From Proposition \ref{locally}, $M_A$ is finitely generated. Then it
is noetherian. Hence $M_\coring{C}$ is noetherian.

(3) Obvious from Proposition \ref{locally} and (2).
\end{proof}

The notation $\cat{C}_f$ stands for the full subcategory of a
Grothendieck category $\cat{C}$ whose objects are the finitely
generated ones. The special case of the following lemma where
$\cat{C}$ is the category of right modules over a ring $R$, recovers
\cite[Exercise 4, p. 109]{Stenstrom:1975}.

\begin{lemma}\label{10}
Let $\cat{C}$ be a locally finitely generated category.
\begin{enumerate}[(1)]
\item\label{additive} The category $\cat{C}_f$ is additive.
\item\label{cokernels} The category $\cat{C}_f$ has cokernels, and
every monomorphism in $\cat{C}_f$ is a monomorphism in $\cat{C}$.
\item\label{locnoeth} The following statements are equivalent:
  \begin{enumerate}[(a)]
   \item\label{haskernels} The category $\cat{C}_f$ has kernels;
   \item\label{locnoeth2} $\cat{C}$ is locally noetherian;
   \item\label{abelian} $\cat{C}_f$ is abelian;
   \item\label{abeliansub} $\cat{C}_f$ is an abelian subcategory of
$\cat{C}$.
  \end{enumerate}
\end{enumerate}
\end{lemma}

\begin{proof}
(1) Obvious from Lemma \ref{abeliansubcategory} and Proposition
\ref{fgfpnoeth}(1).

 (2) That $\mathbf{C}_{f}$ has cokernels is
straightforward from Lemma \ref{V.3.1} (i). Now, let $f:M\to N$ be a
monomorphism in $\mathbf{C}_{f}$ and $\xi:X\rightarrow M$ be a
morphism in $\mathbf{C}$ such that $f\xi=0.$ Suppose that
$X=\bigcup_{i\in I}X_i$, where $X_i\in\mathbf{C}_{f}$, and
$\iota_i:X_i\rightarrow X$, $i\in I$, the canonical injections. Then
$f\xi\iota_i=0$, and $\xi\iota_i=0$,
for every $i$, and by the definition of the inductive limit, $\xi=0.$%

(3) $(b)\Rightarrow(a)$ Straightforward from Proposition
\ref{V.4.1}.

 $(d)\Rightarrow(c)$ and $(c)\Rightarrow(a)$ are trivial.

$(a)\Rightarrow(b)$ Let $M\in\mathbf{C}_{f}$, and $K$ be a subobject
of $M.$ Let $\iota :L\rightarrow M$ be the kernel of the canonical
morphism $f:M\rightarrow M/K$ in $\mathbf{C}_{f}.$ Suppose that
$K=\bigcup_{i\in I}K_i$, where $K_i\in\mathbf{C}_{f}$, for every
$i\in I.$ By the universal property of kernel, there exist a unique
morphism $\alpha:L\rightarrow K,$ and a unique morphism
$\beta_i:K_i\rightarrow L$, for every $i\in I$, making commutative
the diagrams
\[%
\xymatrix{
&   K_i \ar[d] \ar[dl]_{\beta_i} & \\
L \ar[dr]_{\alpha} \ar[r]^{\iota} & M \ar[r]^{f} & M/K \\
   & K \ar[u] & .}
\]
By (2), $\iota$ is a monomorphism in $\mathbf{C}$, then for every
$K_i\subset K_{j}$, the diagram
\[%
\xymatrix{
K_i  \ar[d] \ar[r]^{\beta_i} & L\\
\ar[ur]_{\beta_{j}}  K_{j} &   }
\]
commutes. Therefore we have the commutative diagram
\[%
\xymatrix{
& K \ar[d] \ar[dl]_{{\lim \atop {\longrightarrow \atop I}}\beta_i}\\
L \ar[r]^{\iota} & M.}
\]
Then $K\simeq L$, and hence $K\in\mathbf{C}_{f}.$ Finally, by
Proposition \ref{V.4.1}, $M$ is noetherian in $\mathbf{C.}$

$(b)\Rightarrow(d)$ Straightforward from Lemma
\ref{abeliansubcategory}.
\end{proof}

The following generalization of \cite[Proposition
3.1]{Castano/Gomez/Nastasescu:1999} will allow us to give an
alternative proof to the equivalence ``$(1)\Leftrightarrow(4)$'' of
Theorem \ref{leftrightFrobeniusduality} (see Remark
\ref{dacomodulos} (1)).

\begin{proposition}\label{11}
Let $\mathbf{C}$ and $\mathbf{D}$ be two locally noetherian
categories. Then
\begin{enumerate}[(1)]
\item If $F:\mathbf{C}\rightarrow\mathbf{D}$ is a Frobenius
functor, then its restriction $F_f:\mathbf{C}_{f}\rightarrow
\mathbf{D}_{f}$ is a Frobenius functor.
\item If $H:\mathbf{C}%
_{f}\rightarrow\mathbf{D}_{f}$ is a Frobenius functor, then $H$ can
be
uniquely extended to a Frobenius functor $\overline H:\mathbf{C}%
\to\mathbf{D}.$ \item The assignment $F\mapsto F_f$ defines a
bijective correspondence (up to natural isomorphisms) between
Frobenius functors from $\mathbf{C}$ to $\mathbf{D}$ and Frobenius
functors from $\mathbf{C}_{f}$ to $\mathbf{D}_f.$
\item In particular, if
$\mathbf{C}=\mathcal{M}^{\coring{C}}$ and $\mathbf{D}=\mathcal{M}%
^\coring{D}$ are locally noetherian such that $_A\coring{C}$ and
$_B\coring{D}$ are flat, then
$F:\mathcal{M}^{\coring{C}}\rightarrow\mathcal{M}^\coring{D}$ is a
Frobenius functor if and only if it preserves direct limits and
comodules which are finitely generated as right $A$-modules,
and the restriction functor $F_f:\mathcal{M}_{f}^{\coring{C}}%
\rightarrow\mathcal{M}_{f}^{\coring{D}}$ is a Frobenius functor.
\end{enumerate}
\end{proposition}

\begin{proof}
The proofs of \cite[Proposition 3.1 and Remark
3.2]{Castano/Gomez/Nastasescu:1999} remain valid for our situation,
but with some minor modifications: to prove that $\overline{H}$ is
well-defined, we use Lemma \ref{10}. In the proof of the statements
(1), (2) and (3) we use the Grothendieck AB 5) condition.
\end{proof}

In order to generalize \cite[Proposition A.2.1]{Takeuchi:1977ca} and
its proof, we need the following lemma.

\begin{lemma}\label{12}
\begin{enumerate}[(1)] \item Let $\mathbf{C}$ be a locally noetherian
category, let $\mathbf{D}$ be an arbitrary Grothendieck category,
$F:\mathbf{C}\rightarrow\mathbf{D}$ be an
arbitrary functor which preserves direct limits, and $F_f:\mathbf{C}%
_{f}\rightarrow\mathbf{D}$ be its restriction to $\mathbf{C}_{f}.$
Then $F$ is exact (faithfully exact, resp. left, right exact) if and
only if $F_f$ is exact (faithfully exact, resp. left, right exact).

 In particular, an object $M$ in $\cat{C}_f$ is projective
(resp. projective generator) if and only if it is projective (resp.
projective generator) in $\cat{C}$.

\item Let $\mathbf{C}$ be a locally noetherian category. For every
object $M$ of $\mathbf{C}$, the following conditions are equivalent
\begin{enumerate}[(a)]
\item $M$ is injective (resp. an injective cogenerator); \item the
contravariant functor $\hom{\mathbf{C}}{-}{M} :\mathbf{C}\rightarrow
\mathbf{Ab}$ is exact (resp. faithfully exact); \item the
contravariant
functor $\hom{\mathbf{C}}{-}{M}_{f}:\mathbf{C}_{f}%
\rightarrow\mathbf{Ab}$ is exact (resp. faithfully exact).
\end{enumerate}
 In particular, an object $M$ in $\cat{C}_f$ is injective (resp.
injective cogenerator) if and only if it is injective (resp.
injective cogenerator) in $\cat{C}$.
\end{enumerate}
\end{lemma}

\begin{proof}
(1) The ``only if'' part is straightforward from the fact that the
injection functor $\mathbf{C}_f\rightarrow\mathbf{C}$ is faithfully
exact.

 For the ``if'' part, suppose
that $F_f$ is left exact. Let $f:M\rightarrow N$ be a morphism in
$\mathbf{C}$. Put $M=\bigcup_{i\in I}M_i$ and $N=\bigcup_{j\in
J}N_{j},$ as direct union of directed families of
finitely generated subobjects. For $(i,j)\in I\times J,$ let $M_{i,j}%
=M_i \cap f^{-1}(N_j)$, and $f_{i,j}:M_{i,j}\rightarrow N_{j}$ be
the restriction of $f$ to $M_{i,j}.$ We have
$f=\underset{\underset{I\times J}{\longrightarrow}}{\lim}f_{i,j}$
and then $F(f)=\underset {\underset{I\times
J}{\longrightarrow}}{\lim}F_f(f_{i,j}).$ Hence

 $\ker F(f)=\ker\underset{\underset{I\times
J}{\longrightarrow}}{\lim}F_f(f_{i,j})=\underset {\underset{I\times
J}{\longrightarrow}}{\lim}\ker F_f( f_{i,j})
=\underset{\underset{I\times J}{\longrightarrow}}{\lim}F_f(\ker
f_{i,j})=\underset{\underset{I\times
J}{\longrightarrow}}{\lim}F(\ker f_{i,j})$ (by Lemma \ref{10})
$=F(\underset{\underset{I\times
J}%
{\longrightarrow}}{\lim}\ker f_{i,j})=F(\ker f).$

 Finally $F$ is left exact. Analogously, it can be proved that $F_f$ is
right exact implies that $F$ is also right exact. Now, suppose that
$F_f$ is faithfully exact. We have already proved that $F$ is exact.
It remains to prove that $F$ is faithful. For this, let $0\neq
M=\bigcup_{i\in I}M_i$ be an object of $\mathbf{C}$, where $M_i$ is
finitely generated for every $i\in I.$ We have
\[
F(M)=\underset{\underset{I}{\longrightarrow}}{\lim}F_f(M_i)
\simeq\sum_iF_f(M_i)
\]
(since $F$ is exact). Since $M\neq0$, there exists some $i_{0}\in I$
such that $M_{i_{0}}\neq0.$ By Proposition \ref{faithful-exact},
$F_f(M_{i_{0}})\neq0,$ hence $F(M)\neq0.$ Also by Proposition
\ref{faithful-exact}, $F$ is faithful.

(2) $(a)\Leftrightarrow(b)$ Obvious.

 $(b)\Rightarrow(c)$
Analogous to that of the ``only if'' part of (1).

$(c)\Rightarrow(a)$ That $M$ is injective is a consequence of
Propositions \ref{Baer}, \ref{V.4.1}. Now, suppose moreover that
$\hom{\mathbf{C}}{-}{M}_{f}$ is faithful. Let $L$ be a nonzero
object of $\mathbf{C}$, and $K$ be a nonzero finitely generated
subobject of $L.$ By Proposition \ref{faithful-exact}, there exists
a nonzero morphism $K\to M$. Since $M$ is injective, there exists a
nonzero morphism $L\rightarrow M$ making commutative the following
diagram
\[%
\xymatrix{
&   M   & \\
0  \ar[r] & K \ar[r] \ar[u] & L. \ar[ul]}
\]
From the dual of Proposition \ref{proj-gen}, it follows that $M$ is
a cogenerator.
\end{proof}

\begin{lemma}\emph{\cite[19.19 (1)]{Brzezinski/Wisbauer:2003}}\label{comod.f.p.mod}
\\ Let $M$ be a right $\coring{C}$-module. If
\begin{enumerate}[(a)]
\item $M_A$ is finitely generated and projective, or
\item $\coring{C}_A$ is flat and $M_A$ is finitely presented, or
\item $\coring{C}_A$ is finitely generated and projective,
\end{enumerate}
then the dual module $\rdual{M}=\hom{A}{M}{A}$ is a left
$\coring{C}$-module via the structure map
$$\lambda_{\rdual{M}}:\rdual{M}\to \hom{A}{M}{\coring{C}}\simeq
\coring{C}\tensor{A}\rdual{M}, \quad f\mapsto
(f\tensor{A}\coring{C})\circ\rho_M.$$ (See Lemma
\ref{hom-tens-relations}(6).)
\end{lemma}

If $\coring{C}_A$ is flat and $M \in \rcomod{\coring{C}}$ is
finitely presented as a right $A$-module, then, by Lemma
\ref{comod.f.p.mod}, the dual $A$-module $\rdual{M}$ has a left
$\coring{C}$-comodule structure. Now, if ${}_A\rdual{M}$ turns out
to be finitely presented and ${}_A\coring{C}$ is flat, then
$\ldual{(\rdual{M})} = \hom{A}{\rdual{M}}{A}$ is a right
$\coring{C}$-comodule and the canonical map $\sigma_M : M
\rightarrow \ldual{(\rdual{M})}$ is a homomorphism in
$\rcomod{\coring{C}}$. This construction leads to a duality
\[
\rDual: \rcomod{\coring{C}}_0 \leftrightarrows
\lcomod{\coring{C}}_0: \lDual
\]
between the full subcategories $\rcomod{\coring{C}}_0$ and
$\lcomod{\coring{C}}_0$ of $\rcomod{\coring{C}}$ and
$\lcomod{\coring{C}}$ whose objects are the comodules which are
finitely generated and projective over $A$ on the corresponding side
(this holds even without flatness assumptions of $\coring{C}$). Call
it the \emph{basic duality} (details may be found in
\cite{Caenepeel/DeGroot/Vercruysse:unp}). Of course, in the case
that $A$ is semisimple (e.g. for coalgebras over fields) these
categories are that of finitely generated comodules, and this basic
duality plays a remarkable role in the study of several notions in
the coalgebra setting (e.g. Morita equivalence \cite{Takeuchi:1977},
semiperfect coalgebras \cite{Lin:1977}, Morita duality
\cite{Gomez/Nastasescu:1995}, \cite{Gomez/Nastasescu:1996}, or
Frobenius functors \cite{Castano/Gomez/Nastasescu:1999}). It would
be interesting to know, in the coring setting, to what extent the
basic duality can be extended to the subcategories
$\rcomod{\coring{C}}_f$ and $\lcomod{\coring{C}}_f$, since, as we
will try to show in this section, this allows to obtain better
results. Of course, this is the underlying idea when the ground ring
$A$ is assumed to be quasi-Frobenius (see
\cite{ElKaoutit/Gomez:2002unp}).

Consider contravariant functors between Grothendieck categories $ H:
\cat{A} \leftrightarrows \cat{A}': H', $ together with natural
transformations $\tau : 1_{\cat{A}} \rightarrow H' \circ H$ and
$\tau' : 1_{\cat{A}'} \rightarrow H \circ H'$, satisfying the
condition $H(\tau_X) \circ \tau'_{H(X)} = 1_{H(X)}$ and
$H'(\tau'_{X'}) \circ \tau_{H'(X')} = 1_{H'(X')}$ for $X \in
\cat{A}$ and $X' \in \cat{A}'$. Following \cite{Colby/Fuller:1983},
this situation is called a \emph{right adjoint pair}.

\begin{proposition}\label{lnoethrightadj}
Let $\coring{C}$ be an $A$-coring such that ${}_A\coring{C}$ and
$\coring{C}_A$ are flat. Assume that $\rcomod{\coring{C}}$ and
$\lcomod{\coring{C}}$ are locally noetherian categories. If
${}_A\rdual{M}$ and $\ldual{N}_A$ are finitely generated modules for
every $M \in \rcomod{\coring{C}}_f$ and $N \in
\lcomod{\coring{C}}_f$, then the basic duality extends to a right
adjoint pair $ \rDual: \rcomod{\coring{C}}_f \leftrightarrows
\lcomod{\coring{C}}_f: \lDual$.
\end{proposition}

\begin{proof}
If $M \in \rcomod{\coring{C}}_f$ then, since $\rcomod{\coring{C}}$
is locally noetherian, $M_{\coring{C}}$ is finitely presented. By
Proposition \ref{locally}, $M_A$ is finitely presented and the left
$\coring{C}$-comodule $\rdual{M}$ makes sense. Now, the assumption
${}_A\rdual{M}$ finitely generated implies, by Lemma \ref{locally},
that $\rdual{M} \in \lcomod{\coring{C}}_f$. We have then the functor
$\rDual : \rcomod{\coring{C}}_f \rightarrow \lcomod{\coring{C}}_f$.
The functor $\lDual$ is analogously defined, and the rest of the
proof consists of straightforward verifications.
\end{proof}

\begin{example}
The hypotheses are fulfilled if $_A\coring{C}$ and $\coring{C}_A$
are locally projective and $A$ is left and right noetherian. But
there are situations in which no finiteness condition need to be
required to $A$: this is the case, for instance, of cosemisimple
corings (see Proposition \ref{cosemisimple}). In particular, if an
arbitrary ring $A$ contains a division ring $B$, then, by
Proposition \ref{cosemisimple}, the canonical coring $A \tensor{B}
A$ satisfies all hypotheses in Proposition \ref{lnoethrightadj}.
\end{example}

\begin{definition}
Let $\coring{C}$ be a coring over $A$ satisfying the assumptions of
Proposition \ref{lnoethrightadj}. We will say that $\coring{C}$
\emph{has a duality} if the basic duality extends to a duality $$
\rDual: \rcomod{\coring{C}}_f \leftrightarrows
\lcomod{\coring{C}}_f: \lDual.$$
\end{definition}

We have the following examples of a coring which has a
duality:
\begin{enumerate}[(i)]
\item $\coring{C}$ is a coring over a QF
ring $A$ such that ${}_A\coring{C}$ and $\coring{C}_A$ are flat (and
hence projective);\item $\coring{C}$ is a cosemisimple coring, where
$\rcomod{\coring{C}}_f$ and $\lcomod{\coring{C}}_f$ are equal to
$\rcomod{\coring{C}}_0$ and $\lcomod{\coring{C}}_0$, respectively;
\item $\coring{C}$ is a coring over $A$ such that ${}_A\coring{C}$
and $\coring{C}_A$ are flat and semisimple, $\rcomod{\coring{C}}$
and $\lcomod{\coring{C}}$ are locally noetherian categories, and the
dual of every simple right (resp. left) $A$-module in the
decomposition of $\coring{C}_A$ (resp. $_A\coring{C}$) as a direct
sum of simple $A$-modules is finitely generated and $A$-reflexive
(in fact, every right (resp. left) $\coring{C}$-comodule $M$ becomes
a submodule of the semisimple right (resp. left) $A$-module
$M\otimes_A\coring{C}$, and hence $M_A$ (resp. $_AM$) is also
semisimple).
\end{enumerate}

Let $\coring{C}$ be an $A$-coring such that $_A\coring{C}$ is flat,
and $M\in\rcomod{\coring{C}}$. A comodule $Q\in\rcomod{\coring{C}}$
is said to be $M$-\emph{injective} if for every subcomodule $M'$ of
$M$, and every morphism in $\rcomod{\coring{C}}$, $f:M'\to Q$, there
is $g:M\to Q$ such that the diagram
$$\xymatrix{&&Q\\
0\ar[r] & M' \ar[ur]^f\ar[r] & M\ar@{-->}[u]^g}$$ commutes.
Obviously, $Q\in\rcomod{\coring{C}}$ is injective if and only if it
is $M$-injective for all $M\in\rcomod{\coring{C}}$.

\medskip
The following is a generalization of \cite[Proposition
A.2.1]{Takeuchi:1977ca} and
$(i)\Leftrightarrow(ii)\Leftrightarrow(iv)$ of \cite[Theorem
2.4.17]{Dascalescu/Nastasescu/Raianu: 2001}.

\begin{proposition}\label{13}
Suppose that the coring $\coring{C}$ has a duality. Let
$Q\in\rcomod{\coring{C}}$ such that $Q_A$ is flat. The following are
equivalent
\begin{enumerate}[(1)]
\item $Q$ is coflat (resp. faithfully coflat);
\item $Q$ is $M$-injective for all finitely generated right
$\coring{C}$-comodule $M$ (resp. the last statement hold and
moreover $0\neq M\in \rcomod{\coring{C}}$ is finitely generated
implies $\hom{\coring{C}}{M}{Q}\neq 0$;
\item $Q$ is injective (resp. an injective cogenerator).
\end{enumerate}
\end{proposition}

\begin{proof}
Let $Q\in\rcomod{\coring{C}}$ and
$N\in{}^{\coring{C}}\mathcal{M}_{f}.$ We have the following
commutative diagram (in $\mathcal{M}_{k}$)
\[%
\xymatrix{  0 \ar[r] & Q\square_{\coring{C}}N  \ar[r] &
Q\otimes_AN \ar[rrr]^{\rho_Q\otimes_AN-Q\otimes_A\lambda_N}
\ar[d] & & & Q\otimes_A\coring{C}\otimes_AN \ar[d] \\
 0 \ar[r] & \hom{\coring{C}}{N^*}{Q}\ar[r] &
\hom{A}{N^*}{Q}\ar[rrr]^{f\mapsto\rho_Qf-(f\otimes_A\coring{C})\rho_{N^*}
} & & & \hom{A}{N^*}{Q\otimes_A\coring{C}},}
\]
where the vertical maps are the canonical maps. By the universal
property of kernel, there is a unique morphism
$\eta_{Q,N}:Q\square_{\coring{C}}N\to\hom{\coring{C}}{N^*}{Q}$
making commutative the above diagram. By Lemma \ref{cube}, $\eta$ is
a natural transformation of bifunctors. For example we will verify
that $\eta$ is natural in $N$. For this let $f:N_1\to N_2$ be a
morphism in $^{\coring{C}}\mathcal{M}_{f}$. It is enough to consider
the diagram
$$\xymatrix{
 && Q\cotensor{\coring{C}}N_1 \ar[rrr]^{\eta_{Q,N_1}}
 \ar[dll]^(.3){Q\cotensor{\coring{C}}f}\ar'[d][ddd]&&&
 \hom{\coring{C}}{\rdual{N_1}}{Q}
\ar[dll]^{\hom{\coring{C}}{\rdual{f}}{Q}}\ar[ddd]
\\
Q\cotensor{\coring{C}}N_2 \ar[rrr]^(.35){\eta_{Q,N_2}}\ar[ddd]
&&&\hom{\coring{C}}{\rdual{N_2}}{Q} \ar[ddd]
\\
&&&&&&&&&
\\
 && Q\tensor{A}N_1 \ar'[r][rrr] \ar[dll]^{Q\tensor{A}f} &&&
 \hom{A}{\rdual{N_1}}{Q}\ar[dll]^{\hom{A}{\rdual{f}}{Q}}
\\
Q\tensor{A}N_2 \ar[rrr]  &&& \hom{A}{\rdual{N_2}}{Q}.}$$
If $Q_A$ is
flat then $\eta_{Q,N}$ is an isomorphism for every
$N\in{}^{\coring{C}}\mathcal{M}_{f}.$ We have
\[
Q\square_{\coring{C}}-\simeq \hom{\coring{C}}{-}{Q}_{f} \circ
(-)^*:{}^{\coring{C}}\mathcal{M}_{f}\rightarrow \mathcal{M}_{k}.
\]
Then, by Lemma \ref{12}, $Q_{\coring{C}}$ is coflat (resp.
faithfully coflat) iff
$Q\square_{\coring{C}}-:{}^{\coring{C}}\mathcal{M}_{f}\to
\mathcal{M}_{k}$ is exact (resp. faithfully exact)  iff
$\hom{\coring{C}}{-}{Q}_{f}:\mathcal{M}_{f}^{\coring{C}%
}\rightarrow\mathcal{M}_{k}$ is exact (resp. faithfully exact) iff
$Q_{\coring{C}}$ is injective (resp. an injective cogenerator).
\end{proof}

The particular case of the second part of the following result for
coalgebras over a commutative ring is given in
\cite{AlTakhman:2002}.

\begin{corollary}\label{injectiveinjector}
Let $N\in{}^{\coring{C}}\mathcal{M}^{\coring{D}}$ be a bicomodule.
Suppose that $A$ is a QF ring. \begin{enumerate}[(a)]\item If $N$ is
an injector as a right $\coring{D}$-comodule then $N$ is injective
in $\mathcal{M}%
^{\coring{D}}$. \item If $\coring{D}$ has a duality, $N_B$ is flat
and $N$ is injective in $\mathcal{M}%
^{\coring{D}}$, then $N$ is an injector as a right
$\coring{D}$-comodule.
\end{enumerate}
\end{corollary}

\begin{proof}
(a) Since $A$ is a QF ring, then $A_A$ is injective. Hence
$N_\coring{D}\simeq (A\tensor{A}N)_\coring{D}$ is injective.

(b) Let $X_A$ be an injective module. Since $A$ is a QF ring, $X_A$
is projective. We have then the natural isomorphism
\[
(X\otimes_AN)\square_{\coring{D}}-\simeq
X\otimes_A(N\square_{\coring{D}}-)
:{}^{\coring{D}}\mathcal{M}\rightarrow \mathcal{M}_{k}.
\]
By Proposition \ref{13}, $N_\coring{D}$ is coflat, and then
$X\otimes_AN$ is coflat. Now, since $X\otimes_AN$ is a flat right
$B$-module, and by Proposition \ref{13}, $X\otimes_AN$ is injective
in $\mathcal{M}^{\coring{D}}$.
\end{proof}

 The last two results allow to improve our general statements in
 Section \ref{Frobeniusgeneral} for corings having a duality.

\begin{proposition}\label{adjointpairduality}
Suppose that $\coring{C}$ and $\coring{D}$ have a duality. Consider
the following statements
\begin{enumerate} \item $(-\cotensor{\coring{C}} X,-
\cotensor{\coring{D}}\Lambda)$ is an adjoint pair of functors, with
$_AX$ and $\Lambda_A$ flat;\item $\Lambda$ is quasi-finite injective
as a right ${\coring{C}}$-comodule, with $_AX$ and $\Lambda_A$ flat
and $X\simeq \cohom{\coring{C}}{\Lambda}{\coring{C}}$ in
$^{\coring{C}}\mathcal{M}^{\coring{D}}$.
\end{enumerate}
We have (1) implies (2), and the converse is true if in particular
$B$ is a QF ring.
\end{proposition}

\begin{proof}
$(1)\Rightarrow(2)$ From Proposition \ref{13}, $\coring{D}$ is
injective in $\mathcal{M}^{\coring{D}}$. Since the functor
$-\cotensor{\coring{C}} X$ is exact, $\Lambda\simeq
\coring{D}\cotensor{\coring{D}} \Lambda$ is injective in
$\mathcal{M}^{\coring{C}}$.

 $(2)\Rightarrow(1)$ Assume that $B$ is a QF ring. From Corollary
\ref{injectiveinjector}, $\Lambda$ is quasi-finite injector as a
right ${\coring{C}}$-comodule, and Proposition \ref{19} achieves the
proof.
\end{proof}

We are now in a position to state and prove our main result of this
section.

\begin{theorem}\label{leftrightFrobeniusduality}
Suppose that $\coring{C}$ and $\coring{D}$ have a duality. Let
$X\in{}^{\coring{C}}\mathcal{M}%
^{\coring{D}}$ and
$\Lambda\in{}^{\coring{D}}\mathcal{M}^{\coring{C}}.$ The following
statements are equivalent \begin{enumerate} \item
$(-\cotensor{\coring{C}} X, - \cotensor{\coring{D}}\Lambda)$ is a
Frobenius pair, with $X_B$ and $\Lambda_A$ flat; \item $(\Lambda
\cotensor{\coring{C}} -, X \cotensor{\coring{D}}-)$ is a Frobenius
pair, with $_AX$ and $_B\Lambda$ flat; \item $X$ and $\Lambda$ are
quasi-finite injector on both sides, and $X\simeq
\cohom{\coring{C}}{\Lambda}{\coring{C}}$ in
$^{\coring{C}}\mathcal{M}^{\coring{D}}$ and $\Lambda\simeq
\cohom{\coring{D}}{X}{\coring{D}}$ in
$^{\coring{D}}\mathcal{M}^{\coring{C}}$.

In particular, if $A$ and $B$ are QF rings, then the above
statements are equivalent to \item $X$ and $\Lambda$ are
quasi-finite injective on both sides, and $X\simeq
\cohom{\coring{C}}{\Lambda}{\coring{C}}$ in
$^{\coring{C}}\mathcal{M}^{\coring{D}}$ and $\Lambda\simeq
\cohom{\coring{D}}{X}{\coring{D}}$ in
$^{\coring{D}}\mathcal{M}^{\coring{C}}$.
\end{enumerate}
Finally, suppose that $\coring{C}$ and $\coring{D}$ are cosemisimple
corings. Let
$X\in{}^{\coring{C}}\mathcal{M}%
^{\coring{D}}$ and
$\Lambda\in{}^{\coring{D}}\mathcal{M}^{\coring{C}}.$ The following
statements are equivalent \begin{enumerate} \item
$(-\cotensor{\coring{C}} X, - \cotensor{\coring{D}}\Lambda)$ is a
Frobenius pair; \item $(\Lambda \cotensor{\coring{C}} -, X
\cotensor{\coring{D}}-)$ is a Frobenius pair; \item $X$ and
$\Lambda$ are quasi-finite on both sides, and $X\simeq
\cohom{\coring{C}}{\Lambda}{\coring{C}}$ in
$^{\coring{C}}\mathcal{M}^{\coring{D}}$ and $\Lambda\simeq
\cohom{\coring{D}}{X}{\coring{D}}$ in
$^{\coring{D}}\mathcal{M}^{\coring{C}}$.
\end{enumerate}
\end{theorem}

\begin{proof}
We start by proving the first part. In view of Theorem
\ref{leftrightFrobenius} and Theorem \ref{20a} it suffices to show
that if $(-\cotensor{\coring{C}} X, - \cotensor{\coring{D}}\Lambda)$
is a Frobenius pair, the condition ``$X_{\coring{D}}$ and
$\Lambda_{\coring{C}}$ are coflat'' is equivalent to ``$X_B$ and
$\Lambda_A$  are flat''. Indeed, the first implication is obvious,
for the converse, assume that $X_B$ and $\Lambda_A$ are flat. By
Proposition \ref{adjointpairduality}, $X$ and $\Lambda$ are
injective in $\mathcal{M}^{\coring{D}}$ and
$\mathcal{M}^{\coring{C}}$ respectively, and they are coflat by
Proposition \ref{13}. The particular case is straightforward from
Corollary \ref{injectiveinjector} and the above equivalences.

Now we will show the second part. We know that every comodule
category over a cosemisimple coring is a spectral category (see
Proposition \ref{cosemisimple}). From Proposition \ref{13}, the
bicomodules ${}_{\coring{C}}X_{\coring{D}}$ and
${}_{\coring{D}}\Lambda_{\coring{C}}$ are coflat and injector on
both sides (we can see this directly by using Corollary \ref{split
exactness}). Finally, Theorem \ref{leftrightFrobenius} achieves the
proof.
\end{proof}

\begin{remarks}\label{dacomodulos}
(1) The equivalence ``$(1)\Leftrightarrow(4)$'' of the last theorem
is a generalization of \cite[Theorem
3.3]{Castano/Gomez/Nastasescu:1999}. The proof of \cite[Theorem
3.3]{Castano/Gomez/Nastasescu:1999} gives an alternative proof of
``$(1)\Leftrightarrow(4)$'' of Theorem
\ref{leftrightFrobeniusduality}, using Proposition \ref{11}.

(2) The adjunction of Proposition \ref{19a} and Proposition
\ref{adjointpairduality} generalizes the coalgebra version of
Morita's theorem \cite[Theorem
4.2]{Caenepeel/DeGroot/Militaru:2002}.
\end{remarks}

\begin{example}
Let $A$ be a $k$-algebra. Put $\coring{C}=A$ and $\coring{D}=k.$ The
bicomodule $A\in{}^{\coring{C}}\mathcal{M}^{\coring{D}}$ is
quasi-finite as a right $\coring{D}$-comodule. $A$ is an injector as
a right $\coring{D} $-comodule if and only if the $k$-module $A$ is
flat. If we take $A=k=\mathbb{Z} $, the bicomodule $A$ is
quasi-finite and injector as a right $\coring{D}$-comodule but it is
not injective in $\mathcal{M}^{\coring{D}}.$ Hence, the assertion
``$(-\cotensor{\coring{C}} X, - \cotensor{\coring{D}}\Lambda)$ is an
adjoint pair of functors'' does not imply in general the assertion
``$\Lambda$ is quasi-finite injective as a right
${\coring{C}}$-comodule and $X\simeq
\cohom{\coring{C}}{\Lambda}{\coring{C}}$ in
$^{\coring{C}}\mathcal{M}^{\coring{D}}$'', and the following
statements are not equivalent in general:

(1) $(-\cotensor{\coring{C}} X, - \cotensor{\coring{D}}\Lambda)$ is
a Frobenius pair;

(2) $X$ and $\Lambda$ are quasi-finite injective on both sides, and
$X\simeq \cohom{\coring{C}}{\Lambda}{\coring{C}}$ in
$^{\coring{C}}\mathcal{M}^{\coring{D}}$ and $\Lambda\simeq
\cohom{\coring{D}}{X}{\coring{D}}$ in
$^{\coring{D}}\mathcal{M}^{\coring{C}}$.

On the other hand, there exists a commutative self-injective ring
which is not coherent (see Section \ref{category theory}). By a
theorem of S.U. Chase (see for example \cite[Theorem
19.20]{Anderson/Fuller:1992}), there exists then a $k$-algebra $A$
which is injective, but not flat as $k$-module. Hence, the
bicomodule $A\in{}^{\coring{C}}\mathcal{M}^{\coring{D}}$ is
quasi-finite and injective as a right $\coring{D}$-comodule, but not
an injector as a right $\coring{D}$-comodule.
\end{example}

\section{Applications to induction functors}\label{Frobmor}

We will characterize when the induction functor $-\tensor{A}B:
\rcomod{\coring{C}} \rightarrow \rcomod{\coring{D}}$ defined in
Subsection \ref{IndunctionFunctors} is a Frobenius functor.

\begin{theorem}\label{23}
Let $(\varphi,\rho) :\coring{C}\rightarrow\coring{D}$ be a
homomorphism of corings such that $_A\coring{C}$ and $_B\coring{D}$
are flat. The following statements are
equivalent\begin{enumerate}[(a)] \item
$-\otimes_AB:\rcomod{\coring{C}}\to\rcomod{\coring{D}}$ is a
Frobenius functor;\item the $\coring{C}-\coring{D}$-bicomodule
$\coring{C}\otimes_AB$ is quasi-finite and injector as a right
$\coring{D}$-comodule and there exists an isomorphism of
$\coring{D}-\coring{C}$-bicomodules
$\cohom{\coring{D}}{\coring{C}\otimes_AB}{\coring{D}}\simeq B\otimes
_A\coring{C}.$
\end{enumerate}
Moreover, if $\coring{C}$ and $\coring{D}$ are coseparable, then the
condition ``injector'' in (b) can be deleted.
\end{theorem}

\begin{proof}
First observe that $-\otimes_AB$ is a Frobenius functor if and only
if $(-\otimes_AB,-\square_{\coring{D}}(B\otimes _A\coring{C}))$ is a
Frobenius pair (by Proposition \ref{BW 24.11}). Let $M\in
\rcomod{\coring{C}}$.  The map $\rho_M \tensor{A}B: M\tensor{A}B \to
M \tensor{A}\coring{C}\tensor{A}B$ is a morphism of
$\coring{D}$-comodules. As in the proof of Lemma \ref{caniso}, and
using Proposition \ref{bcomodcotensor} ($_B\coring{D}$ is flat), we
have a commutative diagram in $\rcomod{\coring{D}}$ with exact row
\[
\xymatrix{ 0 \ar[r] & M \cotensor{\coring{C}} (\coring{C} \tensor{A}
B) \ar^{\iota}[r] & M \tensor{A} \coring{C} \tensor{A} B
\ar^-{\omega_{M,\coring{C} \tensor{A} B}}[rr] & & M \tensor{A}
\coring{C} \tensor{A} \coring{C} \tensor{A} B \\
 & M \tensor{A} B \ar^{\psi_M}[u] \ar_-{\rho_M \tensor{A} B}[ur] & &
 &},
\]
where $\psi_M$ is the isomorphism defined by the universal property
of kernel. Now, we will verify that $\psi$ is a natural isomorphism.
For this, let $f:M\to M'$ be a morphism in $\rcomod{\coring{C}}$.
Consider the diagram
$$\xymatrix{
 && M\tensor{A}B \ar[rrr]^{\psi_M}\ar[dll]^{f\tensor{A}B}
 \ar[ddrrr]^(.3){\rho_M\tensor{A}B} &&&
 M\cotensor{\coring{C}}(\coring{C}\tensor{A}B)
 \ar[dll]^{f\cotensor{\coring{C}}(\coring{C}\tensor{A}B)}\ar[dd]
\\
M'\tensor{A}B \ar[rrr]^{\psi_{M'}}
\ar[ddrrr]^{\rho_{M'}\tensor{A}B}&&&
M'\cotensor{\coring{C}}(\coring{C}\tensor{A}B) \ar[dd]
\\
 &&&&&  M\tensor{A}\coring{C}\tensor{A}B
 \ar[dll]^{f\tensor{A}\coring{C}\tensor{A}B}
\\
  &&& M'\tensor{A}\coring{C}\tensor{A}B.}$$
All sides, except perhaps the top one, are commutative, then the top
side is so. The equivalence between (a) and (b) is then obvious from
Propositions \ref{19}, \ref{19a}.
\end{proof}

Let us consider the particular case where $\coring{C}=A$ and
$\coring{D}=B$ are the trivial corings (which are coseparable). Then
Theorem \ref{23} gives functorial Morita's characterization of
Frobenius ring extensions given in \cite[Theorem 5.1]{Morita:1965}:
In the case $\coring{C} =A$, $\coring{D} = B$ we have that (by
Example \ref{quasifinite}), $A\tensor{A}B \simeq B$ is quasi-finite
as a right $B$-comodule if and only if $_AB$ is finitely generated
an projective, and, in this case, $\cohom{B}{B}{-} \simeq
-\tensor{B}\hom{A}{_AB}{A}$.

Theorem \ref{24} generalizes the characterization of Frobenius
extension of coalgebras over fields \cite[Theorem
3.5]{Castano/Gomez/Nastasescu:1999}. It is then reasonable to give
the following definition.

\begin{definition}\label{Frobenius extension}
Let $(\varphi,\rho) :\coring{C}\rightarrow\coring{D}$ be a
homomorphism of corings. It is said to be a \emph{right Frobenius}
morphism of corings if $ -\tensor{A} B : \rcomod{\coring{C}}
\rightarrow \rcomod{\coring{D}}$ is a Frobenius functor.
\end{definition}

\begin{theorem}\label{24}
Suppose that $A$ and $B$ are QF rings. Let
$(\varphi,\rho):\coring{C}\rightarrow\coring{D}$ be a homomorphism
of corings such that the modules $_A\coring{C}$, $\coring{C}_A$ and
$_B\coring{D}$ are projective. Then the following statements are
equivalent
\begin{enumerate}[(a)]
\item $-\otimes_AB:\mathcal{M}^{\coring{C}}\rightarrow
\mathcal{M}^{\coring{D}}$ is a Frobenius functor;\item the
$\coring{C}-\coring{D}$-bicomodule $\coring{C}\otimes_AB$ is
quasi-finite as a right $\coring{D}$-comodule,
$(\coring{C}\otimes_AB)_{\coring{D}}$ is injective and there exists
an isomorphism of
$\coring{D}%
-\coring{C}$-bicomodules
$\cohom{\coring{D}}{\coring{C}\otimes_AB}{\coring{D}}\simeq
B\otimes_A\coring{C}.$
\end{enumerate}
\end{theorem}

\begin{proof}
Obvious from Theorem \ref{23} and Corollary \ref{injectiveinjector}.
\end{proof}

\begin{lemma}
Let $\coring{C}$ be an $A$-coring. If $-\tensor{A}\coring{C}$ is a
left adjoint functor to the forgetful functor, then $_A\coring{C}$
and $\coring{C}_A$ are finitely generated and projective.
\end{lemma}

\begin{proof}
Suppose that the forgetful functor $\rcomod{\coring{C}}\to\rmod{A}$
is a Frobenius functor. Then the functor
$-\otimes_A\coring{C}:\rmod{A}\to\rmod{A}$ is also a Frobenius
functor (since it is a composition of two Frobenius functors) and
$_A\coring{C}$ is finitely generated projective. On the other hand,
since $-\tensor{A}\coring{C}:\rmod{A}\to\rcomod{\coring{C}}$ is a
left adjoint to
$\hom{\coring{C}}{\coring{C}}{-}:\rcomod{\coring{C}}\to\rmod{A}.$
Then $\hom{\coring{C}}{\coring{C}}{-}$ is a Frobenius functor.
Therefore, $\coring{C}$ is finitely generated projective in
$\rcomod{\coring{C}},$ and hence in $\rmod{A}$.
\end{proof}

\begin{lemma}\emph{\cite[Lemma 4.3]{Brzezinski:2002}}\label{B 4.3}
\\ Let $\coring{C}$ be an $A$-coring and let $R$ the opposite
algebra of $\rdual{\coring{C}}$. If $_A\coring{C}$ is finitely
generated and projective, then the categories $\rcomod{\coring{C}}$
and $\rmod{R}$ are isomorphic to each other. Indeed, if $M\in
\rcomod{\coring{C}}$ then $M$ is a right $R$ by $m.r=\sum
m_{(0)}r(m_{(1)})$. This yields an isomorphism of categories
$\rcomod{\coring{C}}\to\rmod{R}$.
\end{lemma}

\begin{lemma}\label{25}
Let $R$ be the opposite algebra of $^*\coring{C}$.
\begin{enumerate}[(1)]
\item $\coring{C}\in{}^{\coring{C}}\mathcal{M}^{A}$ is
quasi-finite (resp. quasi-finite and injector) as a right
$A$-comodule if and only if $_A\coring{C}$ is finitely generated
projective (resp. $_A\coring{C}$ is finitely generated projective
and $_AR$ is flat). Let $\cohom{A}{\coring{C}}{-}
=-\otimes_AR:\mathcal{M}^{A}\rightarrow\mathcal{M}^{\coring{C}}$ be
the cohom functor.\item If $_A\coring{C}$ is finitely generated
projective and $_AR$ is flat, then
\[
_A\cohom{A}{\coring{C}}{A}_\coring{C}\simeq{}_A%
R_{\coring{C}},
\]
where the right $\coring{C}$-comodule structure of $R$ is defined as
in Lemma \ref{B 4.3}
\end{enumerate}
\end{lemma}

\begin{proof}
(1) Straightforward.

(2) From Lemma \ref{B 4.3}, the forgetful functor
$\mathcal{M}^{\coring{C}}\rightarrow\mathcal{M}_A$ is the
composition of functors $\mathcal{M}^{\coring{C}}\rightarrow
\mathcal{M}_{R}\rightarrow\mathcal{M}_A.$ By Proposition
\ref{adjointcohomcotens}, $\cohom{A}{\coring{C}}{-}$ is a left
adjoint to
$-\square_{\coring{C}}\coring{C}:\mathcal{M}^{\coring{C}}\rightarrow
\mathcal{M}_A$ which is isomorphic to the forgetful functor
$\rcomod{\coring{C}}\rightarrow\mathcal{M}_A.$ Then
$\cohom{A}{\coring{C}}{-}$ is isomorphic to the composition of
functors
\[
\xymatrix@1{\mathcal{M}_A\ar[r]^{-\otimes
_AR}&\mathcal{M}_{R}\ar[r]& \mathcal{M}^{\coring{C}}.}
\]
In particular, $_A\cohom{A}{\coring{C}}{A}_{\coring{C}}\simeq {}_A(
A\otimes_AR)_{\coring{C}}\simeq{}_AR_{\coring{C}}.$
\end{proof}

\begin{corollary}\label{26}\emph{\cite[27.10]{Brzezinski/Wisbauer:2003}}
\newline Let $\coring{C}$ be an $A$-coring and let
$R$ be the opposite algebra of $^*\coring{C}.$ Then the following
statements are equivalent
\begin{enumerate}[(a)]
\item The forgetful functor $F:\mathcal{M}^{\coring{C}%
}\rightarrow\mathcal{M}_A$ is a Frobenius functor;
\item
$_A\coring{C}$ is finitely generated projective and
$\coring{C}\simeq R$ as $(A,R)$-bimodules, such that $\coring{C}$ is
a right $R$-module by $c.r=c_{(1)}.r(c_{(2)})$, for every
$c\in\coring{C}$ and $r\in R.$
\end{enumerate}
\end{corollary}

\begin{proof}
Straightforward from Theorem \ref{23} and Lemma \ref{25}.
\end{proof}

The following proposition gives sufficient conditions to have that a
morphism of corings is right Frobenius if and only if it is left
Frobenius. Note that it says in particular that the notion of
Frobenius homomorphism of coalgebras over fields (by (b)) or of
rings (by (d)) is independent on the side. Of course, the latter is
well known.

\begin{proposition}\label{27}
Let $(\varphi,\rho):\coring{C}\rightarrow\coring{D}$ be a
homomorphism of corings such that $_A\coring{C}$, $_B\coring{D}$,
$\coring{C}_A$ and $\coring{D}_B$ are flat. Assume that at least one
of the following holds
\begin{enumerate}[(a)] \item $\coring{C}$ and $\coring{D}$ have a
duality, and $_AB$ and $B_A$ are flat;\item $A$ and $B$ are von
Neumann regular
rings;\item $B\otimes_A\coring{C}$ is coflat in $^{\coring{D}%
}\mathcal{M}$ and $\coring{C}\otimes_AB$ is coflat in $\mathcal{M}%
^{\coring{D}}$ and $_AB $ and $B_A$ are flat;\item $\coring{C}$ and
$\coring{D}$ are coseparable corings.
\end{enumerate} Then the following statements are equivalent
\begin{enumerate}
\item $-\otimes_AB:\mathcal{M}^{\coring{C}%
}\rightarrow\mathcal{M}^{\coring{D}}$ is a Frobenius functor;\item
$B\otimes_A-:{}^{\coring{C}}\mathcal{M}\rightarrow{}^{\coring{D}%
}\mathcal{M}$ is a Frobenius functor.
\end{enumerate}
\end{proposition}

\begin{proof}
Obvious from Theorem \ref{leftrightFrobenius} and Theorem
\ref{leftrightFrobeniusduality}.
\end{proof}

Let us finally show how to derive from our results a remarkable
characterization of the so called \emph{Frobenius corings}.

\begin{corollary}\label{28}\emph{\cite[27.8]{Brzezinski/Wisbauer:2003}}
\newline The following statements are equivalent
\begin{enumerate}[(a)] \item the forgetful
functor $\mathcal{M}^{\coring{C}}\rightarrow\mathcal{M}_A$ is a
Frobenius functor;\item the forgetful functor
$^{\coring{C}}\mathcal{M}\rightarrow{}_A\mathcal{M}$ is a Frobenius
functor;
\item there exist an $(A,A) $-bimodule map
$\eta:A\rightarrow \coring{C}$ and a $(\coring{C},\coring{C})
$-bicomodule map
$\pi:\coring{C}\otimes_A\coring{C}\rightarrow\coring{C}$ such that
$\pi(\coring{C}\otimes_A\eta)
=\coring{C}=\pi(\eta\otimes_A\coring{C}) .$
\end{enumerate}
\end{corollary}

\begin{proof}
The proof of ``$(1)\Leftrightarrow(4)$'' in Proposition \ref{19a}
for $X=\coring{C}\in{}^{A}\mathcal{M}^{\coring{C}}$ and
$\Lambda=\coring{C}\in{}^{\coring{C}}\mathcal{M}^{A}$ remains valid
for our situation. Finally, notice that the condition (4) in this
case is exactly the condition (c).
\end{proof}

In such a case we say that $\coring{C}$ is a \emph{Frobenius
coring}.

\section{Applications to entwined modules}\label{entwined}

In this section we particularize some of our results in Section
\ref{Frobmor} to the category of entwined modules. We start with
some useful results.

\begin{proposition}\label{Takeuchi-Two-sided}
Consider a right-right entwining structure
$(A,C,\psi)\in\mathbb{E}_\bullet^\bullet(k)$
 and a left-left entwining structure $(B,D,\varphi)\in{}_\bullet
 ^\bullet\mathbb{E}(k)$. The category of two-sided entwined modules
$_B^{D}\mathcal{M}(\varphi,\psi)_A^{C}$ (defined in Subsection
\ref{CategoryEntwined}) is isomorphic to the category of bicomodules
$\bcomod{D\otimes B}{A\otimes C}$ over the associated corings.
\end{proposition}

\begin{proof}
Let $\coring{C}=A\otimes C$ and $\coring{D}=D\otimes B$ the
associated $A$-coring and $B$-coring respectively (see Theorems
\ref{Takeuchi} and \ref{Takeuchi'}). Let $M\in{}_B^{D}\mathcal{M}
(\varphi,\psi)_A^{C}$ with $\rho_M:M\to M\otimes C$ and
$\lambda_M:M\to D\otimes M$ as coaction maps and
$\psi^r_M:M\otimes A\to M$ and $\psi^l_M:B\otimes M\to M$ as action
maps. Let $(M,\rho_M')\in\rcomod{\coring{C}}$ and
$(M,\lambda_M')\in\lcomod{\coring{D}}$ the associated right
$\coring{C}$-comodule and left $\coring{D}$-comodule respectively
(see Theorems \ref{Takeuchi} and \ref{Takeuchi'}).

First, we have, $\rho_M'$ is left $B$-linear if and only if
$\psi^l_M$ is right $C$-colinear, and $\lambda_M'$ is right
$A$-linear if and only if $\psi^r_M$ is left $D$-colinear.

Now, consider the diagram
$$\xymatrix{\ar@{}[drr]|{\textbf{(1)}}
M\ar[rr]^-{\rho_M}\ar[d]_{\lambda_M} &&
\ar@{}[dr]|{\textbf{(2)}}M\otimes C\ar[r]^-\simeq
\ar[d]^{\lambda_M\otimes C} & M\tensor{A}(A\otimes C)
\ar[d]^{\lambda_M\tensor{A}(A\otimes C)}
\\ \ar@{}[drr]|{\textbf{(3)}}
D\otimes M \ar[rr]^-{D\otimes\rho_M}\ar[d]_\simeq
&& \ar@{}[dr]|{\textbf{(4)}}D\otimes M\otimes C \ar[r]^-\simeq
\ar[d]^\simeq & (D\otimes M)\tensor{A}(A\otimes C)\ar[d]^\simeq
\\ (D\otimes B)\tensor{B}M\ar[rr]^-{(D\otimes B)\tensor{B}\rho_M} &&
(D\otimes B)\tensor{B}(M\otimes C)\ar[r]^-\simeq & (D\otimes
B)\tensor{B}M\tensor{A}(A\otimes C).}$$ The diagrams (2), (3) and
(4) are commutative. Hence, the greatest rectangle is commutative if
and only if the rectangle (1) is commutative. Finally, from Theorems
\ref{Takeuchi} and \ref{Takeuchi'},
$_B^{D}\mathcal{M}(\varphi,\psi)_A^{C}\simeq
\bcomod{\coring{D}}{\coring{C}}.$
\end{proof}

\begin{proposition}\emph{\cite[Proposition
34]{Caenepeel/Militaru/Zhu:2002}}\label{CMZ 34}
\\If $(A,C,\psi)$ belongs to
$\mathbb{E}_\bullet^{\bullet}(k)$ and is such that $\psi$ is an
isomorphism, then $\psi$ is an isomorphism of corings.
\end{proposition}

\begin{lemma}\label{E0}
Let $\varphi:\coring{C}\to\coring{C}'$ be an isomorphism of
$A$-corings and $\coring{D}$ a $B$-coring. Then there is a canonical
isomorphism of categories $$F:\bcomod{\coring{C}}{\coring{D}}\to
\bcomod{\coring{C}'}{\coring{D}}.$$ Let $M\in
\bcomod{\coring{C}}{\coring{D}}$. Set $F(M)=M$ as
$(A,\coring{D})$-bicomodules. The left $\coring{C}'$-coaction of
$F(M)$ is
$$\xymatrix{M\ar[r]^-{\lambda_M} & \coring{C}\tensor{A}M
\ar[rr]^-{\varphi\tensor{A}M} && \coring{C}'\tensor{A}M.}$$
\end{lemma}

\begin{proof}
Straightforward.
\end{proof}

\begin{corollary}
If $(A,C,\psi)$ and $(A',C',\psi')$ belong to
$\mathbb{E}_\bullet^{\bullet}(k)$ and are such that $\psi$ is an
isomorphism. Then
\begin{enumerate}[(a)]
\item We have
\[
^{A\otimes C}\mathcal{M}^{A'\otimes
C'}\simeq{}_A^{C}\mathcal{M}(\psi^{-1},\psi')_{A'}^{C'}.
\]
\item If the coalgebra $C$ is flat as a $k$-module, then the modules
$_A(A\otimes C)$ and $(A\otimes C)_A$ are flat.
\end{enumerate}
\end{corollary}

\begin{proof}
Follows immediately from Propositions \ref{CMZ 34},
\ref{Takeuchi-Two-sided}, and Lemma \ref{E0}.
\end{proof}

Now, let $(\alpha,\gamma):(A,C,\psi)\rightarrow(A',C',\psi')$ be a
morphism in $\mathbb{E}_\bullet^{\bullet }(k)$. We know that
$(\alpha\otimes\gamma,\alpha):A\otimes{C}\rightarrow{A'}\otimes{C'}$
is a morphism of corings. From \cite[Lemma
8]{Caenepeel/Militaru/Zhu:2002}, we have a functor
$F:\mathcal{M}(\psi)_A^{C}\to \mathcal{M}(\psi')_{A'}^{C'}$ defined
as follows: Let $M\in \mathcal{M}(\psi)_A^{C}$.
$F(M)=M\tensor{A}A'\in\rmod{A'}$. The right $C'$-coaction is defined
by
\begin{equation}
\rho_{M\tensor{A}A'}(m\otimes a')=(m_{(0)}\otimes a'_{\psi'})
\otimes \gamma(m_{(1)})^{\psi'},
\end{equation}
$m\in M,a'\in A'$. For every morphism $f:M\to N$ in
$\mathcal{M}(\psi)_A^{C}$, $F(f)=f\otimes A'$.

\begin{proposition}
 We have the commutative diagram
\[
\xymatrix{\mathcal{M}^{A\otimes
C}\ar[rr]^{-\otimes_AA'}\ar[d]^\simeq & &
\mathcal{M}^{A'\otimes C'}\ar[d]^\simeq \\
\mathcal{M}(\psi)_A^{C}\ar[rr]_{F=-\otimes_AA'} & &
\mathcal{M}(\psi')_{A'}^{C'},}
\]
where $-\otimes_AA':\mathcal{M}^{A\otimes C}\rightarrow
\mathcal{M}^{A'\otimes C'}$ is the induction functor.
\end{proposition}

\begin{proof}
Let $M\in \mathcal{M}(\psi)_A^{C}$. Then
\begin{eqnarray*}
 \rho_{M\tensor{A}A'}'(m\tensor{A}a') &=&
 m_{(0)}\tensor{A} 1_{A'}\tensor{A'}
 (\alpha\otimes\gamma) (1_A\otimes m_{(1)})a'\\
 &=& m_{(0)}\tensor{A} 1_{A'}\tensor{A'}
 \big(1_{A'}\otimes \gamma(m_{(1)})\big)a'\\
 &=& m_{(0)}\tensor{A} 1_{A'}\tensor{A'}
 \big(a'_{\psi'}\otimes \gamma(m_{(1)})^{\psi'}\big)\\
 &=& m_{(0)}\tensor{A} 1_{A'}\tensor{A'}
 a'_{\psi'}\big(1_{A'}\otimes \gamma(m_{(1)})^{\psi'}\big)\\
 &=&m_{(0)}\tensor{A} a'_{\psi'}\tensor{A'}
 \big(1_{A'}\otimes \gamma(m_{(1)})^{\psi'}\big).
\end{eqnarray*}

\end{proof}

We obtain the following result concerning the category of entwined
modules.

\begin{theorem}\label{entwined modules}
Let $(\alpha,\gamma):(A,C,\psi)\rightarrow(A',C',\psi')$ be a
morphism in $\mathbb{E}_\bullet^{\bullet }(k)$, such that $_{k}C$
and $_{k}D$ are flat.
\begin{enumerate} \item The following statements are
equivalent\begin{enumerate}[(a)] \item The functor
$F=-\otimes_AA':\mathcal{M}(\psi)_A^C\rightarrow
\mathcal{M}(\psi')_{A'}^{C'}$ defined above is a Frobenius functor;
\item the $ A\otimes{C}-A'\otimes C'$-bicomodule
$(A\otimes C)\otimes_AA'$ is quasi-finite injector as a right
$A'\otimes C'$-comodule and there exists an isomorphism of
$A'\otimes C'-A\otimes C$-bicomodules $\cohom{A'\otimes
C'}{(A\otimes{C})\otimes_AA'}{A'\otimes C'}\simeq
A'\otimes_A(A\otimes{C}).$
\end{enumerate}
Moreover, if $A\otimes C$ and $A'\otimes C'$ are coseparable
corings, then the condition ``injector'' in (b) can be deleted.

 \item If $A$ and $A'$ are QF rings and the module
$(A'\otimes C')_{A'}$ is projective, then the following are
equivalent\begin{enumerate}[(a)] \item The functor
$F=-\otimes_AA':\mathcal{M}(\psi)_A^{C}\rightarrow\mathcal{M}
(\psi')_{A'}^{C'}$ defined above is a Frobenius functor;
\item the $ A\otimes C-A'\otimes C'$-bicomodule
$(A\otimes{C})\otimes_AA'$ is quasi-finite and injective as a right
$A'\otimes C'$-comodule and there exists an isomorphism of
$A'\otimes C'-A\otimes C$-bicomodules $\cohom{A'\otimes
C'}{(A\otimes{C})\otimes_AA'}{A'\otimes C'}\simeq
A'\otimes_A(A\otimes C).$
\end{enumerate}
\end{enumerate}
\end{theorem}

\begin{proof}
Follows from Theorem \ref{23} and Theorem \ref{24}.
\end{proof}

\begin{remark}
Let a right-right entwining structure
$(A,C,\psi)\in\mathbb{E}_\bullet^\bullet(k)$. The coseparability of
the coring ${A}\tensor{}{C}$ is characterized in \cite[Theorem 38
(1)]{Caenepeel/Militaru/Zhu:2002} (see also Proposition
\ref{cosepcoring}).
\end{remark}

\section{Applications to graded ring theory}\label{graded}

In this section we apply our results in the previous sections to the
category of graded modules by a $G$-set. Let $G$ be a group, $A$ a
$G$-graded $k$-algebra, and $X$ a left $G$-set. We begin by giving
some useful lemmas.

\subsection{Some useful lemmas}

Let $C=kX$ and $C'=kX'$ be two grouplike coalgebras, where $X$ and
$X'$ are arbitrary nonempty sets. We know (see Example \ref{graded
modules}) that the category $\rcomod{C}$ is isomorphic to the
category of $X$-graded modules. Moreover we have the following:

\begin{lemma}\label{graded-A1}
For a $k$-module $M$ which is both a $X$-graded and a $X'$-graded
module, the following are equivalent \begin{enumerate}[(a)]
\item $M\in{}^{C'%
}\mathcal{M}^{C}$; \item for every $m\in M,x\in X,x'\in X', \quad
{}_{x'}(m_{x}) \in M_{x}$; \item for every $m\in M,x\in X,x'\in
{}X', \quad (_{x'}m)_{x}\in{}_{x'}M$;
\item for every $m\in M,x\in X,x'\in {}X', \quad _{x'}(m_{x})
=(_{x'}m)_{x}$.
\end{enumerate}
\end{lemma}

\begin{proof}
Let $M=\bigoplus_{x\in X}M_{x}=\bigoplus_{x'\in X'}{}_{x'}M$. At
first observe that the condition (a) is equivalent to the fact that
the diagram
\[
 \xymatrix{M \ar[rr]^{\rho_M} \ar[d]_{\lambda_M}& & M\otimes kX
\ar[d]^{\lambda_M \otimes kX}
\\kX'\otimes M \ar[rr]_-{kX'\otimes \rho_M} & & kX'\otimes M\otimes
kX\simeq M^{(X'\times X)}}
\]
is commutative.

$(a)\Rightarrow(d)$ Let $m\in M$. We have, $\rho_M(m)=\sum_{x\in
X}m_x \otimes x$ and $\lambda_M(m)=\sum_{x'\in X'}x'\otimes
{_{x'}m}$. From the commutativity of the above diagram,
$$\sum_{x\in X}\sum_{x'\in X'}x'\otimes {_{x'}(m_{x})}\otimes
x=\sum_{x'\in X'}\sum_{x\in X}x'\otimes(_{x'}m) _{x}\otimes x.$$
Hence, $_{x'}(m_{x}) =(_{x'}m) _{x}$ for all $x\in X$, and $x'\in
{}X'$.

$(d)\Rightarrow(b)$ Trivial.

$(b)\Rightarrow(a)$  Let $m\in M$. We have, $$(\lambda_M\otimes
kX)\rho_M(m)=\sum_{x\in X}\sum_{x'\in X'}x'\otimes
{_{x'}(m_{x})}\otimes x,$$ and $$(kX'\otimes
\rho_M)\lambda_M(m)=(kX'\otimes \rho_M)\big(\sum_{x\in X}\sum_{x'\in
X'}x'\otimes {_{x'}(m_{x})}\big)=\sum_{x\in X}\sum_{x'\in
X'}x'\otimes {_{x'}(m_{x})}\otimes x.$$ Hence (a) follows.

$(a)\Leftrightarrow(d)\Leftrightarrow(c)$ Follows by symmetry.
\end{proof}

 Now let $G$ and $G'$ be two groups, $A$ a $G$-graded
$k$-algebra, $A'$ a $G'$-graded $k$-algebra, $X$ a right $G$-set,
and $X'$ a left $G'$-set. Let $kG$ and $kG'$ be the canonical Hopf
algebras.

 Let $\psi:kX\otimes A\to A\otimes kX,$ be the map
defined by $x\otimes a_{g}\mapsto a_{g}\otimes xg$, and
$\psi':A'\otimes kX'\to kX'\otimes A'$, be the map defined by
$a_{g'}'\otimes x'\mapsto g'x'\otimes a_{g'}'.$

From Subsection \ref{CategoryGraded}, we have $(kG,A,kX)
\in\mathbb{DK}_\bullet^\bullet(k) ,$ $(kG',A',kX')
\in{}_\bullet^\bullet\mathbb{DK}(k) ,$ $(A,kX,\psi)
\in\mathbb{E}_\bullet^\bullet(k) ,$ $(A',kX',\psi')
\in{}_\bullet^\bullet\mathbb{E}( k) ,$ and $\mathcal{M}(kG)_A%
^{kX}\simeq gr-(A,X,G),{}$ $_{A'}^{kX'%
}\mathcal{M}(kG')\simeq(G',X',A')-gr.$

\begin{lemma}\label{graded-A2}
\begin{enumerate}[(1)] \item Let $M$ be a $k$-module having the
structure of an $X$-graded right $A$-module and an $X'$-graded left
$A'$-module. The following are equivalent
\begin{enumerate}[(a)]\item
$M\in{}_{A'}^{kX'}\mathcal{M}(\psi',\psi)_A^{kX}$; \item the
following conditions hold \begin{enumerate}[(i)]\item $M$ is a
$(A',A)$-bimodule,\item for every $m\in M,x\in X,x'\in X',
\quad_{x'}(m_{x})\in M_{x}$ (or $(_{x'}m)_{x}\in{}_{x'}M$, or
$_{x'}(m_{x}) =(_{x'}m)_{x}$),\item for every $x\in X,$ $M_{x}$ is a
submodule of $_{A'}M$, \item for every $x'\in X',$ $_{x'}M$ is a
submodule of $M_A.$
\end{enumerate}
\end{enumerate}
\item $M\in{}_{A'}^{kX'}\mathcal{M}(\psi',\psi)_A^{kX}$ if and
only if $M$ is an $X'\times X$-graded $(A',A)$-bimodule (see
Subsection \ref{CategoryGraded}).
\end{enumerate}
\end{lemma}

\begin{proof}
 (1) $(a)\Leftrightarrow(b)$ We will use the definition of an object
in the category of two-sided entwined modules
${}_{A'}^{kX'}\mathcal{M}(\psi',\psi)_A^{kX}$. By Lemma
\ref{graded-A1}, the condition ``$M\in{}^{kX'}\mathcal{M}^{kX}$'' is
equivalent to the condition (ii). We have moreover that the left
$A'$-action on $M$ is $kX$-colinear if and only if for every $x\in
X,m\in M_x$, $\rho_M(a'm)=(a'm)\otimes x$, if and only if for every
$x\in X,m\in M_x$, $a'm\in M_x$, if and only if (iii) holds. By
symmetry, the condition ``the right $A$-action on $M$ is
$kX'$-colinear'' is equivalent to the condition (iv).

(2) The ``if'' part is clear. For the ``only if'' part, put
$M=\bigoplus_{(x',x)\in {X'\times X}}M_{(x',x)}$, where
$M_{(x',x)}={}_{x'}(M_{x}) =(_{x'}M)_{x}$. Let $a'\in A'_{g'}, m\in
M$. Since $a'.{}_{x'}(m_x)\in M_x$ and $a'.{}_{x'}(m_x)\in
{}_{g'x'}(M)$, $a'.{}_{x'}(m_x)\in {}_{a'x'}(M_{x})$. Therefore,
$A'_{g'}.{}_{x'}(M_x)\subset {}_{a'x'}(M_{x})$. By symmetry we
obtain $A'_{g'}M_{(x',x)}A_g\subset M_{g'(x',x)g}=M_{(g'x',xg)}$
($g\in G,g'\in G',x\in X, x'\in X'$).
\end{proof}

\subsection{Adjoint pairs and Frobenius pairs of functors between
categories of graded modules over $G$-sets}

Throughout this subsection, $G$ and $G'$ will be two groups, $A$ a
$G$-graded $k$-algebra, $A'$ a $G'$-graded $k$-algebra, $X$ a right
$G$-set, and $X'$ a right $G'$-set. Let $kG$ and $kG'$ be the
canonical Hopf algebras. Let $\psi:kX\otimes A\to A\otimes kX$ be
the map defined by $x\otimes a_{g}\mapsto a_{g}\otimes xg$.
Analogously we define the map $\psi':kX'\otimes A'\to A'\otimes
kX'$. The comultiplication and the counit maps of the coring
$A\otimes kX$ are defined by:
$$\Delta(a\otimes x)=(a\otimes x)
\tensor{A}(1_A\otimes x), \quad \epsilon(a\otimes x)=a \quad (a\in
A,x\in X).$$

\begin{proposition}\label{coseparable}
The corings $A\otimes kX$ and $A'\otimes kX'$ are coseparable.
\end{proposition}

\begin{proof}
First, \cite[Proposition 101]{Caenepeel/Militaru/Zhu:2002} states
that the forgetful functor
$$\mathcal{M}(kG)^{kX}_A\to \rmod{A}$$ is separable. By
Proposition \ref{cosepcoring}, $A\otimes kX$ is a coseparable
coring.

The following direct proof in the setting of corings is shown to me
by T. Brzezi\'nski. The map $\delta:(A\otimes kX)\tensor{A}(A\otimes
kX)\to A$ defined by $\delta(a\otimes x\otimes y)=a\delta_{x,y}$
(Kronecker's delta) for $a\in A,x,y\in X$, is a cointegral in the
coring $A\otimes kX$ (see Definition \ref{cointegral}). Indeed,
$\delta$ is clearly left $A$-linear. On the other hand, for all
$a\in A,x,y\in X,g\in G, b_g\in A_g$,
\begin{eqnarray*}
 \delta\big((a\otimes x)\tensor{A}(1_A\otimes y)b_g\big)
 &=& \delta\big((a\otimes x)\tensor{A}\psi(y\otimes b_g)\big)\\
 &=& \delta\big((a\otimes x)\tensor{A}(b_g\otimes yg)\big)\\
 &=& \delta\big(a\psi(x\otimes b_g)\tensor{A}(1_A\otimes yg)\big)\\
 &=& \delta\big((ab_g\otimes xg)\tensor{A}(1_A\otimes yg)\big)\\
 &=& ab_g\delta_{xg,yg}\\
 &=& ab_g\delta_{x,y}\\
 &=&\delta\big((a\otimes x)\tensor{A}(1_A\otimes y)\big)b_g.
\end{eqnarray*}
Then $\delta$ is right $A$-linear. Moreover, it is clear that
$\delta\circ\Delta=\epsilon$. Finally, set $\coring{C}=A\otimes kX$.
Then, we have for all $a\in A,x,y\in X$,
\begin{eqnarray*}
(\coring{C}\tensor{A}\delta)\circ(\Delta\tensor{A}\coring{C})
\big((a\otimes x)\tensor{A}(1_A\otimes y)\big) &=&
(\coring{C}\tensor{A}\delta)\big((a\otimes x)\tensor{A}(1_A\otimes
x)\tensor{A}(1_A\otimes y)\big)\\
 &=& (a\otimes x)\delta_{x,y},
\end{eqnarray*}
and
\begin{eqnarray*}
(\delta\tensor{A}\coring{C})\circ(\coring{C}\tensor{A}\Delta)
\big((a\otimes x)\tensor{A}(1_A\otimes y)\big) &=&
(\delta\tensor{A}\coring{C})\big((a\otimes x)\tensor{A}(1_A\otimes
y)\tensor{A}(1_A\otimes y)\big)\\
 &=& a\delta_{x,y}(1_A\otimes y)\\
 &=& a\delta_{x,y}\otimes y.
\end{eqnarray*}
Hence,
$(\coring{C}\tensor{A}\delta)\circ(\Delta\tensor{A}\coring{C})=
(\delta\tensor{A}\coring{C})\circ(\coring{C}\tensor{A}\Delta).$ By
Proposition \ref{cosepcoring}, $A\otimes kX$ is a coseparable
coring.
\end{proof}

\begin{lemma}\label{graded-A3}
\begin{enumerate}[(1)]\item $\psi$ is bijective, $(A,kX,\psi^{-1})
\in{}_\bullet^\bullet\mathbb{E}(k) $, and $$_A%
^{kX}\mathcal{M}(kG):={}_A^{kX}\mathcal{M}(\psi^{-1})
\simeq(G,X,A)-gr,$$ where the structure of left $G$-set on $X$ is
given by $g.x=xg^{-1}$ $(g\in G, x\in X)$.
\item Every object of the category $^{A'\otimes
kX'}\mathcal{M}^{A\otimes kX}\simeq{}_{A'}^{kX'}\mathcal{M}((\psi')
^{-1},\psi)_A^{kX}$ can be identified with an $X'\times X$-graded
$(A',A)$-bimodule.
\end{enumerate}
\end{lemma}

\begin{proof}
(1) From Proposition \ref{antipode}, and since the antipode of
$H=kG$ is $S(g)=g^{-1}$, for all $g\in G$, then $S\circ S=1_{H}$ and
$\overline{S}=S^{-1}=S$ is a twisted antipode of $H.$ Hence $\psi$
is bijective and $(A,kX,\psi^{-1})\in{}_{\bullet
}^\bullet\mathbb{E}(k)  $ (see Subsection \ref{CategoryEntwined}),
and $_A^{kX}\mathcal{M}( kG)  :={}_A^{kX}\mathcal{M}(\psi^{-1})
\simeq(G,X,A)-gr$, since $\psi^{-1}:A\otimes kX\rightarrow kX\otimes
A,$ $a_{g}\otimes x\mapsto g.x\otimes a_{g}$, where $g.x=xg^{-1}.$

(2) It follows from (1) and Lemma \ref{graded-A2}.
\end{proof}

\begin{lemma}\label{graded-A4}
\begin{enumerate}[(1)]
\item Let $M\in gr-(A,X,G)$, and $N\in (G,X,A)-gr$. We know that $M\in
\mathcal{M}^{A\otimes kX}$ by the coaction: $ \rho_M:M\rightarrow
M\tensor{A}{(A\otimes kX)}$, where $\rho_M(m_x)=m_x\tensor{A}
{(1_A\otimes x)}$, and $N\in {}^{A\otimes kX}\mathcal{M}$ by the
coaction: $\lambda_N:N\rightarrow (A\otimes kX)\tensor{A}{N}$, where
$\lambda_N({}_x n)=(1_A\otimes x) \tensor{A}{{}_x n}$. We have
$M\cotensor{(A\otimes kX)}{N}=M\widehat{\otimes}_AN$, where
$M\widehat{\otimes}_AN$ is the additive subgroup of $M\otimes_A N$
generated by the elements  $m\otimes_A n$ where $x\in X,m\in M_x$,
and $n\in {}_xN$ (see Subsection \ref{CategoryGraded}).

\item Let $P$ be an $X\times X'$-graded $(A,A')$-bimodule. We have the
commutative diagram:
\[%
\xymatrix{gr-(A,X,G)\ar[d]^\simeq \ar[rr]^{-\widehat{\otimes}_AP} &
&
gr-( A',X',G') \ar[d]^\simeq  \\
\mathcal{M}^{A\otimes kX} \ar[rr]^{-\cotensor{(A\otimes kX)}{P}} & &
\mathcal{M}^{A'\otimes kX'},}
\]
where $-\widehat{\otimes}_AP: gr-(A,X,G)\rightarrow gr-(A',X',G')$
is the functor defined in Subsection \ref{CategoryGraded}.
\end{enumerate}
\end{lemma}

\begin{proof}
(1) At first we will show that the right $A$-module $A\otimes kX$ is
free. More precisely that the family $\lbrace1_A\otimes x\mid x\in
X\rbrace$ is a right basis of it. It is clear that $a_g\otimes
x=(1_A\otimes xg^{-1})a_g$ for $a_g\in A_g$ and $x\in X$. Now
suppose that $\sum_i (1_A\otimes x_i)a_i=0$, where $x_i\in X, a_i\in
A$. Then  $\sum_i \psi(x_i\otimes a_i)=0$. Since $\psi$ is
bijective, we obtain $\sum_i x_i\otimes a_i=0$. Therefore $a_i=0$
for all $i$. Hence the above mentioned family is a basis of the
right $A$-module $A\otimes kX$. (There is a shorter and indirect
proof of this fact by using that $\psi$ is an isomorphism of
$A$-bimodules, and $\{x\otimes 1_A|x\in X\}$ is a basis of the right
$A$-module $kX\otimes A$.)

Finally, suppose that  $\sum_i m_{x_i}\otimes_A {}_{y_i}n\in
M\cotensor{(A\otimes kX)}{N}$. Then
$$ \sum_i m_{x_i}\otimes_A 1_A\otimes x_i\otimes_A {}_{y_i}n=
\sum_i m_{x_i}\otimes_A 1_A\otimes y_i\otimes_A {}_{y_i}n.$$ Hence
(the right $A$-module $A\otimes kX$ is free), $x_i=y_i$ for all $i$,
and $M\cotensor{(A\otimes kX)}{N}\subset M \widehat{\otimes}_AN$.
The other inclusion is obvious.

(2) It suffices to show that the map $M\widehat{\otimes}_AP
\rightarrow M\widehat{\otimes}_AP \tensor{A'}{(A'\otimes kX')}$
defined by $$\sum_{x\in F}m_x \otimes_A {}_x p\longmapsto
\sum_{x'\in X'}\bigg(\sum_{x\in F}m_x \otimes_A{}_x p\bigg)_{x'}
\tensor{A'}{(1_{A'}\otimes x')},$$ where $F$ is a finite subset of
$X$, makes commutative the following diagram
\[
 \xymatrix{0\ar[r] & M\widehat{\otimes}_AP \ar[rr]^{i} \ar[d]
 & & M\tensor{A}{P} \ar[d]
\\ 0\ar[r] & M\widehat{\otimes}_AP\tensor{A'}{(A'\otimes kX')}
\ar[rr]^{i\tensor{A'}{(A'\otimes kX')}} & &
M\tensor{A}{P}\tensor{A'}{(A'\otimes kX')}.}
\]
That is clear since $(\sum_{x\in F}m_x\otimes_A{}_x
p)_{x'}=\sum_{x\in F}m_x\otimes_A{}(_x p)_{x'}=\sum_{x\in F}m_x
\otimes_A{}p_{(x,x')}$, where $p=\sum_{x\in F}{}_x p$.
\end{proof}

Let $\widehat{A}=A\otimes kX$ be the $X\times X$-graded
$(A,A)$-bimodule associated to the $(A\otimes kX)-(A\otimes
kX)$-bicomodule $A\otimes kX$. It is clear that
$-\widehat{\otimes}_A\widehat{A}= -\cotensor{(A\otimes
kX)}{(A\otimes kX)} \simeq 1_{gr-(A,G,X)}$ and then $\widehat{A}$ is
isomorphic as a bigraded bimodule to Del R\'io's ``$\widehat{A}$''
(see Subsection \ref{CategoryGraded}). The gradings are
$\widehat{A}_x=A\otimes kx$, and ${}_x\widehat{A}=\{\sum_ia_i\otimes
x_i\mid x_ig^{-1}=x, \forall i,\forall g\in G:(a_i)_g\neq 0\}$
($x\in X$).

\begin{proposition}\emph{\cite[Proposition 1.3, Corollary 1.4]{Menini:1993}}
\label{graded-A5}
\begin{enumerate}[(1)]\item The following statements are equivalent for
a $k$-linear functor $F:gr-(A,X,G) \rightarrow gr-(A',X',G')$.
\begin{enumerate}[(a)]
\item $F$ has a right adjoint;
\item $F$ is right exact and preserves coproducts;
\item $F\simeq -\widehat{\otimes}_AP$ for some $X\times
X'$-graded $(A,A')$-bimodule $P$.
\end{enumerate}
\item A $k$-linear functor $G:gr-(A',X',G')\rightarrow gr-(A,X,G)$ has
a left adjoint if and only if $G\simeq \h{P_{A'}}{-}$ for some
$X\times X'$-graded $(A,A')$-bimodule $P$.
\end{enumerate}
\end{proposition}

\begin{proof}
(1) $(a)\Rightarrow (b)$ Clear. $(b)\Rightarrow (c)$ It follows from
Theorem \ref{3} and Lemma \ref{graded-A4}. $(c)\Rightarrow (a)$
Obvious from the above mentioned result.

(2) It follows from (1) and Proposition \ref{Menini 1.2}.
\end{proof}

\begin{lemma}\label{graded-A6}
Let $P$ be an $X\times X'$-graded $(A,A')$-bimodule.
\begin{enumerate}[(1)]\item $\h{P_{A'}}{-}$ is right exact and
preserves direct limits if and only if $_{x}P$ is finitely generated
projective in $\mathcal{M}_{A'}$ for every $x\in X$. \item Suppose
that $_{x}P$ is finitely generated projective in $\mathcal{M}_{A'}$
for every $x\in X$. \begin{enumerate}[(a)] \item For every
$k$-algebra $T$,
$$\xymatrix{\Upsilon_{Z,M}:Z\tensor{T}{\h{P_{A'}}{M}} \ar[r]^-\simeq
& \h{P_{A'}}{Z\tensor{T}{M}}}$$ defined by
$$\Upsilon_{Z,M}(z\tensor{T}{f})(p)=\sum_{x\in
X}z\tensor{T}{f_x(_xp)}$$ $(z\in Z,f=\sum_{x\in X}f_x\in
\h{P_{A'}}{M},p=\sum_{x\in X}{}_xp\in P)$, is the natural
isomorphism associated to the functor (see Section
\ref{Frobeniusgeneral})
$$\h{P_{A'}}{-}:gr-(A',X',G')\rightarrow gr-(A,X,G).$$
\item Moreover we have the natural isomorphism
$$\xymatrix{\eta_N:N \widehat{\otimes}_{A'}\h{P_{A'}}{\widehat{A'}}
\ar[r]^-\simeq & \h{P_{A'}}{N}}$$ defined by
$$\eta_N(n_{x'}\tensor{A'}{}_{x'}f)(p)=\sum_{x\in
X}n_{x'}\delta\big((1_{A'}\otimes
x')\tensor{A'}f_{(x',x)}(_xp)\big)$$ $(n_{x'}\in N_{x'},{}_{x'}f\in
{}_{x'}\h{P_{A'}}{\widehat{A'}},p=\sum_{x\in X}{}_xp\in P)$, where
$\delta$ is the cointegral in the coring $A'\otimes kX'$ defined in
the beginning of this subsection. The left grading on
$\h{P_{A'}}{\widehat{A'}}$ is given by
$$_{x'}\h{P_{A'}}{\widehat{A'}}=\bigg\{f\in
\h{P_{A'}}{\widehat{A'}}\mid \Delta_{(A'\otimes
kX')}\big(f(p)\big)=(1_{A'}\otimes x')\tensor{A'}{\sum_{x\in
X}f_x(_{}xp)},\, \forall p\in P \bigg\}$$ $(x'\in X')$.
\end{enumerate}
\end{enumerate}
\end{lemma}

\begin{proof}
(1) We have
\[
\xymatrix{\h{P_{A'}}{-} \simeq\bigoplus_{x\in
X}\hom{gr-(A',X',G')}{_xP}{-} :gr-(A',X',G') \ar[r] & \mathbf{Ab.}}
\]
Hence, $\h{P_{A'}}{-}$ is right exact and preserves direct limits if
and only if $\hom{gr-(A',X',G')}{_xP}{-}$ is right exact and
preserves direct limits for every $x\in X$ if and only if $_xP$ is
finitely generated projective in $\mathcal{M}_{A'}$ for every $x\in
X$ (by Proposition \ref{V.3.4} and Lemma \ref{fgfpnoeth} (3)).

(2) (a) Let $T$ be a $k$-algebra. Let $M$ be an $X_0\times
X'$-graded $T-A'$-bimodule, where $X_0$ is a singleton, $Z$ be a
right $T$-module, and $x\in X$. We have a sequence of $T$-submodules
$$\h{P_{A'}}{M}_x\leq \h{P_{A'}}{M}\leq \hom{gr-(A',X',G')}{P}{M}\leq
{}_T\hom{A'}{P_{A'}}{M},$$ and the induced structure of left
$T$-module on $\h{P_{A'}}{M}$ is the same structure (see
Section~\ref{Frobeniusgeneral}) of left $T$-module associated to the
functor $\h{P_{A'}}{-}$ on it. Moreover, we have the isomorphism of
left $T$-modules $\h{P_{A'}}{M}_x\simeq
\hom{gr-(A',X',G')}{_xP}{M}$. From Lemma \ref{hom-tens-relations} (4),
for every $x\in X$, there is an isomorphism
\begin{equation} \label{eta x}
\xymatrix{\eta_x:Z\tensor{T}{\hom{A'}{_xP}{M}} \ar[r]^-\simeq &
\hom{A'}{_xP}{Z\tensor{T}{M}}}
\end{equation}
defined by $\eta_x(z\tensor{T}{\gamma_x}):{}_xp\mapsto
z\tensor{T}\gamma_x(_xp)$ $(z\in Z,\gamma_x\in \hom{A'}{_xP}{M})$.
For every $x\in X$, $\eta_x$ induces an isomorphism
\begin{equation} \label{eta'x}
\xymatrix{\eta'_x:Z\tensor{T}{\hom{gr-(A',X',G')}{_xP}{M}}
\ar[r]^-\simeq & \hom{gr-(A',X',G')}{_xP}{Z\tensor{T}{M}}}.
\end{equation}
Now let us consider the isomorphism
$$\xymatrix{\eta'_{Z,M}:=\bigoplus_{x\in
X}\eta'_x:Z\tensor{T}{\h{P_{A'}}{M}} \ar[r]^-\simeq &
\h{P_{A'}}{Z\tensor{T}{M}}}.$$ We have
$\eta'_{Z,M}(z\tensor{T}{f})(p)=\sum_{x\in X}z\tensor{T}{f_x(_xp)}$
$(z\in Z,f=\sum_{x\in X}f_x\in \h{P_{A'}}{M},p=\sum_{x\in X}{}_xp\in
P)$. We can verify easily that $\eta'_{Z,M}$ is a morphism of right
$A$-modules. Hence it is a morphism in $gr-(A,X,G)$. It is clear
that $\eta'_{Z,M}$ is natural in $Z$, and $\eta'_{T,M}$ makes
commutative the diagram \eqref{T} (see Section
\ref{Frobeniusgeneral}). Finally, by Mitchell's Theorem
\ref{Mitchell-Th}, $\eta'_{Z,M}=\Upsilon_{Z,M}.$

(b) To prove the first statement it suffices to use the proof of
Theorem \ref{3}, the property (a), and the fact that every comodule
over an $A$-coseparable coring is $A$-relative injective comodule
(see Proposition \ref{cosepcoring}).

Finally, we will prove the last statement. We know that
$\lambda_{\h{P_{A'}}{\widehat{A'}}}$ is defined to be the unique
$A'$-linear map making commutative the following diagram
\[
\xymatrix{\h{P_{A'}}{\widehat{A'}}
\ar[rr]^{\lambda_{\h{P_{A'}}{\widehat{A'}}}}
\ar[dr]_{\h{P_{A'}}{\Delta_{(A'\otimes kX')}}} & & \widehat{A'}
\tensor{A'} \h{P_{A'}}{\widehat{A'}}
\ar[dl]^{\Upsilon_{\widehat{A'},\widehat{A'}}}
\\  & \h{P_{A'}}{\widehat{A'}\tensor{A'}\widehat{A'}} &
}.
\]
On the other hand,  for each $x'\in X'$,
$_{x'}\h{P_{A'}}{\widehat{A'}}=\{f\in \h{P_{A'}}{\widehat{A'}}\mid
\lambda_{\h{P_{A'}}{\widehat{A'}}}(f)=(1_{A'}\otimes
x')\tensor{A'}f\}$. Since $\Upsilon_{\widehat{A'},\widehat{A'}}$ is
an isomorphism, $_{x'}\h{P_{A'}}{\widehat{A'}}=\{f\in
\h{P_{A'}}{\widehat{A'}}\mid \Delta_{(A'\otimes
kX')}(f(p))=(1_{A'}\otimes x')\tensor{A'}{\sum_{x\in X}f_x({_xp})}
\textrm{ for all }p\in P \}$.
\end{proof}

\begin{remark}\label{graded-A7}
Let $\Sigma$ be a finitely generated projective right $A$-module,
and $M$ be a right $A$-module. It easy to verify that the map
$\alpha_{\Sigma,M}:\hom{A}{\Sigma}{M}\to M\tensor{A}{\Sigma^*}$
defined by
$\alpha_{\Sigma,M}(\varphi)=\sum_i\varphi(e_i)\tensor{A}{e_i^*}$,
where $\{e_i,e_i^*\}_i$ is a dual basis of $\Sigma$, is an
isomorphism, with $\alpha_{\Sigma,M}^{-1}(m\tensor{A}{f})(u)=mf(u)$
$(m\in M,f\in \Sigma^*,u\in \Sigma)$.
\\ Now we assume that the condition of Lemma \ref{graded-A6} (2) holds. It
is easy to verify that for every $x\in X$,
$\eta_x=\alpha^{-1}_{_xP,Z\tensor{T}{M}}\circ(Z\tensor{T}{\alpha_{{}_xP,M}})$,
where $\eta_x$ is the isomorphism \eqref{eta x}. Therefore
$\eta_x^{-1}=(Z\tensor{T}{\alpha^{-1}_{{}_xP,M}})\circ
\alpha_{_xP,Z\tensor{T}{M}}.$
Hence $\Upsilon_{Z,M}^{-1}(g)=\sum_{x\in X}\sum_{i\in
I_x}(Z\tensor{T}{\alpha^{-1}_{{}_xP,M}})(g_x(e_{x,i})\tensor{A'
}{e_{x,i}^*})$, where $\{e_{x,i},e_{x,i}^*\}_{i\in I_x}$ is a dual
basis of $_xP$ $(x\in X)$, and $g=\sum_{x\in X}g_x\in
\h{P_{A'}}{Z\tensor{T}{M}}$. Finally, for each $N\in gr-(A',X',G')$,
$$\xymatrix{\eta^{-1}_N:\h{P_{A'}}{N} \ar[r]^-\simeq & N
\widehat{\otimes}_{A'}\h{P_{A'}}{\widehat{A'}}}$$ is
$\eta^{-1}_N=\Upsilon^{-1}_{N,\widehat{A'}}\circ\h{P_{A'}}{\rho_N}$,
and then $$\eta^{-1}_N(\varphi)=\sum_{x\in X}\sum_{i\in
I_x}\sum_{x'\in
X'}(\varphi_x(e_{x,i}))_{x'}\tensor{A'}\psi'(x'\otimes
e_{x,i}^*(-)),$$ where $\varphi=\sum_{x\in X}\varphi_x\in
\h{P_{A'}}{N}$.

In particular, if $X'=G'=\{e'\}$, then
$$\xymatrix{\eta^{-1}_N:\h{P_{A'}}{N} \ar[r]^-\simeq & N
\otimes_{A'}\h{P_{A'}}{A'}}$$ is defined by
$$\eta^{-1}_N(\varphi)=\sum_{x\in X}\sum_{i\in
I_x}\varphi(e_{x,i})\tensor{A'} e_{x,i}^* \qquad (\varphi\in
\h{P_{A'}}{N}).$$
\end{remark}


\begin{lemma}\label{graded-A8}
\begin{enumerate}[(1)] \item Let $N\in {}_{A'}\mathcal{M}^{A\otimes
kX}.$ Then $N$ is quasi-finite as a right $A\otimes kX$-comodule if
and only if $N_x$ is finitely generated projective in $\lmod{A'}$,
for every $x\in X$. In this case, the cohom functor $\cohom{A\otimes
kX}{N}{-}$ is the composite
\[
  \xymatrix{\mathcal{M}^{A\otimes kX}\ar[r]^-\simeq & gr-(A,X,G)
  \ar[r]^-{-\widehat{\otimes}_A P} & \mathcal{M}_{A'}},
\]
where $P$ is the $X\times X'_0$-graded $A-A'$-bimodule
$\h{_{A'}N}{A'}$ with $X'_0$ is a singleton.

\item Now suppose that $N\in{}^{A'\otimes kX'}\mathcal{M}^{A\otimes
kX}$, and $N_x$ is finitely generated projective in $\lmod{A'}$, for
every $x\in X$. Let $\{e_{x,i},e_{x,i}^*\}_{i\in I_x}$ be a dual
basis of $N_x$ $(x\in X)$. Let $P$ be the bigraded bimodule defined
as above, $\theta:1_{\rcomod{A\otimes kX}}\rightarrow
-\widehat{\otimes}_A P\tensor{A'}N$ be the unit of the adjunction
$(-\widehat{\otimes}_A P,-\tensor{A'}N)$, and let $M\in gr-(A,X,G)$.

We can endow $P$ with a structure of $X\times X'$-graded
$A-A'$-bimodule such that $$\cohom{A\otimes kX}{N}{-}\simeq
-\widehat{\otimes}_A P:gr-(A,X,G)\to gr-(A',X',G').$$ The $A'\otimes
kX'$-comodule structure on $P$, $$\rho_P:P\rightarrow
P\tensor{A'}(A'\otimes kX'),$$ is defined to be the unique
$A'$-linear map satisfying the condition:
\begin{equation}\label{gradedqf}
\sum_{i\in I_x}\rho_P(ae_{x,i}^*)\tensor{A'}e_{x,i}= \sum_{i\in
I_x}\sum_{x'\in X'}ae_{x,i}^*\tensor{A'}(1_{A'}\otimes
x')\tensor{A'}(e_{x,i})_{x'},
\end{equation}
for every $a\in A,x\in X$.
\end{enumerate}
\end{lemma}

\begin{proof}
(1) It follows from Lemma \ref{graded-A4} (2), that the functor
$F:=-\widehat{\otimes}_{A'}N:\rmod{A'}=gr-(A',X_{0}',G')\to
gr-(A,X,G)$, where $X_{0}'$ is a singleton, is the composite
\[
\xymatrix{\rmod{A'}=gr-(A',X_{0}',G') \ar[r]^-{-\otimes_{A'} N} &
\mathcal{M}^{A\otimes kX} \ar[r]^-\simeq & gr-(A,X,G)}.
\]
Then, $N$ is quasi-finite as a right $A\otimes kX$-comodule if and
only if $F$ has a left adjoint, if and only if (by Corollary
\ref{graded-A5} (2)) there exists an $X\times X_{0}'$-graded
$A-A'$-bimodule $P$ such that $F\simeq \h{P_{A'}}{-}$, if and only
if (by Lemma \ref{graded-A6}, Lemma \ref{$T$-objects}, Theorem
\ref{3}) there exists an $X\times X_{0}'$-graded $A-A'$-bimodule $P$
such that $_xP$ is finitely generated projective in
$\mathcal{M}_{A'}$ for every $x\in X,$ and $N\simeq \h{P_{A'}}{A'}$
in $^{A'}\mathcal{M}^{A\otimes kX}.$

Now let us consider the $X\times X'_0$-graded $A-A'$-bimodule
$P:=\h{_{A'}N}{A'}$, and the $X'\times X$-graded $(A',A)$-bimodule
$M:=\h{P_{A'}}{A'}$, where $X'_0$ is a singleton. We have
$$M=\{f\in P^*|f(_xP)=0 \textrm{ for almost all }x\in X \}$$
$$M_x=\big\{f\in P^*|f(_yP)=0 \textrm{ for all }y\in X-\{x\}\big\}\quad
(x\in X),$$ where $P^*=\hom{A'}{P_{A'}}{A'_{A'}}$. The structure of
$A'-A$-bimodule on $P^*$ is given by $(fa)(p)=f(ap),
(a'f)(p)=a'f(p)$ $(f\in P^*, a\in A, a'\in A', p\in P)$. We have
$$M\leq (P^*)_A,\quad  M_x\leq M\leq {}_{A'}(P^*), \quad
\textrm{and} \quad M_x\simeq (_xP)^* \textrm{ in }\lmod{A'}\quad
(x\in X).$$

Analogously,
$$P=\{f\in {}^*N|f(N_x)=0 \textrm{ for almost all }x\in X \}$$
$$_xP=\{f\in {}^*N|f(N_y)=0 \textrm{ for all }y\in X-\{x\}\}\quad (x\in
X),$$ where $^*N=\hom{A'}{_{A'}N}{{}_{A'}A'}$. The structure of
$A-A'$-bimodule on $^*N$ is given by $(af)(n)=f(na),
(fa')(n)=a'f(n)$ $(f\in {}^*N, a\in A, a'\in A', n\in N)$. We have
$$P\leq {}_A(^*N),\quad  _xP\leq P\leq (^*N)_{A'},\quad
\textrm{and} \quad _xP\simeq {}^*(N_x)\textrm{ in } \rmod{A'} \quad
(x\in X).$$ Hence, $N_x$ is finitely generated projective in
$\lmod{A'}$, for every $x\in X$ implies that $_xP$ is finitely
generated projective in $\rmod{A'}$, for every $x\in X$, which
implies that $M_x$ is finitely generated projective in $\lmod{A'}$,
for every $x\in X$.

Finally, let us consider for every $x\in X$ the isomorphism of left
$A'$-modules
$$\xymatrix{H_x:N_x\ar[r]^\simeq & (^*(N_x))^* \ar[r]^\simeq & (_xP)^*
\ar[r]^\simeq & M_x}.$$ We have $H_x(n_x)(\gamma)=\gamma(n_x)$ if
$\gamma\in {}_xP$, and $H_x(n_x)(\gamma)=0$ if $\gamma\in {}_yP$ and
$y\in X-\{x\}$. Set
$$\xymatrix{H=\bigoplus_{x\in X}H_x:N \ar[r]^-\simeq & M}.$$
It can be proved easily that $H$ is a morphism of right $A$-modules.
Hence it is an isomorphism of graded bimodules.

(2) We have $\{e_{x,i}^*,\sigma_x(e_{x,i})\}_{i\in I_x}$ is a dual
basis of $\ldual{(N_x)}\simeq {}_xP$, where $\sigma_x$ is the
evaluation map $(x\in X)$. From the proof of Lemma \ref{graded-A5}
and Remark \ref{graded-A7}, the unit of the adjunction
$(-\widehat{\otimes}_A P,-\tensor{A'}N)$ is $\theta_M:M\rightarrow
M\widehat{\otimes}_A P\tensor{A'}N$, $\theta_M(m)=\sum_{x\in
X}\sum_{i\in I_x}m_x\tensor{A}e_{x,i}^*\tensor{A'}e_{x,i}$ $(m\in
M)$. The coaction on $\cohom{A\otimes kX}{N}{M}=M\widehat{\otimes}_A
P$: $\rho_{M\widehat{\otimes}_A P}:M\widehat{\otimes}_A P\rightarrow
M\widehat{\otimes}_A P\tensor{A'}(A'\otimes kX')$ is the unique
$A'$-linear map satisfying the commutativity of the diagram
\[
\xymatrix{M\ar[rr]^{\theta_M} \ar[d]^{\theta_M} & &
M\widehat{\otimes}_A P\tensor{A'}N
\ar[d]^{\rho_{M\widehat{\otimes}_A P}\tensor{A'}N}
\\ M\widehat{\otimes}_A
P\tensor{A'}N \ar[rr]^-{M\widehat{\otimes}_A P\tensor{A'}\lambda_N }
& &  M\widehat{\otimes}_A P\tensor{A'}(A'\otimes kX')\tensor{A'}N .}
\]
(See Section \ref{Cohom functors}.) Since $(A\otimes
kX)\widehat{\otimes}_A P\simeq P$ as $(A\otimes kX)-A'$-bicomodules,
$P$ can be endowed with a structure of $(A\otimes kX)-(A'\otimes
kX')$-bicomodule such that $(A\times kX)\widehat{\otimes}_A P\simeq
P$ as $(A\otimes kX)-(A'\otimes kX')$-bicomodules. Moreover, the
right coaction on $P$, $\rho_P:P\rightarrow P\tensor{A'}(A'\otimes
kX')$, is the unique $A'$-linear map satisfying the commutativity of
the diagram
\[
\xymatrix{\coring{C} \ar[rr]^{\theta_{\coring{C}}}
\ar[d]^{\theta_{\coring{C}}}&& \coring{C}\widehat{\otimes}_A
P\tensor{A'}N \ar[r]^{\simeq} \ar[d]^{\rho_{(\coring{C}
\widehat{\otimes}_A P)}\tensor{A'}N} &
P\tensor{A'}N\ar[d]^{\rho_P\tensor{A'}N}\\
\coring{C}\widehat{\otimes}_A P\tensor{A'}N\ar[rr]^-{\coring{C}
\widehat{\otimes}_A P\tensor{A'}\lambda_N} && \coring{C}
\widehat{\otimes}_A P\tensor{A'}\coring{C}'\tensor{A'}N
\ar[r]^-{\simeq}& P\tensor{A'}\coring{C}'\tensor{A'}N ,}
\]
where $\coring{C}=A\otimes kX$ and $\coring{C}'=A'\otimes kX'$.
The commutativity of the above diagram is equivalent to the
condition \eqref{gradedqf}. Finally, By Theorem \ref{3},
$\cohom{A\otimes kX}{N}{-}\simeq -\widehat{\otimes}_A((A\otimes
kX)\widehat{\otimes}_A P)\simeq -\widehat{\otimes}_A P.$
\end{proof}

Now we are in a position to state and prove the main results of this
section which characterize adjoint pairs and Frobenius pairs of
functors between categories of graded modules over $G$-sets.

\begin{theorem}
Let $M$ be an $X\times X'$-graded $(A,A')$-bimodule and $N$ an
$X'\times X$-graded $(A',A)$-bimodule. Then the following are
equivalent
\begin{enumerate}[(1)]
\item $(-\widehat{\otimes}_AM,-\widehat{\otimes}_{A'}N)$ is an adjoint
pair;
\item $(N\widehat{\otimes}_A-,M\widehat{\otimes}_{A'}-)$ is an adjoint
pair;
\item $N_x$ is finitely generated projective in $\lmod{A'}$, for every
$x\in X$, and $\h{_{A'}N}{A'} \simeq M$ as $X\times X'$-graded
$A-A'$-bimodules;
\item
there exist bigraded maps
\[
\psi:\widehat{A}\rightarrow M\widehat{\otimes}_{A'}N \text{ and
 }\omega:N \widehat{\otimes}_A M\rightarrow\widehat{A'},
\]
 such that
\begin{equation}
(\omega\widehat{\otimes}_{A'}N)\circ(N \widehat{\otimes}_A\psi)=
\Lambda\text{ and }(M\widehat{\otimes}_{A'}\omega) \circ(\psi
\widehat{\otimes}_AM)=M.
\end{equation}
\end{enumerate}
In particular, $(-\widehat{\otimes}_AM,-\widehat{\otimes}_{A'}N)$ is
a Frobenius pair if and only if
$(M\widehat{\otimes}_{A'}-,N\widehat{\otimes}_A-)$ is a Frobenius
pair.
\end{theorem}

\begin{proof}
We know that the corings $A\otimes kX$ and $A'\otimes kX'$ are
coseparable. Proposition \ref{19a} then achieves the proof.
\end{proof}

\begin{theorem}
Let $G$ and $G'$ be two groups, $A$ a $G$-graded $k$-algebra, $A'$ a
$G'$-graded $k$-algebra,  $X$ a right $G$-set, and $X'$ a right
$G'$-set. For a pair of $k$-linear functors $F:gr-(A,X,G) \to
gr-(A',X',G')$ and $G:gr-(A',X',G') \to gr-(A,X,G)$, the following
statements are equivalent\begin{enumerate}[(a)] \item $(F,G)$ is a
Frobenius pair;\item there exist an $X\times X'$-graded
$(A,A')$-bimodule $M$, and an $X'\times X$-graded $(A',A)$-bimodule
$N$, with the following properties
\begin{enumerate}[(1)] \item $_xM$ and $N_x$ is finitely generated
projective in $\rmod{A'}$ and $\lmod{A'}$ respectively, for every
$x\in X$, and $M_{x'}$ and $_{x'}N$ is finitely generated projective
in $\lmod{A}$ and $\rmod{A}$ respectively, for every $x'\in X'$,
\item $\h{_AM}{A} \simeq N$ as $X'\times X$-graded $(A',A)$-bimodules,
and

$\h{_{A'}N}{A'} \simeq M$ as $X\times X'$-graded $(A,A')$-bimodules,
\item $F\simeq -\widehat{\otimes}_AM$ and
$G\simeq -\widehat{\otimes}_{A'}N$.
\end{enumerate}
\end{enumerate}
\end{theorem}

\begin{proof}
Straightforward from Theorem \ref{leftrightFrobenius}.
\end{proof}

\subsection{When is the induction functor $T^*$ Frobenius?}\label{$T^*$ Frob}

Finally, let $f:G\to G'$ be a morphism of groups, $X$ a right
$G$-set, $X'$ a right $G'$-set, $\varphi:X\to X'$ a map such that
$\varphi(xg) =\varphi(x)f(g) $ for every $g\in G,$ $x\in X.$ Let $A$
be a $G$-graded $k$-algebra, $A'$ a $G'$-graded $k$-algebra, and
$\alpha:A\to A'$ a morphism of algebras such that
$\alpha(A_{g})\subset A_{f(g)}'$ for every $g\in G.$

At first we recall some functors from \cite{Menini:1993} and
\cite[Theorem 2.27]{Menini/Nastasescu:1994}.

Let $T_*:gr-(A',X',G')\to gr-(A,X,G)$ be the functor defined by
$$T_*(M)=\bigoplus_{x\in X}M_{\varphi(x)}^x,$$ for $M\in
gr-(A',X',G')$, where $M_{\varphi(x)}^x=M_{\varphi(x)}$ for every
$x\in X$.

Set $m^x=m$ for $m\in M_{\varphi(x)}$ and $x\in X$. The right
$A$-action on $T_*(M)$ is defined by
$$m^x.a_g=(m\alpha(a_g))^{xg}$$ for all $g\in G,x\in X,a_g\in A_g,
m\in M_{\varphi(x)}.$

Moreover, we have for every $x\in X$,
$$(T_*(M))_x\simeq \hom{gr-(A',X',G')}{A'(\varphi(x))}{M}.$$
Hence,
$$T_*(M)\simeq \bigoplus_{x\in X}\hom{gr-(A',X',G')}{A'(\varphi(x))}{M}.$$
The functor $T_*$ is a covariant exact functor and has a left
adjoint $T^*$. We will recall the definition of this functor later.

Now, we define another functor $\widetilde T:gr-(A,X,G)\to
gr-(A',X',G')$. For every $L\in gr-(A,X,G)$, $$\widetilde
T(L)=\bigoplus_{x'\in X'}(\widetilde T(L))_{x'},$$ where
$$(\widetilde T(L))_{x'}=\hom{gr-(A,X,G)}{T_*(A'(x'))}{L}.$$
Let $g'\in g',a'_{g'}\in A'_{g'},x'\in X',\xi\in (\widetilde
T(L))_{x'}$. Set $$\xi.a'_{g'}=\xi\circ T_*(\lambda_{a'_{g'}}),$$
where $\lambda_{a'_{g'}}:A'(x'g')\to A'(x')$ is the left
multiplication by $a'_{g'}$ on $A'$.

The functor $\widetilde T$ is a right adjoint of $T_*$. Then
$\widetilde T$ is left exact.

For more details we refer to \cite{Menini:1993}.

\begin{theorem}\emph{\cite[Theorem 2.27]{Menini/Nastasescu:1994}}
\label{MN 2.27}
\\The functors $T^*$ and $\widetilde T$ are isomorphic if and only
if the following conditions hold:
\begin{enumerate}[(1)]
\item for every $x'\in X'$, $T_*(A'(x'))$ is finitely generated and
projective in $\rmod{A}$;
\item for every $x\in X$, there exists an isomorphism in
$gr-(A',X',G')$
$$\xymatrix{\theta_x:A'(\varphi(x))\ar[r]^-\sim & \widetilde T(A(x))}$$
such that
$$\big((\theta_{x_1}(1))(a')\big).a=(\theta_{x_2}(1))(a'.\alpha(a))$$ for
every $x_1,x_2,y\in X,a\in A(x_2)_{x_1},a'\in
A'(\varphi(x_1))_{\varphi(y)}$.
\end{enumerate}
\end{theorem}

We have, $\gamma:kX\to kX'$ such that $\gamma(x) =\varphi(x)$ for
each $x\in X$, is a morphism of coalgebras, and $(\alpha ,\gamma)
:(A,kX,\psi) \to(A',kX',\psi') $ is a morphism in
$\mathbb{E}_\bullet^\bullet(k).$

 Let $-\otimes_AA': gr-(A,X,G)\to gr-(A',X',G')$ be
the functor making commutative the following diagram
\[
\xymatrix{gr-(A,X,G)\ar[d]^\simeq \ar[rr]^{-\otimes_AA'} & &
gr-( A',X',G') \ar[d]^\simeq  \\
\mathcal{M}(kG)_A^{kX}
\ar[rr]^{-\otimes_AA'%
} & & \mathcal{M}( kG')_{A'}^{kX'}.}
\]
Let $M\in gr-(A,X,G).$ We have $M\otimes_A A'$ is a right
$A'$-module, and
$$\rho( m_{x}\otimes _Aa_{g'}')
=(m_{x}\otimes_Aa_{g'%
}')  \otimes\varphi( x)  g'.$$ Therefore, $(M\otimes_AA') _{x'}$ is
the $k$-module of $M\otimes_AA'$ spanned by the elements of the form
$m_{x}\otimes_Aa_{g'}'$ where $x\in X,$ $g'\in G',$ $\varphi(x)
g'=x',$ $m_{x}\in M_{x},$ $a_{g'}'\in A_{g'}',$ for every $x'\in
X'.$

 Therefore, the functor $-\otimes_AA':gr-(A,X,G) \rightarrow
gr-(A',X',G')$ is exactly the induction functor $T^*$ defined in
\cite{Menini:1993} (or \cite[p. 531]{Menini/Nastasescu:1994}). Hence
$T^*$ is a Frobenius functor if and only if the induction functor
$-\otimes_AA':\mathcal{M}^{A\otimes kX}\to \mathcal{M}^{A'\otimes
kX'}$ is a Frobenius functor (see Section~\ref{entwined}).

 Moreover, we have the commutativity of the following diagram
\[%
\xymatrix{(G,X,A)-gr \ar[d]^\simeq \ar[rr]^{(T^*)'=A'\otimes_A-}
& & (G',X',A')-gr \ar[d]^\simeq\\
{}_A^{kX}\mathcal{M}(\psi^{-1})
\ar[d]^\simeq\ar[rr]^{A'%
\otimes_A-} & & {}_{A'}^{kX'}\mathcal{M}%
((\psi')^{-1}) \ar[d]^\simeq \\
{}^{kX\otimes A}\mathcal{M} \ar[d]^\simeq
\ar[rr]^{A'\otimes_A-}%
 & & {}^{kX'\otimes A'}\mathcal{M} \ar[d]^\simeq
\\
{}^{A\otimes kX}\mathcal{M} \ar[rr]^{A'\otimes_A-}%
 & & {}^{A'\otimes kX'}\mathcal{M}.}
\]

The following consequence of Theorem \ref{entwined modules} and
Theorem \ref{MN 2.27} give two different characterizations when the
induction functor $T^*$ is a Frobenius functor.

\begin{theorem}
The following statements are equivalent
\begin{enumerate} [(a)]
\item the functor $T^*:gr-(A,X,G)\to gr-(A',X',G')$
is a Frobenius functor;
\item $(T^*(\widehat{A}))_{x'}$ is finitely
generated projective in $\lmod{A}$, for every $x'\in X'$, and there
exists an isomorphism of $X'\times X$ graded $(A',A)$-bimodules
\[
\h{_AT^*(\widehat{A})}{A}\simeq (T^*)'(\widehat{A}).
\]
\end{enumerate}
\end{theorem}

As an immediate consequence of Proposition \ref{27}, we obtain

\begin{proposition}
The following are equivalent
\begin{enumerate} [(a)]\item the functor
$T^*:gr-(A,X,G) \to gr-( A',X',G')$ is a Frobenius functor;
\item the functor $(T^*)':(G,X,A)-gr\to (G',X',A')-gr$ is
a Frobenius functor.
\end{enumerate}
\end{proposition}

\chapter{Comatrix Coring Generalized and Equivalences of Categories of Comodules}

\section{Comatrix coring generalized}\label{comatrixgeneralized}

In this section we generalize the concept of a comatrix coring
defined in \cite{Brzezinski/Gomez:2003} and
\cite{ElKaoutit/Gomez:2003} to the case of a quasi-finite comodule,
and generalize some of its properties.

\medskip
The conjunction of Proposition \ref{19a} and the
following result generalizes \cite[Theorem
2.4]{Brzezinski/Gomez:2003}.

\begin{theorem}\label{Comatrixcoring}
Let $X\in{}^{\coring{C}}\mathcal{M}^{\coring{D}}$ and
$\Lambda\in{}^{\coring{D}}\mathcal{M}^\coring{C}$. Assume that at
least one of the following conditions holds
\begin{enumerate}[(1)] \item
\begin{enumerate}[(a)] \item $_A\coring{C}$,
$_B\coring{D}$, $\coring{C}_A$ and $\coring{D}_B$ are flat, \item
$_{\coring{C}}X$ and $_{\coring{D}}\Lambda$ are coflat; or
\end{enumerate}
\item $A$ and $B$ are von Neumann regular rings; or
\item $_A\coring{C}$,
$_B\coring{D}$ are flat, and $\coring{C}$ and $\coring{D}$ are
coseparable corings.
\end{enumerate}
If there exist bicolinear maps
\[
\psi:\coring{C}\rightarrow X\square_{\coring{D}}\Lambda\;\textrm{and
}\;\omega:\Lambda\square_{\coring{C}}X\rightarrow\coring{D}
\]
in $^{\coring{C}}\mathcal{M}^\coring{C}$ and
$^{\coring{D}}\mathcal{M}^\coring{D}$ respectively, such that the
diagrams
\begin{equation}\label{unitcounit}
\xymatrix{\Lambda\ar[rr]^\simeq \ar[d]_\simeq &&
\Lambda\cotensor{\coring{C}}\coring{C}
\ar[d]^{\Lambda\cotensor{\coring{C}}\psi}
\\\coring{D}\cotensor{\coring{D}}\Lambda &&
\Lambda\cotensor{\coring{C}}X\cotensor{\coring{D}}\Lambda
\ar[ll]^{\omega\cotensor{\coring{D}}\Lambda}}
\xymatrix{X\ar[rr]^\simeq \ar[d]_\simeq &&
\coring{C}\cotensor{\coring{C}}X \ar[d]^{\psi\cotensor{\coring{C}}X}
\\X\cotensor{\coring{D}}\coring{D} &&
X\cotensor{\coring{D}}\Lambda\cotensor{\coring{C}}X
\ar[ll]^{X\cotensor{\coring{D}}\omega}}
\end{equation}
commute, then $\Lambda\square_{\coring{C}}X$ is a $B$-coring with
coproduct
\[
\xymatrix{\Delta:\Lambda\cotensor{\coring{C}}X \ar[r]^-\simeq &
\Lambda\cotensor{\coring{C}}(
\coring{C}\cotensor{\coring{C}}X)\ar[rr]^-{\Lambda\cotensor{\coring{C}}
(\psi\cotensor{\coring{C}}X)} & &
\Lambda\cotensor{\coring{C}}((X\cotensor{\coring{D}}\Lambda)
\cotensor{\coring{C}}X)\ar[r]^-\simeq &
(\Lambda\cotensor{\coring{C}}X)\cotensor{\coring{D}}(
\Lambda\cotensor{\coring{C}}X) \ar@{^{(}->}[d] \\ & & & &
(\Lambda\cotensor{\coring{C}}X)\tensor{B}(\Lambda
\cotensor{\coring{C}}X),}
 \]
 and counit
$\xymatrix{\epsilon=\epsilon_{\coring{D}}\circ\omega
:\Lambda\square_{\coring{C}}X\ar[r] & B}$.
\end{theorem}

\begin{proof}
By the bicolinearity of $\psi$, we have the commutativity of the two
diagrams
\begin{equation}\label{symetrydiagram1}
\xymatrix{\coring{C}\ar[r]^{\Delta_{\coring{C}}}
\ar[ddrr]^{(\psi\square_{\coring{C}}\psi) \circ\Delta_{\coring{C}}}
\ar[d]_{\psi} & \coring{C}\square_{\coring{C}}\coring{C}
\ar[r]^-{\psi\square_{\coring{C}}\coring{C}} &
(X\square_{\coring{D}}\Lambda)\square_{\coring{C}}\coring{C}
\ar[dd]^{ (X\square_{\coring{D}}\Lambda)\square_{\coring{C}}\psi}
\\ (X\square_{\coring{D}}\Lambda)\ar[d]_\simeq & &
\\
(X\square_{\coring{D}}\Lambda)\square_{\coring{C}}\coring{C}
\ar[rr]_{(X\square_{\coring{D}}\Lambda) \square_{\coring{C}}\psi} &
& (X\square_{\coring{D}}\Lambda)
\square_{\coring{C}}(X\square_{\coring{D}}\Lambda)}
\end{equation}
\begin{equation}\label{symetrydiagram2}
\xymatrix{\coring{C}\ar[r]^{\Delta_{\coring{C}}}
\ar[ddrr]^{(\psi\square_{\coring{C}}\psi) \circ\Delta_{\coring{C}}}
\ar[d]_{\psi} & \coring{C}\square_{\coring{C}}\coring{C}
\ar[r]^-{\coring{C}\square_{\coring{C}}\psi} &
\coring{C}\square_{\coring{C}}(X\square_{\coring{D}}\Lambda)
\ar[dd]^{ \psi\square_{\coring{C}}(X\square_{\coring{D}}\Lambda)}
\\ (X\square_{\coring{D}}\Lambda)\ar[d]_\simeq & &
\\
\coring{C}\square_{\coring{C}}(X\square_{\coring{D}}\Lambda)
\ar[rr]_{\psi\square_{\coring{C}} (X\square_{\coring{D}}\Lambda)} &
& (X\square_{\coring{D}}\Lambda)
\square_{\coring{C}}(X\square_{\coring{D}}\Lambda).}
\end{equation}
Therefore, we have the commutativity of the diagram
\[
\xymatrix@C 2pt{\Lambda\square_{\coring{C}}X \ar[r]^\simeq
\ar[d]_\simeq&
\Lambda\square_{\coring{C}}\coring{C}\square_{\coring{C}}X
\ar[rr]^{\Lambda\square_{\coring{C}} \psi\square_{\coring{C}}X} & &
(\Lambda\square_{\coring{C}}X)\square_{\coring{D}}
(\Lambda\square_{\coring{C}}X) \ar[d]^\simeq
 \\
\Lambda\square_{\coring{C}}\coring{C} \square_{\coring{C}}X\ar[dd]_
{\Lambda\square_{\coring{C}}\psi \square_{\coring{C}}X} & & &
(\Lambda\square_{\coring{C}}\coring{C}
\square_{\coring{C}}X)\square_{\coring{D}}
(\Lambda\square_{\coring{C}}X)
\ar[dd]|{(\Lambda\square_{\coring{C}}\psi
\square_{\coring{C}}X)\square_{\coring{D}}
(\Lambda\square_{\coring{C}}X)}
\\ & & & &
\\
(\Lambda\square_{\coring{C}}X)\square_{\coring{D}}
(\Lambda\square_{\coring{C}}X) \ar[dr]_\simeq& & &
(\Lambda\square_{\coring{C}}X)\square_{\coring{D}}
(\Lambda\square_{\coring{C}}X)
\square_{\coring{D}}(\Lambda\square_{\coring{C}}X)
\\ &
(\Lambda\square_{\coring{C}}X)\square_{\coring{D}}
(\Lambda\square_{\coring{C}}\coring{C}\square_{\coring{C}}X).
\ar[urr]|{(\Lambda\square_{\coring{C}}X)\square_{\coring{D}}
(\Lambda\square_{\coring{C}}\psi\square_{\coring{C}}X)} & & &}
\]
Hence the coassociative property of $\Delta$ follows. On the other
hand, if we put
$$i:(\Lambda\square_{\coring{C}}X)\square_{\coring{D}}
(\Lambda\square_{\coring{C}}X)
\hookrightarrow(\Lambda\square_{\coring{C}}X)\otimes_B(\Lambda
\square_{\coring{C}}X)$$ the canonical injection, we have,
\begin{multline}\label{6}
[( \epsilon_{\coring{D}}\circ
\omega\otimes_B(\Lambda\square_{\coring{C}}X)] \circ
i\circ(\Lambda\square_{\coring{C}}\mathcal{\psi}
\square_{\coring{C}}X)\circ[\xymatrix@1{\Lambda\square_{\coring{C}
}X\ar[r]^-\simeq &
\Lambda\square_{\coring{C}}\coring{C}\square_{\coring{C}}X}]
\\
 =[\epsilon_{\coring{D}}
\otimes_B(\Lambda\square_{\coring{C}}X)]\circ
i'\circ(\omega\square_{\coring{D}}(\Lambda
\square_{\coring{C}}X))\circ(\Lambda\square
_{\coring{C}}\mathcal{\psi}\square_{\coring{C}}X)\circ[\xymatrix@1{
\Lambda\square_{\coring{C}}X \ar[r]^-\simeq & \Lambda
\square_{\coring{C}}\coring{C}\square_{\coring{C}}X}],
\end{multline}
where $i':\coring{D}\square_{\coring{D}}(\Lambda\square
_{\coring{C}}X)
\hookrightarrow\coring{D}\otimes_B(\Lambda\square_{\coring{C}}X) $
is the canonical injection.

The first diagram of \eqref{unitcounit} is commutative means that
the diagram
\[
\xymatrix{\Lambda\square_{\coring{C}}
\coring{C}\ar[rr]^-{\Lambda\square _{\coring{C}}\mathcal{\psi}}
\ar[drr]_\simeq & & \Lambda\square
_{\coring{C}}(X\square_{\coring{D}}\Lambda)
\ar[rr]^-{\omega\square_{\coring{D}}\Lambda}
& & \coring{D}\square_{\coring{D}}\Lambda \\
 &  & \Lambda\ar[urr]_\simeq &  & }
\]
commutes. The composition of morphisms
\[
\xymatrix@1{\Lambda\ar[r]^-\simeq &
\Lambda\square_{\coring{C}}\coring{C}\ar[r]^-\simeq &
\Lambda\ar[r]^-\simeq &
\coring{D}\square_{\coring{D}}\Lambda\ar@{^{(}->}[r] &
\coring{D}\otimes_B\Lambda
\ar[rr]^{\epsilon_{\coring{D}}\otimes_B\Lambda} & &
B\otimes_B\Lambda}
\]
is exactly the morphism $[\lambda\mapsto1_B\otimes_B\lambda].$ Then,
\[
(\ref{6})=\Big[\sum_i\lambda_i\otimes_Ax_i\mapsto1_B\otimes
_B(\sum_i\lambda_i\otimes_Ax_i)\Big]  .
\]
Analogously, by using the commutativity of the second diagram of
\eqref{unitcounit}, we complete the proof of the counit property.
\end{proof}

\medskip

Let $N\in{}^{\coring{C}}\mathcal{M}^\coring{D}$ be a bicomodule,
quasi-finite as a right $\coring{D}$-comodule, such that
$_B\coring{D}$ is flat or $\coring{C}$ is a coseparable $A$-coring.
Let $T$ be a $k$-algebra.
\\ We will give an other characterization of the
natural isomorphism (see Section \ref{Frobeniusgeneral}):
\[
\Upsilon_{-,-}:-\otimes_{T}\cohom{\coring{D}}{N}{-}\rightarrow
\cohom{\coring{D}}{N}{-\otimes_{T}-}
\]
associated to the cohom functor $\cohom{\coring{D}}{N}{-}
:\mathcal{M}^\coring{D}\rightarrow\mathcal{M}^\coring{C}.$
\\ Let
$\theta:1_{\mathcal{M}%
^{\coring{D}}}\rightarrow h_{\coring{D}}(N,-)\otimes _AN$ be the
unit of the adjunction $(\cohom{\coring{D}}{N}{-},-\otimes_AN)$.
\\ Let $M\in{}^{T}\mathcal{M}^\coring{D}$ and
$W\in\mathcal{M}^{T}.$ Since the functors
$\cohom{\coring{D}}{N}{-}:\mathcal{M}%
^{\coring{D}}\rightarrow\mathcal{M}^{A}$ and
$-\otimes_AN:\mathcal{M}%
^{A}\rightarrow\mathcal{M}^\coring{D}$ are $k$-linear and preserve
inductive limits, then the functor $\cohom{\coring{D}}{N}{-}
\otimes_AN:\mathcal{M}^\coring{D}\rightarrow\mathcal{M}^{\coring{D}%
}$ is also $k$-linear and preserves inductive limits. By Lemma
\ref{G2} (1), $\theta_M$ is a morphism in
$^{T}\mathcal{M}^\coring{D}.$ Define
\[
\Xi_{W,M}:\cohom{\coring{D}}{N}{W\otimes_{T}M}\rightarrow
W\otimes_{T}\cohom{\coring{D}}{N}{M}
\]
to be the unique morphism in $\mathcal{M}^{A}$ satisfying $(\Xi
_{W,M}\otimes_AN) \theta_{W\otimes_{T}M}=W\otimes_T\theta_M$. We
have
\begin{multline}
\Big[\rho_{W\otimes_T\cohom{\coring{D}}{N}{M}}\otimes_AN-\big(W
\otimes_T\cohom{\coring{D}}{N}{M}\big)
\otimes_A\lambda_N\Big](W\otimes_{T}\theta_M)   \\
=W\otimes_{T}\Big[\big(\rho_{\cohom{\coring{D}}{N}{M}
}\otimes_AN-\cohom{\coring{D}}{N}{M}\otimes_A\lambda_N\big)
\theta_M\Big]=0.
\end{multline}
Then $\operatorname{Im}(W\otimes_T\theta_M)\subset
(W\otimes_T\cohom{\coring{D}}{N}{M}\square_{\coring{C}}N,$ and by
Lemma \ref{cohomind}, $\Xi_{W,M}$ is a morphism in
$\mathcal{M}^\coring{C}$. Now we will verify that $\Xi_{-,M}$ is a
natural transformation. For this, let $f:W\rightarrow W'$ be a
morphism in $\mathcal{M}^{T}.$ We will verify the commutativity of
the diagram
\[%
\xymatrix{\cohom{\coring{D}}{N}{W\otimes_TM}\ar[d]
\ar[rr]^{\Xi_{W,M}
} & & W\otimes_T\cohom{\coring{D}}{N}{M} \ar[d] \\
 \cohom{\coring{D}}{N}{W'\otimes_TM}
\ar[rr]^{\Xi_{W',M}} & &
 W'\otimes_T\cohom{\coring{D}}{N}{M}.}
\]
Since $\theta$ is a natural transformation, the diagram
\[%
 \xymatrix{W\otimes_TM \ar[d]
\ar[rr]^-{\theta_{W\otimes_TM}} & &
\cohom{\coring{D}}{N}{W\otimes_TM}\otimes_AN \ar[d]\\
W'\otimes_TM \ar[rr]^-{\theta_{W'\otimes_TM}} & &
\cohom{\coring{D}}{N}{W'\otimes_TM}\otimes_AN}
\]
commutes. Then, we have
\begin{eqnarray*}
 \Big[\big(\Xi_{W',M}\cohom{\coring{D}}{N}{f\otimes_TM}
\big)\otimes_AN\Big]\theta_{W\otimes_TM} &=& (\Xi_{W',M}\tensor{A}N)
\big(\cohom{\coring{D}}{N}{f\otimes_TM}\tensor{A}N\big)
\theta_{W\otimes_TM}\\
 &=&
 (\Xi_{W',M}\tensor{A}N)\theta_{W'\otimes_TM}(f\tensor{T}M)\\
 &=& (W'\tensor{T}\theta_M)(f\tensor{T}M)\\
 &=& \big(f\tensor{T}\cohom{\coring{D}}{N}{M}\tensor{A}N\big)
 (W\tensor{T}\theta_M)\\
 &=&
 \big(f\tensor{T}\cohom{\coring{D}}{N}{M}\tensor{A}N\big)
 (\Xi_{W,M}\tensor{A}N)\theta_{W\otimes_TM}\\
&=&
 \bigg[\Big(\big(
f\otimes_T\cohom{\coring{D}}{N}{M}\big)\Xi_{W,M}\Big)\otimes_AN\bigg]
\theta_{W\otimes_TM}.
\end{eqnarray*}
By uniqueness, $\Xi_{W',M}\cohom{\coring{D}}{N}{f\otimes_TM}=(
f\otimes_T\cohom{\coring{D}}{N}{M})\Xi_{W,M}$. Finally, we will
verify that $\Upsilon_{T,M}=\Xi_{T,M}^{-1}$, i.e. $\Xi_{T,M}$
satisfies the commutativity of the diagram
\[%
\xymatrix{\cohom{\coring{D}}{N}{T\otimes_TM}\ar[rd]_\simeq
\ar[rr]^{\Xi_{T,M} }&  &
T\otimes_T\cohom{\coring{D}}{N}{M}\ar[ld]^\simeq
\\
 &  \cohom{\coring{D}}{N}{M} & }
\]
(the canonical isomorphism $\xymatrix{T\otimes_TM \ar[r]^-\simeq &
M}$ is an isomorphism in $\mathcal{M}^\coring{D}$). Since $\theta$
is a natural transformation, the diagram
\[%
\xymatrix{T\otimes_TM \ar[d]^\simeq \ar[rr]^-{\theta_{T\otimes_TM}}
& & \cohom{\coring{D}}{N}{T\otimes_TM}\otimes_AN
\ar[d]^\simeq \\
 M \ar[rr]^-{\theta_M} &  &
\cohom{\coring{D}}{N}{M}\otimes_AN}
\]
commutes. It follows from the commutativity of the last diagram and
the fact that $\theta_M$ is $T$-linear, that the following diagram
is also commutative
\[%
\xymatrix{ T\otimes_TM \ar[d]^\simeq \ar[rr]^-{\theta_{T\otimes_TM}}
& & \cohom{\coring{D}}{N}{T\otimes_TM} \otimes_AN
\ar[d]^\simeq\\
 T\otimes_T\cohom{\coring{D}}{N}{M}\otimes_AN
\ar[rr]^\simeq  & & \cohom{\coring{D}}{N}{M}\otimes_AN.}
\]
By Corollary \ref{Mitchell-Cor}, $\Xi_{W,M}$ is an isomorphism, and
by Mitchell's Theorem \ref{Mitchell-Th},
$\Upsilon_{W,M}=\Xi_{W,M}^{-1}$.

\medskip
We obtain then the following generalization of
\cite[1.6]{Takeuchi:1977}:

\begin{proposition}\label{Upsilon}
Let $N\in{}^{\coring{C}}\mathcal{M}^\coring{D}$ be a bicomodule,
quasi-finite as a right $\coring{D}$-comodule, such that
$_B\coring{D}$ is flat or $\coring{C}$ is a coseparable $A$-coring.
Let $\theta:1_{\mathcal{M}^{\coring{D}}}\to
\cohom{\coring{D}}{N}{-}\otimes_AN$ be the unit of the adjunction
$(\cohom{\coring{D}}{N}{-},-\otimes_AN)$, and let $T$ be a
$k$-algebra. If $\Upsilon
_{-,-}:-\otimes_T\cohom{\coring{D}}{N}{-}\rightarrow
\cohom{\coring{D}}{N}{-\otimes_T-}$ is the natural isomorphism
associated to the cohom functor
$\cohom{\coring{D}}{N}{-}:\mathcal{M}^\coring{D}\to\mathcal{M}^\coring{C},$
then $\Upsilon_{W,M}=\Xi_{W,M}^{-1},$ where $\Xi_{W,M}$ is defined
as above.
\end{proposition}

\medskip
As a consequence of the last result we get the following
generalization of \cite[1.13]{Takeuchi:1977}.

\begin{proposition}\label{delta}
Let $N\in\bcomod{\coring{C}}{\coring{D}}$ be a bicomodule,
quasi-finite as a right $\coring{D}$-comodule, such that
$_B\coring{D}$ is flat. Suppose that $_A\coring{C}$ is flat (resp.
 $\mathcal{M}^\coring{C}$ is an abelian category). Let
$\theta:1_{\mathcal{M}%
^{\coring{D}}}\rightarrow \cohom{\coring{D}}{N}{-}\otimes_AN$ be the
unit of the adjunction $(\cohom{\coring{D}}{N}{-},-\otimes_AN)$. If
the cohom functor $\cohom{\coring{D}}{N}{-}$ is exact (resp.
$\coring{D}$ is a coseparable $B$-coring), then the natural
isomorphism $\delta$ (see Corollary \ref{5}):
\[
\cohom{\coring{D}}{N}{-}
\simeq-\square_{\coring{D}}\cohom{\coring{D}}{N}{\coring{D}}:
\mathcal{M}^\coring{D}\rightarrow \mathcal{M}^\coring{C}
\]
satisfy the property:

For every $M\in\mathcal{M}^\coring{D}$,
$\delta_M:\cohom{\coring{D}}{N}{M}\rightarrow
M\square_{\coring{D}}\cohom{\coring{D}}{N}{\coring{D}}$ is the
unique right $A$-linear map satisfying
\[
(\delta_M\otimes_AN)\theta_M=(M\square_{\coring{D}}
\theta_{\coring{D}})\rho_M.
\]
\end{proposition}

\begin{proof}
 Let $M\in\mathcal{M}^\coring{D}$, we have
 $(\Xi_{M,\coring{D}}\otimes_AN)\theta_{M\otimes_B\coring{D}}=
 M\otimes_B\theta_{\coring{D}}$. On the other hand,
since $\theta$ is a natural transformation,
$(\Xi_{M,\coring{D}}\otimes_AN)
 \theta_{M\otimes_B\coring{D}}\rho_M=(\Xi_{M,\coring{D}}\otimes_AN)
 (\cohom{\coring{D}}{N}{\rho_M}\otimes_AN)
 \theta_M=(i\delta_M\otimes_AN)\theta_M$, where
 $i:M\square_{\coring{D}}\cohom{\coring{D}}{N}{\coring{D}}
 \hookrightarrow
M\otimes_B\cohom{\coring{D}}{N}{\coring{D}}$ is the canonical
injection. Therefore,
$(i\otimes_AN)(\delta_M\otimes_AN)\theta_M=(M\otimes_B
\theta_{\coring{D}})\rho_M$.
 Since $i\otimes_AN$ is a monomorphism ($_AN$ is
flat), $(\delta_M\otimes_AN)\theta_M
 =(M\square_{\coring{D}}\theta_{\coring{D}})\rho_M$.
\end{proof}

\medskip

Now we will give a series of properties concerning our comatrix
coring.

\begin{proposition}\label{comatrixproperty1}
Let $\Lambda\in \bcomod{\coring{D}}{\coring{C}}$ be a bicomodule,
quasi-finite as a right $\coring{C}$-comodule, such that
$_A\coring{C}$ and $_B\coring{D}$ are flat. Set
$X=\cohom{\coring{C}}{\Lambda}{\coring{C}}\in
\bcomod{\coring{D}}{\coring{C}}$.
 If
 \begin{enumerate}[(a)]\item $\coring{C}_A$ and $\coring{D}_B$ are flat,
 the cohom functor $\cohom{\coring{C}}{\Lambda}{-}$ is exact and
$_{\coring{D}}{\Lambda}$ is coflat, or
\item $A$ and $B$ are von Neumann regular rings and the cohom functor
$\cohom{\coring{C}}{\Lambda}{-}$ is exact, or
\item $\coring{C}$ and
$\coring{D}$ are coseparable corings,
\end{enumerate}
then we have
\begin{enumerate} [(1)]\item
$\xymatrix{\delta_{\Lambda}:\e{\coring{C}}{\Lambda}\ar[r] &
\Lambda\cotensor{\coring{C}}X}$ is an isomorphism of $B$-corings;
\item
$\xymatrix{\chi_{\coring{D}}\circ
\cohom{\coring{C}}{\Lambda}{\lambda_{\Lambda}}=
\omega\circ\delta_{\Lambda}: \e{\coring{C}}{\Lambda}\ar[r] &
\coring{D}}$ and
$\xymatrix{\omega:\Lambda\cotensor{\coring{C}}X\ar[r] &\coring{D}}$
is a homomorphism of $B$-corings, where
$\lambda_{\Lambda}:\Lambda\to\coring{D}\cotensor{\coring{D}}{\Lambda}$
is the left coaction of $\Lambda$, $\chi$ the counit of the
adjunction
$(\cohom{\coring{C}}{\Lambda}{-},-\square_{\coring{D}}\Lambda)$, and
$\omega$ is the unique $\coring{D}-\coring{D}$-bicolinear map such
that \eqref{counitcohom} is the counit of the adjunction
$(-\cotensor{\coring{C}}{\cohom{\coring{C}}{\Lambda}{\coring{C}}},
-\cotensor{\coring{D}} {\Lambda})$.
\end{enumerate}
\end{proposition}

\begin{proof}
Let $\theta:1_{\mathcal{M}^{\coring{C}}}\to
\cohom{\coring{C}}{\Lambda}{-}\tensor{B}{\Lambda}$ be the unit of
the adjunction
$(\cohom{\coring{C}}{\Lambda}{-},-\tensor{B}{\Lambda})$.

(1) By Theorem \ref{3}, $\cohom{\coring{C}}{\Lambda}{-}\simeq
-\square_{\coring{C}}\cohom{\coring{C}}{\Lambda}{\coring{C}}$. In
the first case, $\cohom{\coring{C}}{\Lambda}{\coring{C}}$ is coflat
as a left $\coring{C}$-comodule. By Proposition \ref{G4},
$\theta_{\coring{C}}:\coring{C}\rightarrow
\cohom{\coring{C}}{\Lambda}{\coring{C}} \square_{\coring{D}}\Lambda$
is a $\coring{C}$-bicolinear map.

From Proposition \ref{delta}, for every
$M\in{}\mathcal{M}^\coring{C}$, we have a commutative diagram
\[
\xymatrix{M\ar[d]_{\theta_M}\ar[rr]^\simeq & &
M\square_{\coring{C}}\coring{C}\ar[d]^{M
\square_{\coring{C}}\theta_{\coring{C}}}
\\
\cohom{\coring{C}}{\Lambda}{M}\square_{\coring{D}}\Lambda
\ar[rr]_-{\delta_M \square_{\coring{D}}\Lambda} & &
M\square_{\coring{C}}\cohom{\coring{C}}{\Lambda}
{\coring{C}}\square_{\coring{D}}\Lambda.}
\]
From Lemma \ref{18}, the unit and the counit of the adjunction
$(-\cotensor{\coring{C}}{\cohom{\coring{C}}
{\Lambda}{\coring{C}}},-\cotensor{\coring{D}} {\Lambda})$ are
\begin{equation}\label{unitcohom}
\xymatrix@1{\eta:1_{\mathcal{M}^\coring{C}}\ar[r]^\simeq &
-\square_{\coring{C}}\coring{C}\ar[rr]^-
{-\square_{\coring{C}}\theta_{\coring{C}}} & &
-\square_{\coring{C}}\cohom{\coring{C}}{\Lambda}{\coring{C}}
\square_{\coring{D}}\Lambda}
\end{equation}
 and
\begin{equation}\label{counitcohom}
\xymatrix@1{\varepsilon:-\square_{\coring{D}}
\Lambda\square_{\coring{C}} \cohom{\coring{C}}{\Lambda}{\coring{C}}
\ar[rr]^-{-\square_{\coring{D}}\omega} & &
-\square_{\coring{D}}\coring{D}\ar[r]^\simeq &
1_{\mathcal{M}^\coring{D}}.}
\end{equation}
To show that $\delta_\Lambda$ is a homomorphism of $B$-corings, it
suffices to prove (see Proposition \ref{BW 23.8}):
\begin{multline}
(\delta_{\Lambda}\tensor{B}{\delta_{\Lambda}}\tensor{B}{\Lambda})
\circ(\e{\coring{C}}{\Lambda}\tensor{B}{\theta_{\Lambda}})
\circ\theta_{\Lambda}
\\
=(\Lambda\cotensor{\coring{C}}{\theta_\coring{C}}
\cotensor{\coring{C}}X\tensor{B}{\Lambda})
\circ[\xymatrix{\Lambda\cotensor{\coring{C}}X\tensor{B}{\Lambda}
\ar[r]^-\simeq&
\Lambda\cotensor{\coring{C}}{\coring{C}}\cotensor{\coring{C}}X
\tensor{B}{\Lambda}}]
\circ(\delta_{\Lambda}\tensor{B}{\Lambda})\circ\theta_{\Lambda},
\end{multline}
and
\begin{multline}
(\epsilon_{\coring{D}}\tensor{B}{\Lambda})\circ
(\omega\tensor{B}{\Lambda})
\circ(\delta_{\Lambda}\tensor{B}{\Lambda}) \circ\theta_{\Lambda}=
[\xymatrix{\Lambda\ar[r]^-\simeq& B\tensor{B}{\Lambda}}].
\end{multline}

Now we will prove $(8)$. We have,
\begin{multline}
\begin{split}
(\delta_{\Lambda}\tensor{B}{\delta_{\Lambda}} \tensor{B}{\Lambda})
\circ(\e{\coring{C}}{\Lambda}\tensor{B} {\theta_{\Lambda}})
\circ\theta_{\Lambda} &
=\Big[\delta_{\Lambda}\tensor{B}{\big((\delta_{\Lambda}
\tensor{B}{\Lambda})\circ\theta_{\Lambda}\big)}\Big]
\circ\theta_{\Lambda}\\
& =\Big[(\Lambda\cotensor{\coring{C}}X)\tensor{B}
{\big((\delta_{\Lambda}\tensor{B}{\Lambda})
\circ\theta_{\Lambda}\big)\Big]}\circ(\delta_{\Lambda}
\tensor{B}{\Lambda})\circ\theta_{\Lambda}\\
& =\Big[(\Lambda\cotensor{\coring{C}}X)\tensor{B}
{\big((\Lambda\cotensor{\coring{C}}
{\theta_{\coring{C}}})\circ\rho_{\Lambda})\Big]}
\circ(\Lambda\cotensor{\coring{C}}
{\theta_{\coring{C}}})\circ\rho_{\Lambda}\\ & \qquad\textrm{(using
Proposition \ref{delta})}\\
& =(\Lambda\cotensor{\coring{C}}X\tensor{B}
{\Lambda\cotensor{\coring{C}} {\theta_{\coring{C}}}})\circ
(\Lambda\cotensor{\coring{C}}X\tensor{B}
{\rho_{\Lambda}})\circ(\Lambda\cotensor{\coring{C}}
{\theta_{\coring{C}}})\circ\rho_{\Lambda}.
\end{split}
\end{multline}
On the other hand, by Proposition \ref{delta},
\begin{multline}
(\Lambda\cotensor{\coring{C}}{\theta_\coring{C}}
\cotensor{\coring{C}}X\tensor{B}{\Lambda})
\circ[\xymatrix{\Lambda\cotensor{\coring{C}}X \tensor{B}{\Lambda}
\ar[r]^-\simeq& \Lambda\cotensor{\coring{C}}{\coring{C}}
\cotensor{\coring{C}}X\tensor{B}{\Lambda}}]
\circ(\delta_{\Lambda}\tensor{B}{\Lambda})
\circ\theta_{\Lambda}\\
=(\Lambda\cotensor{\coring{C}}{\theta_\coring{C}}
\cotensor{\coring{C}}X\tensor{B}{\Lambda})
\circ[\xymatrix{\Lambda\cotensor{\coring{C}}X\tensor{B}
{\Lambda}\ar[r]^-\simeq & \Lambda\cotensor{\coring{C}}
{\coring{C}}\cotensor{\coring{C}}X\tensor{B}{\Lambda}}]
\circ(\Lambda\cotensor{\coring{C}}{\theta_{\coring{C}}})
\circ\rho_{\Lambda}.
\end{multline}
 Now, by Proposition \ref{G4},
$\theta_{\coring{C}}:\coring{C}\rightarrow
\cohom{\coring{C}}{\Lambda}{\coring{C}}\tensor{B}{\Lambda}$ is
 $\coring{C}$-bicolinear. Then
$\theta_{\coring{C}}$
 makes commutative the diagrams
\eqref{symetrydiagram1} and \eqref{symetrydiagram2}
 (by replacing $\psi$ by $\theta_{\coring{C}}$). Hence
$(10)=(11)$.

Finally, we will prove $(9)$.
\begin{eqnarray*}
(\epsilon_{\coring{D}}\tensor{B}{\Lambda})
\circ(\omega\tensor{B}{\Lambda})
\circ(\delta_{\Lambda}\tensor{B}{\Lambda}) \circ\theta_{\Lambda}
&=& (\epsilon_{\coring{D}}\tensor{B}{\Lambda})
\circ(\omega\tensor{B}{\Lambda}) \circ(\Lambda\cotensor{\coring{C}}
{\theta_{\coring{C}}})\circ\rho_{\Lambda} \\
&& \textrm{(using Proposition \ref{delta})}\\
&=&(\epsilon_{\coring{D}}\tensor{B}{\Lambda}) \circ\lambda_{\Lambda} \\
&& \textrm{(using the first equality of \eqref{unitcounit})} \\
&=& [\xymatrix{\Lambda \ar[r]^-\simeq& B\tensor{B}{\Lambda}}].
\end{eqnarray*}

(2) Since $\delta$ is a natural isomorphism, the diagram
\[
\xymatrix{\cohom{\coring{C}}{\Lambda}{\Lambda}
\ar[d]_{\cohom{\coring{C}}{\Lambda}{\lambda_{\Lambda}}}
\ar[rr]^{\delta_{\Lambda}} & & {\Lambda}\cotensor{\coring{C}}X
\ar[d]^{\lambda_{\Lambda}\cotensor{\coring{C}}X}
\\
\cohom{\coring{C}}{\Lambda}{\coring{D}\cotensor{\coring{D}}{\Lambda}}
\ar[rr]_-{\delta_{\coring{D}\cotensor{\coring{D}}{\Lambda}}} & &
(\coring{D}\cotensor{\coring{D}}{\Lambda})\cotensor{\coring{C}}X}
\]
is commutative. The counit of the adjunction
$(\cohom{\coring{C}}{\Lambda}{-},-\cotensor{\coring{D}}{\Lambda})$
is \begin{equation}
\xymatrix@1{\cohom{\coring{C}}{\Lambda}{-\cotensor{\coring{D}}
{\Lambda}}\ar[rr]^{\delta_{-\cotensor{\coring{D}}{\Lambda}}} & &
-\cotensor{\coring{D}}{\Lambda}
\cotensor{\coring{C}}X\ar[rr]^-{-\cotensor{\coring{D}}{\omega}} & &
-\cotensor{\coring{D}}{\coring{D}}\ar[r]^\simeq &
1_{\mathcal{M}^\coring{D}}.}
\end{equation}
Then,
$\chi_{\coring{D}}\circ\cohom{\coring{C}}{\Lambda}{\lambda_{\Lambda}}=
[\xymatrix{\coring{D}\cotensor{\coring{D}}{\coring{D}}\ar[r]^\simeq
& \coring{D}}]\circ(\coring{D}\cotensor{\coring{D}}{\omega})
\circ(\lambda_{\Lambda}\cotensor{\coring{C}}X)\circ\delta_{\Lambda}$.\\
Since $\omega$ is left $\coring{D}$-colinear, and
$\lambda_{\Lambda\cotensor{\coring{C}}X}=\lambda_{\Lambda}
\cotensor{\coring{C}}X$, and from the commutativity of the following
diagram
\[
\xymatrix{\Lambda\cotensor{\coring{C}}X \ar[r]^-\simeq
\ar[dr]_{\lambda_{\Lambda\cotensor{\coring{C}}X}} &
\coring{D}\cotensor{\coring{D}}(\Lambda\cotensor{\coring{C}}X)
\ar[rr]^{\coring{D}\cotensor{\coring{D}}{\omega}} \ar@{^{(}->}[d] &
&
\coring{D}\cotensor{\coring{D}}\coring{D}\ar[r]^-\simeq\ar@{^{(}->}[d]
 & \coring{D} \ar[dl]^{\Delta_{\coring{D}}}
 \\ &
\coring{D}\tensor{B}(\Lambda\cotensor{\coring{C}}X)
\ar[rr]_{\coring{D}\tensor{B}{\omega}}
 & & \coring{D}\tensor{B}{\coring{D}} & ,}
\]
$\chi_{\coring{D}}\circ\cohom{\coring{C}}{\Lambda}{\lambda_{\Lambda}}=
\omega\circ\delta_{\Lambda}$. Finally, by Remark \ref{KG 5.2'}, it
is a homomorphism of $B$-corings. Hence, $\omega$ is also a
homomorphism of $B$-corings.
\end{proof}

\begin{corollary}\label{coendcoalgebra}
Let $C$ be a coalgebra over a field $k$, and let
$\Lambda\in\rcomod{C}$ be a quasi-finite and injective comodule.
Then $\e{C}{\Lambda}\simeq \Lambda\cotensor{C}\cohom{C}{\Lambda}{C}$
as coalgebras.
\end{corollary}

\begin{example}
Let $_A\coring{C}$ and $_B\coring{D}$ be flat. Let
$(\varphi,\rho):\coring{C}\rightarrow\coring{D}$ be a homomorphism
of corings. Suppose that $\coring{C}_A$ and $\coring{D}_B$ are flat
and $_{\coring{D}}(B\tensor{A}{\coring{C}})$ is coflat, or
$\coring{C}$ and $\coring{D}$ are coseparable corings. Then, we have
$(-\tensor{A}{B},-\cotensor{\coring{D}}(B\tensor{A}{\coring{C}}))$
is an adjoint pair, and
$-\tensor{A}{B}\simeq-\cotensor{\coring{C}}{(\coring{C}\tensor{A}{B})}$
(see the proof of Theorem \ref{23}). By Proposition \ref{BW 23.9},
the comatrix coring
$(B\tensor{A}{\coring{C}})\cotensor{\coring{C}}(\coring{C}\tensor{A}{B})$
($\simeq\e{\coring{C}}{B\tensor{A}{\coring{C}}}$) is isomorphic (as
corings) to the coring $B\coring{C}B$ (see Example \ref{examples
corings} (3)).
\end{example}

Now, we need the left-handed version of Corollary \ref{26}.

\begin{lemma}\label{26'}\emph{\cite[27.13]{Brzezinski/Wisbauer:2003}}
\\Let $\coring{C}$ be an $A$-coring and $T=\coring{C}^*.$
Then the following statements are equivalent
\begin{enumerate}[(a)]
\item $\coring{C}$ is a Frobenius coring;
\item $\coring{C}_A$ is finitely generated projective and
$\coring{C}\simeq T$ as $(A^\circ,T)$-bimodules, such that
$\coring{C}$ is a right $T$-module by $c.t=t(c_{(1)}).c_{(2)}$, for
every $c\in\coring{C}$ and $t\in T.$
\end{enumerate}
\end{lemma}

\begin{theorem}\label{comatrixproperty2}
Let $_A\coring{C}$ be flat. Let $\Lambda
\in{}_B\mathcal{M}^\coring{C}$ be a bicomodule, quasi-finite as a
right $\coring{C}$-comodule. Set
$X=\cohom{\coring{C}}{\Lambda}{\coring{C}}\in{}^{\coring{C}}\mathcal{M}_B$.
 If $\coring{C}_A$ is flat and the cohom functor
$\cohom{\coring{C}}{\Lambda}{-}$ is exact, or if $\coring{C}$ is a
coseparable coring, then we have
\begin{enumerate} [(1)]\item The functor
$-\tensor{B}{\Lambda}:
\mathcal{M}_B\rightarrow\mathcal{M}^\coring{C}$ is separable if and
only if the comatrix coring $\Lambda\cotensor{\coring{C}}X$ is a
cosplit $B$-coring. \item If the cohom functor
$\cohom{\coring{C}}{\Lambda}{-}$ is separable then the comatrix
coring $\Lambda\cotensor{\coring{C}}X$ is a coseparable $B$-coring.
\item If $-\tensor{B}{\Lambda}:
\mathcal{M}_B\rightarrow\mathcal{M}^\coring{C}$ is a Frobenius
functor, either $X_B$ is flat and $\Lambda_{\coring{C}}$ is coflat
or $\coring{C}$ is coseparable, and if for all
$f\in(\Lambda\cotensor{\coring{C}}X)^*$ and all $b\in B$,
$bH(f)=H(f)b$, where $\xymatrix{H:(\Lambda\cotensor{\coring{C}}X)^*
\ar[r] & {\Lambda\cotensor{\coring{C}}X}}$ is the isomorphism
defined by \eqref{frobeniuscomatrix} (this is the case if
$\Lambda\cotensor{\coring{C}}X$ is a coalgebra (i.e. $B=k$)), then
the comatrix coring $\Lambda\cotensor{\coring{C}}X$ is a Frobenius
$B$-coring.
\end{enumerate}
\end{theorem}

\begin{proof}
(1) Let $\omega:\Lambda\cotensor{\coring{C}}X\rightarrow B$ the
unique $B$-bilinear map such that
\begin{equation}
\xymatrix@1{\varepsilon:-\otimes_B\Lambda\square_{\coring{C}}
\cohom{\coring{C}}{\Lambda}{\coring{C}} \ar[rr]^-{-\otimes_B\omega}
& & -\otimes_B\coring{D}\ar[r]^\simeq & 1_{\mathcal{M}_B}}
\end{equation}
is the counit of the adjunction
$(-\cotensor{\coring{C}}{\cohom{\coring{C}}{\Lambda}
{\coring{C}}},-\tensor{B}{\Lambda})$. We have,
$\epsilon_{\Lambda\cotensor{\coring{C}}X}=\omega$. By Lemma \ref{18}
and Rafael's Proposition \ref{Rafael}, the functor
$-\tensor{B}{\Lambda}:
\mathcal{M}_B\rightarrow\mathcal{M}^\coring{C}$ is separable if and
only if there exists a $B$-bilinear map
$\omega':B\rightarrow\Lambda\cotensor{\coring{C}}X$
 satisfying $\omega\circ\omega'=1_B$. Hence, (1)
 follows from the definition of a cosplit coring.

(2) We know that the unit of the adjunction
$(-\cotensor{\coring{C}}{\cohom{\coring{C}}{\Lambda}
{\coring{C}}},-\tensor{B} {\Lambda})$ is
\begin{equation}
\xymatrix@1{\eta:1_{\mathcal{M}^\coring{C}}\ar[r]^\simeq &
-\square_{\coring{C}}\coring{C}\ar[rr]^-{-\square_{\coring{C}}
\theta_{\coring{C}}} & &
-\square_{\coring{C}}\cohom{\coring{C}}{\Lambda}
{\coring{C}}\otimes_B\Lambda}.
\end{equation}
By Lemma \ref{18} and Rafael's Proposition \ref{Rafael}, the cohom
functor is separable if and only if there exists a
$\coring{C}-\coring{C}$-bicolinear map
$\psi':\cohom{\coring{C}}{\Lambda}{\coring{C}}\otimes_B\Lambda\to
\coring{C}$ satisfying
$\psi'\circ\theta_{\coring{C}}=1_{\coring{C}}$. We have,
$\Delta_{\Lambda\cotensor{\coring{C}}X}=(\Lambda
\cotensor{\coring{C}}{\theta_{\coring{C}}}
\cotensor{\coring{C}}X)\circ
[\xymatrix{\Lambda\cotensor{\coring{C}}X \ar[r]^-\simeq &
\Lambda\cotensor{\coring{C}}{\coring{C}}\cotensor
{\coring{C}}X}]$.\\
Let $f=[\xymatrix{\Lambda\cotensor{\coring{C}}{\coring{C}}
\cotensor{\coring{C}}X\ar[r]^-\simeq &
\Lambda\cotensor{\coring{C}}X}]\circ
(\Lambda\cotensor{\coring{C}}{\psi'}\cotensor{\coring{C}}X)$. We
have, $f\circ\Delta_{\Lambda\cotensor{\coring{C}}X}=1_{\Lambda
\cotensor{\coring{C}}X}$,
\begin{multline*}
\begin{split}
 \Delta_{\Lambda\cotensor{\coring{C}}X}\circ f &
=\big[[\xymatrix{\Lambda\cotensor{\coring{C}}
{\coring{C}}\ar[r]^-\simeq &
\Lambda}]\circ(\Lambda\cotensor{\coring{C}}{\psi'})
\big]\cotensor{\coring{C}} \big[(\psi\cotensor{\coring{C}}X)
\circ[\xymatrix{X\ar[r]^-\simeq &
\coring{C}\cotensor{\coring{C}}X}]\big]\\ &
=(f\tensor{B}{\Lambda\cotensor{\coring{C}}X})\circ
(\Lambda\cotensor{\coring{C}}X\tensor{B}
{\Delta_{\Lambda\cotensor{\coring{C}}X}}),
\end{split}
\end{multline*}
and
\begin{multline*}
\begin{split}
 \Delta_{\Lambda\cotensor{\coring{C}}X}\circ f &
=\big[(\Lambda\cotensor{\coring{C}}{\psi})
\circ[\xymatrix{\Lambda\ar[r]^-\simeq
 & \Lambda\cotensor{\coring{C}}{\coring{C}}}]
 \big]\cotensor{\coring{C}}
\big[[\xymatrix{\coring{C}\cotensor{\coring{C}}X \ar[r]^-\simeq &
 X}]\circ(\psi'\cotensor{\coring{C}}X)\big] \\
 & =(\Lambda\cotensor{\coring{C}}X\tensor{B}{f})\circ
(\Delta_{\Lambda\cotensor{\coring{C}}X}\tensor{B}
{\Lambda\cotensor{\coring{C}}X}).
\end{split}
\end{multline*}
From the definition of a coseparable coring,
$\Lambda\cotensor{\coring{C}}X$ is a coseparable coring.

(3) We will prove it using Lemma \ref{26'}. Set
$T=(\Lambda\cotensor{\coring{C}}X)^*$,
$\Delta=\Delta_{\Lambda\cotensor{\coring{C}}X}$ and
$\epsilon=\epsilon_{\Lambda\cotensor{\coring{C}}X}$. At first from
Theorem \ref{leftrightFrobenius},
$(X\tensor{B}{-},\Lambda\cotensor{\coring{C}}{-})$ is a Frobenius
pair, and then $\Lambda\cotensor{\coring{C}}X\tensor{B}{-}$ is a
Frobenius functor. Hence $(\Lambda\cotensor{\coring{C}}X)_B$ is
finitely generated projective.\\
Let us consider the isomorphism
$$\xymatrix{H:T \ar[r]^-\simeq &
{\e{\coring{C}}{\Lambda}^*} \ar[r]^-{\phi_{\Lambda,B}}&
{\operatorname{End}_{\coring{C}}(\Lambda)} \ar[r]^-\simeq &
{\hom{\coring{C}}{B\tensor{B}{\Lambda}}{\Lambda}} \ar[r]^-\simeq &
{\hom{B}{B}{\Lambda\cotensor{\coring{C}}X}} \ar[r]^-\simeq &
{\Lambda\cotensor{\coring{C}}X}},$$ where
$\xymatrix{\phi_{\Lambda,B}:\e{\coring{C}}{\Lambda}^*\ar[r] &
{\operatorname{End}_{\coring{C}}(\Lambda)}}$ is the canonical
anti-isomorphism of rings defined in Proposition \ref{BW 23.8}.
Then, (by using the following consequence of Proposition
\ref{delta}; $\Delta=\big[(\delta_{\Lambda}\tensor{B}{\Lambda})
\theta_{\Lambda}\big]\cotensor{\coring{C}}X$), for $f\in T$,
\begin{multline}\label{frobeniuscomatrix}
H(f)=\Big[\big(\phi_{\Lambda,B}(f\delta_{\Lambda})
[\xymatrix{B\tensor{B}{\Lambda}\ar[r]^-\simeq
& {\Lambda}}]\big)\cotensor{\coring{C}}X\Big]\xi_B(1_B) \\
={[\xymatrix{B\tensor{B}{\Lambda}
\cotensor{\coring{C}}X\ar[r]^-\simeq &
{\Lambda\cotensor{\coring{C}}X}}]
(f\tensor{B}{\Lambda}\cotensor{\coring{C}}X)
\Delta[\xymatrix{B\tensor{B}{\Lambda}
\cotensor{\coring{C}}X\ar[r]^-\simeq &
{\Lambda\cotensor{\coring{C}}X}}]}\xi_B(1_B),
\end{multline}
where $\xi:1_{\mathcal{M}_B}\rightarrow-\tensor{B}
{\Lambda}\cotensor{\coring{C}}X$ is the unit of the adjunction
$(-\tensor{B}{\Lambda},-\cotensor{\coring{C}}X)$.

Now, let $f,f'\in T$, $\phi_{\Lambda,B}(f*^{r}f')=\phi_{\Lambda,B}
(f'\delta_{\Lambda}) \phi_{\Lambda,B}(f'\delta_{\Lambda})$.
\begin{multline*}
\begin{split}
H(f*^{r}f') & =\big(\phi_{\Lambda,B}(f'\delta_{\Lambda})
\cotensor{\coring{C}}X\big)\big(H(f)\big)
\\
& =[\xymatrix{B\tensor{B}{\Lambda}\cotensor{\coring{C}}
X\ar[r]^-\simeq &
{\Lambda\cotensor{\coring{C}}X}}](f'\tensor{B}{\Lambda}
\cotensor{\coring{C}}X)
\Delta\big(H(f)\big)\\
 \\
& =H(f).f'.
\end{split}
\end{multline*}
Finally, let $f\in T$, $b\in B$ and let $i_R:B\to T$ be the
anti-morphism of rings defined by $i_R(b)=\epsilon(b-)$ for $B\in B$
(see Proposition \ref{BW 17.8} (1)).
\begin{multline*}
\begin{split}
H(b.f) & =H(i_R(b)*^{r}f) &
=\bigg[\Big(\phi_{\Lambda,B}(f\delta_{\Lambda})
\phi_{\Lambda,B}\big(i_R(b)\delta_{\Lambda}\big) [\xymatrix{\Lambda
\ar[r]^-\simeq & {B\tensor{B}{\Lambda}}}]\Big)\cotensor{\coring{C}}
X\bigg]\xi_B(1_B).
\end{split}
\end{multline*}
On the other hand,
\begin{multline*}
\begin{split}
\phi_{\Lambda,B}(i_R(b)\delta_{\Lambda})\cotensor {\coring{C}}X & =
[\xymatrix{B\tensor{B}{\Lambda}\cotensor{\coring{C}}X\ar[r]^-\simeq
&
{\Lambda\cotensor{\coring{C}}X}}]\big(\epsilon(b-)\tensor{B}{\Lambda}
\cotensor{\coring{C}}X\big)
\Delta \\
& =b1_{\Lambda\cotensor{\coring{C}}X}.
\end{split}
\end{multline*}
Hence $H(b.f)=bH(f)$.
\end{proof}

\medskip

The last result is a generalization of \cite[Theorem 3.2(1), Theorem
3.5(1)]{Brzezinski/Gomez:2003}.

\begin{corollary}\label{comatrix}
Let $_BM_A$ be a $(B,A)$-bimodule such that $M_A$ is finitely
generated projective. Then \begin{enumerate}[(1)]
\item $_AM^*_B$ is
a separable bimodule if and only if the comatrix coring
$M^*\tensor{B}M$ is a cosplit $A$-coring.
\item If $M$ is a separable bimodule, then the
comatrix coring $M^*\tensor{B}M$ is a coseparable $A$-coring. \item
If $M$ is a Frobenius bimodule, then the comatrix coring
$M^*\tensor{B}M$ is a Frobenius $A$-coring.
\end{enumerate}
\end{corollary}

\begin{proof}
 Define $\psi:B\rightarrow M\tensor{A}M^*$ by
$b\mapsto
 \sum_{i=1}^{n}be_i\otimes
e_i^*=\sum_{i=1}^{n}e_i\otimes e_i^*b$, where
$\{e_i,e_i^*\}_{i\in\{1,\dots,n\}}$ is a dual basis basis of $M_A$,
and $\omega:M^*\tensor{B}M\rightarrow A$ by $\varphi\otimes
m\mapsto\varphi(m)$ (the evaluation map). Then
$(A,B,{_AM^*_B},{_BM_A},\omega,\psi)$ is a comatrix coring context
(see Example \ref{examples corings} (5)). From Proposition
\ref{19a}, we have the adjunction
$(-\tensor{B}{M},-\tensor{A}{M^*})$, in particular,
$\Lambda=M^*\in{_A\mathcal{M}_B}$ is quasi-finite as a right
$B$-(co)module, and the cohom functor is $-\tensor{B}{M}$.

(1) It is clear from the definition of a separable bimodule (see
Section \ref{category theory}), that $_AM^*_B$ is a separable
bimodule if and only if $\omega$ is an $A$-bimodule retraction.

(2) Let $S=\operatorname{End}_A(M_A)$. Define $\lambda:B\rightarrow
S$ by $b\mapsto\lambda_B:M\rightarrow M,[m\mapsto bm]$, and
$\phi:S\rightarrow M\tensor{A}{M^*}$ by $s\mapsto
\sum_is(e_i)\otimes e_i^*$, where
$\{e_i,e_i^*\}_{i\in\{1,\dots,n\}}$ is a dual basis of $M_A$. $\phi$
is an isomorphism of $B$-bimodule, with inverse map $m\otimes
\varphi\mapsto[x\mapsto m\varphi(x)]$. It is easy to see that
$\psi=\phi\circ\lambda$. From Sugano's Proposition \ref{Sugano}, the
bimodule $M$ is separable if and only if $\lambda$ is a split
extention, i.e., $\lambda$ is a $B$-bimodule section. However that
is equivalent to $\psi$ is a $B$-bimodule section.

(3) This is \cite[Theorem 3.7(1)]{Brzezinski/Gomez:2003}.
\end{proof}

\section{Applications to equivalences of categories of comodules}\label{equivalences}

In this section we will generalize and improve the main results
concerning equivalences between categories of comodules given in
\cite{Takeuchi:1977}, \cite{AlTakhman:2002} and
\cite{Brzezinski/Wisbauer:2003}. We also give new characterizations
of equivalences between categories of comodules over separable
corings or corings with a duality.

\medskip

We recall from \cite[p. 93]{MacLane:1998} that an \emph{adjoint
equivalence of categories} is an adjunction $(S,T,\eta,\varepsilon)$
in which both the unit $\eta:1\rightarrow TS$ and the counit
$\varepsilon:ST\rightarrow 1$ are natural isomorphisms.

\begin{proposition}\label{adjoint-equivalence}
Let $S:\cat{C}\to \cat{D}$ and $T:\cat{D}\to \cat{C}$ be functors.
Then $(S,T)$ is a pair of inverse equivalences if and only if $S$
and $T$ are part of an adjoint equivalence $(S,T,\eta,\varepsilon)$.
\end{proposition}

\begin{proof}
The proof is analogous to that of \cite[Theorem
IV.4.1]{MacLane:1998}. The ``if'' part is trivial. For the ``only
if'' part, let $\eta:1\rightarrow TS$ be a natural equivalence, then
$\varphi_{C,D}(\alpha)= T(\alpha).\eta_{C}$ for
$\alpha:S(C)\rightarrow D$, is a natural transformation in $C$ and
$D$. Since $T$ is faithful and full, $\varphi$ is a natural
equivalence, and $T$ is a right adjoint to $S$. Therefore the counit
of this adjunction $\varepsilon:ST\rightarrow 1$ satisfies
$T(\varepsilon_{D}).\eta_{T(D)}=1_{T(D)}$ for every $D\in\cat{D}$,
and consequently $T(\varepsilon_{D})=(\eta_{T(D)})^{-1}$ is
invertible. Since $T$ is faithful and full, $\varepsilon_{D}$ is
also invertible, and $\varepsilon$ is a natural equivalence (to show
that $\varepsilon$ is a natural equivalence we can also use
Proposition \ref{Maclane-adj-fun}).
\end{proof}

\medskip

The following result generalizes the case of the category of modules
\cite[II (2.4)]{Bass:1968}.

\begin{proposition}\label{equivalence1}
Suppose that $_A\coring{C}$, $\coring{C}_A$, $_B\coring{D}$ and
$\coring{D}_B$ are flat. Let
$X\in{}^{\coring{C}}\mathcal{M}^\coring{D}$ and
$\Lambda\in{}^{\coring{D}}\mathcal{M}^\coring{C}$.
\item The
following statements are equivalent:
\begin{enumerate}[(1)]
\item $(-\square
_{\coring{C}}X,-\square_{\coring{D}}\Lambda)$ is a pair of inverse
equivalences;
\item there exist bicomodule isomorphisms
\[
 f:X\square_{\coring{D}}\Lambda
\rightarrow\coring{C}\quad \text{and}\quad
g:\Lambda\square_{\coring{C}}X\rightarrow\coring{D}
\]
in $\bcomod{\coring{C}}{\coring{C}}$ and
$\bcomod{\coring{D}}{\coring{D}}$ respectively, such that
\begin{enumerate}[(a)] \item $_A X$ and $_B \Lambda$
are flat, and
$\omega_{X,\Lambda}=\rho{}_X\otimes_B\Lambda-X\otimes_A\rho{}_\Lambda$
is pure in $\lmod{A}$ and
$\omega_{\Lambda,X}=\rho{}_\Lambda\otimes_AX-\Lambda\otimes_B\rho{}_X$
is pure in $\lmod{B}$, or \item ${}_{\coring{C}}X$ and
${}_{\coring{D}}\Lambda$ are coflat.
\end{enumerate}
\end{enumerate}
In such a case the diagrams
\begin{equation}
\xymatrix{\Lambda\cotensor{\coring{C}}X\cotensor{\coring{D}}\Lambda
\ar[rr]^{\Lambda\cotensor{\coring{C}}f}
\ar[d]_{g\cotensor{\coring{D}}\Lambda} &&
\Lambda\cotensor{\coring{C}}\coring{C} \ar[d]^\simeq
\\\coring{D}\cotensor{\coring{D}}\Lambda
\ar[rr]^\simeq && \Lambda}
\qquad
\xymatrix{X\cotensor{\coring{D}}\Lambda\cotensor{\coring{C}}X
\ar[rr]^{f\cotensor{\coring{C}}X} \ar[d]_{X\cotensor{\coring{D}}g}&&
\coring{C}\cotensor{\coring{C}}X\ar[d]^\simeq
\\X\cotensor{\coring{D}}\coring{D} \ar[rr]^\simeq &&
X}
\end{equation}
commute.

If $A$ and $B$ are von Neumann regular rings, or if $\coring{C}$ and
$\coring{D}$ are coseparable corings (without $\coring{C}_A$ and
$\coring{D}_B$ are flat), the conditions (a) and (b) can be deleted.
\end{proposition}

\begin{proof}
Clear from Lemma \ref{18} and Propositions \ref{19a},
\ref{adjoint-equivalence}.
\end{proof}

\begin{definition}\label{7}
A bicomodule $N\in \bcomod{\coring{C}}{\coring{D}}$ is called a
\emph{cogenerator preserving} (resp. \emph{an injector-cogenerator})
as a right $\coring{D}$-comodule if the functor
$-\otimes_AN:\rmod{A}\to\rcomod{\coring{D}}$ preserves cogenerator
(resp. injective cogenerator) objects.
\end{definition}

The version for module categories of the first part of the following
lemma is given in \cite[Exercise 20.8]{Anderson/Fuller:1992}.

\begin{lemma}\label{injector-cogenerator}
Let $\cat{C}$ and $\cat{D}$ be two abelian categories and
$S:\cat{C}\rightarrow \cat{D}$ a functor. Let
$T:\cat{D}\rightarrow \cat{C}$ be a right adjoint of $S$. We have
the following properties:
\begin{enumerate}[(1)]
\item If $N\in\cat{D}$ is cogenerator and $T(N)\in\cat{C}$ is
injective, then $S$ is exact.
\item If $N\in\cat{D}$ is cogenerator and $S$ is
faithful, then $T(N)\in\cat{C}$ is a cogenerator. \item If
$T(N)\in\cat{C}$ is a cogenerator for some $N\in\cat{D}$, then $S$
is faithful.

If moreover $\cat{D}$ is an AB 3 category with generators and has
enough injectives (e.g. if $\cat{D}$ is a Grothendieck category),
then the following properties hold: \item $S$ is faithful if and
only if $T$ preserves cogenerator objects.\item $S$ is faithfully
exact if and only if $T$ preserves injective cogenerator objects.
\end{enumerate}
\end{lemma}

\begin{proof}
Let $N\in\cat{D}$. Then we have a natural isomorphism
$$\eta_{-,N}:\hom{\cat{D}}{S(-)}{N}\to \hom{\cat{C}}{-}{T(N)}.$$
The functor $\hom{\cat{D}}{S(-)}{N}$ is the composition of covariant
functors
$$\xymatrix{\cat{C}^\circ\ar[r]^{S^\circ}& \cat{D}^\circ
\ar[rr]^{\hom{\cat{D}}{-}{N}}&& \textbf{Ab}},$$ where $S^\circ$ is
the dual functor of $S$. Note that $S^\circ$ is exact (resp. faithful)
if and only if $S$ is exact (resp. faithful).

(1) By assumptions, $\hom{\cat{D}}{-}{T(N)}$ is exact and
$\hom{\cat{C}}{-}{N}$ is faithful. By Theorem \ref{faithful},
$\hom{\cat{C}}{-}{N}$ reflects exact sequences. Then $S^\circ$ is
exact. Hence $S$ is so.

(2) Follows directly from the fact that the composition of two
faithful functors is also faithful.

(3) Follows directly from the fact that if the composition of two
functors $GF$ is faithful then $F$ is also faithful.

(4) is an immediate consequence of (2), (3) and the fact that
$\cat{D}$ has a cogenerator object (see Proposition
\ref{inj.cogen.obj}).

 (5) The ``only if'' part is obvious from Lemma
\ref{right-adjoint-inj-obj} and (4). The ``if'' part is an immediate
consequence of Proposition \ref{inj.cogen.obj} and (1) and (3).
\end{proof}

\begin{corollary}
Let $N\in\bcomod{\coring{C}}{\coring{D}}$ be a bicomodule,
quasi-finite as a right $\coring{D}$-comodule, such that
$_A\coring{C}$ and $_B\coring{D}$ are flat. $N$ is a cogenerator
preserving (resp. an injector-cogenerator) as a right
$\coring{D}$-comodule if and only if the cohom functor
$h_{\coring{D}}(N,-)$ is faithful (resp. faithfully exact).
\end{corollary}

\medskip
The first part of the following is
\cite[23.10]{Brzezinski/Wisbauer:2003}. We think that the proof we
give here is more clear.

\begin{proposition}\label{equivalence2}
Let $_A\coring{C}$ be flat and let
$\Lambda\in{}_B\mathcal{M}^\coring{C}$ be quasi-finite as a right
$\coring{C}$-comodule. Let $\cohom{\coring{C}}{\Lambda}{-}$ and
$\e{\coring{C}}{\Lambda}$ be the cohom functor and the
coendomorphism coring of $\Lambda$ respectively. Suppose that
$\e{\coring{C}}{\Lambda}$ is flat as left $B$-module.
\begin{enumerate}[(1)] \item If $\coring{C}_A$ is flat
and $\e{\coring{C}}{\Lambda}$ is flat as right $B$-module, and
$\Lambda$ is an injector-cogenerator as a right
$\coring{C}$-comodule, then the functors
\begin{equation}\label{functors}
-\square_{\e{\coring{C}}{\Lambda}}\Lambda:
\mathcal{M}^{\e{\coring{C}}{\Lambda}}\to
\mathcal{M}^\coring{C},\qquad
\cohom{\coring{C}}{\Lambda}{-}:\mathcal{M}^\coring{C}\rightarrow
\mathcal{M}^{\e{\coring{C}}{\Lambda}},
\end{equation}
are inverse equivalences.
\item If $\coring{C}$ and
$\e{\coring{C}}{\Lambda}$ are coseparable corings, and $\Lambda$ is
a cogenerator preserving as a right $\coring{C}$-comodule, then the
functors of \eqref{functors} are inverse equivalences.
\end{enumerate}
\end{proposition}

\begin{proof}
We will prove at the same time the two statements. Let
$\theta:1_{\mathcal{M}^\coring{C}}\rightarrow
\cohom{\coring{C}}{\Lambda}{-}\tensor{B}{\Lambda}$ be the unit of
the adjunction
$(\cohom{\coring{C}}{\Lambda}{-},-\tensor{B}{\Lambda})$. By Theorem
\ref{3}, $\cohom{\coring{C}}{\Lambda}{-}\simeq
-\square_{\coring{C}}\cohom{\coring{C}}{\Lambda}{\coring{C}}$. In
the case (1), $\cohom{\coring{C}}{\Lambda}{\coring{C}}$ is coflat as
a left $\coring{C}$-comodule. We put
$\coring{D}=\e{\coring{C}}{\Lambda}$. By Proposition \ref{G4},
$\delta_{\Lambda}:\coring{D}\rightarrow
\Lambda\square_{\coring{C}}\cohom{\coring{C}} {\Lambda}{\coring{C}}$
is $\coring{D}$-bicolinear. Therefore,
$$-\square_{\coring{D}}\Lambda\square_{\coring{C}}
\cohom{\coring{C}}{\Lambda}{\coring{C}}
\simeq-\square_{\coring{D}}\coring{D}$$ is exact. Hence, in the case
(1), the functor $-\square_{\coring{D}}\Lambda$ is exact (since
$\cohom{\coring{C}}{\Lambda}{-}$ is faithful). From the proof of
Proposition \ref{comatrixproperty1}(1), the unit of the adjunction
$(-\square_{\coring{C}}\cohom{\coring{C}}{\Lambda}{\coring{C}},
-\square_{\coring{D}}\Lambda)$ is
$$\xymatrix@1{\eta
:1_{\mathcal{M}^\coring{C}}\ar[r]^\simeq & -\square
_{\coring{C}}\coring{C}
\ar[rr]^-{-\square_{\coring{C}}\theta_{\coring{C}}} & &
-\square_{\coring{C}}\cohom{\coring{C}}{\Lambda}{\coring{C}}
\square_{\coring{D}}\Lambda}.$$ From Proposition \ref{delta},
$\delta_{\Lambda}$ is the unique map making the following diagram
commutative
\[
\xymatrix{\Lambda \ar[d]_\simeq\ar[rr]^\simeq & &
\Lambda\square_{\coring{C}}\coring{C}\ar[d]^{\Lambda
\square_{\coring{C}}\theta_{\coring{C}}}
\\
\coring{D}\square_{\coring{D}}\Lambda
\ar[rr]_-{\delta_{\Lambda}\square_{\coring{D}}\Lambda} & &
\Lambda\square_{\coring{C}}(\cohom{\coring{C}}{\Lambda}
{\coring{C}}\square_{\coring{D}}\Lambda).}
\]
By Lemma \ref{18}, Poposition \ref{19a}, the counit of the
adjunction
$(-\square_{\coring{C}}\cohom{\coring{C}}{\Lambda}{\coring{C}},
-\square_{\coring{D}}\Lambda)$ is
$$\xymatrix@1{\varepsilon:-\square_{\coring{D}}
\Lambda\square_{\coring{C}} \cohom{\coring{C}}{\Lambda}{\coring{C}}
\ar[rr]^-{-\square_{\coring{D}}\delta^{-1}_{\Lambda}} & &
-\square_{\coring{D}}\coring{D}\ar[r]^\simeq &
1_{\mathcal{M}^\coring{D}},} $$ and we have the commutative diagram
\[
\xymatrix{\cohom{\coring{C}}{\Lambda}{\coring{C}}
\ar[d]_\simeq\ar[rr]^\simeq & &
\coring{C}\square_{\coring{C}}\cohom{\coring{C}}{\Lambda}{\coring{C}}
\ar[d]^{\theta_{\coring{C}}\square_{\coring{C}}\cohom{\coring{C}}
{\Lambda}{\coring{C}}}
\\
\cohom{\coring{C}}{\Lambda}{\coring{C}}\square_{\coring{D}}\coring{D}
& & \cohom{\coring{C}}{\Lambda}{\coring{C}}\square_{\coring{D}}
\Lambda\square_{\coring{C}} \cohom{\coring{C}}{\Lambda}{\coring{C}}
\ar[ll]_-{\cohom{\coring{C}}{\Lambda}{\coring{C}}\square_{\coring{D}}
\delta^{-1}_{\Lambda}}.}
\]
Then, $\theta_{\coring{C}}\square_{\coring{C}}
\cohom{\coring{C}}{\Lambda}{\coring{C}}$ is an isomorphism. Since
$-\square_{\coring{C}}\cohom{\coring{C}}{\Lambda}{\coring{C}}$ is
faithful, $\theta_{\coring{C}}$ is also an isomorphism. Finally, the
unit and the counit of the adjunction
$(-\square_{\coring{C}}\cohom{\coring{C}}{\Lambda}{\coring{C}},
-\square_{\coring{D}}\Lambda)$ are natural isomorphisms.
\end{proof}

\medskip

The first part of the following is contained in \cite[Theorem
3.10]{ElKaoutit/Gomez:2003}.

\begin{corollary}\label{KG 3.10}
Let ${}_BM_A$ be a $(B,A)$-bimodule such that $M_A$ is finitely
generated projective. \begin{enumerate}[(1)] \item The following
statements are equivalent \begin{enumerate}[(a)]\item
${}_A(M^*\otimes_BM)$ is flat and
$-\otimes_BM:\mathcal{M}_B\rightarrow\mathcal{M}^{M^*\otimes_BM}$ is
an equivalence of categories; \item $_B{}M$ is faithfully flat.
\end{enumerate}
\item If $\mathcal{M}^{M^*\otimes_BM}$ is an abelian
category, and $M^*\otimes_BM$ is an $A$-coseparable coring, then the
following statements are equivalent
\begin{enumerate}[(a)]
\item $-\otimes_BM:\mathcal{M}_B\to\mathcal{M}^{M^*\otimes_BM}$ is an
equivalence of categories;
\item $_BM$ is completely faithful.
\end{enumerate}
\end{enumerate}
\end{corollary}

\begin{proof}
It suffices to take $\Lambda=M^*\in{}_A\mathcal{M}_B$ in Proposition
\ref{equivalence2}.
\end{proof}

The first part of the following is contained in
\cite[23.12]{Brzezinski/Wisbauer:2003}. We think that the proof we
give here is more clear.

\begin{proposition}\label{equivalence3}
Suppose that $_A\coring{C}$ and $_B\coring{D}$ are flat, and let
$X\in{}^{\coring{C}}\mathcal{M}^\coring{D}$ and
$\Lambda\in{}^{\coring{D}}\mathcal{M}^\coring{C}$.

If $\coring{C}_A$ and $\coring{D}_B$ are flat, then the following
statements are equivalent
\begin{enumerate}[(1)]
\item $(-\square_{\coring{C}}X,-\square_{\coring{D}}\Lambda)$ is a
pair of inverse equivalences;
\item $\Lambda$ is quasi-finite injector-cogenerator as a right
$\coring{C}$-comodule, $\e{\coring{C}}{\Lambda}\simeq \coring{D}$ as
corings and $X\simeq \cohom{\coring{C}}{\Lambda}{\coring{C}}$ in
$^{\coring{C}}\mathcal{M}^\coring{D}$;
\item $X$ is quasi-finite injector-cogenerator as a right
$\coring{D}$-comodule, $\e{\coring{D}}{X}\simeq \coring{C}$ as
corings and $\Lambda\simeq \cohom{\coring{D}}{X}{\coring{D}}$ in
$^{\coring{D}}\mathcal{M}^\coring{C}$.

If moreover $\coring{C}$ and $\coring{D}$ are coseparable or
cosemisimple, then (1) is equivalent to
\item $\Lambda$ is quasi-finite cogenerator preserving
as a right
$\coring{C}%
$-comodule, $\e{\coring{C}}{\Lambda}\simeq \coring{D}$ as corings
and $X\simeq \cohom{\coring{C}}{\Lambda}{\coring{C}}$ in
$^{\coring{C}}\mathcal{M}^\coring{D}$;\item $X$ is quasi-finite
cogenerator preserving as a right $\coring{D}%
$-comodule, $\e{\coring{D}}{X}\simeq \coring{C}$ as corings and
$\Lambda\simeq \cohom{\coring{D}}{X}{\coring{D}}$ in
$^{\coring{D}}\mathcal{M}^\coring{C}$.
\end{enumerate}
\end{proposition}

\begin{proof}
At first we will prove the first part.

$(1)\Rightarrow(2)$ By Proposition \ref{19}, $\Lambda$ is
quasi-finite injector as a right $\coring{C}$-comodule, and $X\simeq
\cohom{\coring{C}}{\Lambda}{\coring{C}}$ in
$^{\coring{C}}\mathcal{M}^\coring{D}$. Therefore, $-\square
_{\coring{C}}X\simeq \cohom{\coring{C}}{\Lambda}{-}$ is faithful. By
Lemma \ref{injector-cogenerator}, $\Lambda$ is quasi-finite
injector-cogenerator as a right $\coring{C}$-comodule. Finally, by
Proposition \ref{comatrixproperty1}, $\e{\coring{C}}{\Lambda}\simeq
\coring{D}$ as corings.

$(2)\Rightarrow(1)$ By Proposition \ref{19}, $-\square
_{\coring{C}}X\simeq \cohom{\coring{C}}{\Lambda}{-}$. Hence (1)
follows obviously from Proposition \ref{equivalence2}(1).

$(1)\Leftrightarrow(3)$ Follows by symmetry.

 Now we will prove the second part. The case ``cosemisimple''
is obvious from the first part. It suffices to show the
``coseparable'' case.

$(1)\Rightarrow(4)$ Obvious from the first part.

$(4)\Rightarrow(1)$ By Proposition \ref{19a}, $-\square
_{\coring{C}}X\simeq \cohom{\coring{C}}{\Lambda}{-}$. Hence (1)
follows obviously from Proposition \ref{equivalence2}(2).

$(1)\Leftrightarrow(5)$ Follows by symmetry.
\end{proof}

\medskip

As an immediate consequence of Proposition \ref{19a} and Proposition
\ref{equivalence3}, we get the two following theorems.

\begin{theorem}\label{equivalence4}
Suppose that $_A\coring{C}$, $\coring{C}_A$, $_B\coring{D}$ and
$\coring{D}_B$ are flat, and let
$X\in{}^{\coring{C}}\mathcal{M}^{\coring{D}}$ and
$\Lambda\in{}^{\coring{D}}\mathcal{M}^\coring{C}$. The following
statements are equivalent:
\begin{enumerate}[(1)]
\item $(-\square_{\coring{C}}X,-\square_{\coring{D}}\Lambda)$ is a
pair of inverse equivalences with $X_{\coring{D}}$ and
$\Lambda_{\coring{C}}$ coflat;
\item $(\Lambda\cotensor{\coring{C}}-,X \cotensor{\coring{D}}-)$ is
a pair of inverse equivalences with $_{\coring{C}}X$ and
$_{\coring{D}}\Lambda$ coflat;
\item $\Lambda$ is quasi-finite injector-cogenerator as a right
$\coring{C}$-comodule with $X_{\coring{D}}$ and
$\Lambda_{\coring{C}}$ coflat, $\e{\coring{C}}{\Lambda}\simeq
\coring{D}$ as corings and $X\simeq
\cohom{\coring{C}}{\Lambda}{\coring{C}}$ in
$^{\coring{C}}\mathcal{M}^\coring{D}$;
\item $X$ is quasi-finite injector-cogenerator as a right
$\coring{D}$-comodule with $X_{\coring{D}}$ and
$\Lambda_{\coring{C}}$ coflat, $\e{\coring{D}}{X}\simeq \coring{C}$
as corings and $\Lambda\simeq \cohom{\coring{D}}{X}{\coring{D}}$ in
$^{\coring{D}}\mathcal{M}^\coring{C}$;
\item $\Lambda$ is quasi-finite injector-cogenerator coflat on both
sides, $\e{\coring{C}}{\Lambda}\simeq \coring{D}$ as corings, and
$X\simeq \cohom{\coring{C}}{\Lambda}{\coring{C}}$ in
$^{\coring{C}}\mathcal{M}^\coring{D}$;\item $X$ is quasi-finite
injector-cogenerator coflat on both sides, $\e{\coring{D}}{X}\simeq
\coring{C}$ as corings, and $\Lambda\simeq
\cohom{\coring{D}}{X}{\coring{D}}$ in
$^{\coring{D}}\mathcal{M}^\coring{C}.$
\end{enumerate}
\end{theorem}

\begin{theorem}\label{equivalence5}
Suppose that $_A\coring{C}$, $\coring{C}_A$, $_B\coring{D}$ and
$\coring{D}_B$ are flat, and let
$X\in{}^{\coring{C}}\mathcal{M}%
^{\coring{D}}$ and
$\Lambda\in{}^{\coring{D}}\mathcal{M}^\coring{C}$. If $\coring{C}$
and $\coring{D}$ are coseparable or cosemisimple (resp. $A$ and $B$
are von Neumann regular ring), then the following statements are
equivalent:
\begin{enumerate}[(1)] \item $(-\square
_{\coring{C}}X,-\square_{\coring{D}}\Lambda)$ is a pair of inverse
equivalences; \item $(\Lambda \cotensor{\coring{C}}-,X
\cotensor{\coring{D}}-)$ is a pair of inverse equivalences; \item
$\Lambda$ is quasi-finite cogenerator preserving (resp.
injector-cogenerator) as a right
$\coring{C}%
$-comodule, $\e{\coring{C}}{\Lambda}\simeq \coring{D}$ as corings
and $X\simeq \cohom{\coring{C}}{\Lambda}{\coring{C}}$ in
$^{\coring{C}}\mathcal{M}^\coring{D}$;\item $X$ is quasi-finite
cogenerator preserving (resp. injector-cogenerator) as a right
$\coring{D}%
$-comodule, $\e{\coring{D}}{X}\simeq \coring{C}$ as corings and
$\Lambda\simeq \cohom{\coring{D}}{X}{\coring{D}}$ in
$^{\coring{D}}\mathcal{M}^\coring{C}$;\item $\Lambda$ is
quasi-finite cogenerator preserving (resp. injector-cogenerator) on
both sides, $\e{\coring{C}}{\Lambda}\simeq \coring{D}$ as corings,
and $X\simeq \cohom{\coring{C}}{\Lambda}{\coring{C}}$ in
$^{\coring{C}}\mathcal{M}^\coring{D}$;\item $X$ is quasi-finite
cogenerator preserving (resp. injector-cogenerator) on both sides,
$\e{\coring{D}}{X}\simeq \coring{C}$ as corings, and $\Lambda\simeq
\cohom{\coring{D}}{X}{\coring{D}}$ in
$^{\coring{D}}\mathcal{M}^\coring{C}.$
\end{enumerate}
\end{theorem}

The following result, which recovers Morita's characterization of
equivalence \cite[Theorem 22.2]{Anderson/Fuller:1992}, is a relevant
consequence of the last theorem.

\begin{corollary}\label{Morita}
Let $A$ and $B$ be $k$-algebras and let
$$F:\rmod{A}\to \rmod{B}\qquad\textrm{and
}\qquad G:\rmod{B}\to \rmod{A}$$ be $k$-linear functors. Then the
following statements are equivalent
\begin{enumerate}[(1)]\item $F$ and $G$ are inverse
equivalences;
\item there exists a bimodule $_AM_B$ such that:
\begin{enumerate}[(a)]\item $_AM$ and $M_B$ are
progenerators (i.e., finitely generated, projective and generators),
\item $_AM_B$
is (faithfully) balanced (see \cite[p. 60]{Anderson/Fuller:1992}),
\item $F\simeq -\tensor{A}{M}\qquad\textrm{and}\qquad
G\simeq \hom{B}{M}{-}$;
\end{enumerate}
\item there exists a bimodule $_AM_B$ such that: $M_B$
is finitely generated projective, $_AM$ is completely faithful, the
evaluation map $M^*\tensor{A}{M}\rightarrow B$ is an isomorphism,
and
\\ $F\simeq -\tensor{A}{M}\qquad\textrm{and}\qquad
G\simeq \hom{B}{M}{-}$.
\end{enumerate}
Moreover, in such a case, $_BM^*$ and $M^*_A$ are progenerators, and
$$F\simeq \hom{A}{M^*}{-}\qquad\textrm{and}\qquad
G\simeq -\tensor{B}{M^*}.$$
\end{corollary}

\begin{proof}
$(1)\Rightarrow(2)$ There exists a bimodule $_AM_B$ such that
$F\simeq -\tensor{A}{M}$ and $G\simeq \hom{B}{M}{-}$. Since $G$ is
an equivalence, $G\simeq -\tensor{B}{M^*}$ (by Eilenberg-Watts
Theorem or by Theorem \ref{3}). Then $M_B$ is finitely generated
and projective. From $G$ faithful, it follows that $M_B$ is a
progenerator. By symmetry, it follows that $_AM$, $_BM^*$ and
$M^*_A$ are also progenerators.

As in the proof of Corollary \ref{comatrix},
$(B,A,{_BM^*_A},{_AM_B},\omega,\psi)$ is a comatrix coring context,
and $\psi=\phi\circ\lambda$, where $\omega$ , $\psi$, $\phi$ are
defined as in that proof, and
$\lambda:A\to\operatorname{End}_B(M_B)$ is defined by
$\lambda(a)=\lambda_a$ with $\lambda_a(m)=am$. In particular we have
an adjunction $(-\tensor{A}{M},-\tensor{B}{M^*})$. Using Proposition
\ref{M IV.3} and since $\phi$ is an isomorphism, we obtain that $F$
is fully faithful if and only if $\psi$ is an isomorphism, if and
only if $\lambda$ is an isomorphism. Then, by symmetry, $_AM_B$ is
faithfully balanced.

 $(2)\Rightarrow(3)$ From Lemma \ref{completely faithful} (2), it
 follows that $_AM$ is completely faithful. In particular it is
 faithful. By symmetry $_AM_B$ is faithfully balanced. Now, we
 will prove that the coring homomorphism
$\epsilon_{M^*\tensor{A}{M}}:M^*\tensor{A}{M}\rightarrow B$ is an
isomorphism. Let $N\in\rmod{B}$ and let us consider the map
$$\xymatrix{\xi_N:\hom{B}{M}{N}\tensor{A}{M}
\ar[r]^-\nu & \hom{B}{\hom{A}{M}{M}}{N} \ar[r]^-\simeq &
\hom{B}{B}{N} \ar[r]^-\simeq & N},$$ where $\nu$
is defined by
$\nu(\gamma\tensor{A}{m}):\alpha\mapsto\gamma(\alpha(m)).$ From
Lemma \ref{hom-tens-relations} (7) ($_AM$ is finitely
generated projective), $\nu$ and $\xi_N$ are isomorphisms. It is
easy to see that $\xi_N(\gamma\tensor{A}{m})=\gamma(m)$, for
every $\gamma$ in $\hom{B}{M}{N}$ and every $m$ in $M$. In
particular ($N=B$), $\xi_B=\epsilon_{M^*\tensor{A}M}$.

$(3)\Rightarrow(1)$ Since $M_B$ is finitely generated projective,
$\Lambda={}_BM^*_A$ is $(B,A)$-quasi-finite, and
the cohom functor is $F\simeq -\tensor{A}M$. $F$ is faithful
since $_AM$ is completely faithful. Moreover the coring
homomorphism $\epsilon_{M^*\tensor{A}M}:M^*\tensor{A}M\to B$ is
an isomorphism. The implication $(3)\Rightarrow(1)$ of Theorem
\ref{equivalence5} then achieves the proof.
\end{proof}

\begin{corollary}
Let $A$ and $B$ be $k$-algebras. $\rmod{A}$ and $\rmod{B}$ are
equivalent if and only if $\lmod{A}$ and $\lmod{B}$ are equivalent.
In such a case we say that $A$ and $B$ are \emph{Morita equivalent}.
\end{corollary}

\begin{proof}
Obvious from Theorem \ref {equivalence5}.
\end{proof}

\medskip
Now, we will study when the category $\rcomod{\coring{C}}$ for some
coring $\coring{C}$ is equivalent to a category of modules
$\rmod{B}$.

Let $\Sigma\in {}_B\rcomod{\coring{C}}$ such that $\Sigma_A$ is
finitely generated projective with a dual basis $\{e_i,e_i^*\}_i$,
and let $M\in \rcomod{\coring{C}}$. By \cite[Proposition
1.4]{Caenepeel/DeGroot/Vercruysse:unp}, the canonical isomorphism
$\hom{A}{\Sigma}{M}\to M\tensor{A}\rdual{\Sigma}$ yields an
isomorphism $\hom{\coring{C}}{\Sigma}{M}\to
M\cotensor{\coring{C}}\rdual{\Sigma}$, where the left coaction on
$\Sigma^*$ is given by
$$ \lambda_{\Sigma^*}:\Sigma^*\to
\coring{C}\tensor{A}\Sigma^*,\; \lambda_{\Sigma^*}(f)=\sum_i
f(e_{i(0)})e_{i(1)}\tensor{A}e_i^*.$$

Moreover, we have an adjoint pair $(F,G)$, where
$$F=-\tensor{B}\Sigma:\rmod{B}\to
\rcomod{\coring{C}},$$ and
$$G=\hom{\coring{C}}{\Sigma}{-}:\rcomod{\coring{C}}\to
\rmod{B}.$$ The unit and the counit of this adjunction are given by:
For $N\in\rmod{B}$,
$$\upsilon_N:N\to
\hom{\coring{C}}{\Sigma}{N\tensor{B}\Sigma},\;
\upsilon_N(n)(u)=n\tensor{B}u,$$ or
$$\upsilon_N:N\to
(N\tensor{B}\Sigma)\cotensor{\coring{C}}\rdual{\Sigma},\;
\upsilon_N(n)=\sum_i(n\tensor{B}e_i)\tensor{A}e_i^*,$$ and for $M\in
\rcomod{\coring{C}}$:
$$\zeta_M:\hom{\coring{C}}{\Sigma}{M}\tensor{B}\Sigma\to
M,\;\zeta_M(\varphi\tensor{B}u)=\varphi(u),$$ or
$$\zeta_M:(M\cotensor{\coring{C}}\rdual{\Sigma})\tensor{B}\Sigma\to
M, \;
\zeta_M\Big((\sum_jm_j\tensor{A}f_j)\tensor{B}u\Big)=\sum_jm_jf_j(u)$$
(see \cite[Proposition 1.5]{Caenepeel/DeGroot/Vercruysse:unp}).

Let us define the map
$$\textsf{can}:\rdual{\Sigma}\tensor{B}\Sigma\to
\coring{C},\; \textsf{can}(f\tensor{B}u)=f(u_{(0)})u_{(1)}.$$ From
\cite[Lemma 3.1]{Caenepeel/DeGroot/Vercruysse:unp}, we have
$\textsf{can}$ is a morphism of corings. It follows from this that
$\textsf{can}$ is a $\coring{C}$-bicolinear map. Moreover, we can
verify easily that for every $M\in \rcomod{\coring{C}}$, the
following diagram is commutative:
\begin{equation}\label{counitmod}
\xymatrix{(M\cotensor{\coring{C}}\rdual{\Sigma})\tensor{B}\Sigma
\ar[rr]^-{\zeta_M} \ar[d]^{\psi_M} & & M\ar[d]^\simeq
\\
M\cotensor{\coring{C}}(\rdual{\Sigma}\tensor{B}\Sigma)\ar[rr]^-
{M\cotensor{\coring{C}}\textsf{can}} & &
M\cotensor{\coring{C}}\coring{C},}
\end{equation}
where $\psi_M$ is the canonical map. We obtain that if $\psi_M$ is
isomorphism, for every $M\in \rcomod{\coring{C}}$ (for example if
$_B\Sigma$ is flat or if $\coring{C}$ is coseparable), then $G$ is
fully faithful if and only if $\textsf{can}$ is an isomorphism (see
Proposition \ref{Maclane-adj-fun}).

Let us consider the map $\upsilon:B\to
\Sigma\cotensor{\coring{C}}\rdual{\Sigma},\;
\upsilon(b)=\sum_ibe_i\tensor{A}e_i^*$ $(b\in B).$ For every $N\in
\rmod{B}$, we have the following commutative diagram
\begin{equation}\label{unitmod}
\xymatrix{N \ar[rr]^-{\upsilon_N}\ar[d]_\simeq & &
(N\tensor{B}\Sigma)\cotensor{\coring{C}}\rdual{\Sigma}
\\ N\tensor{B}B\ar[rr]^-{N\tensor{B}\upsilon}
& &
N\tensor{B}(\Sigma\cotensor{\coring{C}}\rdual{\Sigma})\ar[u]_{\psi_N},}
\end{equation}
where $\psi_N$ is the canonical map. Let $\phi:B\to
\operatorname{End}_\coring{C}(\Sigma)$ be the canonical morphism of
$k$-algebras which define the left action on $\Sigma$
($\phi(b)(u):=bu$) (see Section \ref{Frobeniusgeneral}). We have
$\upsilon:\xymatrix{B\ar[r]^-\phi &
\operatorname{End}_\coring{C}(\Sigma)\ar[r]^-\simeq &
\Sigma\cotensor{\coring{C}}\rdual{\Sigma}}$. We obtain that if
$\psi_N$ is an isomorphism, for every $N\in \rmod{B}$ (for example
if $\Sigma_\coring{C}$ is projective or if $\coring{C}$ is
coseparable), then $F$ is fully faithful if and only if $\phi$ is an
isomorphism (see Proposition \ref{Maclane-adj-fun}).

Furthermore, if $_B\Sigma$ is flat or $\coring{C}$ is coseparable,
then we have the commutative diagram:
\begin{equation}\label{unitcounitmod}
\xymatrix{\Sigma\ar[rr]^\simeq \ar[d]_\simeq & &
B\tensor{B}\Sigma\ar[d]^{\upsilon\tensor{B}\Sigma}\\
\Sigma\cotensor{\coring{C}}\coring{C} & & \Sigma
\cotensor{\coring{C}}\Sigma^*\tensor{B}\Sigma
\ar[ll]^{\Sigma\cotensor{\coring{C}}\textsf{can}}.}
\end{equation}

Now we are ready to state and prove the following theorem. In the
particular case where $B=\operatorname{End}_\coring{C}(\Sigma)$, the
first part of Theorem is known, see \cite[Theorem
3.2]{ElKaoutit/Gomez:2003} and
\cite[18.27]{Brzezinski/Wisbauer:2003}.

\begin{theorem}\label{equivmodcomod}
Let $\Sigma\in {}_B\rcomod{\coring{C}}$ such that $\Sigma_A$ is
finitely generated projective. Let $F=-\tensor{B}\Sigma:\rmod{B}\to
\rcomod{\coring{C}}$ and
$G=\hom{\coring{C}}{\Sigma}{-}:\rcomod{\coring{C}}\to \rmod{B}$. The
following statements are equivalent:
\begin{enumerate}[(1)]
\item$(F,G)$ is a pair of inverse equivalences with
$_A\coring{C}$ flat;
\item $_B\Sigma$ is flat, $\Sigma_\coring{C}$ is
projective, and $\textsf{can}$, and the morphism of $k$-algebras
$\phi:B\to \operatorname{End}_\coring{C}(\Sigma)$ defined by
$\phi(b)(u)=bu$ for $b\in B,u\in \Sigma$, are isomorphisms;
\item $_A\coring{C}$ is flat, $\Sigma_\coring{C}$ is
a projective generator, and the morphism of $k$-algebras $\phi:B\to
\operatorname{End}_\coring{C}(\Sigma)$ defined by $\phi(b)(u)=bu$
for $b\in B,u\in \Sigma$, is an isomorphism;
\item $_B\Sigma$ is faithfully flat, and
$\textsf{can}$ is an isomorphism.
\end{enumerate}
If moreover $\coring{C}$ is coseparable, and $_A\coring{C}$ is
projective, then the following statements are equivalent:
\begin{enumerate}[(1)]\item $(F,G)$ is a pair of
inverse equivalences;
\item $\Sigma_\coring{C}$ is a generator, and the
morphism of $k$-algebras $\phi:B\to
\operatorname{End}_\coring{C}(\Sigma)$ defined by $\phi(b)(u)=bu$
for $b\in B,u\in \Sigma$, is an isomorphism;
\item the morphism of $k$-algebras $\phi:B\to
\operatorname{End}_\coring{C}(\Sigma)$ defined by $\phi(b)(u)=bu$
for $b\in B,u\in \Sigma$, is an isomorphism, and $\textsf{can}$ is a
surjective map;
\item $_B\Sigma$ is completely faithful, and
$\textsf{can}$ is a bijective map.
\end{enumerate}
\end{theorem}

\begin{proof}
First we will prove the first statement. It is obvious that the
condition (1) implies the other conditions.

$(2)\Rightarrow(1)$ That $_A\coring{C}$ is flat follows from
$_B\Sigma$ is flat and $\textsf{can}$ is an isomorphism of
$A$-bimodules. To prove that $(F,G)$ is a pair of equivalences, it
is enough to use the commutativity of the diagrams \eqref{counitmod}
and \eqref{unitmod}.

$(3)\Rightarrow(1)$ Follows from the Gabriel-Popescu Theorem
\ref{Gabriel-Popescu}, and the commutativity of the diagram
\eqref{unitmod}.

$(4)\Rightarrow(1)$ From $\textsf{can}$ is an isomorphism, it
follows that $-\cotensor{\coring{C}}(\Sigma^*\tensor{B}\Sigma)\simeq
-\cotensor{\coring{C}}\coring{C}$, and the following diagram is
commutative and each of its morphisms is an isomorphism:
$$
\xymatrix{N\tensor{B}(\Sigma\cotensor{\coring{C}}
\Sigma^*)\tensor{B}\Sigma \ar[rr]^{\psi_N\tensor{B}\Sigma}
\ar[d]^\simeq &&
\big((N\tensor{B}\Sigma)\cotensor{\coring{C}}\Sigma^*\big)
\tensor{B}\Sigma
\ar[d]^\simeq\\
N\tensor{B}\big(\Sigma\cotensor{\coring{C}}(\Sigma^*\tensor{B}
\Sigma)\big)\ar[rr]^\simeq &&
(N\tensor{B}\Sigma)\cotensor{\coring{C}}(\Sigma^*\tensor{B}\Sigma).}
$$
Then $\psi_N$ is an isomorphism (since $-\tensor{B}\Sigma$ is
faithful). From the commutativity of the diagram
\eqref{unitcounitmod}, $\textsf{can}$ is an isomorphism, and
$-\tensor{B}\Sigma$ is faithful, we have $\upsilon$ is an
isomorphism. Finally from the commutativity of the diagram
\eqref{unitmod}, $\upsilon_N$ is an isomorphism for every
$N\in\rmod{B}$.

Now we will prove the second statement. Obviously the condition (1)
implies the condition (4).

$(2)\Rightarrow(1)$ By Proposition \ref{19a},
$(\Sigma^*\tensor{B}-,\Sigma\cotensor{\coring{C}}-)$ is an adjoint
pair, and furthermore, $\Sigma_\coring{C}$ is generator if and only
if $\Sigma\cotensor{\coring{C}}-$ is faithful. From the
commutativity of the diagram \eqref{unitcounitmod}, and
$\Sigma\cotensor{\coring{C}}-$ is faithful, it follows that can is
an isomorphism. Hence (2) follows.

$(3)\Rightarrow(2)$ Since $_A\coring{C}$ is projective, and
$\textsf{can}$ is surjective, $\textsf{can}$ is a retraction in
$\lmod{A}$. Since $\coring{C}$ is coseparable, it follows that
$\textsf{can}$ is a retraction in $\lcomod{\coring{C}}$. Then for
every $M\in \rcomod{\coring{C}}$,
$M\cotensor{\coring{C}}\mathsf{can}$ is a retraction in $\rmod{k}$,
and from the commutativity of the diagram \eqref{counitmod},
$\zeta_M$ is surjective for every $M\in \rcomod{\coring{C}}$. By
Proposition \ref{Maclane-adj-fun}, $\Sigma_\coring{C}$ is a
generator.

$(4)\Rightarrow(3)$ We have that
$\Sigma\cotensor{\coring{C}}\mathsf{can}$ is bijective. From the
commutativity of the diagram \eqref{unitcounitmod},
$\upsilon\tensor{B}\Sigma$ is bijective. Since $_B\Sigma$ is
completely faithful, $\upsilon$ and $\phi$ are also bijective.
\end{proof}

\begin{remark}
In \cite[Proposition 5.6]{Caenepeel/DeGroot/Vercruysse:unp}, the
authors have a similar version of our second statement. They state
that if $\coring{C}_A$ is projective,
$B=\operatorname{End}_\coring{C}(\Sigma)$, and $\coring{C}$ is
coseparable, then $(F,G)$ defined as above is a pair of inverse
equivalences.
\end{remark}

In order to give a generalization of \cite[Theorem
3.5]{Takeuchi:1977} and \cite[Corollary 7.6]{AlTakhman:2002}, we
need the following result.

\begin{lemma}\label{injectivecogenerator}
Let $N\in{}^\coring{C}\mathcal{M}^\coring{D}$ be a bicomodule.
Suppose that $A$ is a QF ring. \begin{enumerate}[(a)] \item If $N$
is an injector-cogenerator as a right $\coring{D}$-comodule, then
$N$ is an injective cogenerator in $\mathcal{M}^\coring{D}$.
\item If $\coring{D}$ has a duality, and $N$
is an injective cogenerator in $\mathcal{M}%
^\coring{D}$ such that $N_B$ is flat, then $N$ is an
injector-cogenerator as a right $\coring{D}$-comodule.
\end{enumerate}
\end{lemma}

\begin{proof}
(a) Since $A$ is a QF ring, then $A_A$ is an injective cogenerator.
Hence $N_\coring{D}\simeq (A\tensor{A}N)_\coring{D}$ is an injective
cogenerator.

(b) Let $X_A$ be an injective cogenerator module. Since $A$ is a QF
ring, $X_A$ is projective. We have then the natural isomorphism
\[
(X\otimes_AN)\square_{\coring{D}}-\simeq
X\otimes_A(N\square_{\coring{D}}-)
:{}^{\coring{D}}\mathcal{M}\rightarrow \mathcal{M}_{k}.
\]
By Proposition \ref{13}, $N_{\coring{D}}$ and $X_A$ are faithfully
coflat, and then $X\otimes_AN$ is faithfully coflat. Once again by
Proposition \ref{13} ($X\otimes_AN$ is a flat right $B$-module),
$X\otimes_AN$ is injective cogenerator in $\mathcal{M}^\coring{D}$.
\end{proof}

\begin{theorem}\label{equivalence6}
Suppose that $_A\coring{C}$, $\coring{C}_A$, $_B\coring{D}$ and
$\coring{D}_B$ are flat, and let
$X\in{}^{\coring{C}}\mathcal{M}^{\coring{D}}$ and
$\Lambda\in{}^{\coring{D}}\mathcal{M}^\coring{C}$. If $\coring{C}$
and $\coring{D}$ have a duality, then the following statements are
equivalent:
\begin{enumerate}[(1)] \item $(-\square
_{\coring{C}}X,-\square_{\coring{D}}\Lambda)$ is a pair of inverse
equivalences with $X_B$ and $\Lambda_A$ are flat;
\item $(\Lambda
\cotensor{\coring{C}}-,X \cotensor{\coring{D}}-)$ is a pair of
inverse equivalences with $_AX$ and $_B\Lambda$ are flat;

If in particular $A$ and $B$ are QF rings, then (1) and (2) are
equivalent to \item $\Lambda$ is quasi-finite injective cogenerator
as a right $\coring{C}$-comodule with $X_B$ and $\Lambda_A$ are
flat, $\e{\coring{C}}{\Lambda}\simeq \coring{D}$ as corings and
$X\simeq \cohom{\coring{C}}{\Lambda}{\coring{C}}$ in
$^{\coring{C}}\mathcal{M}^\coring{D}$;\item $X$ is quasi-finite
injective cogenerator as a right $\coring{D}%
$-comodule with $X_B$ and $\Lambda_A$ are flat,
$\e{\coring{D}}{X}\simeq \coring{C}$ as corings and $\Lambda\simeq
\cohom{\coring{D}}{X}{\coring{D}}$ in
$^{\coring{D}}\mathcal{M}^\coring{C}$;\item $\Lambda$ is
quasi-finite injective cogenerator on both sides,
$\e{\coring{C}}{\Lambda}\simeq \coring{D}$ as corings, and $X\simeq
\cohom{\coring{C}}{\Lambda}{\coring{C}}$ in
$^{\coring{C}}\mathcal{M}^\coring{D}$;\item $X$ is quasi-finite
injective cogenerator on both sides, $\e{\coring{D}}{X}\simeq
\coring{C}$ as corings, and $\Lambda\simeq
\cohom{\coring{D}}{X}{\coring{D}}$ in
$^{\coring{D}}\mathcal{M}^\coring{C}.$
\end{enumerate}
\end{theorem}

\begin{proof}
The equivalence between (1) and (2) follows from Theorem
\ref{equivalence4}, and the fact that if $(-\square
_{\coring{C}}X,-\square_{\coring{D}}\Lambda)$ is a Frobenius pair,
$X_{\coring{D}}$ and $\Lambda_{\coring{C}}$ are coflat if and only
if $X_B$ and $\Lambda_A$ are flat (see the proof of Theorem
\ref{leftrightFrobeniusduality}).

$(1)\Leftrightarrow(3)$ Obvious from Proposition
\ref{equivalence3}$(I)$, and Lemma \ref{injectivecogenerator}.

$(1)\Leftrightarrow(4)$ The proof is analogous to that of
``$(1)\Leftrightarrow(3)$''.

$(5)\Rightarrow(3)$ Trivial.

 $(1)\Rightarrow(5)$ Obvious from the
equivalence of (1), (2) and (3).

$(1)\Leftrightarrow(6)$ Follows by symmetry.

\end{proof}

\begin{remark}
The second part of \cite[Corollary 7.6]{AlTakhman:2002} (and also
the second part of \cite[12.14]{Brzezinski/Wisbauer:2003}) is true
in a more general context. Let $_A\coring{C}$ be flat. Consider the
statements: \begin{enumerate}[(1)] \item
$\Lambda\in{}_B\mathcal{M}^\coring{C}$ be quasi-finite as a right
$\coring{C}$-comodule; \item $\hom{\coring{C}}{M_{0}}{\Lambda}$ is a
finitely generated left $B$-module, for every finitely generated
comodule $M_{0}\in\mathcal{M}^\coring{C}$; \item
$\hom{\coring{C}}{M_{0}}{\Lambda}$ is a finitely generated
projective left $B$-module, for every finitely generated comodule
$M_{0}\in\mathcal{M}^\coring{C}$. \end{enumerate} We have,
$(1)\Rightarrow(2)$ holds if $B$ is a QF ring, and the category
$\mathcal{M}^\coring{C}$ is locally finitely generated, and the
converse implication holds if in particular $B$ is a semisimple
ring, and the category $\mathcal{M}^\coring{C}$ is locally finitely
generated. $(1)\Leftrightarrow(3)$ holds if $\coring{C}$ is a
cosemisimple coring. Moreover, for the two cases, if (1) holds, then
for every comodule $M\in\mathcal{M}^\coring{C},$
$$\cohom{\coring{C}}{\Lambda}{M}\simeq
\underset{\underset{I}{\longrightarrow}}{\lim}
\hom{\coring{C}}{M_i}{\Lambda}^*,$$ where $(M_i)_{i\in I}$ is the
family of all finitely generated subcomodules of $M$.

We will prove at the same time $(1)\Rightarrow(2)$ and
$(1)\Rightarrow(3)$. Let $M_{0}\in\mathcal{M}^\coring{C}$ be a
finitely generated comodule. We have,
$\cohom{\coring{C}}{\Lambda}{M_{0}}^*=\hom{B}{\cohom{\coring{C}}
{\Lambda}{M_{0}}}{B}\simeq \hom{\coring{C}}{M_{0}}{\Lambda}$. Since
the functor $-\tensor{B}{\Lambda}$ is exact and preserves
coproducts, the cohom functor $\cohom{\coring{C}}{\Lambda}{-}$
preserves finitely generated (resp. finitely generated projective)
objects. In particular, $\cohom{\coring{C}}{\Lambda}{M_{0}}$ is a
finitely generated (resp. finitely generated projective) right
$B$-module, and $\hom{\coring{C}}{M_{0}}{\Lambda}$ so is. Therefore,
$\cohom{\coring{C}}{\Lambda}{M_{0}}\simeq\cohom{\coring{C}}
{\Lambda}{M_{0}}^{**}\simeq \hom{\coring{C}}{M_{0}}{\Lambda}^*$.
Hence, $$\cohom{\coring{C}}{\Lambda}{M}\simeq
\underset{\underset{I}{\longrightarrow}}{\lim}\cohom{\coring{C}}
{\Lambda}{M_i}\simeq
\underset{\underset{I}{\longrightarrow}}{\lim}\hom{\coring{C}}
{M_i}{\Lambda}^*,$$ where $(M_i)_{i\in I}$ is the family of all
finitely generated subcomodules of $M$ (since the cohom functor
preserves inductive limits).

 The proof of the implication $(i)\Rightarrow(ii)$ of
\cite[Proposition 1.3]{Takeuchi:1977}, remains valid to prove
$(2)\Rightarrow(1)$ and $(3)\Rightarrow(1)$.
\end{remark}

\section{Applications to induction functors}\label{induction}

In this section we particularize our results in the previous
sections to induction functors defined in Subsection
\ref{IndunctionFunctors}.

\begin{theorem}\label{equivmor}
Let $(\varphi,\rho) :\coring{C}\rightarrow\coring{D}$ be a
homomorphism of corings such that $_A\coring{C}$ and $_B\coring{D} $
are flat. If $\coring{C}_A$ and $\coring{D}_B$ are flat (resp.
$\coring{C}$ and $\coring{D}$ are coseparable), then the following
statements are equivalent
\begin{enumerate}[(a)] \item
$(-\otimes_AB,-\square_{\coring{D}}(B\tensor{A}\coring{C}))$ is a
pair of inverse equivalences;\item the functor $-\otimes_AB$ is
exact and faithful (resp. faithful), and there exists an isomorphism
of $B$-corings $\coring{D}\simeq B\coring{C}B.$
\end{enumerate}
\end{theorem}

\begin{proof}
From the proof of Theorem \ref{23}, $-\otimes_AB\simeq
-\cotensor{\coring{C}}{(\coring{C}\otimes_AB)}$. The use of
Proposition \ref{equivalence3}, Proposition \ref{comatrixproperty1},
and Proposition \ref{BW 23.9} achieves the proof.

\end{proof}

\begin{corollary}
Let $\coring{C}$ be an $A$-coring. Then the following statements are
equivalent
\begin{enumerate}[(a)]
\item The forgetful functor
$U_r:\mathcal{M}^{\coring{C}%
}\rightarrow\mathcal{M}_A$ is an equivalence of categories;
\item The forgetful functor
$U_l:{}^{\coring{C}}\mathcal{M}\rightarrow{}_A\mathcal{M}$ is an
equivalence of categories;\item there exists an isomorphism of
$A$-corings $A\simeq \coring{C}.$
\end{enumerate}
\end{corollary}

\begin{proof}
It is enough to apply the last theorem to the particular
homomorphism of corings
$(\epsilon_{\coring{C}},1_A):\coring{C}\rightarrow A$ which gives
the well known adjunction $(U_r,-\tensor{A}{\coring{C}})$, and
observing that the map $A\coring{C}A\to \coring{C}$ defined by
$a\otimes c\otimes a'\mapsto aca'$, is an isomorphism of corings.
\end{proof}

Finally, given a homomorphism of corings, we give sufficient
conditions to have that the right induction functor is an
equivalence if and only if the left induction functor so is. Note
that for the case of coalgebras over fields (by (b)), or of rings
(well known) (by (d)), we have the left right symmetry.

\begin{proposition}\label{equivmorsym}
Let $(\varphi,\rho):\coring{C}\rightarrow\coring{D}$ be a
homomorphism of corings such that $_A\coring{C}$, $_B\coring{D}$,
$\coring{C}_A$ and $\coring{D}_B$ are flat. Assume that at least one
of the following holds
\begin{enumerate}[(a)] \item $\coring{C}$ and
$\coring{D}$ have a duality, and $_AB$ and $B_A$ are flat;\item $A$
and $B$ are von Neumann regular rings;\item $B\otimes_A\coring{C}$
is coflat in $^{\coring{D}}\mathcal{M}$ and $\coring{C}\otimes_AB$
is coflat in $\mathcal{M}^{\coring{D}}$ and $_AB $ and $B_A$ are
flat;\item $\coring{C}$ and $\coring{D}$ are coseparable corings.
\end{enumerate} Then the following statements are
equivalent
\begin{enumerate}
\item $-\otimes_AB:\mathcal{M}^{\coring{C}}\to\mathcal{M}^\coring{D}$
is an equivalence of categories;
\item $B\otimes_A-:{}^{\coring{C}}\mathcal{M}\to{}^{\coring{D}}\mathcal{M}$
is an equivalence of categories.
\end{enumerate}
\end{proposition}

\begin{proof}
Obvious from Theorem \ref{equivalence4}, Theorem \ref{equivalence5}
and Theorem \ref{equivalence6}.
\end{proof}

\section{Applications to entwined modules and graded
ring theory}\label{entwined'}

In this section we particularize our results to corings associated
to entwined structures and in particular those associated to a
$G$-graded algebra and a right $G$-set, where $G$ is group.

\medskip
 We obtain the following result concerning the category of
entwined modules.

\begin{theorem}
Let $(\alpha,\gamma):(A,C,\psi)\to (A',C',\psi')$ be a morphism in
$\mathbb{E}_\bullet^\bullet(k)$, such that $_{k}C$ and $_{k}D$ are
flat. If either $(A\otimes C)_A$ and $(A'\otimes C')_{A'}$ are flat
(e.g., if $\psi$ and $\psi'$ are isomorphisms), (resp. $A\otimes C$
and $A'\otimes C'$ are coseparable), then the following statements
are equivalent
\begin{enumerate}[(a)] \item The functor
$-\otimes_AA':\mathcal{M}(\psi)_A^C\rightarrow
\mathcal{M}(\psi')_{A'}^{C'}$ defined in Section \ref{entwined} is
an equivalence;
\item the functor
$-\otimes_AA'$ is exact and faithful (resp. faithful), and there
exists an isomorphism of $A'$-corings $A'\otimes C'\simeq
A'(A\otimes C)A'.$
\end{enumerate}
\end{theorem}

\begin{proof}
Follows immediately from Theorem \ref{equivmor}.
\end{proof}

\begin{corollary}\emph{\label{graded1}\cite[Proposition 2.1]{DelRio:1992}}
\\Let $G$ and $G'$ be two groups, $A$ a $G$-graded
$k$-algebra, $A'$ a $G'$-graded $k$-algebra, $X$ a right $G$-set,
and $X'$ a right $G'$-set. The following statements are equivalent
\begin{enumerate}[(1)] \item the categories
$gr-(A,X,G)$ and $gr-(A',X',G')$ are equivalent;
\item there are an $X\times X'$-graded $(A,A')$-bimodule
$P$ and an $X'\times X$-graded $(A',A)$-bimodule $Q$ such that
$$P\widehat{\otimes}_{A'}Q\simeq \widehat{A}\quad
\textrm{and}\quad
 Q\widehat{\otimes}_AP\simeq \widehat{A'}.$$
\end{enumerate}
Moreover, if $P$ and $Q$ satisfy the condition (2), then
$-\widehat{\otimes}_AP$ and $-\widehat{\otimes}_{A'}Q$ are inverse
equivalences.
\end{corollary}

\begin{proof}
Clear from Proposition \ref{equivalence1} (using the fact that the
corings $A\otimes kX$ and $A'\otimes kX'$ are coseparable) and
Proposition \ref{graded-A5}.
\end{proof}

The following result is similar to \cite[Corollary
2.4]{DelRio:1992}.

\begin{corollary}
Let $A$ be a $G$-graded $k$-algebra, $X$ a right $G$-set, and $B$ a
$k$-algebra. The following statements are equivalent
\begin{enumerate}[(1)] \item the category $gr-(A,X,G)$
is equivalent to $\rmod{B}$;
\item the category $(G,X,A)-gr$ is equivalent to
$\lmod{B}$;
\item there exists an $X_0\times X$-graded
$B-A$-bimodule $P$, such that $X_0$ is a singleton, $P$ is finitely
generated projective in $\rmod{A}$, and generator in $gr-(A,X,G)$,
and the morphism of $k$-algebras $\phi:B\to
\operatorname{End}_{gr-(A,X,G)}(P)$ defined by $\phi(b)(p)=bp$ for
$b\in B,p\in P$, is an isomorphism.
\end{enumerate}
\end{corollary}

\begin{proof}
Follows immediately from Theorem \ref{equivmodcomod}.
\end{proof}

In \cite[Theorem 2.3]{DelRio:1992}, A. Del R\'{\i}o gave a
characterization when $(-\widehat{\otimes}_AP,\h{P_{A'}}{-})$ is a
pair of inverse equivalences. We think that the following result
gives a simple characterization of it. Our result is also a
generalization of Morita's characterization of equivalence Remark
\ref{Morita}.

\begin{theorem}\label{graded2}
Let $P$ be an $X\times X'$-graded $(A,A')$-bimodule. Then the
following are equivalent
\begin{enumerate}[(1)] \item
$(-\widehat{\otimes}_AP,\h{P_{A'}}{-})$ is a pair of inverse
equivalences;
\item \begin{enumerate}[(a)]\item $_xP$ is finitely
generated projective in $\rmod{A'}$ for every $x\in X$, and $P_{x'}$
is finitely generated projective in $\lmod{A}$ for every $x'\in X'$
(resp. and $P$ is a generator in both $gr-(A',X',G')$ and
$(G,X,A)-gr$),
\item the following bigraded bimodules maps:
 $\psi:\widehat{A}\rightarrow \h{P_{A'}}{P}$ defined
by $\psi(a\otimes x)(p)=a({}_xp)\in
 P$ $(a\in A,x\in X,p\in P)$, and
$\psi':\widehat{A'}\rightarrow \h{_AP}{P}$ defined by
 $\psi'(a'\otimes x')(p)=(p_{x'})a'\in
 P$ $(a'\in A',x'\in X',p\in P)$, are isomorphisms
(resp. are surjective maps).
\end{enumerate}

\item \begin{enumerate}[(a)]\item $_xP$ is finitely
generated projective in $\rmod{A'}$ for every $x\in X$, \item the
evaluation map
$$\varepsilon_{\widehat{A'}}:\h{P_{A'}}{\widehat{A'}}\widehat{\otimes}_AP\to
\widehat{A'},\; \varepsilon_{\widehat{A'}}(f\tensor{A}{P})=f(p)$$
($f\in \h{P_{A'}}{\widehat{A'}}_x, p\in {}_xP, x\in X$) is an
isomorphism,
\item the functor $-\widehat{\otimes}_AP$ is faithful.
\end{enumerate}
\end{enumerate}
\end{theorem}

\begin{proof}
From \cite[Proposition 1.2]{Menini:1993}, the unit and the counit
of the adjunction \\$(-\widehat{\otimes}_AP,\h{P_{A'}}{-})$ are
given respectively by $\eta_M:M\rightarrow
\h{P_{A'}}{M\widehat{\otimes}_AP},$ $\eta_M(m)(p)=\sum_{x\in
X}m_x\tensor{A}{_xp}$ ($m=\sum_{x\in X}m_x\in M,p=\sum_{x\in
X}{}_xp\in P$), and
$\varepsilon_N:\h{P_{A'}}{N}\widehat{\otimes}_AP\rightarrow N$,
$\varepsilon_N(f\tensor{A}{P})=f(p)$ ($f\in \h{P_{A'}}{N}_x,p\in
{}_xP, x\in X$). By Lemma \ref{graded-A6} (1), the functor
$\h{P_{A'}}{-})$ preserves inductive limits if and only if the
condition $(3)(a)$ holds. Then by Theorem \ref{3} (the coring
$A'\otimes kX'$ is coseparable), there is a natural isomorphism
$\xymatrix{\delta:\h{P_{A'}}{-}\ar[r]^-\simeq &
-\widehat{\otimes}_{A'}\h{P_{A'}}{\widehat{A'}}}$. Set
$Q=\h{P_{A'}}{\widehat{A'}}$.

$(1)\Leftrightarrow(3)$ It follows from Proposition
\ref{equivalence3}(II) and Proposition \ref{comatrixproperty1}(2).

$(1)\Leftrightarrow(2)$ We can suppose that the condition $(2)(a)$
holds. We have $\eta_{\widehat{A}}:M\rightarrow
\h{P_{A'}}{\widehat{A}\widehat{\otimes}_AP},$
$\eta_{\widehat{A}}(a\otimes x)(p)=(a\otimes
x)\tensor{A}{}_xp=a(1_A\otimes
x)\tensor{A}{}_xp=a\lambda_P({}_xp)=\lambda_P(a({}_xp))=
\lambda_P(\psi(a\otimes x)(p)),$ $(a\in A,x\in X,p\in P)$. Since
$\delta$ is a natural isomorphism, the following diagram is
commutative:
$$
\xymatrix{\widehat{A}\ar[ddr]^{\psi}\ar[r]^-{\eta_{\widehat{A}}} &
\h{P_{A'}}{\widehat{A}\widehat{\otimes}_AP}
\ar[rr]^{\delta_{\widehat{A}\widehat{\otimes}_AP}} & &
\widehat{A}\widehat{\otimes}_AP\widehat{\otimes}_{A'}Q
\ar[r]^-\simeq & P\widehat{\otimes}_{A'}Q \\
& & & & \\
& \h{P_{A'}}{P}\ar[uu]^\simeq_{\h{P_{A'}}{\lambda_P}}
\ar[uurrr]^{\delta_P}& & & .}
$$
Therefore the unit of the adjunction
$(-\widehat{\otimes}_AP,-\widehat{\otimes}_{A'}Q)$ is
$$\xymatrix@1{1_{gr-(A,X,G)}\ar[r]^-\simeq &
-\widehat{\otimes}_A\widehat{A}\ar[rr]^-{-\widehat{\otimes}_A\psi_0}
& & -\widehat{\otimes}_AP\widehat{\otimes}_{A'}Q},$$ where $\psi_0$
is the map $\psi_0:\xymatrix{\widehat{A}\ar[r]^-{\psi} &
\h{P_{A'}}{P}\ar[r]^-{\delta_P} & P\widehat{\otimes}_{A'}Q}$. Hence
$-\widehat{\otimes}_AP$ is fully faithful if and only if $\psi$ is
an isomorphism. By Proposition \ref{19a},
$(Q\widehat{\otimes}_A-,P\widehat{\otimes}_{A'}-)$ is an adjoint
pair. Moreover $-\widehat{\otimes}_{A'}Q$ is fully faithful if and
only if $P\widehat{\otimes}_{A'}-$ is fully faithful, if and only if
$\psi'$ is an isomorphism (see Proposition \ref{Maclane-adj-fun}).

Finally, if the maps $\psi$ and $\psi'$ defined in $(2)(b)$ are
surjective, then, $P$ is a generator in $gr-(A',X',G')$ if and only
if $-\widehat{\otimes}_{A'}Q$ is faithful, if and only if
$P\widehat{\otimes}_{A'}-$ is faithful, if and only if $\psi'$ is an
injective map (see Proposition \ref{Maclane-adj-fun}). By symmetry,
$P$ is a generator in $(G,X,A)-gr$ if and only if $\psi$ is an
injective map.
\end{proof}

Finally, let $f:G\to G'$ be a morphism of groups, $X$ a right
$G$-set, $X'$ a right $G'$-set, $\varphi:X\to X'$ a map such that
$\varphi(xg) =\varphi(x)f(g) $ for every $g\in G,$ $x\in X.$ Let $A$
be a $G$-graded $k$-algebra, $A'$ a $G'$-graded $k$-algebra, and
$\alpha:A\to A'$ a morphism of algebras such that $\alpha( A_{g})
\subset A_{f(g)}'$ for every $g\in G.$

 We have, $\gamma:kX\to kX'$ such that
$\gamma(x)=\varphi(x)$ for each $x\in X$, is a morphism of
coalgebras, and $(\alpha ,\gamma):(A,kX,\psi)\to (A',kX',\psi')$
is a morphism in $\mathbb{E}_\bullet^\bullet(k).$

 Let $T^*=-\otimes_AA': gr-( A,X,G)\to
gr-(A',X',G')$ be the functor defined in Subsection \ref{$T^*$ Frob}.
We know that $T^*$ makes commutative the following diagram
\[
\xymatrix{gr-(A,X,G)\ar[d]^\simeq \ar[rr]^{T^*} & &
gr-(A',X',G') \ar[d]^\simeq  \\
\mathcal{M}(kG)_A^{kX}
\ar[rr]^{-\otimes_AA'} & & \mathcal{M}(kG')_{A'
}^{kX'}.}
\]
Moreover, we have the commutativity of the following diagram
\[%
\xymatrix{(G,X,A)-gr \ar[d]^\simeq \ar[rr]^{(T^*)'=A'\otimes_A-}
& & (G',X',A')-gr \ar[d]^\simeq\\
{}_A^{kX}\mathcal{M}(\psi^{-1})
\ar[d]^\simeq\ar[rr]^{A'%
\otimes_A-} & & {}_{A'}^{kX'}\mathcal{M}%
((\psi')^{-1}) \ar[d]^\simeq \\
{}^{kX\otimes A}\mathcal{M} \ar[d]^\simeq
\ar[rr]^{A'\otimes_A-}%
 & & {}^{kX'\otimes A'}\mathcal{M} \ar[d]^\simeq
\\
{}^{A\otimes kX}\mathcal{M} \ar[rr]^{A'\otimes_A-}%
 & & {}^{A'\otimes kX'}\mathcal{M}.}
\]

A. Del R\'{\i}o gave in \cite[Example 2.6]{DelRio:1992} an
interesting characterization when $T^*$ is an equivalence. Our
result gives an other characterization of it.

\begin{theorem}\label{graded3}
The following statements are equivalent
\begin{enumerate} [(1)]
\item the functor $T^*:gr-(A,X,G)\to
gr-(A',X',G')$ is an equivalence;
\item $T^*$ is faithful, and the map
$$A'\tensor{A}\widehat{A}\tensor{A}A'\to
\widehat{A'}, \quad a'\tensor{A}(a\otimes x)\tensor{A}a''\mapsto
a'\alpha(a)a''\otimes \varphi(x)g',$$ where $a\in A,a'\in A',x\in X,
g'\in G',a''\in A'_{g'}$, is bijective.
\end{enumerate}
\end{theorem}

\begin{proof}
Obvious from the commutativity of the diagram \eqref{counit-ind-fun-diag}
and Theorem \ref{equivmor}.
\end{proof}

Finally we give the following consequence of Proposition
\ref{equivmorsym}:

\begin{proposition}\label{graded4}
The following are equivalent
\begin{enumerate} [(1)]\item the functor
$T^*:gr-(A,X,G) \to gr-( A',X',G')$ is an equivalence;
\item the functor $(T^*)':(G,X,A)-gr\to
(G',X',A')-gr$ is an equivalence.
\end{enumerate}
\end{proposition}

\chapter{The Picard Group of Corings}

In this chapter, all the considered corings are flat on both sides
over their base rings.

\section{The Picard group of corings}\label{Picard}

Corings over $k$-algebras and their morphisms form a category. We
denote it by $\mathbf{Crg}_k$. Notice that the category of
$k$-algebras $\mathbf{Alg}_k$, and the category of $k$-coalgebras
$\mathbf{Coalg}_k$, are full subcategories of $\mathbf{Crg}_k$.

\begin{lemma}
A morphism $(\varphi,\rho)\in \mathbf{Crg}_k$ is an isomorphism if
and only if $\varphi$ and $\rho$ are bijective.
\end{lemma}

\begin{proof}
Straightforward.
\end{proof}

\medskip
Now we recall from \cite{Bass:1968} the definition of the Picard
group of a $k$-abelian category.

Let $\cat{C}$ be a $k$-abelian category, $\pic{\cat{C}}$ is the
group of isomorphism classes $(T)$ of $k$-equivalences
$T:\cat{C}\to\cat{C}$ with the composition law $(T)(S)=(ST)$.

\begin{definition}
We say that a $(\coring{C},\coring{D})$-bicomodule $X$ is
\emph{right invertible} if the functor
$-\cotensor{\coring{C}}X:\rcomod{\coring{C}}\to\rcomod{\coring{D}}$
is an equivalence. We say that $\coring{C}$ and $\coring{D}$ are
\emph{Morita-Takeuchi right equivalent} if there is a right
invertible $(\coring{C},\coring{D})$-bicomodule.
\end{definition}

The isomorphism classes $(X)$ of right invertible
$\coring{C}$-bicomodules with the composition law
$$(X_1)(X_2)=(X_1\cotensor{\coring{C}}X_2),$$
form the \emph{right Picard group} of $\coring{C}$. We denote it by
$\rpic{\coring{C}}$. This law is well defined since
$-\cotensor{\coring{C}}(X_1\cotensor{\coring{C}}X_2)\simeq
(-\cotensor{\coring{C}}X_1)\cotensor{\coring{C}}X_2$ (From
Propositions \ref{assocotens1}, \ref{assocotens2} ($_\coring{C}X_2$
is coflat)). The associativity of this law follows from the
associativity of the cotensor product ($_\coring{C}X_1$ and
$_\coring{C}X_2$ are coflat). The isomorphism class $(\coring{C})$
is the identity element. $(X)^{-1}=(\Lambda)$, where $\Lambda$ is
such that $(-\cotensor{\coring{C}}X,-\cotensor{\coring{C}}\Lambda)$
is a pair of inverse equivalences (see Proposition
\ref{equivalence1}).

It follows from Theorems \ref{equivalence5}, \ref{equivalence6} that
if $\coring{C}$ is coseparable (this case includes that of
algebras), or cosemisimple, or $A$ is von Neumann regular ring, or
if $\coring{C}$ is a coalgebra over a QF ring, then a
$\coring{C}$-bicomodule $X$ is right invertible if and only if it is
left invertible. Hence, for each case we have
$\rpic{\coring{C}}=\lpic{\coring{C}}$. We will denote this last by
$\pic{\coring{C}}$.

\medskip

The following proposition is the motivation of the study of the
Picard group of corings. The proof is straightforward using Theorem
\ref{3}.

\begin{proposition}\label{picard-cat-crg}
There are inverse isomorphisms of groups
$$\xymatrix{\rpic{\coring{C}}\ar@<2pt>[r]^-\alpha &
\pic{\rcomod{\coring{C}}}\ar@<2pt>[l]^-\beta},$$
$\alpha((X))=(-\cotensor{\coring{C}}X)$ and
$\beta((T))=(T(\coring{C}))$.
\end{proposition}

\medskip
Now we introduce two categories, $\mathfrak{Crg}_k$ and
$\mathfrak{MT-Crg}_k^r$. The objects of both are the corings. The
morphisms of $\mathfrak{Crg}_k$ are the isomorphisms of corings. The
morphism of $\mathfrak{MT-Crg}_k^r$ are the Morita-Takeuchi right
equivalences. A right \emph{Morita-Takeuchi equivalence}
$\coring{C}\sim\coring{D}$ is simply an isomorphism class $(M)$
where $_\coring{D}M_\coring{C}$ is a right invertible
$(\coring{D},\coring{C})$-bicomodule. If $(N)$ is a Morita-Takeuchi
right equivalence $\coring{D}\sim\coring{C}'$, the composite is
given by $(N\cotensor{\coring{D}}M)$.

The categories $\mathfrak{Crg}_k$ and $\mathfrak{MT-Crg}_k^r$ are
groupoids (i.e., categories in which every morphism is an
isomorphism). Then $\aut{\coring{C}}$ and $\rpic{\coring{C}}$ are
respectively the endomorphism group of the object $\coring{C}$ in
$\mathfrak{Crg}_k$ and in $\mathfrak{MT-Crg}_k$.

Let $f:(\coring{C}:A)\to(\coring{D}:B)$ and
$g:(\coring{C}_1:A_1)\to(\coring{D}_1:B_1)$ be two morphisms of
corings, and let $_{\coring{C}_1}X_{\coring{C}}$ be a bicomodule.

Consider the right induction functor $-\tensor{A}B:
\rcomod{\coring{C}}\to \rcomod{\coring{D}}$ and the left induction
functor $B_1\tensor{A_1}-: \lcomod{\coring{C_1}}\to
\lcomod{\coring{D_1}}$.

Now, consider $_gX_f:=B_1\tensor{A_1}X\tensor{A}B$ which is a
$(\coring{D}_1,\coring{D})$-bicomodule, and $X_f:=X\tensor{A}B$
which is a $(\coring{C}_1,\coring{D})$-bicomodule. It can be showed
easily that $_1X_f\simeq X_f$ as bicomodules, and that if
$\xymatrix{\coring{C}\ar[rr]^{f=(\varphi,\rho)} && \coring{D}
\ar[rr]^{f'=(\varphi',\rho')}&& \coring{C}'}$ are morphisms of
corings, then $(X_f)_{f'}\simeq X_{f'f}$ as bicomodules. Moreover,
from Theorem \ref{3}, $X_f\simeq X\cotensor{\coring{C}}\coring{C}_f$
as $(\coring{C}_1,\coring{D})$-bicomodules.

Now, assume that $f$ is an isomorphism. Define $X_f'=X_f$ as
$\coring{C}_1$-comodules, with right $B$-action and right
$\coring{D}$-coaction given by
$$xb=x\rho^{-1}(b) \textrm{ and }
\rho_{X_f'}(x)= x_{(0)}\tensor{B}\varphi(x_{(1)}),$$ where $x\in
X,b\in B$ and $\rho_X(x)= x_{(0)}\tensor{A}x_{(1)}$. $X_f$ and
$X_f'$ are isomorphic as $(\coring{C},\coring{D})$-bicomodules. The
connecting isomorphism is the canonical isomorphism of abelian
groups $\alpha :X_f\to X_f'$ defined by
$\alpha(x\tensor{A}b)=x\rho^{-1}(b)$ whose inverse map is
 defined by $\beta(x)=x\tensor{A}1$. Similarly, we define $_fX'$.

\begin{lemma}\label{functor}
We have a functor $$\Omega:\mathfrak{Crg}_k \to
\mathfrak{MT-Crg}_k^r,$$ which is the identity on objects and
associates with an isomorphism of corings
$f:(\coring{C}:A)\to(\coring{D}:B)$ the isomorphism class of the
right invertible $(\coring{D},\coring{C})$-bicomodule $_f\coring{C}$.
\end{lemma}

\begin{proof}
Let $f=(\varphi,\rho):(\coring{C}:A)\to(\coring{D}:B)$ be an
isomorphism of corings. Since
$B\tensor{A}\coring{C}\tensor{A}B\to\coring{D},\;
b\tensor{A}c\tensor{A}b\mapsto b\varphi(c)b$ is a morphism of
$B$-corings, it follows from Theorem \ref{equivmor} that
$_f\coring{C}$ is a right invertible
$(\coring{D},\coring{C})$-bicomodule. Now let
$\xymatrix{\coring{C}\ar[r]^f & \coring{D} \ar[r]^g& \coring{C}'}$
be two morphisms in $\mathfrak{Crg}_k$. Then
$_{gf}\coring{C}\simeq{}_g(_f\coring{C})\simeq {}_g\coring{D}
\cotensor{\coring{D}}{}_f\coring{C}$ is an isomorphism of
$(\coring{C'},\coring{C})$-bicomodules.
\end{proof}

\begin{lemma}\emph{\cite[18.12 (1)]{Brzezinski/Wisbauer:2003}}\label{BW 18.12}
\\The map $\operatorname{End}_\coring{C}(\coring{C})
\to\rdual{\coring{C}}$, $f\mapsto \epsilon f$, is an algebra
anti-morphism with inverse map $h\mapsto
(h\tensor{A}\coring{C})\Delta$.
\end{lemma}

Now we are ready to state and prove the main result of this section.

\begin{theorem}\label{coringexseq}
Let $\coring{C}$ be an $A$-coring such that $_A\coring{C}$ and
$\coring{C}_A$ are flat. \begin{enumerate}[(1)]
\item We have an exact sequence
$$\xymatrix{1\ar[r] & \rinn{\coring{C}}\ar[r]&
\aut{\coring{C}}\ar[r]^\Omega & \rpic{\coring{C}}},$$ where
$\rinn{\coring{C}}$ is the set of $f=(\varphi,\rho)\in
\aut{\coring{C}}$ such that there is $p\in \coring{C}^*$ invertible
with
\begin{equation}\label{coringexseq_cond}
\varphi(c_{(1)})p(c_{(2)})= p(c_{(1)})c_{(2)},
\end{equation}
for every $c\in \coring{C}$. As for algebras and coalgebras over
fields, we call $\rinn{\coring{C}}$ the \emph{right group of inner
automorphisms} of $\coring{C}$

In particular, there is a monomorphism from the quotient group
$$\rout{\coring{C}}:=\aut{\coring{C}}/\rinn{\coring{C}}$$ to
$\rpic{\coring{C}}$. As for algebras and coalgebras over fields, we
call $\rout{\coring{C}}$ the right \emph{outer automorphism group}
of $\coring{C}$.
\item If $\coring{C}$ and $\coring{D}$ are Morita-Takeuchi right
equivalent, then $\rpic{\coring{C}}\simeq\rpic{\coring{D}}$.
\end{enumerate}
\end{theorem}

\begin{proof}
(1) Let $f=(\varphi,\rho)\in \aut{\coring{C}}$ and
$_f\coring{C}\simeq \coring{C}$ as $\coring{C}$-bicomodules.
 We have $_f\coring{C}\simeq {}_f\coring{C}'$ as
$\coring{C}$-bicomodules. Now, let $h:{}_f\coring{C}'\to \coring{C}$
be an
 isomorphism of $\coring{C}$-bicomodules.
$h$ is a morphism of right $\coring{C}$-bicomodules is equivalent to
the existence of $p\in \coring{C}^*$ invertible such that $h(c)=
p(c_{(1)})c_{(2)},$ for every $c\in \coring{C}$ (see Lemma \ref{BW
18.12}). Then, $p=\epsilon\circ h$. That $h$ is a morphism of left
$\coring{C}$-bicomodules means that for every $a\in A,
c\in\coring{C}$
\begin{equation}\label{e1}
h(ac)=\rho(a)h(c),
\end{equation}
and for every $c\in \coring{C}$,
\begin{equation}\label{e2}
\varphi(c_{(1)})\tensor{A}h(c_{(2)})=
h(c)_{(1)}\tensor{A}h(c)_{(2)}.
\end{equation}
The condition \eqref{e2} is equivalent to
$$\varphi(c_{(1)})\tensor{A}p(c_{(2)})c_{(3)}=
p(c_{(1)})c_{(2)}\tensor{A}c_{(3)},$$ for every $c\in \coring{C}$.
Hence the condition \eqref{coringexseq_cond} follows. To prove the
converse, take $p\in \coring{C}^*$ invertible satisfying
\eqref{coringexseq_cond} and set $h(c)=p(c_{(1)})c_{(2)},$ for every
$c\in \coring{C}$. Using the left $A$-linearity of $\varphi$, the
condition \eqref{e1} follows. On the other hand, the condition
\eqref{e2} follows easily.

(2) It follows immediately from that $\coring{C}$ and $\coring{D}$
are isomorphic to each other in the category
$\mathfrak{MT-Crg}_k^r$. More explicitly, let
$_\coring{C}M_\coring{D}$ be a right invertible bicomodule. The map
$$(X)\mapsto (M^{-1}\cotensor{\coring{C}}X\cotensor{\coring{C}}M)$$
is an isomorphism from $\rpic{\coring{C}}$ to $\rpic{\coring{D}}$.
\end{proof}

If we apply Theorem \ref{coringexseq} in the case $\coring{C}=A$,
then we obtain \cite[Proposition II (5.2)(3)]{Bass:1968} (see also
\cite[Theorems 55.9, 55.11]{Curtis/Reiner:1987}). The special case
where $\coring{C}$ is a $k$-coalgebra ($A=k$) with $\coring{C}$ is
flat as a $k$-module, recovers \cite[Theorem
2.7]{Torrecillas/Zhang:1996}.

\begin{corollary}
\begin{enumerate}[(1)]
\item For a $k$-algebra $A$, we have an exact sequence
$$\xymatrix{1\ar[r] & \inn{A}\ar[r]& \aut{A}\ar[r]^\Omega & \pic{A}},$$
where $\inn{A}:=\{a\mapsto bab^{-1}\mid b \;\textrm{invertible in}
\;A\},$ the \emph{group of inner automorphisms} of $A$.
\item  For a $k$-coalgebra $C$ such that $_kC$ is flat, we have an
exact sequence
$$\xymatrix{1\ar[r] & \inn{C}\ar[r]& \aut{C}\ar[r]^\Omega &
\rpic{C}},$$ where $\inn{C}$ (the \emph{group of inner
automorphisms} of $C$) is the set of $\varphi\in \aut{C}$ such that
there is $p\in C^*$ invertible with $\varphi(c)=
p(c_{(1)})c_{(2)}p^{-1}(c_{(3)}),$ for every $c\in C.$
\end{enumerate}
\end{corollary}

\begin{proof}
(1) It suffices to see that for $p\in \coring{C}^*$ where
$\coring{C}=A$, $p$ is invertible if and only if
$p=\lambda_b:a\mapsto ba$ and $b$ is invertible in $A$.

(2) Suppose that $\varphi\in \aut{C}$ and $
\varphi(c_{(1)})p(c_{(2)})=p(c_{(1)})c_{(2)},$ for every $c\in C$.
Then $\varphi(c)=\varphi(c_{(1)})\epsilon_C(c_{(2)})=
\varphi(c_{(1)})p(c_{(2)})p^{-1}(c_{(3)})=
p(c_{(1)})c_{(2)}p^{-1}(c_{(3)}),$ for every $c\in C$. The converse
is obvious.
\end{proof}

\medskip
In order to consider the one-sided comodule structure of invertible
bicomodules we need

\begin{proposition}\label{oneside}
Let $(M),(N)\in \rpic{\coring{C}}$. Then
\begin{enumerate}[(1)] \item $(N)\in \mathrm{Im}\Omega.(M)$  if and
only if $N\simeq {}_fM$ as bicomodules for some $f\in
\aut{\coring{C}}$. \item $_AM_\coring{C}\simeq {}_AN_\coring{C}$ if
and only if $N\simeq {}_fM$ as bicomodules for some
$f=(\varphi,1_A)\in \aut{\coring{C}}$.
\end{enumerate}
\end{proposition}

\begin{proof}
(1) Straightforward from the fact that for every $f\in
\aut{\coring{C}}$, there is an isomorphism of
$\coring{C}$-bicomodules
$_fM\simeq{}_f\coring{C}\cotensor{\coring{C}}M$.

(2) $(\Leftarrow)$ Obvious.  $(\Rightarrow)$ Let
$h:{}_AM_\coring{C}\to {}_AN_\coring{C}$ be a bicomodule
isomorphism. Since $M$ is a $(A,\coring{C})$-quasi-finite comodule
(see Proposition \ref{equivalence3}), it has a structure of
$(\e{\coring{C}}{M},\coring{C})$-bicomodule. Then $N$ has a
structure of $(\e{\coring{C}}{M},\coring{C})$-bicomodule induced by
$h$, and $h$ is an $(\e{\coring{C}}{M},\coring{C})$-bicomodule
isomorphism. On the other hand, $N$ is a
$(A,\coring{C})$-quasi-finite comodule. Let
$F=:\cohom{\coring{C}}{N}{-}:\rcomod{\coring{C}}\to
\rcomod{\e{\coring{C}}{M}}$ be the cohom functor, and let
$\theta:1_{\rcomod{\coring{C}}}\to
F(-)\cotensor{\e{\coring{C}}{M}}N$ and
$\chi:F(-\cotensor{\e{\coring{C}}{M}}N)\to
1_{\rcomod{\e{\coring{C}}{M}}}$ be respectively the unit and the
counit of this adjunction. From Proposition \ref{KG 5.2} and Remark
\ref{KG 5.2'}, the map $\varphi:=\chi_{\e{\coring{C}}{M}}\circ
F(\lambda_N^h):\e{\coring{C}}{N}\to \e{\coring{C}}{M}$, where
$\lambda_N^h:N\to \e{\coring{C}}{M}\cotensor{\e{\coring{C}}{M}}N$ is
the left comodule structure map on $N$, is a morphism of
$A$-corings, and making commutative the diagram
$$\xymatrix{N
\ar[ddrr]^{\lambda_N^h}\ar[rr]^{\theta_N} &&
\e{\coring{C}}{N}\tensor{A}N\ar[dd]^{\varphi\tensor{A}1_N}
\\ &&&\\
& & \e{\coring{C}}{M}\tensor{A}N.}$$  Since $F$ is a left adjoint to
$-\cotensor{\e{\coring{C}}{M}}N$ which is an equivalence,
$\chi_{\e{\coring{C}}{M}}$ and $\varphi$ are isomorphisms. Since $M$
and $N$ are right invertible,
$\e{\coring{C}}{M}\simeq\coring{C}\simeq\e{\coring{C}}{N}$ (see
Proposition \ref{equivalence3}). We identify $\e{\coring{C}}{M}$ and
$\e{\coring{C}}{N}$ with $\coring{C}$. Then $\varphi$ is a coring
endomorphism of $\coring{C}$. Hence $M\simeq {}_fN$ as
$\coring{C}$-bicomodules, where $f=(\varphi,1_A)$.
\end{proof}

The special case where $\coring{C}$ is a $k$-coalgebra ($A=k$) with
$\coring{C}$ is flat as a $k$-module, recovers \cite[Proposition
2.8]{Torrecillas/Zhang:1996}. Its version for algebras is given in
\cite[Proposition II.5.2(4)]{Bass:1968} or \cite[Theorem
55.12]{Curtis/Reiner:1987}.

\begin{corollary}
Let $C$ be a $k$-coalgebra such that $_kC$ is flat, and let
$(M),(N)\in \rpic{C}$. Then $M_C\simeq N_C$ if and only if $(N)\in
\mathrm{Im}\Omega.(M)$, that is, $N\simeq {}_fM$ as bicomodules for
some $f\in \aut{C}$.
\end{corollary}

\begin{proposition}\emph{[Bass]}
\\Let $A$ be a $k$-algebra
and let $(M),(N)\in \pic{A}$. Then $M_A\simeq N_A$ if and only if
$(N)\in \mathrm{Im}\Omega.(M)$, that is, $N\simeq {}_fM$ as
bicomodules for some $f\in \aut{A}$.
\end{proposition}

\section{The Aut-Pic Property}\label{Aut-Pic}

Following \cite{Bolla:1984} and
\cite{Cuadra/Garcia/Torrecillas:2000}, we give the following

\begin{definition}
We say that a coring $\coring{C}$ has the \emph{right Aut-Pic
property} if the group morphism $\Omega$ of Theorem
\ref{coringexseq} is surjective. In such a case,
$\rout{\coring{C}}\simeq \rpic{\coring{C}}$.
\end{definition}

Bolla proved (\cite[p. 264]{Bolla:1984}) that every ring such that
all right progenerators (=finitely generated projective generators)
are free, has Aut-Pic. In particular, local rings (by using a
theorem of Kaplansky), principal right ideal domains (by
\cite[Corollary 2.27]{Lam:1999}), and polynomial rings
$k[X_1,\dots,X_n]$, where $k$ is a field (by Quillen-Suslin
Theorem), have Aut-Pic. He also proved that every basic
(semiperfect) ring has Aut-Pic (see \cite[Proposition
3.8]{Bolla:1984}). Moreover, it is well known that a semiperfect
ring $A$ is Morita equivalent to its basic ring $eAe$ (see
\cite[Proposition 27.14]{Anderson/Fuller:1992}), and by Theorem
\ref{coringexseq}, $\pic{A}\simeq \out{eAe}$. In
\cite{Cuadra/Garcia/Torrecillas:2000}, the authors gave several
interesting examples of coalgebras over fields having Aut-Pic. For
instance, they proved that every basic coalgebra has Aut-Pic. On the
other hand, we know \cite[Corollary 2.2]{Chin/Montgomery:1995} that
given a coalgebra over a field $C$ is Morita-Takeuchi equivalent to
a basic coalgebra $C_0$, and by Theorem \ref{coringexseq},
$\pic{C}\simeq \out{C_0}$.

\medskip
Of course all of the examples of rings and coalgebras over fields
having Aut-Pic mentioned in \cite{Bolla:1984} and
\cite{Cuadra/Garcia/Torrecillas:2000} are examples of corings having
Aut-Pic. In order to give other examples of corings satisfying the
Aut-Pic property we need the next lemma which is a generalization of
\cite[Proposition 4.1]{Nastasescu/Torrecillas/VanOystaeyen:2000}.

\medskip
We recall from \cite{Nastasescu/Torrecillas/VanOystaeyen:2000} that
an object $M$ in an additive category $\cat{C}$ has the
\emph{invariant basis number property} (IBN for short) if $M^n\simeq
M^m$ implies $n=m$. For example every nonzero finitely generated
projective module over a semiperfect ring has IBN (from
\cite[Theorem 27.11]{Anderson/Fuller:1992}).

\begin{lemma}\label{IBN}
Let $\coring{C}$ be an $A$-coring such that $_A\coring{C}$ is flat.
If $0\neq M\in{}_B\rcomod{\coring{C}}$ is
$(B,\coring{C})$-quasi-finite such that
\begin{enumerate}[(a)]
\item $B$ is semisimple and $\rcomod{\coring{C}}$ is locally finitely
generated; or
\item $B$ is semiperfect and there is a finitely generated
projective $M_0$ in $\rcomod{\coring{C}}$ such that
$\hom{\coring{C}}{M_{0}}{M}\neq 0$. (the last condition is fulfilled
if $\coring{C}$ is cosemisimple, or (by Proposition
\ref{semiperfect}) if $A$ is right artinian, $_A\coring{C}$ is
projective, and $\coring{C}$ is right semiperfect.);
\end{enumerate}
then $M$ has IBN as a $(B,\coring{C})$-bicomodule.

In particular, if $A$ is semisimple, or $A$ is semiperfect and there
is a nonzero finitely generated projective in $\rcomod{\coring{C}}$,
then $0\neq \coring{C}$ has IBN as an $(A,\coring{C})$-bicomodule.
\end{lemma}

\begin{proof}
We prove at the same time the two statements. Let $M_0$ be a
finitely generated subcomodule of $M$ (resp. a finitely generated
projective in $\rcomod{\coring{C}}$ such that
$\hom{\coring{C}}{M_{0}}{M}\neq 0$). We have
$\cohom{\coring{C}}{M}{M_{0}}^*=\hom{B}{\cohom{\coring{C}}{M}{M_{0}}}{B}\simeq
\hom{\coring{C}}{M_{0}}{M}$ as left $B$-modules. Since the functor
$-\tensor{B}{M}:\rmod{B}\to\rcomod{\coring{C}}$ is exact and
preserves coproducts, the cohom functor
$\cohom{\coring{C}}{\Lambda}{-}:\rcomod{\coring{C}}\to\rmod{B}$
preserves finitely generated (resp. finitely generated projective)
objects. In particular, $\cohom{\coring{C}}{M}{M_{0}}$ is a finitely
generated (resp. finitely generated projective) right $B$-module.
Hence the left $B$-module $\hom{\coring{C}}{M_{0}}{M}$ is so.

Now suppose that $M^m\simeq M^n$. Since the functor
$\hom{\coring{C}}{M_{0}}{-}:{}_B\rcomod{\coring{C}}\to {}\lmod{B}$
is $k$-linear, then $\hom{\coring{C}}{M_{0}}{M^m}\simeq
\hom{\coring{C}}{M_{0}}{M^n}$ as left $B$-bimodules. Hence
$\hom{\coring{C}}{M_{0}}{M}^m \\\simeq \hom{\coring{C}}{M_{0}}{M}^n$
as left $B$-bimodules. Finally, by \cite[Theorem
27.11]{Anderson/Fuller:1992}, $m=n$.

For the particular case we take $M=\coring{C}$.
\end{proof}

\begin{corollary}
Let $C$ be a $k$-coalgebra such that $_kC$ is flat. If $0\neq
M\in{}\rcomod{C}$ is quasi-finite, such that
\begin{enumerate}[(a)]
\item $k$ is semisimple; or
\item $k$ is semiperfect and there is a finitely generated
projective $M_0$ in $\rcomod{C}$ such that $\hom{C}{M_{0}}{M}\neq
0$;
\end{enumerate}
then $M$ has IBN as a ${C}$-comodule.

In particular, if $k$ is semisimple, or $k$ is semiperfect and there
is a nonzero finitely generated projective in $\rcomod{C}$, then
$0\neq C$ has IBN as an $C$-bicomodule.
\end{corollary}

\begin{proposition}
Let $\coring{C}\neq0$ be an $A$-coring such that $A$ is semisimple,
or $A$ is semiperfect and there is a nonzero finitely generated
projective in $\rcomod{\coring{C}}$. If every
$\coring{C}$-bicomodule which is $(A,\coring{C})$-injector, is
isomorphic to $\coring{C}^{(I)}$ as $(A,\coring{C})$-bicomodules for
some set $I$, then the coring $\coring{C}$ has right Aut-Pic.
\end{proposition}

\begin{proof}
The proof is analogous to that of \cite[Proposition
2.4]{Cuadra/Garcia/Torrecillas:2000} and the proof of Bolla
(\cite[p. 264]{Bolla:1984}) of the fact that every ring such that
all left progenerators are free has Aut-Pic.

Let $M$ be a right invertible $\coring{C}$-bicomodule. By
assumptions, $M\simeq \coring{C}^{(I)}$ as
$(A,\coring{C})$-bicomodules for some set $I$. Let $M_0$ be a
finitely generated (resp. finitely generated projective) right
comodule. Since $M_0$ is a small object, then the $k$-linear functor
$\hom{\coring{C}}{M_{0}}{-}:\rcomod{\coring{C}}\to {}\rmod{k}$
preserves coproducts. Therefore, $\hom{\coring{C}}{M_{0}}{M}\simeq
\hom{\coring{C}}{M_{0}}{\coring{C}}^{(I)}$ as left $A$-modules.
Hence I is a finite set, and $M\simeq \coring{C}^{(n)}$ as
$(A,\coring{C})$-bicomodules for some $n\geq1$. Let $N$ be a
bicomodule such that $(N)$ is the inverse of $(M)$ in
$\rpic{\coring{C}}$, then
$$\coring{C}\simeq M\cotensor{\coring{C}}N\simeq
\coring{C}^{(n)}\cotensor{\coring{C}}N\simeq N^{(n)}$$ as
$(A,\coring{C})$-bicomodules. On the other hand, $N\simeq
\coring{C}^{(m)}$ as $(A,\coring{C})$-bicomodules for some set
$m\geq1$. It follows that $\coring{C}\simeq \coring{C}^{(nm)}$ as
$(A,\coring{C})$-bicomodules. By Lemma \ref{IBN}, $nm=1$ and then
$M\simeq \coring{C}$ as $(A,\coring{C})$-bicomodules. Finally, from
Proposition \ref{oneside}, $M\simeq {}_f\coring{C}$ as bicomodules
for some $f=(\varphi,1_A)\in \aut{\coring{C}}$.
\end{proof}

\begin{corollary}
Let $C$ be a $k$-coalgebra such that $k$ is semisimple, or $k$ is a
QF ring and there is a nonzero finitely generated projective in
$\rcomod{C}$. If every right injective $C$-comodule is free, then
the coalgebra $C$ has right Aut-Pic.
\end{corollary}

The following result allows to simplify the computation of the
Picard group of some interesting corings.

\begin{proposition}\label{examples}
\begin{enumerate}[(1)] \item Let $_B\Sigma_A$ be a bimodule
such that $\Sigma_A$ is finitely generated and projective and
$_B\Sigma$ is faithfully flat, and
$(\Sigma^*\tensor{B}\Sigma)_A$ is flat. Then
$\rpic{\Sigma^*\tensor{B}\Sigma}\simeq \pic{B}$, where
$\Sigma^*\tensor{B}\Sigma$ is the comatrix coring associated to
$\Sigma$. If moreover $B$ has Aut-Pic then
$\rpic{\Sigma^*\tensor{B}\Sigma}\simeq \out{B}$.
\\ In particular, For a $k$-algebra $A$ and $n\in \mathbb{N}$, we
have $\rpic{M^c_n(A)}\simeq \pic{A}$, where $M^c_n(A)$ is the
$(n,n)$-matrix coring over $A$ defined in Example \ref{examples
corings} (6). If moreover $A$ has Aut-Pic then
$\rpic{M^c_n(A)}\simeq \out{A}$.
\item Let $\coring{C}$ be an $A$-coring such that $\coring{C}_A$ is
flat, and $R$ the opposite algebra of $^*\coring{C}$. If
$_A\coring{C}$ is finitely generated projective (e.g. if
$\coring{C}$ is a Frobenius ring), then $\rpic{\coring{C}}\simeq
\pic{R}$.
 \item If $\coring{C}$ is an $A$-coring such that $_A\coring{C}$
 and $\coring{C}_A$ are flat and the category $\rcomod{\coring{C}}$
 has a small projective generator $U$, then $\rpic{\coring{C}}\simeq
\pic{\operatorname{End}_\coring{C}(U)}$.
\end{enumerate}
\end{proposition}

\begin{proof}
(1) From Corollary \ref{KG 3.10} (1), we have an equivalence of
categories $-\tensor{B}\Sigma:{}\rmod{B}\to
\rcomod{\Sigma^*\tensor{B}\Sigma}$. Theorem \ref{coringexseq}
achieves then the proof. Now we prove the particular case. Let $A$
be a $k$-algebra and $n\in \mathbb{N}$. If we take $\Sigma=A^n$, the
comatrix coring $\Sigma^*\tensor{B}\Sigma$ can be identified with
$M^c_n(A)$.

(2) From Lemma \ref{B 4.3}, the categories $\rcomod{\coring{C}}$
and $\rmod{R}$ are isomorphic to each other. It is enough to apply
Theorem \ref{coringexseq}.

(3) This follows immediately from Theorems \ref{Harada} and
\ref{coringexseq}.
\end{proof}

\section{Application to the Picard group of $gr-(A,X,G)$} \label{Picard graded}

We begin by giving two useful lemmas whose proofs are
straightforward.

\begin{lemma} A morphism $(\alpha,\gamma)\in
\mathbb{E}_\bullet^\bullet(k)$ is an isomorphism if and only if
$\alpha$ and $\gamma$ are bijective.
\end{lemma}

\begin{lemma}
 If
$(\alpha,\gamma):(A,C,\psi)\rightarrow(A',C',\psi')$ is a morphism
of entwined structures, then $(\alpha\otimes\gamma,\alpha):A\otimes
C\to A\otimes C$ is a morphism of corings. Hence we have a functor
$$F:\mathbb{E}_\bullet^\bullet(k)\to\mathbf{Crg}_k.$$
\end{lemma}

By Theorem \ref{coringexseq}, we obtain

\begin{proposition}\label{entwining}
For an entwining structure $(A,C,\psi)$ there is an exact sequence
$$\xymatrix{1\ar[r] & \operatorname{Ker}(\Omega\circ F)\ar[r]&
\aut{(A,C,\psi)} \ar[r]^-{\Omega\circ F} & \rpic{A\otimes C}}\simeq
\rpic{\mathcal{M}(\psi)_A^C},$$ and $\operatorname{Ker}(\Omega\circ
F)$ is the set of $(\alpha,\gamma)\in \aut{(A,C,\psi)}$ such that
there is $p\in (A\otimes C)^*$ satisfying $\sum
\big(\alpha(a)\otimes\gamma(c_{(1)})\big)p(1_A\otimes c_{(2)})=\sum
p(a\otimes c_{(1)})\otimes c_{(2)}$ for all $a\in A,c\in C$.
\end{proposition}

\begin{lemma}
A morphism $(\hbar,\alpha,\gamma)\in \mathbb{DK}_\bullet^\bullet(k)$
is an isomorphism if and only if $\hbar$, $\alpha$ and $\gamma$ are
bijective.
\end{lemma}

\begin{proof}
Straightforward.
\end{proof}

Now consider the functor $G:\mathbb{DK}_\bullet^\bullet(k)\to
\mathbb{E}_\bullet^\bullet(k)$ defined in Proposition \ref{CMZ 17}.
Hence Proposition \ref{entwining} yields the following

\begin{proposition}\label{DK}
For a DK structure $(H,A,C)$ there is an exact sequence
$$\xymatrix{1\ar[r] & \operatorname{Ker}(\Omega\circ F\circ G)\ar[r]&
\aut{(H,A,C)}\ar[rr]^-{\Omega\circ F\circ G} && \rpic{A\otimes
C}\simeq\rpic{\mathcal{M}(H)_A^C}},$$ and
$\operatorname{Ker}(\Omega\circ F\circ G)$ is the set of
$(\hbar,\alpha,\gamma)\in \aut{(H,A,C)}$ such that there is $p\in
(A\otimes C)^*$ satisfying $\sum
\big(\alpha(a)\otimes\gamma(c_{(1)})\big)p(1_A\otimes c_{(2)})=\sum
p(a\otimes c_{(1)})\otimes c_{(2)}$ for all $a\in A,c\in C$.
\end{proposition}

Finally, we consider the category of right modules graded by a
$G$-set, $gr-(A,X,G)$, where $G$ is a group and $X$ is a right
$G$-set. This category is equivalent to a category of modules over a
ring if $X=G$ and $A$ is strongly graded (by Dade's theorem
\cite[Theorem 3.1.1]{Nastasescu/VanOystaeyen:2004}), or if $X$ is a
finite set (by \cite[Theorem
2.13]{Nastasescu/Raianu/VanOystaeyen:1990}). But there are examples
of a category of graded modules which is not equivalent to a
category of modules, see \cite[Remark 2.4]{Menini/Nastasescu:1988}.

Let $G$ be a group, $X$ a right $G$-set, and $A$ be a $G$-graded
$k$-algebra. We know that $(kG,A,C)$ is a DK structure, with $kG$ is
a Hopf algebra. Then we have, $(A,kX,\psi)$ is an entwining
structure where $\psi:kX\otimes A\rightarrow A\otimes kX$ is the map
defined by $\psi(x\otimes a_{g})=a_{g}\otimes xg$ for all $x\in
X,g\in G,a_g\in A_g$. $A\otimes kX$ is an $A$-coring with the
comultiplication and the counit maps:
$$\Delta_{A\otimes kX}(a\otimes x)=(a\otimes
x)\tensor{A}(1_A\otimes x), \quad \epsilon_{A\otimes kX}(a\otimes
x)=a \quad (a\in A,x\in X).$$ We know that this coring is
coseparable (see Proposition \ref{coseparable}). From Subsection
\ref{CategoryGraded}, $gr-(A,X,G)\simeq \mathcal{M}(kG)_A^{kX}$.
Moreover, we know that the family $\{1_A\otimes x\mid x\in X\}$ is a
basis of both the left and the right $A$-module $A\otimes kX$ (see
the proof of Lemma \ref{graded-A4}). Each $p\in (A\otimes kX)^*$ is
entirely determined by the data of $p(1_A\otimes x)$, for all $x\in
X$. The same thing holds for $\varphi$, where $(\varphi,\rho)\in
\aut{A\otimes kX}$.

\medskip
By Theorem \ref{coringexseq} we get

\begin{proposition}
\begin{enumerate}[(a)]\item We have, $p\in (A\otimes kX)^*$ is
invertible
 if and only if there exists $q\in (A\otimes kX)^*$ such that
 $\sum_{h\in G}q(1_A\otimes xh^{-1})p(1_A\otimes
x)_h=1_A=\sum_{h\in G}p(1_A\otimes xh^{-1})q(1_A\otimes x)_h$ for
every $x\in X$.
\item We have an exact sequence
$$\xymatrix{1\ar[r] & \operatorname{Ker}(\Omega)\ar[r]&
\aut{A\otimes kX}\ar[r]^-\Omega & \pic{A\otimes kX}\simeq
\pic{gr-(A,X,G)}},$$ and $(\varphi,\rho)\in
\operatorname{Ker}(\Omega)$ if and only if there exists $p\in
(A\otimes kX)^*$ invertible such that
\begin{enumerate}[(i)] \item for
all $x\in X$, $a_y^xp(1_A\otimes x)_h=0$ for all $y\in X,h\in G$
such that $yh\neq x$, where $\varphi(1_A\otimes x)=\sum_{y\in
X}a_y^x\otimes y$, and \item for all $x\in X,a\in A$, $p(a\otimes
x)=\rho(a)p(1_A\otimes x).$
\end{enumerate}
\end{enumerate}
\end{proposition}

\begin{proof}
We have, $p\in (A\otimes kX)^*$ is invertible if and only if there
exists $q\in (A\otimes kX)^*$ such that
\begin{eqnarray}
q\big(p(a\otimes x)(1_A\otimes x)\big)=\epsilon(a\otimes x)=a
\label{G-Inv-1}\\
p\big(q(a\otimes x)(1_A\otimes x)\big)=\epsilon(a\otimes x)=a
\label{G-Inv-2},
\end{eqnarray}
for all $a\in A, x\in X$. (from Proposition \ref{BW 17.8}.) On the
other hand,
\begin{eqnarray*}
q\big(p(a_g\otimes x)(1_A\otimes x)\big) &=& q\big(p(1_A\otimes
xg^{-1})(a_g\otimes x)\big) \\
&=& q\big(\sum_{h\in G}(b_ha_g)\otimes
x\big),\quad\textrm{where}\quad p(1_A\otimes xg^{-1})=\sum_{h\in
G}b_h \\
&=& \sum_{h\in G}q\big(1_A\otimes x(hg)^{-1}\big)(b_ha_g), \quad
\textrm{since}\quad b_ha_g\in A_{hg} \\
&=& \sum_{h\in G}q(1_A\otimes xg^{-1}h^{-1})p(1_A\otimes
xg^{-1})_ha_g,
\end{eqnarray*}
for all $x\in X,g\in G, a_g\in A_g$. Hence
\begin{eqnarray*}
 \eqref{G-Inv-1}&\Longleftrightarrow
 & q\big(p(a_g\otimes x)(1_A\otimes x)\big)=
a_g\quad \textrm{for all}\, x\in X,g\in G, a_g\in A_g \\
&\Longleftrightarrow& \sum_{h\in G}q(1_A\otimes xh^{-1})p(1_A\otimes
x)_h=1_A \quad \textrm{for all}\, x\in X.
\end{eqnarray*}
By symmetry, $$ \eqref{G-Inv-2}\Longleftrightarrow \sum_{h\in
G}p(1_A\otimes xh^{-1})q(1_A\otimes x)_h=1_A \quad \textrm{for
all}\, x\in X.$$

By Theorem \ref{coringexseq}, $(\varphi,\rho)\in
\operatorname{Ker}(\Omega)$ if and only if there exists $p\in
(A\otimes kX)^*$ invertible with $\varphi(a\otimes x)p(1_A\otimes
x)=p(a\otimes x)\otimes x$ for all $a\in A,x\in X$. On the other
hand,
\begin{eqnarray*}
\varphi(a\otimes x)p(1_A\otimes x)&=& \rho(a)\varphi(1_A\otimes
x)p(1_A\otimes x) \\
&=& \rho(a)\sum_{y\in X}\sum_{h\in G}(a_y^x\otimes
y)b_h, \\
&& \textrm{where}\quad \varphi(1_A\otimes x)=\sum_{y\in
X}a_y^x\otimes y,\quad p(1_A\otimes x)=\sum_{h\in G}b_h, \\
&=& \rho(a)\sum_{y\in X}\sum_{h\in G}a_y^xb_h\otimes yh.
\end{eqnarray*}
Since $\sum_{y\in X}a_y^x=\epsilon\circ\varphi(1_A\otimes
x)=\rho\circ\epsilon(1_A\otimes x)=1_A$, (b) follows.
\end{proof}

Now, let $f:G\to G'$ be a morphism of groups, $X$ a right $G$-set,
$X'$ a right $G'$-set, and $\varphi:X\rightarrow X'$ a map. Let $A$
be a $G$-graded $k$-algebra, $A'$ a $G'$-graded $k$-algebra, and
$\alpha:A\rightarrow A'$ a morphism of algebras. We have,
$\hbar:kG\rightarrow kG'$ defined by $\hbar(g) =f(g)$ for each $g\in
G$, is a morphism of Hopf algebras, and $\gamma:kX\rightarrow kX'$
defined by $\gamma(x) =\varphi(x)$ for each $x\in X$, is a morphism
of coalgebras. It is easy to show that $(\hbar,\alpha,\gamma)$ is a
morphism of DK structures if and only if
\begin{eqnarray}
\varphi(xg) =\varphi(x)f(g)\quad  \textrm{for all}\, g\in G,x\in X
\label{DK-G-1}\\
\alpha( A_{g}) \subset A_{f(g)}' \quad\textrm{for all}\, g\in G
\label{DK-G-2}.
\end{eqnarray}

We define a category as follows. The objects are the triples
$(G,X,A)$, where $G$ is a group, $X$ is a right $G$-set, and $A$ is
a graded $k$-algebra. The morphisms are the triples
$(f,\varphi,\alpha)$, where $f:G\rightarrow G'$ is a morphism of
groups, $\varphi:X\rightarrow X'$ is a map, and $\alpha:A\rightarrow
A'$ is a morphism of algebras such that \eqref{DK-G-1} and
\eqref{DK-G-2} hold. We denote this category by $\mathbb{G}^r(k)$.
There is a faithful functor
$$H:\mathbb{G}^r(k)\to \mathbb{DK}_\bullet^\bullet(k)$$ defined by
$H((G,X,A))=(kG,A,kX)$ and
$H((f,\varphi,\alpha))=(\hbar,\alpha,\gamma)$. We obtain then the
next result. We think that our subgroup of $\pic{gr-(A,X,G)}$ is a
natural subgroup and it is more simple than the Beattie-Del R\'io
subgroup (see \cite[\S2]{Beattie/Delrio:1996}).

\begin{proposition}
There is an exact sequence

$\xymatrix{1\ar[r] & \operatorname{Ker}(\Omega\circ F\circ G\circ
H)\ar[r]& \aut{(G,X,A)}\ar[rr]^-{\Omega\circ F\circ G\circ H} & &
\pic{A\otimes kX}}$

$\xymatrix{&&&&&&&&&&& \simeq \pic{gr-(A,X,G)}},$ and
$\operatorname{Ker}(\Omega\circ F\circ G\circ H)$ is the set of
$(f,\varphi,\alpha)\in \aut{(G,X,A)}$ such that there is $p\in
(A\otimes kX)^*$ invertible satisfying
\begin{enumerate}[1)]
\item $p(a\otimes x)=\alpha(a)p(1_A\otimes x)$ for all $x\in X,a\in
A$, and
\item $p(1_A\otimes x)_h=0$ for all $x\in X, h\in G$ such
that $\varphi(x)h\neq x$.
\end{enumerate}
\end{proposition}

\end{document}